\documentclass[mnsc,nonblindrev]{informs3-VVM-noheader}

\SingleSpacedXI
\usepackage{natbib}
 \bibpunct[, ]{(}{)}{,}{a}{}{,}%
 \def\bibfont{\small}%

\usepackage{amsmath}
\usepackage{amssymb}
\usepackage{enumerate}
\usepackage{graphicx}

\usepackage{algorithm}
\usepackage{algorithmic}

\usepackage{subcaption}
\usepackage{color}

\usepackage[colorinlistoftodos, textwidth=\marginparwidth, textsize = footnotesize]{todonotes}
\usepackage{ctable}
\usepackage{forest}

\def\xb{\mathbf{x}}
\def\Xcal{\mathcal{X}}
\def\Ibb{\mathbb{I}}

\def\Acal{\mathcal{A}}

\def\stop{\mathbf{stop}}
\def\go{\mathbf{go}}
\def\Exp{\mathbb{E}}

\def\Rbb{\mathbb{R}}

\def\zb{\mathbf{z}}
\def\yb{\mathbf{y}}
\def\tipi{\tilde{\pi}}
\def\Fcal{\mathcal{F}}

\def\Lcal{\mathcal{L}}
\def\ib{\mathbf{i}}
\def\Ibb{\mathbb{I}}

\def\Tcal{\mathcal{T}}
\def\leftchild{\mathbf{leftchild}}
\def\rightchild{\mathbf{rightchild}}
\def\Ncal{\mathcal{N}}
\def\leaves{\mathbf{leaves}}
\def\splits{\mathbf{splits}}
\def\dir{D}

\def\vb{\mathbf{v}}
\def\thetab{\boldsymbol \theta}
\def\ab{\mathbf{a}}

\def\Bcal{\mathcal{B}}

\newcommand{\defeq}{\ensuremath{\triangleq}}

\DeclareMathOperator{\st}{s.t.}
\DeclareMathOperator{\as}{a.s.}

\TheoremsNumberedThrough     %
\ECRepeatTheorems

\EquationsNumberedThrough    %
\MANUSCRIPTNO{--}

\begin{document}
\RUNAUTHOR{Ciocan and Mi\v{s}i\'{c}}

\RUNTITLE{Interpretable Optimal Stopping}

\TITLE{Interpretable Optimal Stopping}

\ARTICLEAUTHORS{%
\AUTHOR{Dragos Florin Ciocan}
\AFF{INSEAD; Boulevard de Constance 77305, Fontainebleau, France, \EMAIL{florin.ciocan@insead.edu}} %
\AUTHOR{Velibor V. Mi\v{s}i\'c}
\AFF{Anderson School of Management, University of California, Los Angeles; 110 Westwood Plaza, Los Angeles, CA 90095, USA, \EMAIL{velibor.misic@anderson.ucla.edu}}
} %

\ABSTRACT{%
Optimal stopping is the problem of deciding when to stop a stochastic system to obtain the greatest reward, arising in numerous application areas such as finance, healthcare and marketing. State-of-the-art methods for high-dimensional optimal stopping involve approximating the value function or the continuation value, and then using that approximation within a greedy policy. Although such policies can perform very well, they are generally not guaranteed to be interpretable; that is, a decision maker may not be able to easily see the link between the current system state and the policy's action. In this paper, we propose a new approach to optimal stopping, wherein the policy is represented as a binary tree, in the spirit of naturally interpretable tree models commonly used in machine learning. We show that the class of tree policies is rich enough to approximate the optimal policy. We formulate the problem of learning such policies from observed trajectories of the stochastic system as a sample average approximation (SAA) problem. We prove that the SAA problem converges under mild conditions as the sample size increases, but that computationally even immediate simplifications of the SAA problem are theoretically intractable. We thus propose a tractable heuristic for approximately solving the SAA problem, by greedily constructing the tree from the top down. We demonstrate the value of our approach by applying it to the canonical problem of option pricing, using both synthetic instances and instances using real S\&P-500 data. Our method obtains policies that (1) outperform state-of-the-art non-interpretable methods, based on simulation-regression and martingale duality, and (2) possess a remarkably simple and intuitive structure. }%


\maketitle

\section{Introduction}
\label{sec:introduction}

We consider the problem of optimal stopping, which can be described as follows: a system evolves stochastically from one state to another in discrete time steps. At each decision epoch, a decision maker (DM) chooses whether to stop the system or allow it to continue for one more time step. If the DM chooses to stop the system, she garners a reward that is dependent on the current state of the system; if she chooses to allow it to continue, she does not receive any reward in the current period, but can potentially stop it at a future time to obtain a higher reward. The DM must specify a policy, which prescribes the action to be taken (stop/continue) for each state the system may enter. The optimal stopping problem is to find the policy that achieves the highest possible reward in expectation.

The optimal stopping problem is a key problem in stochastic control and arises in many important applications; we name a few below:
\begin{enumerate}
\item \textbf{Option pricing}. One of the most important applications of optimal stopping is to the pricing of financial options that allow for early exercise, such as American and Bermudan options. The system is the collection of underlying securities (typically stocks) that the option is written on. The prices of the securities comprise the system state. The decision to stop the system corresponds to exercising the option and receiving the corresponding payoff; thus, the problem of obtaining the highest expected payoff from a given option corresponds to finding an optimal stopping policy. The highest expected payoff attained by such an optimal stopping (optimal exercise) policy is then the price that the option writer should charge for the option. 
\item \textbf{Healthcare}. Consider an organ transplant patient waiting on a transplant list, who may be periodically offered an organ for transplantation. In this context, the system corresponds to the patient and the currently available organ, and the system state describes the patient's health and the attributes of the organ. The decision to stop corresponds to the patient accepting the organ, where the reward is the estimated quality-adjusted life years (QALYs) that the patient will garner upon transplantation. The problem is then to find a policy that prescribes for a given patient and a given available organ whether the patient should accept the organ, or wait for the next available organ, so as to maximize the QALYs gained from the transplant.%
\item \textbf{Marketing}. %
Consider a retailer selling a finite inventory of products over some finite time horizon. The system state describes the remaining inventory of the products, which evolves over time as customers buy the products from period to period. The action of stopping the system corresponds to starting a price promotion that will last until the end of the finite time horizon. The problem is to decide at each period, based on the remaining inventory, whether to commence the promotion, or to wait one more period, in order to maximize the total revenue garnered by the end of the horizon.
\end{enumerate}

Optimal stopping problems are solved via dynamic programming (DP); the paradigm here is to solve the DP and obtain a value function that, for a given system state, specifies the best possible expected reward that can be attained when one starts in that state. With this value function in hand, one can obtain a good policy by considering the greedy policy with respect to the value function. Due to the well known curse of dimensionality, large-scale optimal stopping problems that occur in practice are typically solved by approximate dynamic programming (ADP) methods. The goal in such ADP methods is to replace the true value function with an approximate value function that can be computed more tractably.

In either DP or ADP methods, the value function may provide some insight into which states of the system state space are more desirable. However, the policy that one obtains by being greedy with respect to the this value function need not have any readily identifiable structure and indeed, the DP paradigm is to use the mapping from states to value functions and then actions in a ``black-box'' fashion. This is disadvantageous, because in many optimal stopping problems, we are not only interested in policies that attain high expected reward, but also policies that are \emph{interpretable}. A policy that is interpretable is one where we can directly see how the state of the system maps to the recommended action, and the relation between state and action is sufficiently transparent.

Interpretability is desirable for three reasons. First, in modeling a real world system, it is useful to obtain some insight about what aspects of the system state are important for controlling it optimally or near optimally. Second, a complex policy, such as a policy that is obtained via an ADP method, may not be operationally feasible in many real life contexts. Lastly -- and most importantly -- in many application domains where the decision maker is legally responsible or accountable for the action taken by a policy, interpretability is not merely a desirable feature, but a requirement: in such settings a decision maker will simply not adopt the policy without the ability to explain the policy's mechanism of action. There is moreover a regulatory push to increase the transparency and interpretability of customer facing data-driven algorithms; as an example, General Data Protection Regulation rules set by the EU \citep{doshi2017towards} dictate that algorithms which can differentiate between users must provide explanations for their decisions if such queries arise.

In this paper, we consider the problem of constructing interpretable optimal stopping policies from data. Our approach to interpretability is to consider policies that are representable as a binary tree. In such a tree, each leaf node is an action (stop or go), and each non-leaf node (also called a \emph{split} node) is a logical condition in terms of the state variables that determines whether we proceed to the left or the right child of the current node. To determine the action we should take, we take the current state of the system, run it down the tree until we reach a leaf, and take the action prescribed in that leaf. An example of such a tree-based policy is presented in Figure~\ref{figure:example_tree_policy}. Policies of this kind are simple, and allow the decision maker to directly see the link between the current system state and the action.

Before delving into our results, we comment on how one might compare such an aforementioned interpretable policy with an optimal or, possibly, near-optimal heuristic policy. Our informal goal is to look for the best policy within the constraints of an interpretable policy architecture, which in our case is tree-based policies. However, without any a priori knowledge that a given optimal stopping problem is ``simple enough'', one would expect that enforcing that the policy architecture be interpretable carries a performance price; in other words, interpretable policies should generally not perform as well as a state-of-the-art heuristic. On the other hand, one can hope that there exist stopping problems where the price of interpretability is low, in that interpretable policies, while sub-optimal, do not carry a large optimality gap. The present paper is an attempt to (a) exhibit stopping problems of practical interest for which tree-based policies attain near-optimal  performance, along with being interpretable and (b) provide an algorithm to find such interpretable policies directly from data. 

We make the following specific contributions:
\begin{enumerate}
\item \textbf{Sample average approximation.} We formulate the problem of learning an interpretable tree-based policy from a sample of trajectories of the system (sequences of states and payoffs) as a sample average approximation (SAA) problem. To the best of our knowledge, the problem of directly learning a policy in the form of a tree for optimal stopping problems (and Markov decision processes more generally) has not been studied before. We show that under mild conditions, the tree policy SAA problem defined using a finite sample of trajectories converges almost surely in objective value to the tree policy problem defined with the underlying stochastic process, as the number of trajectories available as data grows large. We also prove that one can approximate the optimal policy to an arbitrary precision by a tree policy of sufficient depth.  %
\item \textbf{Computational tractability.} From a computational complexity standpoint, we establish that three natural simplifications of the SAA problem are NP-Hard and thus finding good solutions to the SAA problem is challenging. 

In response, we present a computationally tractable methodology for solving the learning problem. Our method is a construction algorithm: starting from a degenerate tree consisting of a single leaf, the algorithm grows the tree in each iteration by splitting a single leaf into two new child leaves. The split is chosen greedily, and the algorithm stops when there is no longer a sufficient improvement in the sample-based reward. Key to the procedure is determining the optimal split point at a candidate split; we show that this problem can be solved in a computationally efficient manner. In this way, the overall algorithm is fully data-driven, in that the split points are directly chosen from the data and are not artificially restricted to a collection of split points chosen a priori. While the algorithm resembles top-down tree induction methods from classification/regression, several important differences arising from the temporal nature of the problem make this construction procedure algorithmically nontrivial. %
\item \textbf{Practical performance versus state-of-the-art ADP methods.} Using both synthetic and real S\&P-500 data, we numerically demonstrate the value of our methodology by applying it to the problem of pricing a Bermudan option, which is a canonical optimal stopping problem in finance. We show that our tree policies outperform two state-of-the-art approaches, namely the simulation-regression approach of \cite{longstaff2001valuing} and the martingale duality-based approach of \cite{desai2012pathwise}. At the same time, we also show that the tree policies produced by our approach are remarkably simple and intuitive. We further investigate the performance of our approach by testing it on a stylized one-dimensional optimal stopping problem (not drawn from option pricing), where the exact optimal policy can be computed; in general, our policy is either optimal or very close to optimal. 
\end{enumerate}

The rest of this paper is organized as follows. In Section~\ref{sec:literature_review}, we discuss the relevant literature in ADP, machine learning and interpretable decision making. In Section~\ref{sec:problem}, we formally define our optimal stopping problem and its sample-based approximation, we define the problem of finding a tree policy from sample data, and theoretically analyze this problem. In Section~\ref{sec:construction}, we present a heuristic procedure for greedily constructing a tree directly from data. In Section~\ref{sec:results}, we present an extensive computational study in option pricing comparing our algorithm to alternate approaches. In Section~\ref{sec:results_1D_Uniform}, we evaluate our algorithm on the aforementioned one-dimensional problem. Finally, we conclude in Section~\ref{sec:conclusion}.

\section{Literature review}
\label{sec:literature_review}

Our paper relates to three different broad streams of research: the optimal stopping and ADP literature; the machine learning literature; and the growing literature on interpretable decision making. We survey each of these below. \\

\noindent \textbf{Approximate dynamic programming (ADP)}. ADP has been extensively studied in the operations research community since the mid-1990s as a solution technique for Markov decision processes \citep{powell2007approximate,van2002neuro}. In the last fifteen years, there has been significant interest in solving MDPs by approximating the linear optimization (LO) model; at a high level, one formulates the MDP as a LO problem, reduces the number of variables and constraints in a tractable manner, and solves the more accessible problem to obtain a value function approximation. Some examples of this approach include \cite{defarias2003linear}, \cite{adelman2008relaxations}, \cite{desai2012approximate} and \cite{bertsimas2016decomposable}. 

A large subset of the ADP research literature, originating in both the operations research and finance communities, has specifically studied optimal stopping problems. The seminal papers of \cite{carriere1996valuation}, \cite{longstaff2001valuing} and \cite{tsitsiklis2001regression} propose simulation-regression approaches, where one simulates trajectories of the system state and uses regression to compute an approximation to the optimal continuation value at each step. Later research has considered the use of martingale duality techniques. The idea in such approaches is to relax the non-anticipativity requirement of the policy by allowing the policy to use future information, but to then penalize policies that use this future information, in the spirit of Lagrangean duality. This approach yields upper bounds on the optimal value and can also be used to derive high quality stopping policies. Examples include \cite{rogers2002monte}, \cite{andersen2004primal}, \cite{haugh2004pricing}, \cite{chen2007additive}, \cite{brown2010information}, \cite{desai2012pathwise} and \cite{goldberg2018beating}.  

Our approach differs from these generic ADP approaches and optimal stopping-specific approaches in two key ways. First, general purpose ADP approaches, as well as those specifically designed for optimal stopping, are focused on obtaining an approximate value function or an upper bound on the value function, which is then used in a greedy manner. In contrast, our approach involves optimizing over a policy directly, without computing/optimizing over a value function. Second, our approach is designed with interpretability in mind, and produces a policy that can be easily visualized as a binary tree. In contrast, previously proposed ADP methods are not guaranteed to result in policies that are interpretable. 

Lastly, we note that some methods in finance for pricing options with early exercise involve tree representations; examples include the binomial lattice approach \citep{cox1979option} and the random tree method \citep{broadie1997pricing}. However, the trees found in these methods represent discretizations of the sample paths of the underlying asset prices, which provide a tractable way to perform scenario analysis. In contrast, the trees in our paper represent policies, not sample paths. As such, our approach is unrelated to this prior body of work. \\

\noindent \textbf{Machine learning}. Interpretability has been a goal of major interest in the machine learning community, starting with the development of decision trees in the 1980s \citep{breiman1984classification,quinlan1986induction,quinlan1993c4}. A stream of research has considered interpretable scoring rules for classification and risk prediction; recent examples include \cite{ustun2015supersparse,ustun2016learning} and \cite{zeng2017interpretable}. Another stream of research considers the design of disjunctive rules and rule lists; recent examples include \cite{wang2015or,wang2015learning, wang2017bayesian}, \cite{letham2015interpretable}, \cite{angelino2017learning} and \cite{lakkaraju2016interpretable}. Other research has also considered how to extract interpretable models from complicated black-box models \citep{bastani2018interpreting}. %

The algorithm we will present is closest in spirit to classical tree algorithms like CART and ID3. Our algorithm differs from these prior methods in that it is concerned with optimal stopping, which is a stochastic control problem that involves making a decision over time, and is fundamentally different from classification and regression. In particular, a key part of estimating a classification or regression tree is determining the leaf labels, which in general is a computationally simple task. As we will see in Section~\ref{subsec:problem_complexity}, the analogous problem in the optimal stopping realm is NP-Hard. As a result, this leads to some important differences in how the tree construction must be done to account for the temporal nature of the problem; we comment on these in more detail in Section~\ref{subsec:construction_comparison_to_classification}. \\

\noindent \textbf{Interpretable decision making}. In the operations research community, there is growing interest in interpretability as it pertains to dynamic decision making; we provide some recent examples here. 
With regard to dynamic problems, \cite{bertsimas2013fairness} considers the problem of designing a dynamic allocation policy for allocating deceased donor kidneys to patients requiring a kidney transplant that maximizes efficiency while respecting certain fairness constraints. To obtain a policy that is sufficiently interpretable to policy makers, \cite{bertsimas2013fairness} further propose using ordinary least squares regression to find a scoring rule that predicts the fairness-adjusted match quality as a function of patient and donor characteristics. In more recent work, \cite{azizi2018designing} consider a related approach for dynamically allocating housing resources to homeless youth. The paper proposes a general mixed-integer optimization framework for selecting an interpretable scoring rule (specifically, linear scoring rules, decision tree rules with axis-aligned or oblique splits, or combinations of both linear and decision tree rules) for prioritizing youth on the waiting list. %
Lastly, the paper of \cite{bravo2017discovering} considers the use of machine learning for analyzing optimal policies to MDPs. The key idea of the paper is to solve instances of a given MDP to optimality, and then to use the optimal policies as inputs to machine learning methods. The paper applies this methodology to classical problems such as inventory replenishment, admission control and multi-armed bandits. This approach differs from ours, in that we do not have access to the optimal policy; in fact, the goal of our method is to directly obtain a near-optimal policy. Stated differently, the paper of \cite{bravo2017discovering} seeks to understand how interpretable optimal policies are, whereas our paper seeks to design interpretable policies that deliver good performance. %

\section{Problem definition}
\label{sec:problem}

We begin by defining the optimal stopping problem in Section~\ref{subsec:problem_model}, and its sample-based counterpart in Section~\ref{subsec:problem_SAA}. We then define the tree policy sample-average approximation (SAA) problem, where the class of policies is restricted to those that can be represented by a binary tree, in Section~\ref{subsec:problem_trees}. We show that the tree policy SAA problem converges in Section~\ref{subsec:problem_convergence}. In Section~\ref{subsec:problem_complexity}, we show that the tree policy SAA problem is NP-Hard when one considers three specific simplifications. Finally, in Section~\ref{subsec:tree_optimality}, we establish that under mild conditions, any optimal policy can be approximated to an arbitrary precision by a tree policy of sufficient depth.

\subsection{Optimal stopping model}
\label{subsec:problem_model}

Consider a system with state given by $\xb = (x_1, \dots, x_n) \in \Xcal_1 \times \dots \times \Xcal_n \defeq \Xcal \subseteq \mathbb{R}^n$. We let $\xb(t)$ denote the state of the system at time $t$, which evolves according to some stochastic process. We let $\Acal = \{ \stop, \go \}$ denote the action space of the problem; we may either stop the system ($\stop$) or allow it to continue for one more time step ($\go$). We assume a finite horizon problem with $T$ periods, starting at period $t = 1$. We let $g(t,\xb)$ denote the reward or payoff from stopping the system when the current period is $t$ and the current system state is $\xb$. We assume that all rewards are discounted by a factor of $\beta$ for each period. 

We define a policy $\pi$ as a mapping from the state space $\Xcal$ to the action space $\Acal$. We let $\Pi = \{ \pi \, | \, \pi: [T] \times \Xcal \to \Acal \}$ be the set of all possible policies, where we use the notation $[N] = \{1,\dots, N\}$ for any integer $N$. For a given realization of the process $\{ \xb(t) \}_{t=1}^T$, we define the stopping time $\tau_{\pi}$ as the first time at which the policy $\pi$ prescribes the action $\stop$: 
\begin{equation}
\tau_{\pi} = \min \{ t \in [T] \, | \,  \pi(t, \xb(t) ) = \stop \},
\end{equation}
where we take the minimum to be $+\infty$ if the set is empty. Our goal is to find the policy that maximizes the expected discounted reward over the finite horizon, which can be represented as the following optimization problem:
\begin{equation}
\underset{\pi \in \Pi}{\text{maximize}} \ \Exp \left[  \beta^{ \tau_\pi - 1} \cdot g( \tau_{\pi}, \xb( \tau_\pi ) ) \ | \ \xb(1) = \xb \right].  \label{prob:optimal_stopping_canonical}
\end{equation}
For any policy $\pi$, we let $J^\pi(\xb) \defeq \Exp \left[  \beta^{ \tau_\pi - 1} \cdot g( \tau_{\pi}, \xb( \tau_\pi ) ) \ | \ \xb(1) = \xb \right]$. For simplicity, we assume that there exists a starting state $\bar\xb$ at which the system is started, such that $\xb(1) = \bar\xb$  always. %

\subsection{Sample average approximation} 
\label{subsec:problem_SAA}

In order to solve problem~\eqref{prob:optimal_stopping_canonical}, we need to have a full specification of the stochastic process $\{ \xb(t) \}_{t=1}^T$. In practice, we may not have this specification or it may be too difficult to work with directly. Instead of this specification, we may instead have data, that is, we may have access to specific realizations or trajectories of the process $\{ \xb(t) \}_{t=1}^T$. In this section, we describe a \emph{sample-average approximation} (SAA) formulation of the optimal stopping problem~\eqref{prob:optimal_stopping_canonical} that will allow us to design a policy directly from these trajectories, as opposed to a probabilistic definition of the stochastic process. 

We will assume that we have a set of $\Omega \in  \mathbb{N}^+$ trajectories, indexed by $\omega \in [\Omega]$. We denote the state of the system in trajectory $\omega$ at time $t$ by $\xb(\omega, t)$. Thus, each trajectory $\omega$ corresponds to a sequence of system states $\xb(\omega, 1), \xb(\omega, 2), \dots, \xb(\omega, T)$, with $\xb(\omega, 1) = \bar\xb$. 

We can now define a sample-average approximation (SAA) version of the optimal stopping problem~\eqref{prob:optimal_stopping_canonical}. Let $\tau_{\pi, \omega}$ denote the time at which the given policy $\pi$ recommends that the system be stopped in the trajectory $\omega$; mathematically, it is defined as 
\begin{equation}
\tau_{\pi, \omega} = \min \{ t \in [T] \, | \, \pi( t, \xb(\omega, t)) = \stop \},
\end{equation}
where we again take the minimum to be $+\infty$ if the set is empty. Our SAA optimal stopping problem can now be written as 
\begin{equation}
\underset{ \pi \in \Pi}{ \text{maximize} } \ \frac{1}{\Omega}  \sum_{\omega = 1}^\Omega \beta^{ \tau_{\pi, \omega} - 1} \cdot g(  \tau_{\pi, \omega},  \xb( \omega,  \tau_{\pi, \omega} ) ).   \label{prob:optimal_stopping_SAA}
\end{equation}
Note that, in order to handle the case when $\pi$ does not stop on a given trajectory, we define $\beta^{\tau_{\pi,\omega} - 1}$ to be 0 if $\tau_{\pi,\omega} = + \infty$.
Lastly, we introduce the following short-hand notation for the sample average value of a policy $\pi$:
\begin{equation*}
\hat J^\pi(\bar\xb) \defeq \frac{1}{\Omega} \sum_{\omega=1}^{\Omega} \beta^{ \tau_{\pi, \omega} - 1} \cdot g(  \tau_{\pi, \omega},  \xb( \omega,  \tau_{\pi, \omega} ) ).
\end{equation*}

\subsection{Tree policies}
\label{subsec:problem_trees}

Problem~\eqref{prob:optimal_stopping_SAA} defines an approximation of the original problem~\eqref{prob:optimal_stopping_canonical} that uses data -- specifically, a finite sample of trajectories of the stochastic process $\{ \xb(t) \}_{t=1}^{T}$. However, despite this simplification that brings the problem closer to being solvable in practice, problem~\eqref{prob:optimal_stopping_SAA} is still difficult because it is an optimization problem over the set of all possible stopping policies $\Pi$. Moreover, as discussed in Section~\ref{sec:introduction}, we wish to restrict ourselves to policies that are sufficiently simple and interpretable. 

In this section, we will define the class of \emph{tree policies}. A tree policy is specified by a binary tree that corresponds to a recursive partitioning of the state space $\Xcal$. Each tree consists of two types of nodes: split nodes and leaf nodes. Each split node is associated with a query of the form $x_{i} \leq \theta$; we call the state variable $x_i$ that participates in the query the \emph{split variable}, the index $i$ the \textit{split variable index} and the constant value $\theta$ in the inequality the \emph{split point}. If the query is true, we proceed to the left child of the current split node; otherwise, if it is false, we proceed to the right child. %

We let $\Ncal$ denote the set of all nodes (splits and leaves) in the tree. We use $\splits$ to denote the set of split nodes and $\leaves$ to denote the set of leaf nodes. We define the functions $\leftchild: \splits \to \Ncal$ and $\rightchild: \splits \to \Ncal$ to indicate the left and right child nodes of each split node, i.e., for a given split node $s$, $\leftchild(s)$ is its left child and $\rightchild(s)$ is its right child. We use $\Tcal$ to denote the \emph{topology} of the tree, which we define as the tuple $\Tcal = (\Ncal, \leaves, \splits, \leftchild, \rightchild)$. 

Given the topology $\Tcal$, we use $\vb = \{ v(s) \}_{s \in \splits}$ and $\thetab = \{ \theta(s) \}_{s \in \splits}$ to denote the collection of all split variable indices and split points, respectively, where $v(s)$ is the split variable index and $\theta(s)$ is the split point of split $s$. 
We let $\ab = \{ a(\ell) \}_{\ell \in \leaves}$ denote the collection of leaf actions, where $a(\ell) \in \Acal$ is the action we take if the current state is mapped to leaf $\ell$. A complete tree is therefore specified by the tuple $(\Tcal, \vb, \thetab, \ab)$, which specifies the tree topology, the split variable indices, the split points and the leaf actions. 

Given a complete tree $(\Tcal, \vb, \thetab, \ab)$, we let $\ell( \xb ; \Tcal, \vb, \thetab)$ denote the leaf in $\leaves$ that the system state $\xb \in \Xcal$ is mapped to, and define the stopping policy $\pi( \cdot ; \Tcal, \vb, \thetab, \ab)$ by taking the action of the leaf to which $\xb$ is mapped:
\begin{equation*}
\pi(\xb; \Tcal, \vb, \thetab, \ab) = a( \, \ell( \xb ; \Tcal, \vb, \thetab) \, ). 
\end{equation*}

We provide an example of a tree policy below.

\begin{example}
\label{example:tree_policy}

Consider a system where $\Xcal = \mathbb{R}^3$, for which the policy is the tree given in Figure~\ref{figure:example_tree_policy}. In this example, suppose that we number the nodes with the numbers 1 through 7, from top to bottom, left to right. Then, $\Ncal = \{1,\dots, 6\}$, $\leaves = \{4,5,6,7\}$ and $\splits = \{1,2,3\}$. The $\leftchild$ and $\rightchild$ mappings are:
\begin{align*}
& \leftchild(1) = 2, & & \rightchild(1) = 3, \\
& \leftchild(2) = 4, & & \rightchild(2) = 5, \\
& \leftchild(3) = 6, & & \rightchild(3) = 7.
\end{align*}
The topology with the split and leaf labels is visualized in Figure~\ref{figure:example_tree_topology}. 

The split variable indices and split points for the split nodes $\{1,2,3\}$ are
\begin{align*}
& v(1) = 3, & & \theta(1) = 2.5, \\
& v(2) = 1, & & \theta(2) = 0.9, \\
& v(3) = 2, & & \theta(3) = 1.5;
\end{align*}
and the leaf actions for the leaf nodes $\{4,5,6,7\}$ are
\begin{align*}
& a(4) = \go, & &  a(6) = \go,\\
& a(5) = \stop, & & a(7) = \stop. 
\end{align*}
As an example, suppose that the current state of the system is $\xb = (1.2, 0.8, 2.2)$. To map this observation to an action, we start at the root and check the first query, $x_3 \leq 2.5$. Since $x_3 = 2.2$, the query is true, and we proceed to the left child of the root node. This new node is again a split, so we check its query, $x_1 \leq 0.9$. Since $x_1 = 1.2$, this query is false, so we proceed to its right child node, which is a leaf. The action of this leaf is $\stop$, and thus our policy stops the system. \Halmos

\begin{figure}
\centering
\begin{subfigure}[t]{0.4\textwidth}
                \centering
                \includegraphics[width = 0.65\textwidth]{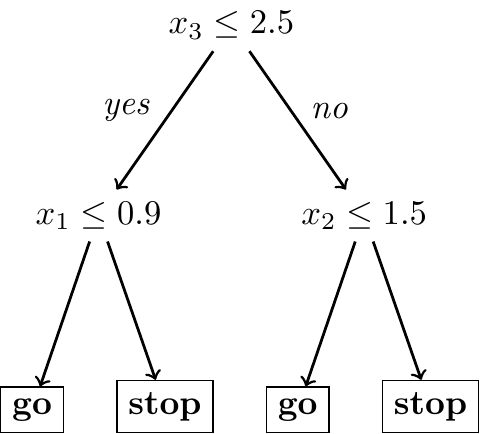}
                \caption{Tree policy.}
                \label{figure:example_tree_policy}
        \end{subfigure}  \qquad \qquad
\begin{subfigure}[t]{0.4\textwidth}
                \centering
                \includegraphics[width = 0.65\textwidth]{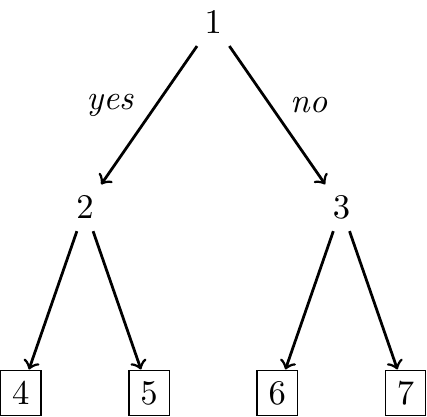}
                \caption{Node indices for tree topology. (Numbers indicate index of node in $\Ncal$.)}
                \label{figure:example_tree_topology}
        \end{subfigure}
\caption{Visualization of tree policy in Example~\ref{example:tree_policy}, for $\Xcal = \mathbb{R}^3$. \label{figure:example_tree}}
\end{figure}
\end{example}

Letting $\Pi_{\text{tree}}$ denote the set of all tree policies specified as above, we wish to find the tree policy that optimizes the sample-average reward:
\begin{equation}
\underset{\pi \in \Pi_{\text{tree}}}{\text{maximize}} \  \frac{1}{\Omega}  \sum_{\omega = 1}^\Omega \beta^{ \tau_{\pi, \omega} - 1} \cdot g(  \tau_{\pi, \omega},  \xb( \omega,  \tau_{\pi, \omega} ) ). 
\label{prob:optimal_stopping_SAA_treepolicy}
\end{equation}
We refer to this problem as the \emph{tree policy SAA problem}. In addition, we also define the counterpart of problem~\eqref{prob:optimal_stopping_canonical} restricted to policies in $\Pi_{\text{tree}}$, which we refer to simply as the \emph{tree policy problem}:
\begin{equation}
\underset{\pi \in \Pi_{\text{tree}}}{\text{maximize}} \  \Exp \left[  \beta^{ \tau_\pi - 1} \cdot g( \tau_{\pi}, \xb( \tau_\pi ) ) \ | \ \xb(1) = \xb \right],\label{prob:optimal_stopping_canonical_treepolicy}
\end{equation}
where $\tau_{\pi}$ is defined as in Section~\ref{subsec:problem_model}. 

We comment on two aspects of this modeling approach. First, the class of tree policies, as defined above, are stationary: the policy's behavior does not change with the period $t$. This turns out to not be a limitation, because it is always possible to augment the state space with an additional state variable to represent the current period $t$. By then allowing the tree to split on $t$, it is possible to obtain a time-dependent policy. We will follow this approach in our numerical experiments with option pricing in Section~\ref{sec:results} and the stylized one-dimensional problem in Section~\ref{sec:results_1D_Uniform}. 

Second, our approach requires access to trajectories that are \emph{complete}, that is, the trajectories are not terminated early/censored by the application of some prior policy, and the system state is known at every $t \in [T]$. In the case of censored trajectories, our methodology can potentially be applied by first fitting a stochastic model to the available trajectories and then simulating the model to fill in the missing data. The extension of our methodology to censored trajectories is beyond the scope of this paper, and left to future research.

\subsection{Convergence of tree policy SAA problem}
\label{subsec:problem_convergence}

A natural expectation for the SAA problem~\eqref{prob:optimal_stopping_SAA_treepolicy} is that its objective approaches the optimal value of \eqref{prob:optimal_stopping_canonical_treepolicy}, as the number of trajectories $\Omega$ available as samples increases. In this section, we show that this is indeed a property of our SAA problem. More precisely, we show that in the limit of $\Omega \to \infty$ and if we restrict the policy class $\Pi_{\text{tree}}$ to only trees of arbitrary \textit{bounded} depth, the optimal value of the SAA problem converges to the value of the optimal tree-based policy for the true problem almost surely.

We note that establishing such a convergence result in the case of tree-based policies is challenging due to two difficulties. First, a tree policy is defined by a tuple $(\mathcal{T}, \vb, \thetab, \ab)$; the space of such tuples could in general be uncountable, and thus we cannot directly invoke the strong law of large numbers to guarantee that almost sure convergence holds simultaneously over all $\pi\in \Pi_{\rm tree}$. Secondly, invoking the strong law of large numbers over all members of a finite cover of this space is also non-trivial. This is due to the fact that $\hat J^\pi(\bar\xb)$ is not necessarily continuous in $\pi$, even as we change the split points $\thetab$ and keep all other tree parameters constant, and as such we cannot rely on Lipschitz continuity arguments.

In order to handle these challenges, we restrict our convergence analysis to the class of tree policies corresponding to trees of depth at most $d$ for some finite parameter $d$. We denote this class of tree policies by $\Pi_{\rm tree}(d) \subseteq \Pi_{\rm tree}$. This restriction limits the set of $(\mathcal{T}, \vb, \ab)$ parameters specifying a tree-based policy to be finite, although the $\thetab$ parameters still lie in a potentially uncountable set. We remark that, because our focus is on interpretable policies, limiting the potential depth of allowable trees that our SAA approach can optimize over is reasonable. 

We also make some relatively mild assumptions regarding the underlying structure of the optimal stopping problem, which facilitate the subsequent analysis. The first two enforce the total boundedness of the state space $\mathcal{X}$ and normalize the absolute magnitude of state variables (this can be done without loss of generality), as well as place a universal upper bound on the magnitude of the cost function $g(\cdot, \cdot)$:
\begin{assumption}
\label{ass:tb}
The state space $\mathcal{X}$ is totally bounded with respect to the $||\cdot||_\infty$ norm and $\max_{\xb, \mathbf{y} \in \Xcal} ||\xb - \mathbf{y}||_\infty \leq 1$. 
\end{assumption}
\begin{assumption}
\label{ass:bounded_cost}
There exists a constant $G$ such that for any $t \in [T]$, $\xb \in \mathcal{X}$, $0 \leq g(t, \xb) \leq G$. \end{assumption}
With regard to Assumption~\ref{ass:tb}, we note that in many optimal stopping problems, the state space may not be bounded. For example, in an option pricing problem involving assets whose prices follow a geometric Brownian motion model, the asset prices at each period may be arbitrarily large. One could potentially circumvent the requirement of a bounded state space by making an additional assumption that the state variable at each period is contained in a bounded set with high probability. We proceed with Assumption~\ref{ass:tb} in order to maintain the simplicity of the analysis that follows.

The third assumption essentially imposes that the distribution of the state variable is sufficiently smooth:
\begin{assumption}
\label{ass:state_smoothness}
 There exists a function $f : [0, \infty) \rightarrow [0,1]$ such that,
\begin{enumerate}
\item Given the Borel measure $\mu$ on $\mathbb{R}^n$, for any subset $A \in \mathbb{R}^n$ and for any $1 < t  \leq T$,
\begin{equation*}\Pr\left[\xb(t) \in A \right] \leq f(\mu(A)).\end{equation*}
\item $f(0)=0$.
\item $f$ is strictly increasing and continuous.
\end{enumerate}
\end{assumption}
To gain some intuition regarding this assumption, it is helpful to interpret it in the case that the state is univariate. Then, the assumption reduces to whether the probability that the state variable $x \in [a,b]$ is bounded by a function which only depends on the width $|b-a|$, and which goes to zero as this width vanishes. This is satisfied, for example, by distributions with probability density functions that can be uniformly upper bounded across their entire domain.

Additionally, we note that, while this assumption seems to potentially rule out a convergence result in the case of categorical state variables, this is in fact not the case. As alluded to in the above, the main challenge in establishing our convergence result is the uncountable number of possible split points; to deal with this uncountability issue, we rely on Assumption \ref{ass:state_smoothness}. If some of the state variables were categorical, we only need to consider a finite number of split points for splits on these variables and only require Assumption \ref{ass:state_smoothness} on the remaining subset of continuous state variables. For simplicity however, we omit this generalization here.

Having set up these assumptions, the following theorem proves almost sure convergence of $\hat J^\pi(\bar\xb)$ to $J^\pi(\bar\xb)$ over all policies $\pi$ in $\Pi_{\rm tree}(d)$. 

\begin{theorem}
\label{thm:pointwise_convergence}
Fix any finite tree depth $d$, initial state $\bar\xb$ and arbitrary $\epsilon > 0$. Then, with probability one, there exists a finite sample size $\Omega_0$ such that for all $\Omega \geq \Omega_0$ and all tree policies $\pi \in \Pi_{\rm tree}(d)$,
\begin{equation*}
\left|
J^{\pi}(\bar\xb) - \hat J^\pi(\bar\xb)
\right|
\leq \epsilon.
\end{equation*}
\end{theorem}

We make a few technical remarks about the above theorem:
\begin{enumerate}[(a)]
\item The proof of this result is driven by the fact that the sample average problem, for a fixed $\pi \in \Pi_{\rm tree}(d)$, converges to its true expectation almost surely by the strong law of large numbers. The rate of convergence of the SAA objective is  of the order $O\left(\sqrt{\log \log \Omega / \Omega}\right)$ by the law of the iterated logarithm. 
\item The restriction to constant depth $d$ tree policies is not an artifact of the proof. In fact, in Section \ref{sec:convergence_can_fail}, we show that convergence fails if we allow for any tree policy in $ \Pi_{\rm tree}$.
\end{enumerate}

Theorem \ref{thm:pointwise_convergence} is the main technical result of this section and enables us to prove our desired convergence of the SAA optimal objective. This is stated in the following corollary:
\begin{corollary}
\label{cor:convergence}
Fix any finite tree depth $d$, initial state $\bar\xb$ and arbitrary $\epsilon > 0$. Then, with probability one, there exists a finite sample size $\Omega_0$ such that for all $\Omega \geq \Omega_0$,
\begin{equation*}
\left|
\sup_{\pi \in \Pi_{\rm tree}(d)} J^{\pi}(\bar\xb) - \sup_{\pi \in \Pi_{\rm tree}(d)} \hat J^\pi(\bar\xb)
\right|
\leq \epsilon.
\end{equation*}
\end{corollary}
The corollary above shows that the decision maker indeed obtains, via solving the sample average problem \eqref{prob:optimal_stopping_SAA_treepolicy}, a policy whose value approximates with arbitrarily high precision the value of the best tree-based policy for the true problem, as long as the decision maker has access to a sufficiently large sample of state trajectories. 

As the proofs of Theorem~\ref{thm:pointwise_convergence} and Corollary~\ref{cor:convergence} are quite involved, we present them in Section \ref{sec:convergence_appendix}.

\subsection{Complexity of tree policy SAA problem}
\label{subsec:problem_complexity}

Having established that the tree policy SAA problem converges to the exact tree policy optimization problem, we now turn our attention to the computational complexity of solving the tree policy SAA problem. In this section, we will consider three simplified versions of problem~\eqref{prob:optimal_stopping_SAA_treepolicy} and show that each of these is theoretically intractable.

To motivate the first of our intractability results, observe that problem~\eqref{prob:optimal_stopping_SAA_treepolicy} allows for the tree policy to be optimized along all four dimensions: the topology $\Tcal$, the split variable indices $\vb$, the split points $\thetab$ and the leaf actions $\ab$. Rather than optimizing over all four variables, let us consider instead a simplified version of problem~\eqref{prob:optimal_stopping_SAA_treepolicy}, where the topology $\Tcal$, split variable indices $\vb$ and split points $\thetab$ are fixed, and the leaf actions $\ab$ are the only decision variables. We use $\Pi( \Tcal, \vb, \thetab)$ to denote the set of all policies with the given $\Tcal$, $\vb$ and $\thetab$, i.e., 
\begin{equation}
\Pi( \Tcal, \vb, \thetab ) = \{ \pi(\cdot; \Tcal, \vb, \thetab, \ab) \, | \, \ab \in \{ \stop, \go \}^{\leaves} \}.
\end{equation}
The \emph{leaf action SAA problem} is to find the policy in $\Pi( \Tcal, \vb, \thetab)$ that optimizes the sample average reward:
\begin{equation}
\underset{\pi \in \Pi(\Tcal, \vb, \thetab) }{\text{maximize}}\  \frac{1}{\Omega}  \sum_{\omega = 1}^\Omega \beta^{ \tau_{\pi, \omega} - 1} \cdot g(  \tau_{\pi, \omega},  \xb( \omega,  \tau_{\pi, \omega} ) ).  
\label{prob:optimal_stopping_SAA_leafdecision}
\end{equation}
While problem~\eqref{prob:optimal_stopping_SAA_treepolicy} is difficult to analyze in generality, the simpler problem~\eqref{prob:optimal_stopping_SAA_leafdecision} is more amenable to analysis. It turns out, perhaps surprisingly, that this simplified problem is already hard to solve: 
\begin{proposition}
The leaf action problem~\eqref{prob:optimal_stopping_SAA_leafdecision} is NP-Hard when $|\leaves|$, $n$ and $\Omega$ are unrestricted.
\label{proposition:leaf_decision_NPHard}
\end{proposition}
This result is significant for two reasons. First, it establishes that even in this extremely simplified case, where we have already selected a tree topology, split variable indices and split points, the resulting problem is theoretically intractable.  

Second, this result points to an important distinction between optimal stopping tree policies and classification trees that are used in machine learning. For binary classification, a classification tree consists of a tree topology, split variable indices, split points and leaf labels. Given $\Tcal$, $\vb$ and $\thetab$, determining the leaf labels $\ab$ that minimize the 0-1 classification error on a training sample is trivial: for each leaf, we simply predict the class that is most frequent among those observations that are mapped to that leaf. More importantly, each leaf's label can be computed independently. The same holds true for the regression setting, where $a(\ell)$ is the continuous prediction of leaf $\ell$: to find the values of $\ab$ that minimize the squared prediction error, we can set each $a(\ell)$ as the average of the dependent variable for all training observations that are mapped to that leaf. 

For optimal stopping, the situation is strikingly different. Determining the leaf actions is much more difficult because a single trajectory has a time dimension: as time progresses, a tree policy may map the current state of that trajectory to many different leaves in the tree. As a result, the decision of whether to stop or not in one leaf (i.e., to set $a(\ell) = \stop$ for a leaf $\ell$) cannot be made independently of the decisions to stop or not in the other leaves: the leaf actions are coupled together. For example, if we set $a(\ell ) = \stop$ for a given leaf $\ell$ and $a(\ell') = \stop$ for a different leaf $\ell'$, then a trajectory that reaches $\ell$ before it reaches $\ell'$ could never stop at $\ell'$: thus, whether we choose to stop at $\ell'$ depends on whether we choose to stop at $\ell$. %

The proof of this result (see Section~\ref{proof:leaf_decision_NPHard}) follows by a reduction from the minimum vertex cover problem. The specific reduction involves an instance where $\Omega$, $n$ and $|\leaves|$ scale with the size of the vertex cover instance. Thus, a natural question is whether the problem remains intractable when any of these quantities are fixed to constant values. To answer this question, we provide the following insights:
\begin{enumerate}
\item If only the number of leaves $|\leaves|$ is fixed, the problem becomes tractable: the number of possible leaf actions $\ab$ is $2^{|\leaves|}$, so enumerating all of the values of $\ab$ and evaluating each $\ab$ on $\Omega$ trajectories of length $T$ will require $O(\Omega \cdot T \cdot 2^{|\leaves|} )$ steps, which is clearly polynomial in the parameters besides $|\leaves|$. 
\item If only the number of trajectories $\Omega$ is fixed, the problem also becomes tractable. As a specific example, let $\Omega = 1$, and define $t_{\ell}$ to be the first time at which the trajectory enters leaf $\ell$ (and define it to be $\infty$ if the trajectory never enters leaf $\ell$). We know that the trajectory can be stopped in at most one period, and that period will correspond to a leaf. The set of possible times it may stop at is therefore exactly $\{ t_{\ell} \}_{\ell \in \leaves}$; since we are optimizing for a single trajectory, the problem simplifies to selecting a time at which to stop. More precisely, in order to stop at $t_{\ell}$, one must set $a(\ell) = \stop$ and must set $a(\ell') = \go$ for all leaves $\ell'$ with $t_{\ell'} < t_{\ell}$. Therefore, while there are $2^{|\leaves|}$ possible values of $\ab$, one only needs to consider $O( |\leaves| )$ possible values of $\ab$ to find the optimal $\ab$. The same type of logic can be extended to a general fixed $\Omega$, allowing one to solve problem~\eqref{prob:optimal_stopping_SAA_leafdecision} through an enumeration scheme in $O(T^{\Omega} \cdot \Omega \cdot |\leaves|)$ time, which is polynomial in the parameters besides $\Omega$. The details of this enumeration scheme are provided in Section~\ref{appendix:leaf_action_problem_fixed_Omega}.
\item The state variable dimension $n$ does not appear to affect the intractability of the problem. The proof in Section~\ref{proof:leaf_decision_NPHard} requires that $n$ is set to the number of vertices in the vertex cover instance. It turns out that there exists an alternate reduction from the minimum vertex cover problem to problem~\eqref{prob:optimal_stopping_SAA_leafdecision} where there is only one state variable ($n = 1$). The details of this alternate reduction are provided in Section~\ref{appendix:leaf_action_problem_NPHard_neq1}.
\end{enumerate}

Proposition~\ref{proposition:leaf_decision_NPHard} establishes that when the topology, split variable indices and split points are fixed, optimizing over the leaf actions is an intractable problem. It turns out that optimizing individually over the split variable indices and the split points, with the rest of the tree policy parameters fixed, is also intractable. Let $\Pi(\Tcal, \thetab, \ab)$ denote the set of all tree policies with the given $\Tcal$, $\thetab$ and $\ab$ (i.e., the split variable indices $\vb$ may be optimized over), and let $\Pi(\Tcal, \vb, \ab)$ denote the set of all tree policies with the given $\Tcal$, $\vb$ and $\ab$ (i.e., the split points $\thetab$ may be optimized over):
\begin{align}
& \Pi( \Tcal, \thetab, \ab ) = \{ \pi(\cdot; \Tcal, \vb, \thetab, \ab) \, | \, \vb \in [n]^{\splits} \}, \\
& \Pi( \Tcal, \vb, \ab ) = \{ \pi(\cdot; \Tcal, \vb, \thetab, \ab) \, | \, \thetab \in \mathbb{R}^{\splits} \}.
\end{align}
Let us define the \emph{split variable index SAA problem} as the problem of finding a policy in $\Pi(\Tcal, \thetab, \ab)$ to maximize the sample-based reward:
\begin{equation}
\underset{\pi \in \Pi(\Tcal, \thetab, \ab) }{\text{maximize}}\  \frac{1}{\Omega}  \sum_{\omega = 1}^\Omega \beta^{ \tau_{\pi, \omega} - 1} \cdot g(  \tau_{\pi, \omega},  \xb( \omega,  \tau_{\pi, \omega} ) ). 
\label{prob:optimal_stopping_SAA_splitvariable}
\end{equation}
Similarly, let us define the \emph{split point SAA problem} as the analogous problem of optimizing over policies in $\Pi(\Tcal, \vb, \ab)$:
\begin{equation}
\underset{\pi \in \Pi(\Tcal, \vb, \ab) }{\text{maximize}}\  \frac{1}{\Omega}  \sum_{\omega = 1}^\Omega  \beta^{ \tau_{\pi, \omega} - 1} \cdot g(  \tau_{\pi, \omega},  \xb( \omega,  \tau_{\pi, \omega} ) ) .
\label{prob:optimal_stopping_SAA_splitpoint}
\end{equation}
We then have the following two intractability results. 
\begin{proposition}
The split variable index SAA problem~\eqref{prob:optimal_stopping_SAA_splitvariable} is NP-Hard. \label{proposition:split_variable_NPHard}
\end{proposition}
\begin{proposition}
The split point SAA problem~\eqref{prob:optimal_stopping_SAA_splitpoint} is NP-Hard. \label{proposition:split_point_NPHard}
\end{proposition}
Like Proposition~\ref{proposition:leaf_decision_NPHard}, the proofs of Propositions~\ref{proposition:split_variable_NPHard} and \ref{proposition:split_point_NPHard}, found in Sections~\ref{proof:split_variable_NPHard} and \ref{proof:split_point_NPHard} respectively, also follow by a reduction from the minimum vertex cover problem. With regard to Proposition~\ref{proposition:split_point_NPHard}, we remark that a special simplified version of the split point SAA problem, in which one optimizes over tree policies that are defined by a collection of thresholds $\theta_1, \dots, \theta_{T-1}$ on the reward $g(t,\cdot)$ at each period $t$ prior to the final period $t = T$, also turns out to be NP-Hard; we refer the reader to Section~\ref{appendix:threshold_opt_nphard} for more detail.

Finally, before concluding this section, we comment on the overall complexity of solving problem~\eqref{prob:optimal_stopping_SAA_treepolicy} when the topology $\Tcal$ is fixed, and the leaf actions $\ab$, split variable indices $\vb$ and split points $\thetab$ are allowed to vary. In this case, enumerating all possible combinations of $(\ab, \vb, \thetab)$ will have a complexity of $O( 2^{|\leaves|} \cdot n^{|\splits|} \cdot (\Omega T + 1)^{|\splits|})$. The first two factors in the complexity bound correspond to the number of possible choices of $\ab$ and $\vb$, respectively, while the last term corresponds to the number of possible choices of $\thetab$. (Note that $\Omega T$ is the maximum number of unique values a state variable may take in the sample, across all trajectories; since all split points between consecutive unique values of a state variable will result in the same stopping behavior, one only needs to look at how many intervals of equivalent split points there are in the data, which is at most $\Omega T + 1$.)

\subsection{Approximate optimality of tree policies}
\label{subsec:tree_optimality}
We end this section by providing additional justification for our choice of policy class. Specifically, we prove that if we assume that the decision boundary of an optimal stopping policy is well-behaved, then this policy can be approximated to arbitrary fidelity by a tree policy of finite, albeit possibly large, depth. This result shows that trees are in general sufficiently expressive to represent optimal policies. %

For the results in this section to hold, we make some additional assumptions:%
\begin{assumption}
\label{ass:compactness}
\begin{enumerate}
\item We augment the state $\xb \in \Xcal$ with a $0$-th additional component $x_0(t) = t/T$ for all times $t\in [T]$.
\item There exists an optimal stopping policy $\pi^*$ solving \eqref{prob:optimal_stopping_canonical} such that, for every $t \in [T] \setminus \{1\}$, the set $\Xcal_t^\stop \defeq \{\xb \in \Xcal \mid \pi^*(t, \xb) = \stop\}$ is compact.
\end{enumerate}
\end{assumption}
We make a few clarifying comments about the need for Assumption \ref{ass:compactness}. The addition of a variable encoding the time into the state space is not surprising: note that the optimal policy $\pi^*$ can be time-dependent; thus, it is natural that a tree policy that approximates it well should also be allowed to split on the time period, which is why we augment the state space in this way.
The second assumption can be thought of as imposing some regularity on the decision boundary of $\pi^*$, which conveniently allows us to approximate this decision boundary with a finite collection of boxes in $\mathbb{R}^n$, and which is then amenable to a representation via trees. 

We have the following theorem which is the main result of this section:
\begin{theorem}
\label{thm:tree_policy_appr_opt}
If Assumption \ref{ass:compactness} holds, then for any $\epsilon > 0$ and starting state $\bar\xb$ there exists a tree policy $\pi \in \Pi_{\rm tree}$ of finite depth such that $J^*(\bar\xb) - J^\pi(\bar\xb) \leq \epsilon$.
\end{theorem}
The proof of this theorem, given in Section~\ref{sec:tree_optimality_appendix}, involves using Assumption~\ref{ass:compactness} to cover the stopping region $\Xcal_t^{\stop}$ by a finite cover of boxes, and then representing this ``approximate'' stopping region through a tree. The main takeaway from this theorem is that there exist tree policies with performance arbitrarily close to optimal. In Sections~\ref{sec:results} and \ref{sec:results_1D_Uniform}, we will see numerically that tree policies can obtain near-optimal performance.

Although Theorem~\ref{thm:tree_policy_appr_opt} provides a guarantee that a tree policy of arbitrary optimality will exist under mild conditions, it does not provide any insight on the depth of such a tree. In Section~\ref{subsec:epsilon_depth_bound_R2} of the electronic companion, we provide an alternate analysis under slightly stronger assumptions, which guarantees the existence of an $\epsilon$-optimal tree policy whose depth scales like $O(T + n \log n + n \log(1 / \epsilon))$ (Theorem~\ref{theorem:ultimate_depth_bound}). This alternate construction relies on discretizing the state space into a collection of hypercubes, and requires the assumption that the reward function $g$ and the optimal continuation value function are Lipschitz continuous. We remark that the analysis in Section~\ref{subsec:epsilon_depth_bound_R2} provides only one approach to obtaining approximately optimal trees with an accompanying depth guarantee. The questions of whether there exist tree policies that have the same degree of suboptimality but with a lower depth under the same assumptions, or whether there are alternate assumptions that yield trees of lower depth, are interesting directions for future research.

\section{Construction algorithm}
\label{sec:construction}
As we saw in Section~\ref{subsec:problem_trees}, three major simplifications of our SAA problem are theoretically intractable. Expanding the scope of the problem to allow for joint optimization over the split points, split variable indices and leaf actions, as well as the topology, renders the problem even more difficult. Given the difficulty of this problem, we now present a practical heuristic algorithm for approximately solving the SAA problem. Our algorithm is a greedy procedure that grows/induces the tree from the top down. In Section~\ref{subsec:construction_description}, we provide a description of the overall algorithm. In Section~\ref{subsec:construction_OptimizeSplitPoint}, we provide a procedure for performing the key step of our construction algorithm, which is finding the optimal split point. In Section~\ref{subsec:construction_complexity}, we comment on the complexity of the overall algorithm. Finally, in Section~\ref{subsec:construction_comparison_to_classification}, we compare our construction algorithm and extant classification tree approaches such as CART \citep{breiman1984classification}. We emphasize that the algorithm presented in this section is generic and does not require any of the assumptions made in Sections \ref{subsec:problem_convergence} and \ref{subsec:tree_optimality}.

\subsection{Algorithm description}
\label{subsec:construction_description}

Our algorithm to heuristically solve \eqref{prob:optimal_stopping_SAA_treepolicy} is a construction procedure where we greedily grow a tree up to the point where we no longer observe an improvement by adding another split point. %
At a high level, our algorithm works as follows. We start from a degenerate tree, consisting of a single leaf with the action of $\go$ as the root node. At this leaf, we consider placing a split. Such a split will result in a left child and a right child, with both child nodes being leaves. For each potential split variable, we find the best possible split point assuming that the left child leaf will be $\stop$ (and the right child will be $\go$), and assuming that the right child leaf will be $\stop$ (and the left child will be $\go$). We find the best possible combination of the split variable, split point and direction (left child is $\stop$ and right child is $\go$, or right child is $\stop$ and left child is $\go$), and add the split to the tree. %
We continue the procedure if the split resulted in a sufficient improvement on the current objective; otherwise, we terminate the procedure. 

\begin{figure}
	\begin{subfigure}[t]{0.2\textwidth}
                \centering
                \includegraphics[width = 1.2\textwidth]{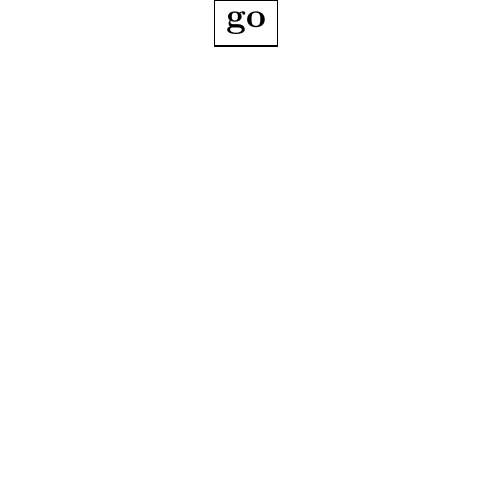} %
                \caption{Iteration 0.}
                \label{fig:gull}
        \end{subfigure} \hfill
        \begin{subfigure}[t]{0.2\textwidth}
                \centering
                \includegraphics[width = 1.2\textwidth]{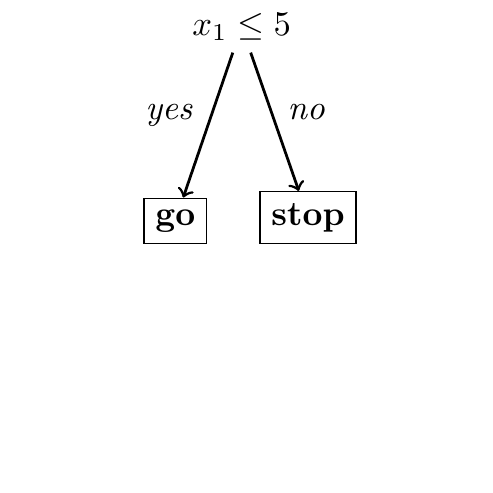} %
                \caption{Iteration 1.}
                \label{fig:gull2}
        \end{subfigure} \hfill 
         \begin{subfigure}[t]{0.2\textwidth}
                \centering
                \includegraphics[width = 1.2\textwidth]{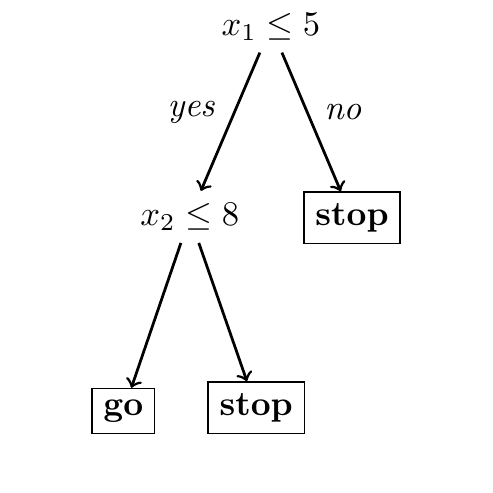} %
                \caption{Iteration 2.}
                \label{fig:gull2}
        \end{subfigure} \hfill
        \begin{subfigure}[t]{0.2\textwidth}
                \centering
                \includegraphics[width = 1.2\textwidth]{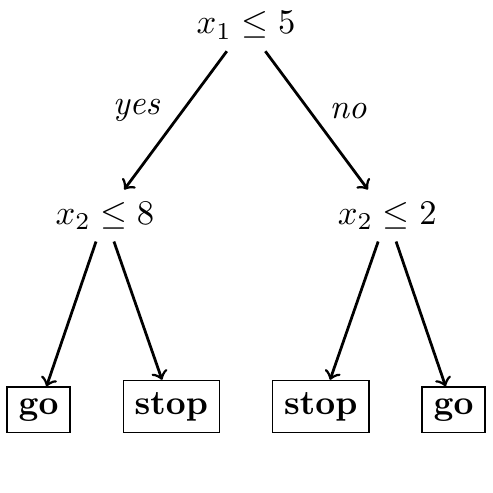} %
                \caption{Iteration 3.}
                \label{fig:gull}
        \end{subfigure} \hfill
\caption{Example of evolution of tree with each iteration of construction procedure (Algorithm~\ref{algorithm:construction}). \label{figure:example_construction_iter} }
\vspace{-1em}
\end{figure}

In the next iteration, we repeat the process, except that now we also optimize over the leaf that we select to split: we compute the best possible split at each leaf in the current tree, and take the best split at the best leaf. We continue in this way until there is no longer sufficient improvement in an iteration. At each iteration, we expand the tree at a single leaf and add two new leaves to the tree. Figure~\ref{figure:example_construction_iter} provides a visualization of several iterations of the algorithm.

The notion of sufficient improvement that we use in our algorithm is that of relative improvement. Specifically, if $Z'$ is the objective value with the best split and $Z$ is the current objective, then we continue running the algorithm if $Z' \geq (1 + \gamma)Z$, where $\gamma \geq 0 $ is a user-specified tolerance on the relative improvement; otherwise, if the relative improvement is lower than $\gamma$, the algorithm is terminated. Lower values of $\gamma$ correspond to trees that are deeper with a larger number of nodes. %

\begin{figure}
\begin{center}
 \begin{subfigure}[t]{0.4\textwidth}
                \centering
                \includegraphics{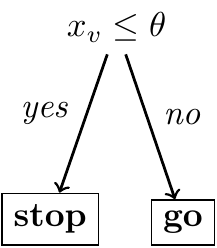}
                \caption{Left-stop subtree.}
                \label{fig:gull}
        \end{subfigure} \quad
        \begin{subfigure}[t]{0.4\textwidth}
                \centering
                \includegraphics{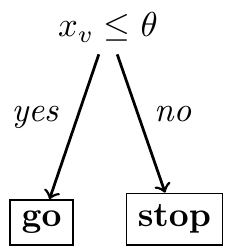}
                \caption{Right-stop subtree.}
        \end{subfigure}
\end{center}
\caption{Visualization of the two different types of subtrees that can be used to split a leaf in the construction algorithm. \label{figure:subtrees}}
\end{figure}

Our construction algorithm is defined formally as Algorithm~\ref{algorithm:construction}. The tree policy is initialized to the degenerate policy described above, which prescribes the action $\go$ at all states. Each iteration of the loop first computes the objective attained from choosing the best split point at each leaf $\ell$ with each variable $v$ assuming that either the left child leaf will be a $\stop$ action or the right child leaf will be a $\stop$ action. In the former case, we call the subtree rooted at node $\ell$ a left-stop subtree, and in the latter case, we call it a right-stop subtree; Figure~\ref{figure:subtrees} shows both of these subtrees. Then, if the best such split achieves an objective greater than the current objective value, we grow the tree in accordance to that split: we add two child nodes to the leaf $\ell^*$ using the \textsc{GrowTree} function, which is defined formally as Algorithm~\ref{algorithm:grow}; we set the actions of those leaves according to $\dir^*$; we set the split variable index of the new split as $v^*$; and finally, we set the split point as $\theta^*_{\ell^*, v^*, \dir^*}$. The \textbf{existsImprovement} flag is used to terminate the algorithm when it is no longer possible to improve on the objective of the current tree by at least a factor of $\gamma$. Section~\ref{sec:toy_example} of the electronic companion presents a small example illustrating a couple of iterations of the algorithm.

\begin{algorithm}
{\SingleSpacedXI
\begin{algorithmic}
\REQUIRE User-specified parameter $\gamma$ (relative improvement tolerance)
\STATE Initialization:
\STATE \quad $\Tcal \leftarrow ( \{1\}, \{1\}, \emptyset, \leftchild, \rightchild )$, $\vb \leftarrow \emptyset$, $\thetab \leftarrow \emptyset$, $a(1) \leftarrow \go$
\STATE \quad $Z \leftarrow 0$
\STATE \quad $\textbf{existsImprovement} \leftarrow \textbf{true}$ \\[1em]
\WHILE{ \textbf{existsImprovement} }
	\FOR{$\ell \in \leaves$, $v \in [n]$, $\dir \in \{ \text{left}, \text{right} \}$ }
				\STATE $Z^*_{\ell,v, \dir }, \theta^*_{\ell,v, \dir} \leftarrow \textsc{OptimizeSplitPoint}(\ell, v, \dir; \Tcal, \vb, \thetab, \ab)$
	\ENDFOR\\[1em]
	$\textbf{existsImprovement} \leftarrow \Ibb \left[\max_{\ell, v, \dir} Z^*_{\ell,v, \dir} \geq (1 + \gamma)  Z \right]$\\[1em]
	\IF{ $\max_{\ell, v, \dir} Z^*_{\ell,v, \dir} >  Z$}
		\STATE ($\ell^*, v^*, \dir^*) \leftarrow \arg \max_{\ell, v, D} Z^*_{\ell,v, \dir }$
		\STATE $\textsc{GrowTree}(\Tcal, \ell^*)$
		\STATE $v(\ell^*) \leftarrow v^*$
		\STATE $\theta(\ell^*) \leftarrow \theta^*_{\ell^*, v^*, \dir^*}$ \\[0.5em]
		\IF{ $\dir^* = \text{left}$ }
		 	\STATE $a( \leftchild(\ell^*) ) \leftarrow \stop$,  $a( \rightchild(\ell^*) ) \leftarrow \go$
		\ELSE
		 	\STATE $a( \leftchild(\ell^*) ) \leftarrow \go$, $a( \rightchild(\ell^*) ) \leftarrow \stop$
		\ENDIF\\
		\STATE $Z \leftarrow Z^*_{\ell^*, v^*, \dir^*}$ 
	\ENDIF
\ENDWHILE
\RETURN Tree policy $(\Tcal, \vb, \thetab, \ab)$.
\end{algorithmic}
}
\caption{Tree construction algorithm. \label{algorithm:construction}}
\end{algorithm}

\begin{algorithm}
{\SingleSpacedXI
\begin{algorithmic}
\REQUIRE Tree topology $\Tcal = (\Ncal, \leaves, \splits, \leftchild, \rightchild)$; target leaf $\ell$. 
\STATE $\Ncal \leftarrow \Ncal \cup \{ |\Ncal| + 1, |\Ncal| + 2 \}$
\STATE $\leaves \leftarrow (\leaves \setminus \{ \ell \} ) \cup \{ |\Ncal| + 1, |\Ncal| + 2 \}$
\STATE $\splits \leftarrow \splits \cup \{ \ell \}$
\STATE $\leftchild(\ell) = |\Ncal| + 1$
\STATE $\rightchild(\ell) = |\Ncal| + 2$
\end{algorithmic}
}
\caption{\textsc{GrowTree} function. \label{algorithm:grow}}
\end{algorithm}

We comment on three important aspects of Algorithm~\ref{algorithm:construction}. The first aspect is the determination of the optimal split point for a given leaf, a given split variable index and the direction of the stop action (left/right). At present, this optimization is encapsulated in the function \textsc{OptimizeSplitPoint} to aid in the exposition of the overall algorithm; we defer the description of this procedure to the next section (Section~\ref{subsec:construction_OptimizeSplitPoint}).

The second aspect is the assumption of how the child actions are set when we split a leaf. At present, we consider two possibilities: the left child is $\stop$ and the right child is $\go$ (this corresponds to $\dir^* = \text{left}$) or the left child is $\go$ and the right child is $\stop$ (this corresponds to $\dir^* = \text{right}$). These are not the only two possibilities, as we can also consider setting both child nodes to $\go$ or to $\stop$. We do not explicitly consider these. First, \textsc{OptimizeSplitPoint} can potentially return a value of $-\infty$ or $+\infty$, effectively resulting in a degenerate split where we always go to the left or to the right. It is straightforward to see that such a split is equivalent to setting both child actions to $\go$ or $\stop$. Second, one of these two possibilities -- $\stop$-$\stop$ or $\go$-$\go$ -- will result in exactly the same behavior as the current tree (for example, if $a(\ell^*) = \go$ in the current tree and we set $a(\leftchild(\ell^*)) = a(\rightchild(\ell^*)) = \go$ in the new tree, the new tree policy will prescribe the same actions as the current tree), and thus cannot result in an improvement. In addition, the other possibility will also have been considered in a previous step (for example, if $a(\ell^*) = \go$ in the current tree, and we considered setting $a(\leftchild(\ell^*)) = a(\rightchild(\ell^*)) = \stop$, this would be the same as simply setting $a(\ell^*) = \stop$ in the current tree, which by the aforementioned point about setting the split point to $-\infty$ or $+\infty$ would have been previously considered).

The third aspect is complexity control. The algorithm in its current form stops building the tree when it is no longer possible to improve the objective value by a factor of at least $\gamma$. The parameter $\gamma$ plays an important role in controlling the complexity of the tree: if $\gamma$ is set to a very low value, the algorithm may produce deep trees with a large number of nodes, which is undesirable because (1) the resulting tree may lack interpretability and (2) the resulting tree may not generalize well to new data. In contrast, higher values of $\gamma$ will cause the algorithm to terminate earlier with smaller trees whose training set performance will be closer to their out-of-sample performance; however, these trees may not be sufficiently complex to lead to good performance. %
We note that this relative improvement criterion is not the only way of terminating the algorithm/controlling complexity. Another possibility is to, for example, terminate when the tree reaches a specific depth or a specific limit on the number of nodes. Our use of the relative improvement criterion and the parameter $\gamma$ is partially inspired by cost complexity pruning in the original CART algorithm, where one effectively terminates the induction procedure when the change in classification accuracy is smaller than a user-specified cost complexity parameter $\alpha$ \citep{breiman1984classification}. The difference between the complexity parameter $\alpha$ in CART and our parameter $\gamma$ is that $\alpha$ controls the \emph{absolute} improvement in accuracy, whereas our parameter $\gamma$ controls the \emph{relative} improvement in reward. 

With regard to how one may tune $\gamma$, in Section \ref{sec:kfold_cross_validation} we present an algorithm to calibrate this parameter based on $k$-fold cross-validation. In that section, we observe that the naive way of doing cross-validation, which is running a grid search for $\gamma$, may overlook some good values of this parameter; in addition, this grid search is computationally wasteful since for a smaller $\gamma$ value one has to repeat the iterations that had already been computed for a larger $\gamma$ value. Instead of this approach, we propose a tailored algorithm that is significantly more efficient by  integrating the search for $\gamma$ into the construction procedure.

\subsection{Finding the optimal split point}
\label{subsec:construction_OptimizeSplitPoint}

The key part of Algorithm~\ref{algorithm:construction} is the \textsc{OptimizeSplitPoint} function, which aims to answer the following question: what split point $\theta$ should we choose for the split at node $\ell$ on variable $v$ so as to maximize the sample-based reward? At first glance, this question appears challenging because we could choose any real number to be $\theta$, and it is not clear how we can easily optimize over $\theta$. Fortunately, it will turn out that the sample average reward of the new tree will be a piecewise-constant function of $\theta$, which we will be able to optimize over easily. The piecewise-constant function is obtained by taking a weighted combination of a collection of trajectory-specific piecewise-constant functions. These trajectory-specific piecewise-constant functions can be computed by carefully analyzing how the split variable $x_v$ changes along each trajectory. 

Assume that the split variable index $v$, the leaf $\ell$ and the subtree direction $\dir$ are fixed. Our first step is to determine the behavior of the trajectories with respect to the leaf $\ell$. For each trajectory, we first find the time at which each trajectory is stopped at a leaf different from $\ell$. We call this the \emph{no-stop time}. We denote it with $\tau_{-\ell, \omega}$ and define it as 
\begin{equation}
\tau_{-\ell, \omega} =  \min \{ t \in [T] \ | \ \ell( \xb(\omega, t) ) \neq \ell \ \text{and} \ a ( \ell(\xb(\omega,t)) = \stop \}, \label{eq:nostop_time}
\end{equation}
where we define the minimum to be $+\infty$ if the set is empty. We can effectively think of this as the time at which trajectory $\omega$ will stop if the new split that we place at $\ell$ does not result in the trajectory being stopped or equivalently, if we were to set the action of leaf $\ell$ to $\go$. We define the \emph{no-stop value}, $f_{\omega,\text{ns}}$, as the reward we garner if we allow the trajectory to stop at the no-stop time: 
\begin{equation}
f_{\omega, \text{ns}} = \left\{ \begin{array}{ll} g(\tau_{-\ell, \omega}, \xb(\omega, \tau_{-\ell, \omega} ) ) & \text{if}\ \tau_{-\ell, \omega} < +\infty, \\
0 & \text{otherwise}. \end{array} \right. 
\label{eq:nostop_value}
\end{equation}
Note that the value is zero if the trajectory is never stopped outside of leaf $\ell$. Lastly, we also determine the set of periods $S_{\omega}$ when the policy maps the system state to the leaf $\ell$. We call these the \emph{in-leaf periods}, and define the set as
\begin{equation}
S_{\omega} = \{ t \in [T] \ | \  \ell ( \xb(\omega, t) ) = \ell \ \text{and} \ t < \tau_{-\ell, \omega} \} .  \label{eq:inleaf_periods}
\end{equation}

The second step is to determine the \emph{permissible stop periods}. Depending on whether we are optimizing for the left-stop subtree or for the right-stop subtree, these are the periods at which the subtree could stop. The right-stop permissible stop periods are defined as
\begin{equation}
P_{\omega} = \left\{ t \in S_{\omega} \ | \ x_v(\omega,t) > \max\{ x_v(\omega, t') \, | \, t' \in S_{\omega} \ \text{and} \ t' < t\} \right\}, \label{eq:permissible_stop_rightstop_subtree}
\end{equation}
where the maximum is defined as $-\infty$ if the corresponding set is empty. The left-stop permissible stop periods are similarly defined as 
\begin{equation}
P_{\omega} = \left\{ t \in S_{\omega} \ | \ x_v(\omega,t) < \min\{ x_v(\omega, t') \, | \, t' \in S_{\omega} \ \text{and} \ t' < t\} \right\},
\label{eq:permissible_stop_leftstop_subtree}
\end{equation}
where the minimum is defined as $+\infty$ if the corresponding set is empty. 

The third step is to construct the trajectory-specific piecewise-constant functions. Each such piecewise constant function will tell us how the reward of a trajectory $\omega$ varies as a function of $\theta$. Let us order the times in $P_{\omega}$ as  $P_{\omega} = \{ t_{\omega,1}, t_{\omega,2}, \dots, t_{\omega, |P_{\omega}|} \}$, 
where $t_{\omega,1} < t_{\omega,2} < \dots < t_{\omega, |P_{\omega}|}$. We use those times to define the corresponding breakpoints of our piecewise constant function; the $i$th breakpoint, $b_{\omega,i}$, is defined as 
\begin{equation*}
b_{\omega,i} = x_v(\omega, t_{\omega,i} ).
\end{equation*}
The function value of the $i$th piece is given by $f_{\omega,i}$, which is defined as the value if we stopped at the $i$th permissible stop period:
\begin{equation*}
f_{\omega,i} = g(t_{\omega,i}, \xb(\omega, t_{\omega,i}) ).
\end{equation*}
For a right-stop subtree, the corresponding piecewise constant function is denoted by $F_{\omega}$, and it is defined as follows:
\begin{equation}
F_{\omega}(\theta) = \left\{  \begin{array}{ll} 
f_{\omega,1} & \text{if} \ \theta < b_{\omega,1}, \\
f_{\omega,2} & \text{if} \ b_{\omega,1} \leq \theta < b_{\omega,2}, \\
f_{\omega,3} & \text{if} \ b_{\omega,2} \leq \theta < b_{\omega,3}, \\
\vdots & \vdots \\
f_{\omega, |P_{\omega}|} & \text{if} \ b_{\omega, |P_{\omega}|-1 } \leq \theta < b_{\omega, |P_{\omega}| }, \\
f_{\omega, \text{ns}} & \text{if} \ \theta \geq b_{\omega, |P_{\omega}| }.
\end{array} \right. \label{eq:fomega_rightstop}
\end{equation}
For a left-stop subtree, the piecewise constant function is defined similarly, as 
\begin{equation}
F_{\omega}(\theta) = \left\{ \begin{array}{ll}
f_{\omega,1} & \text{if} \ \theta \geq b_{\omega,1}, \\
f_{\omega,2} & \text{if} \  b_{\omega,2} \leq \theta < b_{\omega,1}, \\
f_{\omega,3} & \text{if} \  b_{\omega,3} \leq \theta < b_{\omega,2}, \\
\vdots & \vdots \\
f_{\omega,|P_{\omega}|} &  \text{if} \ b_{\omega,|P_{\omega}|} \leq \theta < b_{\omega,|P_{\omega}|-1}, \\
f_{\omega, \text{ns}} & \text{if}\ \theta < b_{\omega,|P_{\omega}|}.
\end{array} \right. \label{eq:fomega_leftstop}
\end{equation}
The function value $F_{\omega}(\theta)$ is exactly the reward that will be garnered from trajectory $\omega$ if we set the split point to $\theta$. By averaging these functions over $\omega$, we obtain the function $F$, which returns the average reward, over the whole training set of trajectories, that ensues from setting the split point to $\theta$:
\begin{equation}
F(\theta) = (1/ \Omega) \sum_{\omega=1}^{\Omega} F_{\omega}(\theta). \label{eq:F_average}
\end{equation}
We wish to find the value of $\theta$ that maximizes the function $F$, i.e., 
\begin{equation*}
\theta^* \in \arg \max_{\theta \in \mathbb{R}} F(\theta).
\end{equation*}
Note that because each $F_{\omega}(\cdot)$ is a piecewise constant function of $\theta$, the overall average $F(\cdot)$ will also be a piecewise constant function of $\theta$, and as a result there is no unique maximizer of $F(\cdot)$. This is analogous to the situation faced in classification, where there is usually an interval of split points that lie between the values of the given coordinate of two points. Let $I =  \arg \max_{\theta \in \mathbb{R}} F(\theta)$ be the interval on which $F(\cdot)$ is maximized. For convenience, consider the interior $\text{int}(I)$ of this interval, which will always be an interval of the form $(-\infty,b)$, $(b, b')$ or $(b, \infty)$ for some values $b, b' \in \mathbb{R}$, with $b < b'$. We can then set the split point as follows:
\begin{equation}
\theta^* = \left\{ \begin{array}{ll} 
-\infty & \text{if} \ \text{int}(I) = (-\infty, b) \ \text{for some} \ b \in \mathbb{R},  \\
(b + b')/2 & \text{if} \ \text{int}(I) = (b, b')\ \text{for some} \ b, b' \in \mathbb{R}, \ b < b', \\
+\infty & \text{if} \ \text{int}(I) = (b, +\infty)\ \text{for some} \ b \in \mathbb{R}.
\end{array} \right. \label{eq:argmax_to_thetastar}
\end{equation}
In the case that $I$ is a bounded interval, we set the split point to the midpoint of the interval. If $I$ is an unbounded interval, then there is no well-defined midpoint of the interval, and we set the split point to be either $+\infty$ or $-\infty$. (Note that we do not impose any restriction on the split point. This is in contrast to our convergence results in Section~\ref{subsec:problem_convergence}, where we required the state space to be totally bounded. In the setting we consider here, the function $F$ is derived from a finite sample and is piecewise constant, allowing it to be easily optimized even if the split point can be chosen from $\mathbb{R}$. In the setting of the convergence result, the need for total boundedness arises because we are optimizing over the split point when the expected reward is with respect to an \emph{infinite} sample, i.e., it is the exact expected value of the policy's reward with respect to the underlying stochastic process.)

We summarize the procedure formally as Algorithm~\ref{algorithm:OptimizeSplitPoint}. We provide an example in Section~\ref{sec:example_optimizesplitpoint} to illustrate the procedure. %

\begin{algorithm}
{\SingleSpacedXI
\begin{algorithmic}
\REQUIRE Coordinate $v$, leaf $\ell$, direction of subtree $\dir$, current tree policy $(\Tcal, \vb, \thetab, \ab)$. 
\FOR{$\omega \in [\Omega]$}
	\STATE Compute no-stop time $\tau_{-\ell,\omega}$ using equation~\eqref{eq:nostop_time}. 
	\STATE Compute no-stop value $f_{\omega, \text{ns}}$ using equation~\eqref{eq:nostop_value}.
	\STATE Compute in-leaf periods $S_{\omega}$ using equation~\eqref{eq:inleaf_periods}. 
	\IF{$\dir = \text{left} $}
		\STATE Compute the permissible stop periods $P_{\omega}$ for left-stop subtree using equation~\eqref{eq:permissible_stop_leftstop_subtree}.
		\STATE Order $P_{\omega}$ as $P_{\omega} = \{ t_{\omega,1}, t_{\omega,2}, \dots, t_{\omega, |P_{\omega}|} \}$, where $t_{\omega,1} < t_{\omega,2} < \dots < t_{\omega, |P_{\omega}|}$. 
		\STATE Compute the trajectory function $F_{\omega}(\cdot)$ using equation~\eqref{eq:fomega_leftstop}.
	\ELSE
		\STATE Compute the permissible stop periods $P_{\omega}$ for right-stop subtree using equation~\eqref{eq:permissible_stop_rightstop_subtree}.
		\STATE Order $P_{\omega}$ as $P_{\omega} = \{ t_{\omega,1}, t_{\omega,2}, \dots, t_{\omega, |P_{\omega}|} \}$, where $t_{\omega,1} < t_{\omega,2} < \dots < t_{\omega, |P_{\omega}|}$. 
		\STATE Compute the trajectory function $F_{\omega}(\cdot)$ using equation~\eqref{eq:fomega_leftstop}.
	\ENDIF
\ENDFOR
\STATE Compute the sample average function $F(\cdot)$ using equation~\eqref{eq:F_average}.
\STATE Compute the maximizer set $I \leftarrow \arg \max_{\theta \in \mathbb{R}} F(\theta)$. 
\STATE Compute $\theta^*$ using equation~\eqref{eq:argmax_to_thetastar}.
\RETURN $\theta^*$.
\end{algorithmic}
}
\caption{ \textsc{OptimizeSplitPoint} function. \label{algorithm:OptimizeSplitPoint}}
\end{algorithm}

\subsection{Algorithm complexity}
\label{subsec:construction_complexity}

We make a few comments about the asymptotic complexity of our construction procedure in the case that the chosen stopping criterion is a maximum depth $d$. The complexity bound of one run of \textsc{OptimizeSplitPoint} (Algorithm~\ref{algorithm:OptimizeSplitPoint}) is driven by the computation of $F_\omega(\cdot)$ which can be done in $O(T)$ time. Then averaging all $F_\omega(\cdot)$ functions into $F(\cdot)$ requires time $O(T \Omega)$, and computing $I$ requires a sort for a total bound of $O(T\Omega (\log T + \log \Omega))$. In turn, Algorithm \ref{algorithm:construction} runs \textsc{OptimizeSplitPoint} at most $n \cdot 2^{d} \cdot 2$ times, that is, all possible split variable indices, for all leaves at the current tree depth up to $d$, and for the two possible subtree orientations. The complexity of one iteration of Algorithm \ref{algorithm:construction} is therefore $O(n2^d T \Omega (\log T + \log \Omega))$. Thus, up to logarithmic factors, the construction procedure is linear in all problem parameters except the maximum tree depth $d$.

\subsection{Comparison to top-down tree induction for classification}
\label{subsec:construction_comparison_to_classification}

We note that Algorithm \ref{algorithm:construction} is inspired by classical top-down classification tree induction algorithms such as CART \citep{breiman1984classification}, ID3 \citep{quinlan1986induction} and C4.5 \citep{quinlan1993c4}. However, there are a number of important differences, in both the problem and the algorithm. We have already discussed in Section~\ref{subsec:problem_complexity} that the problem of finding a stopping policy in the form of a tree is structurally different than finding a classifier in the form of a tree -- in particular, deciding leaf labels in a tree classifier is easy to do, whereas deciding leaf actions in a tree policy for a stopping problem is an NP-Hard problem. 

With regard to the algorithms themselves, there are several important differences. First, in our algorithm, we directly optimize the in-sample reward as we construct the tree. This is in contrast to how classification trees are built for classification, where typically it is not the classification error that is directly minimized, but rather an impurity metric such as the Gini impurity (as in CART) or the information gain/entropy (as in ID3 and C4.5). Second, in the classification tree setting, each leaf can be treated independently; once we determine that we should no longer split a leaf (e.g., because there is no more improvement in the impurity, or the number of points in the leaf is too low), we never need to consider that leaf again. In our algorithm, we must consider every leaf in each iteration, even if that leaf may not have resulted in improvement in the previous iteration; this is because the actions we take in the leaves interact with each other and are not independent of each other. Lastly, as mentioned earlier, determining the optimal split point is much more involved than in the classification tree setting: in classification, the main step is to sort the observations in a leaf by the split variable. In our optimal stopping setting, determining the split point is a more involved calculation that takes into account when each trajectory is in a given leaf and how the cumulative maximum/minimum of the values of a given split variable change with the time $t$.

\section{Application to option pricing}
\label{sec:results}

In this section, we report on the performance of our method in an application drawn from option pricing. We define the option pricing problem in Section~\ref{subsec:results_problemdefinition}. We compare our policies against two benchmarks in terms of out-of-sample performance and computation time in Sections~\ref{subsec:results_OOS_performance} and \ref{subsec:results_computation_time}, respectively. In Section~\ref{subsec:results_policy_structure}, we present our tree policies for this application and discuss their structure. Finally, in Section~\ref{subsec:results_SP500}, we evaluate our policies on instances derived from real S\&P500 stock price data. Additional numerical results are provided in Section~\ref{appendix:results}.

\subsection{Problem definition}
\label{subsec:results_problemdefinition}

High-dimensional option pricing is one of the classical applications of optimal stopping, and there is a wide body of literature devoted to developing good heuristic policies (see \citealt{glasserman2013monte} for a comprehensive survey).  
In this section, we illustrate our tree-based method on a standard family of option pricing problems from \cite{desai2012pathwise}. We consider two benchmarks. Our first benchmark is the least-squares Monte Carlo method from \cite{longstaff2001valuing}, which is a commonly used method for option pricing. Our second benchmark is the pathwise optimization (PO) method from \cite{desai2012pathwise}, which was shown in that paper to yield stronger exercise policies than the method of Longstaff and Schwartz. Both of these methods involve building regression models that estimate the continuation value at each time $t$ using a collection of basis functions of the underlying state, and using them within a greedy policy. 

In each problem, the option is a Bermudan max-call option, written on $n$ underlying assets. The stopping problem is defined over a period of 3 calendar years with $T = 54$ equally spaced exercise opportunities. The price paths of the $n$ assets are generated as a geometric Brownian motion with drift equal to the annualized risk-free rate $r=5\%$ and annualized volatility $\sigma = 20\%$, starting at an initial price $\bar{p}$. The pairwise correlation $\rho_{ij}$ between different assets $i \neq j$ is set to $\bar{\rho} = 0$. The strike price for each option is set at $K = 100$. Each option has a knock-out barrier $B = 170$, meaning that if any of the underlying stock prices exceeds $B$ at some time $t_0$, the option is ``knocked out'' and the option value becomes $0$ at all times $t \geq t_0$. Time is discounted continuously at the risk-free rate, which implies a discrete discount factor of $\beta = \exp( - 0.05 \times 3 / 54) = 0.99723$. 

Our state variable is defined as $\xb(t) = ( t, p_1(t), \dots, p_n(t), y(t), g(t) )$, where $t$ is the index of the period; $p_j(t)$ is the price of asset $j$ at exercise time $t \in \{1,\dots, T\}$; $y(t)$ is a binary variable that is 0 or 1 to indicate whether the option has not been knocked out by time $t$, defined as
\begin{equation}
y(t) = \Ibb\left\{  \max_{1 \leq j \leq n, 1 \leq t' \leq t} p_j(t)  < B \right\};
\end{equation}
and $g(t)$ is the payoff at time $t$, defined as 
\begin{equation}
g(t) = \max\left\{ 0, \max_{1 \leq j \leq n} p_j(t) - K \right\}  \cdot y(t).
\end{equation}

In our implementation of the tree optimization algorithm, we vary the subset of the state variables that the tree model is allowed to use. We use \textsc{time} to denote $t$, \textsc{prices} to denote $p_1(t), \dots, p_n(t)$, \textsc{payoff} to denote $g(t)$, and \textsc{KOind} to denote $y(t)$. We set the relative improvement parameter $\gamma$ to 0.005, requiring that we terminate the construction algorithm when the improvement in in-sample objective becomes lower than 0.5\%. We report on the sensitivity of our algorithm to the parameter $\gamma$ in Section~\ref{appendix:results_sensitivity}. In addition, in Section~\ref{appendix:results_kfoldcv}, we report results on our $k$-fold cross-validation algorithm (defined in Section~\ref{sec:kfold_cross_validation}) for tuning $\gamma$.

In our implementation of the Longstaff-Schwartz (LS) algorithm, we vary the basis functions that are used in the regression. We follow the same notation as for the tree optimization algorithm in denoting different subsets of state variables; we additionally define the following sets of basis functions:
\begin{itemize}
\item \textsc{one}: the constant function, $1$. 
\item \textsc{pricesKO}: the knock-out (KO) adjusted prices defined as $p_i(t) \cdot y(t)$ for $1 \leq i \leq n$.
\item \textsc{maxpriceKO} and \textsc{max2priceKO}: the largest and second largest KO adjusted prices.%
\item \textsc{prices2KO}: the knock-out adjusted second-order price terms, defined as $p_i(t) \cdot p_j(t) \cdot y(t)$ for $1 \leq i \leq j \leq n$. 
\end{itemize}

In our implementation of the PO algorithm, we also vary the basis functions, and follow the same notation as for LS. We use 500 inner samples, as in \cite{desai2012pathwise}.

We vary the number of stocks as $n = 4, 8, 16$ and the initial price of all assets as $\bar{p} = 90, 100, 110$. For each combination of $n$ and $\bar{p}$, we consider ten replications. In each replication, we generate $\Omega = 20,000$ trajectories for training the methods and 100,000 trajectories for out-of-sample testing. In Section~\ref{appendix:results_highdim} we also consider higher dimensional option pricing problems, where the number of stocks is on the order of the number of trajectories. All replications were executed on the Amazon Elastic Compute Cloud (EC2) using a single instance of type \texttt{r4.4xlarge} (Intel Xeon E5-2686 v4 processor with 16 virtual CPUs and 122 GB memory). All methods were implemented in the Julia technical computing language, version 0.6.2 \citep{bezanson2017julia}. All linear optimization problems for the pathwise optimization method were formulated using the JuMP package for Julia \citep{lubin2015computing,dunning2017jump} and solved using Gurobi 8.0 \citep{gurobi}.

\subsection{Out-of-sample performance}
\label{subsec:results_OOS_performance}

Table~\ref{table:OOS_LS_vs_tree_neq8} shows the out-of-sample reward garnered by the LS and PO policies for different basis function architectures and the tree policies for different subsets of the state variables, for different values of the initial price $\bar{p}$. The rewards are averaged over the ten replications, with standard errors reported in parentheses. For ease of exposition, we focus only on the $n = 8$ assets, as the results for $n = 4$ and $n = 16$ are qualitatively similar; for completeness, these results are provided in Section~\ref{appendix:results_rho0_performance}. In addition, Section~\ref{appendix:results_rhoneq0_performance} provides additional results for when the common correlation $\bar{\rho}$ is not equal to zero. 
\begin{table}[ht]
\caption{Comparison of out-of-sample performance between LSM, PO and tree policies for $n = 8$ assets, for different initial prices $\bar{p}$. In each column, the best performance is indicated in bold.}
 \label{table:OOS_LS_vs_tree_neq8}
\centering
\small
\begin{tabular}{lllccc}
  \toprule
$n$  & Method & State variables / Basis functions &  \multicolumn{3}{c}{Initial Price} \\
 &  & & $\bar{p} = 90$ & $\bar{p} = 100$ & $\bar{p} = 110$ \\ 
  \midrule
 8 & LS &\textsc{one} & 33.82\enskip (0.021) & 38.70\enskip (0.023) & 43.13\enskip (0.015) \\ 
    8 & LS & \textsc{prices} & 33.88\enskip (0.019) & 38.59\enskip (0.023) & 43.03\enskip (0.014) \\ 
    8 & LS & \textsc{pricesKO} & 41.45\enskip (0.027) & 49.33\enskip (0.017) & 53.08\enskip (0.009) \\ 
    8 & LS & \textsc{pricesKO}, \textsc{KOind} & 41.86\enskip (0.021) & 49.36\enskip (0.020) & 53.43\enskip (0.012) \\ 
    8 & LS & \textsc{pricesKO}, \textsc{KOind}, \textsc{payoff} & 43.79\enskip (0.022) & 49.86\enskip (0.013) & 53.07\enskip (0.009) \\[0.25em] 
    8 & LS & \textsc{pricesKO}, \textsc{KOind}, \textsc{payoff},  & 43.83\enskip (0.021) & 49.86\enskip (0.013) & 53.07\enskip (0.009) \\
    & & \textsc{maxpriceKO} & \\[0.25em]
    8 & LS & \textsc{pricesKO}, \textsc{KOind}, \textsc{payoff}, & 43.85\enskip (0.022) & 49.87\enskip (0.012) & 53.06\enskip (0.008) \\ 
    & & \textsc{maxpriceKO}, \textsc{max2priceKO} \\ [0.25em]
    8 & LS & \textsc{pricesKO}, \textsc{payoff} & 44.06\enskip (0.013) & 49.61\enskip (0.010) & 52.65\enskip (0.008) \\[0.25em] 
    8 & LS & \textsc{pricesKO}, \textsc{prices2KO},  & 44.07\enskip (0.013) & 49.93\enskip (0.010) & 53.11\enskip (0.010) \\
    & & \textsc{KOind}, \textsc{payoff} \\ [0.5em]
      8 & PO & \textsc{prices} & 40.94\enskip (0.012) & 44.84\enskip (0.016) & 47.48\enskip (0.014) \\ 
    8 & PO & \textsc{pricesKO}, \textsc{KOind}, \textsc{payoff} & 44.01\enskip (0.019) & 50.71\enskip (0.011) & 53.82\enskip (0.009) \\[0.25em]
    8 & PO & \textsc{pricesKO}, \textsc{KOind}, \textsc{payoff},  & 44.07\enskip (0.017) & 50.67\enskip (0.011) & 53.81\enskip (0.011) \\ 
&     & \textsc{maxpriceKO}, \textsc{max2priceKO} & \\[0.25em]
    8 & PO & \textsc{pricesKO}, \textsc{prices2KO}, & 44.66\enskip (0.018) & 50.67\enskip (0.010) & 53.77\enskip (0.008) \\
    & & \textsc{KOind}, \textsc{payoff} \\  [0.5em]
      8 & Tree & \textsc{payoff}, \textsc{time} & \bfseries 45.40\enskip (0.018) & \bfseries 51.28\enskip (0.016) & \bfseries 54.52\enskip (0.006) \\ 
    8 & Tree & \textsc{prices} & 35.86\enskip (0.170) & 43.42\enskip (0.118) & 46.95\enskip (0.112) \\ 
    8 & Tree & \textsc{prices}, \textsc{payoff} & 39.13\enskip (0.018) & 48.37\enskip (0.014) & 53.61\enskip (0.010) \\ 
    8 & Tree & \textsc{prices}, \textsc{time} & 38.21\enskip (0.262) & 40.21\enskip (0.469) & 42.69\enskip (0.133) \\ 
    8 & Tree & \textsc{prices}, \textsc{time}, \textsc{payoff} & \bfseries 45.40\enskip (0.017) & \bfseries 51.28\enskip (0.016) & 54.51\enskip (0.006) \\ 
    8 & Tree & \textsc{prices}, \textsc{time}, \textsc{payoff}, \textsc{KOind} & \bfseries 45.40\enskip (0.017) &  \bfseries51.28\enskip (0.016) & 54.51\enskip (0.006) \\
   \bottomrule
\end{tabular}
\end{table}

From this table, it is important to recognize three key insights. First, for all three values of $\bar{p}$, the best tree policies -- specifically, those tree policies that use the \textsc{time} and \textsc{payoff} state variables -- are able to outperform all of the LS and PO policies. Relative to the best LS policy for each $\bar{p}$, these improvements range from 2.04\% ($\bar{p} = 110$) to 3.02\% ($\bar{p} = 90$), which is substantial given the context of this problem. Relative to the best PO policy for each $\bar{p}$, the improvements range from 1.12\% ($\bar{p} = 100$) to 1.66\% ($\bar{p} = 90$), which is still a remarkable improvement. 

Second, observe that this improvement is attained despite an experimental setup biased in favor of LS and PO. In this experiment, the tree optimization algorithm was only allowed to construct policies using subsets of the primitive state variables. In contrast, both the LS and the PO method were tested with a richer set of basis function architectures that included knock-out adjusted prices, highest and second-highest prices and second-order price terms. From this perspective, it is significant that our tree policies could outperform the best LS policy and the best PO policy. This also highlights an advantage of our tree optimization algorithm, which is its nonparametric nature: if the boundary between $\stop$ and $\go$ in the optimal policy is highly nonlinear with respect to the state variables, then by estimating a tree policy one should (in theory) be able to closely approximate this structure with enough splits (as suggested by Theorem~\ref{thm:tree_policy_appr_opt}). In contrast, the performance of LS and PO is highly dependent on the basis functions used, and requires the DM to specify a basis function architecture. %

Third, with regard to our tree policies specifically, we observe that policies that use \textsc{time} and \textsc{payoff} perform the best. The \textsc{time} state variable is critical to the success of the tree policies because the time horizon is finite: as such, a good policy should behave differently near the end of the horizon from how it behaves at the start of the time horizon. The LS algorithm handles this automatically because it regresses the continuation value from $t = T-1$ to $t = 1$, so the resulting policy is naturally time-dependent. The PO policy is obtained in a similar way, with the difference that one regresses an upper bound from $t = T-1$ to $t = 1$, so the PO policy is also time-dependent. Without time as an explicit state variable, our tree policies will estimate stationary policies, which are unlikely to do well given the nature of the problem. Still, it is interesting to observe some instances in our results where our time-\emph{in}dependent policies outperform time-dependent ones from LS (for example, compare tree policies with \textsc{prices} only to LS with \textsc{prices} or \textsc{one}). With regard to \textsc{payoff}, we note that because the payoff $g(t)$ is a function of the prices $p_1(t), \dots, p_n(t)$, one should in theory be able to replicate splits on \textsc{payoff} using a collection of splits on \textsc{price}; including \textsc{payoff} explicitly helps the construction algorithm recognize such collections of splits through a single split on the payoff variable.

In addition to the reward, it is also interesting to ask the question of how close to optimal our tree policies are. To answer this question, we can use the PO approach, which can produce an upper bound on the optimal expected reward. In Table~\ref{table:PO_UB}, we report on the upper bounds produced by the PO method for different initial prices and different basis function architectures. We remark here that the bound we report is the \emph{biased} upper bound, which is the objective value of the PO linear optimization problem. Ideally, as discussed in \cite{desai2012pathwise}, one would use an unbiased upper bound, which would be obtained by evaluating the solution from the PO linear optimization problem on a new independent sample of trajectories, together with a corresponding collection of inner samples. In our experimentation, we found that computing the unbiased bound was computationally quite prohibitive, and that on small examples, the unbiased and biased bounds were quite close; thus, for simplicity, we report the biased bound. Notwithstanding this bias, we can see that as compared to the tightest upper bounds, our best tree policies result in optimality gaps no greater than 1.5\%, suggesting that our policies are quite close to optimal.

\begin{table}[ht]
\caption{PO upper bound for $n = 8$ assets, for different initial prices $\bar{p}$. (For ease of comparison, the performance of the tree policy with \textsc{payoff} and \textsc{time} from Table~\ref{table:OOS_LS_vs_tree_neq8} is reproduced at the bottom.)  }
 \label{table:PO_UB}
\centering
\small
\begin{tabular}{lllccc}
  \toprule
$n$  & Method & State variables / Basis functions &  \multicolumn{3}{c}{Initial Price} \\
 &  & & $\bar{p} = 90$ & $\bar{p} = 100$ & $\bar{p} = 110$ \\ 
  \midrule
    8 & PO-UB & \textsc{prices} & 51.39\enskip (0.023) & 57.21\enskip (0.009) & 60.32\enskip (0.006) \\ 
    8 & PO-UB &  \textsc{pricesKO}, \textsc{KOind}, \textsc{payoff} & 46.13\enskip (0.022) & 52.04\enskip (0.022) & 55.05\enskip (0.015) \\ 
    8 & PO-UB &  \textsc{pricesKO}, \textsc{KOind}, \textsc{payoff}, & 46.11\enskip (0.025) & 52.03\enskip (0.021) & 55.05\enskip (0.016) \\ 
 &   &\textsc{maxpriceKO}, \textsc{max2priceKO}  & & \\[0.5em]
    8 & PO-UB & \textsc{pricesKO}, \textsc{prices2KO}, & 46.08\enskip (0.022) & 51.97\enskip (0.023) & 55.00\enskip (0.016) \\ 
    & & \textsc{KOind}, \textsc{payoff} \\[0.5em]
  8 &  Tree & \textsc{payoff}, \textsc{time} & 45.40\enskip (0.018) & 51.28\enskip (0.016) & 54.52\enskip (0.006) \\ 
   \bottomrule
\end{tabular}
\end{table}

\subsection{Computation time} 
\label{subsec:results_computation_time}

Table~\ref{table:OOS_LS_vs_tree_neq8_rho0} reports the computation time for the LS, PO and tree policies for $n = 8$ and $\bar{p} \in \{90, 100, 110\}$, for the uncorrelated ($\bar{\rho} = 0$) case. The computation times are averaged over the ten replications for each combination of $n$ and $\bar{p}$. As with the performance results, we focus on $n = 8$ to simplify the exposition; additional timing results for $n = 4$ and $n = 16$ are provided in Section~\ref{appendix:results_rho0_timing}. For LS, the computation time consists of only the time required to perform the regressions from $t = T-1$ to $t = 1$. For PO, the computation time consists of the time required to formulate the linear optimization problem in JuMP, the solution time of this problem in Gurobi, and the time required to perform the regressions  from $t = T-1$ to $t= 1$ (as in Longstaff-Schwartz). For the tree method, the computation consists of the time required to run Algorithm~\ref{algorithm:construction}. 

\begin{table}
\caption{Comparison of estimation time between Longstaff-Schwartz, pathwise optimization and tree policies for $n = 8$ assets, for different initial prices $\bar{p}$ and common correlation $\bar{\rho} = 0$.}
 \label{table:OOS_LS_vs_tree_neq8_rho0}
 \small
  \centering
\begin{tabular}{lllrrr}
  \toprule
$n$  & Method & State variables / Basis functions &  \multicolumn{3}{c}{Initial Price} \\
 &  & & $\bar{p} = 90$ & $\bar{p} = 100$ & $\bar{p} = 110$ \\ 
  \midrule
  8 & LS & \textsc{one} & 1.2\enskip (0.0) & 1.2\enskip (0.0) & 1.2\enskip (0.0) \\ 
    8 & LS & \textsc{prices} & 1.4\enskip (0.0) & 1.4\enskip (0.0) & 1.4\enskip (0.0) \\ 
    8 & LS & \textsc{pricesKO} & 1.4\enskip (0.1) & 1.5\enskip (0.1) & 1.5\enskip (0.1) \\ 
    8 & LS & \textsc{pricesKO}, \textsc{KOind} & 1.6\enskip (0.1) & 1.4\enskip (0.1) & 1.5\enskip (0.1) \\ 
    8 & LS & \textsc{pricesKO}, \textsc{KOind}, \textsc{payoff} & 1.9\enskip (0.2) & 1.9\enskip (0.1) & 1.8\enskip (0.1) \\ 
    8 & LS & \textsc{pricesKO}, \textsc{KOind}, \textsc{payoff},  & 2.5\enskip (0.2) & 2.4\enskip (0.2) & 2.4\enskip (0.2) \\ 
    & & \textsc{maxpriceKO} \\
    8 & LS & \textsc{pricesKO}, \textsc{KOind}, \textsc{payoff}, & 2.7\enskip (0.2) & 2.6\enskip (0.3) & 2.4\enskip (0.2) \\ 
    & & \textsc{maxpriceKO}, \textsc{max2priceKO}  \\
    8 & LS & \textsc{pricesKO}, \textsc{payoff} & 1.7\enskip (0.2) & 1.6\enskip (0.1) & 1.4\enskip (0.1) \\ 
    8 & LS & \textsc{pricesKO}, \textsc{prices2KO}, \textsc{KOind}, \textsc{payoff} & 5.5\enskip (0.4) & 4.4\enskip (0.2) & 4.5\enskip (0.2) \\[0.5em]
 8 & PO & \textsc{prices} & 33.3\enskip (0.7) & 35.5\enskip (0.7) & 32.8\enskip (0.7) \\ 
    8 & PO & \textsc{pricesKO}, \textsc{KOind}, \textsc{payoff} & 76.8\enskip (2.8) & 73.8\enskip (5.0) & 56.9\enskip (2.1) \\ 
    8 & PO & \textsc{pricesKO}, \textsc{KOind}, \textsc{payoff}, & 104.7\enskip (5.8) & 79.6\enskip (3.3) & 66.8\enskip (3.8) \\ 
    & & \textsc{maxpriceKO}, \textsc{max2priceKO}  \\
    8 & PO & \textsc{pricesKO}, \textsc{prices2KO}, \textsc{KOind}, \textsc{payoff} & 221.3\enskip (9.4) & 180.7\enskip (4.1) & 142.0\enskip (4.2) \\[0.5em]
  8 & Tree & \textsc{payoff}, \textsc{time} & 7.6\enskip (0.3) & 3.9\enskip (0.3) & 3.2\enskip (0.1) \\ 
    8 & Tree & \textsc{prices} & 124.7\enskip (8.6) & 125.5\enskip (4.8) & 125.0\enskip (5.8) \\ 
    8 & Tree & \textsc{prices}, \textsc{payoff} & 5.5\enskip (0.1) & 5.3\enskip (0.1) & 5.1\enskip (0.1) \\ 
    8 & Tree & \textsc{prices}, \textsc{time} & 158.8\enskip (12.5) & 101.0\enskip (12.8) & 51.1\enskip (2.7) \\ 
    8 & Tree & \textsc{prices}, \textsc{time}, \textsc{payoff} & 20.4\enskip (1.0) & 10.9\enskip (0.6) & 9.3\enskip (0.1) \\ 
    8 & Tree & \textsc{prices}, \textsc{time}, \textsc{payoff}, \textsc{KOind} & 21.5\enskip (1.5) & 11.2\enskip (0.6) & 9.3\enskip (0.2) \\ \bottomrule
\end{tabular}
\end{table}

From this table, we can see that although our method requires more computation time than LS, the times are in general quite modest: our method requires no more than 2.5 minutes on average in the largest case. (In experiments with $n = 16$, reported in Section~\ref{appendix:results_rho0_timing}, we find that the method requires no more than 5 minutes on average in the largest case.) The computation times of our method also compare quite favorably to the computation times for the PO method. We also remark here that our computation times for the PO method do \emph{not} include the time required to generate the inner paths and to pre-process them in order to formulate the PO linear optimization problem. For $n = 8$, including this additional time increases the computation times by a large amount, ranging from 540 seconds (using \textsc{prices} only; approximately 9 minutes) to 2654 seconds (using \textsc{pricesKO}, \textsc{prices2KO}, \textsc{KOind}, \textsc{payoff}; approximately 44 minutes).

\subsection{Policy structure}

\label{subsec:results_policy_structure}

\begin{figure}
\caption{Examples of tree policies for initial price $\bar{p} = 90$ with state variables \textsc{prices}, \textsc{time}, \textsc{payoff} and \textsc{KOind} for a single replication.}
\label{figure:example_tree_policies_peq90}
\begin{subfigure}[t]{0.3\textwidth}
                \centering
                \includegraphics[width = \textwidth]{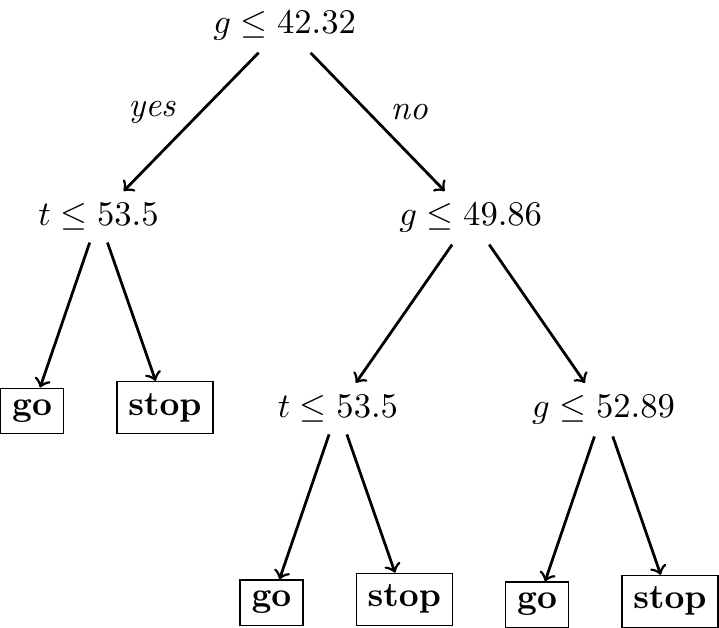}
                \caption{$n = 4$.}
                \label{fig:gull}
        \end{subfigure}  \hfill
\begin{subfigure}[t]{0.3\textwidth}
                \centering
                \includegraphics[width = \textwidth]{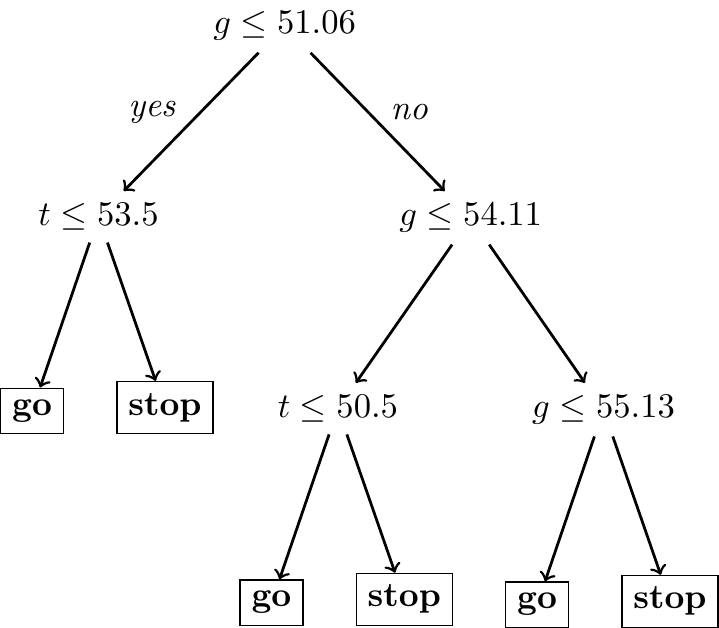}
                \caption{$n = 8$.}
                \label{fig:gull}
        \end{subfigure} \hfill
\begin{subfigure}[t]{0.3\textwidth}
                \centering
                \includegraphics[width = \textwidth]{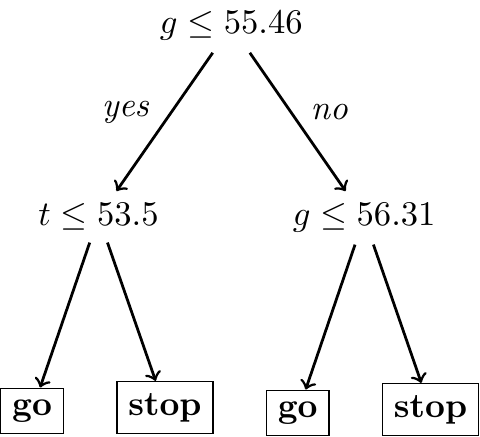}
                \caption{$n = 16$.}
                \label{fig:gull}
        \end{subfigure}
\vspace{-1.5em}
\end{figure}

It is also interesting to examine the structure of the policies that emerge from our tree optimization algorithm. Figure~\ref{figure:example_tree_policies_peq90} shows trees obtained using \textsc{prices}, \textsc{time}, \textsc{payoff} and \textsc{KOind} for one replication with $\bar{p} = 90$, for $n = 4, 8, 16$.
This figure presents a number of important qualitative insights about our algorithm. First, observe that the trees are extremely simple: there are no more than seven splits in any of the trees. The policies themselves are easy to understand and sensible. Taking $n = 8$ as an example, we see that if the payoff is lower than 51.06, it does not stop unless we are in the last period ($t = 54$), because there is still a chance that the payoff will be higher by $t = 54$. If the payoff is greater than 51.06 but less than or equal to 54.11, then the policy does not stop unless we are in the last four periods ($t = $ 51, 52, 53 or 54). If the payoff is greater than 54.11 but less than or equal to 55.13, then we continue. Otherwise, if the payoff is greater than 55.13, then we stop no matter what period we are in; this is likely because when such a payoff is observed, it is large enough and far enough in the horizon that it is unlikely a larger reward will be realized later. In general, as the payoff becomes larger, the policy will recommend stopping earlier in the horizon. Interestingly, the policies do not include any splits on the prices and the KO indicator, despite the construction algorithm being allowed to use these variables: this further underscores the ability of the construction algorithm to produce simple policies. It is also interesting to note that the tree structure is quite consistent across all three values of $n$. To the best of our knowledge, we do not know of any prior work suggesting that simple policies as in Figure~\ref{figure:example_tree_policies_peq90} can perform well against mainstream ADP methods for high-dimensional option pricing. In Section~\ref{appendix:results_sparseKO} of the ecompanion, we consider a more complicated family of instances where the effective dimension of the problem is larger, and in order for a policy to do well, it is necessary to split on additional variables beyond the time and the payoff.

We observe in the $n = 4$ and $n = 8$ trees that there is some redundancy in the splits. For example, for $n = 4$, observe that the left subtree of the split $g \leq 42.32$ is identical to the left subtree of the split $g \leq 49.86$; thus, the policy will take the same action whether $g \in [0, 42.32]$ or $g \in (42.32, 49.86]$, and the entire tree could be simplified by replacing it with the subtree rooted at the split $g \leq 49.86$. The reason for this redundancy is due to the greedy nature of the construction procedure. At the start of the construction procedure, splitting on $g \leq 42.32$ leads to the best improvement in the reward, but as this split and other splits are added, the split on $g \leq 49.86$ becomes more attractive. 

In addition to this redundancy, we also note that in all three trees, there is an interval $(g_1, g_2]$ such that if $g \in (g_1, g_2]$, the policy will continue (for example, in $n =8$, this interval is $(54.11, 55.13]$). This property of the policies is suboptimal because if $g$ is inside that interval, the policy may choose to continue even if $t = 54$, leading to a reward of zero for that trajectory; thus, we could in theory improve the performance of the policy by requiring the policy to stop if $t = 54$ (i.e., adding a right-stop subtree with the split $t \leq 53.5$). This may occur for two reasons: first, due to the sample-based nature of the optimization algorithm, the number of training set trajectories that are inside the interval $(g_1, g_2]$ at $t = 54$ may be small enough that adding the split $t \leq 53.5$ will not improve the overall sample-based objective by a relative factor of more than $1 + \gamma = 1.005$. The second reason is that, even if the sample is sufficiently large, the split may still fail to offer a sufficient improvement to the objective; this could be the case if the policy is such that there is a very small probability that a trajectory makes it to $t = 54$ and that $g(t)$ is in the interval $(g_1, g_2]$ at $t = 54$, so that taking the optimal action has a negligible impact on the objective. %

Finally, it is worth qualitatively comparing the policies that arise from our algorithm to policies derived from LS or the PO method. The trees in Figure~\ref{figure:example_tree_policies_peq90} fully specify the policy: the decision to stop or go can be made by checking at most three logical conditions. Moreover, the trees in Figure~\ref{figure:example_tree_policies_peq90} directly tell us when the policy will stop: either when the reward is very high or when we are in the final period of the horizon. In contrast, a policy derived from the LS or the PO method will consist of regression coefficients for all basis functions for each $t$ in the time horizon; we would not consider this policy interpretable due to the high dimensionality of its specification, as well as the fact that it does not drop, or clearly delineate, which basis functions are non-informative to the policy. From such a specification, it is difficult to immediately understand what the behavior of the policy will be (i.e., at what times and in what part of the state space the policy will stop), and even more so when the basis functions are more complex (such as the second-order price terms). We provide an example of a LS policy for $n = 8$ in Section~\ref{appendix:results_LS_policy_example} to further illustrate this.

\subsection{Out-of-sample performance with S\&P-500 data calibration}

\label{subsec:results_SP500}

In the previous sections, we considered options with artificial parameters specifying the stock price dynamics. In this section, we evaluate the performance of our tree policies on option pricing instances corresponding to real S\&P-500 stocks. An appealing feature of these experiments is that we use the raw S\&P-500 data as trajectories; this is opposed to making a distributional assumption, such as prices following a geometric Brownian motion, and then simulating trajectories that follow this distribution.

To create these instances, we consider stocks in the S\&P-500 from January 3, 2000 to November 17, 2017. This period consists of 4500 trading days. We remove any stocks for which data is missing on any trading day in this period. From the remaining set of 318 stocks, we create an \emph{instance} as follows. We sample four stocks without replacement from the set of 318. We create 100 such instances. The 4500 trading days are divided into sets of 30 consecutive trading days; each such set of 30 days is used to form a trajectory for the four selected stocks. This results in 150 trajectories, of which the first 100 are used for the training set, and the remaining 50 are used for the test set, and each trajectory corresponds to $T = 30$ exercise opportunities.

For each instance, the optimal stopping problem is to find an exercise policy on a max-call option written on the four stocks. To adapt each trajectory to this option pricing problem, we rescale each stock's prices over the 30 day period so that the initial price is \$100, and we define the payoff as $g(t) = \max\{ 0, \max_{1 \leq j \leq 4} p_j(t) - K\}$, where we set the strike price $K$ to \$105. Unlike the instances in Section~\ref{subsec:results_problemdefinition}, we do not define a knock-out barrier. 

We assume that payoffs are discounted at a continuous (annualized) interest rate of 0.02. We also assume that the exercise opportunities are equispaced, with each period corresponding to a day. We make this assumption for simplicity, as each set of 30 days of price data will contain prices that are separated by periods in which no trading occurred, such as weekends or public holidays.

Before continuing to our results, we wish to clarify that the goal of this experiment is to study option pricing problem instances where the stock price dynamics are as realistic as possible. The procedure we have described aims to accomplish this by using daily stock prices, with only the minor adjustment of normalizing the prices so that all stocks start at the same price in each trajectory. Aside from this modification, the trajectories here are completely unlike the ones we use in our prior experiments: the price of a single stock on a given day may exhibit a complicated dependence with the other stock prices on the same day or on earlier days, with the price of the same stock on earlier days or with the day itself. Stated differently, the dynamics of these trajectories need not be consistent with a well-behaved stochastic process such as geometric Brownian motion.

Table~\ref{table:SP500_LS_vs_tree_OOS} shows the average out-of-sample performance of the tree and LS policies with a variety of different state variable sets/basis function architectures. From this table, we can see that the tree policies that involve \textsc{payoff} and \textsc{time} achieve the highest average performance, and their average performance is about 15\% higher than the best performing LS policy.

\begin{table}[!ht]
\centering
\small
\begin{tabular}{llc} \toprule
Method & State variables /  & Out-of-sample  \\ 
& basis functions & performance \\ 
\midrule
Tree & \textsc{payoff}, \textsc{time} & \bfseries 4.71\enskip (0.163) \\ 
   & \textsc{prices} & 2.30\enskip (0.104) \\ 
   & \textsc{prices}, \textsc{payoff} & 2.71\enskip (0.115) \\ 
   & \textsc{prices}, \textsc{time} & 4.65\enskip (0.166) \\ 
   & \textsc{prices}, \textsc{time}, \textsc{payoff} & 4.68\enskip (0.165) \\[0.5em]
  LS & \textsc{one} & 3.97\enskip (0.143) \\ 
   & \textsc{prices} & 3.93\enskip (0.152) \\ 
   & \textsc{prices}, \textsc{one} & 4.11\enskip (0.160) \\ 
   & \textsc{prices}, \textsc{one}, \textsc{payoff} & 3.73\enskip (0.179) \\ 
   & \textsc{prices}, \textsc{one}, \textsc{payoff}, \textsc{maxprice} & 3.50\enskip (0.177) \\ 
   & \textsc{prices}, \textsc{payoff} & 4.01\enskip (0.180) \\ 
   & \textsc{prices}, \textsc{prices2}, \textsc{one}, \textsc{payoff} & 3.06\enskip (0.149) \\ 
   \bottomrule
\end{tabular}

\caption{Average out-of-sample reward of LS and tree policies for real S\&P-500 instances. Values reported are averages over the 100 instances, with standard errors in parentheses; bold is used to indicate the best policy. For LS, \textsc{prices2} denotes second-order price basis functions, and \textsc{maxprice} denotes the maximum price. \label{table:SP500_LS_vs_tree_OOS} }
\end{table}

\begin{figure}
\centering
\includegraphics[width=0.4\textwidth]{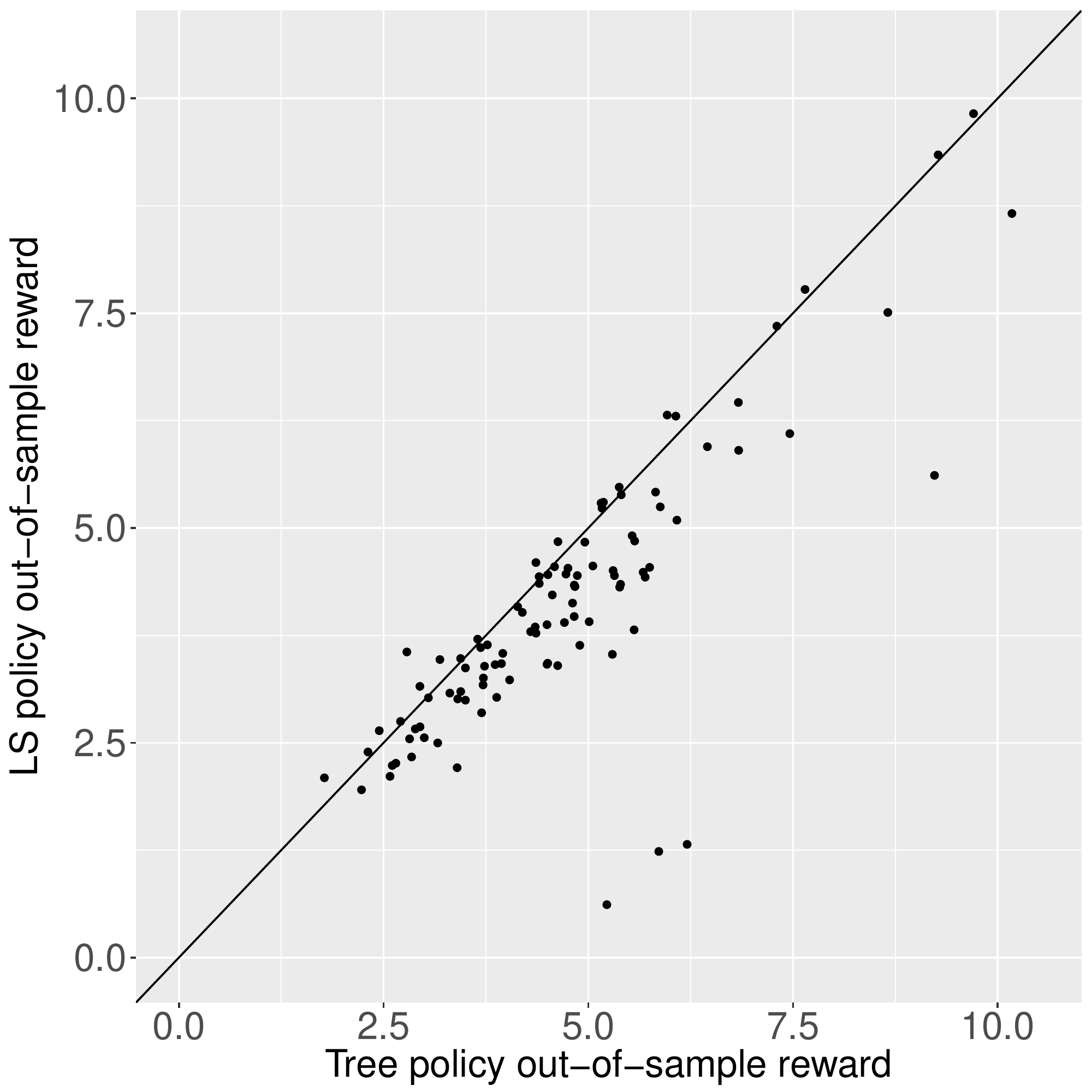} \quad \includegraphics[width=0.4\textwidth]{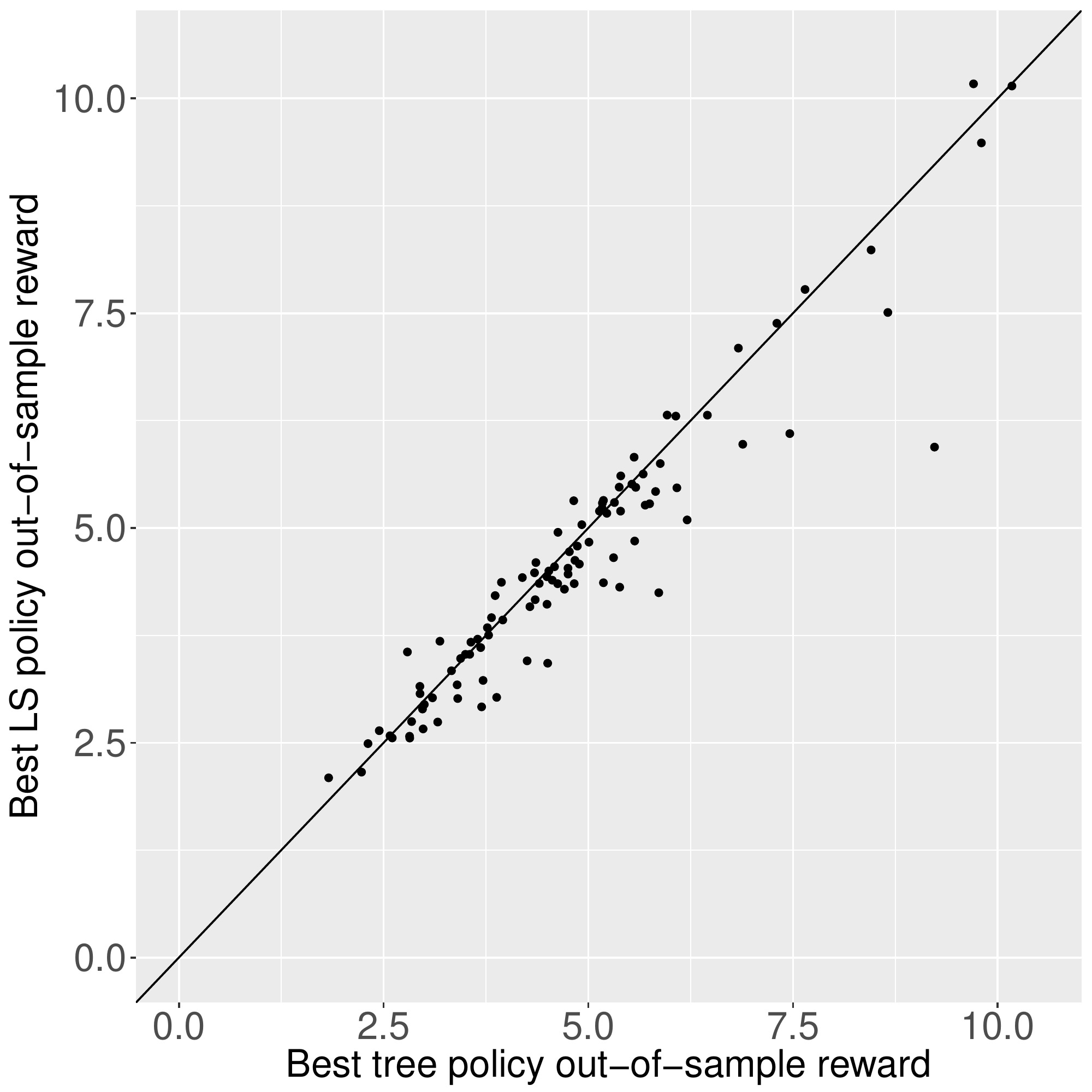}
\caption{Comparison of tree policies and LS policies for S\&P500 instances. The left hand plot compares the tree policy (with \textsc{payoff}, \textsc{time}) out-of-sample performance and  LS policy (with \textsc{one} and \textsc{prices}) out-of-sample performance for the 100 instances. The right hand plot compares, for each instance, the best tree out-of-sample performance and the best LS out-of-sample performance. \label{figure:SP500_LS_vs_tree_both} }
\end{figure}

To further compare the methods, the left hand plot of Figure~\ref{figure:SP500_LS_vs_tree_both} displays, for each of the 100 instances, the out-of-sample performance of the tree policy with \textsc{payoff} and \textsc{time} and the LS policy with the basis functions \textsc{one} and \textsc{prices}. From this plot, we can see that in general, on a per-instance basis, the tree policy outperforms the LS policy, as the majority of the points (approximately 80\%) are below the $y = x$ line. As a further comparison, the right hand plot of Figure~\ref{figure:SP500_LS_vs_tree_both} displays, for each of the 100 instances, the best out-of-sample performance achieved by any tree policy against the best out-of-sample performance achieved by any LS policy. From this plot, we can see that the tree policies are still in general better than the LS policies, although the edge is not as large as in the left hand plot. The best tree policy achieves a higher out-of-sample reward than the best LS policy in about two-thirds of the instances. Overall, these results suggest that our method can naturally learn policies that perform well on trajectories derived from real data, that may exhibit richer and more complicated dynamics than the instances based on geometric Brownian motion.

Before ending this section, we emphasize that our approach here is designed as a proof of concept for how one would solve these option pricing problems using strictly the raw data. This approach has some potential weaknesses compared to, for example, using S\&P-500 time-series data to calibrate geometric Brownian motion price processes for the underlying stocks, and then generating trajectories from these processes. Firstly, it limits the number of trajectories available for training and testing, whereas the classical approach would give control over the number of such trajectories. Secondly, the trajectories are generated by sequentially picking contiguous 30 day intervals over the entire $17$ year horizon. One could then expect that changes in market dynamics over such a long period would make the first and last stock price trajectories structurally different; it is thus conceivable that our $100$ training trajectories are ``stale'' and less informative concerning the more recent $50$ testing trajectories. In light of these challenges regarding the data being fed to the algorithm, it is quite encouraging that our approach remains robust and produces promising results versus incumbent methods.

\section{Application to one-dimensional uniform problem}
\label{sec:results_1D_Uniform}

In this section, we consider the following simple optimal stopping problem: we have a one-dimensional stochastic process, $x(1), x(2), \dots, x(T)$, where at each time $x(t)$ is independently drawn from a continuous $\text{Uniform}(0,1)$ distribution. The payoff at each $t$ is given by $g(t, x) = x$ and rewards are discounted by a discount factor of $\beta$. The rationale for considering this problem is its simplicity: it is small and simple enough that we can solve for the optimal policy directly, and obtain insight from comparing the performance of our tree policies to this optimum.

We compare both our tree policies and Longstaff-Schwartz against the optimal policy. We run the tree construction algorithm with \textsc{payoff} and \textsc{time} as state variables, and we run LS with the constant basis function (1). For the construction algorithm we set $\gamma = 0.005$. Note that for LS, it does not make sense to use $x(t)$ as a predictor in the regression; this is because the state variable $x(t)$ is drawn independently at each time, and so the optimal continuation value does not vary with the current $x(t)$. For this reason, we only consider the basis function architecture consisting of 1. We use 20,000 trajectories to build each model, and 100,000 to perform out-of-sample evaluation. We test values of $\beta$ in $\{0.9, 0.95, 0.97, 0.98, 0.99, 0.995, 0.999, 0.9999, 1.0 \}$. 

\begin{figure}
\centering 
\begin{subfigure}[t]{0.18\textwidth}
                \centering
                \includegraphics[width = 0.7\textwidth]{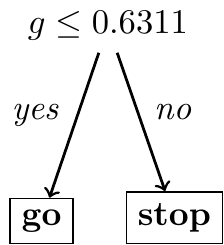}
                \caption{$\beta = 0.9$.}
        \end{subfigure}  \hfill
\begin{subfigure}[t]{0.18\textwidth}
                \centering
                \includegraphics[width = 0.7\textwidth]{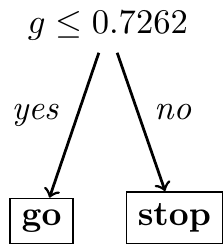}
                \caption{$\beta = 0.95$.}
        \end{subfigure}  \hfill
\begin{subfigure}[t]{0.18\textwidth}
                \centering
                \includegraphics[width = \textwidth]{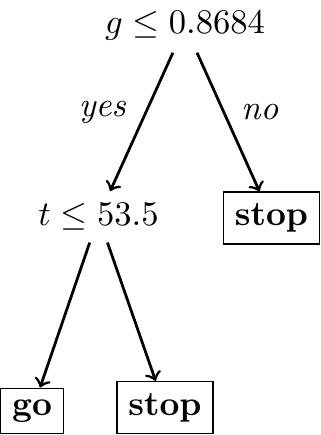}
                \caption{$\beta = 0.99$.}
        \end{subfigure}  \hfill
\begin{subfigure}[t]{0.18\textwidth}
                \centering
                \includegraphics[width = \textwidth]{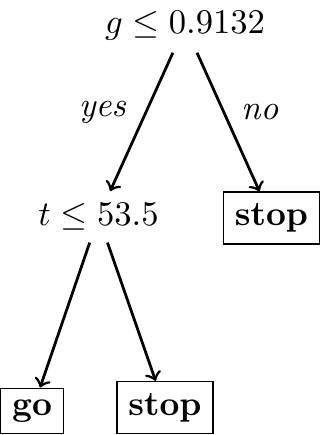}
                \caption{$\beta = 1.0$.}
        \end{subfigure} 
        
\caption{Examples of tree policies for $\beta = 0.9, 0.95, 0.99, 1.0$. \label{figure:1D_uniform_example_tree}}
\end{figure}

\begin{figure}
\begin{subfigure}[t]{0.5\textwidth}
                \centering
                \includegraphics[width = 0.85\textwidth]{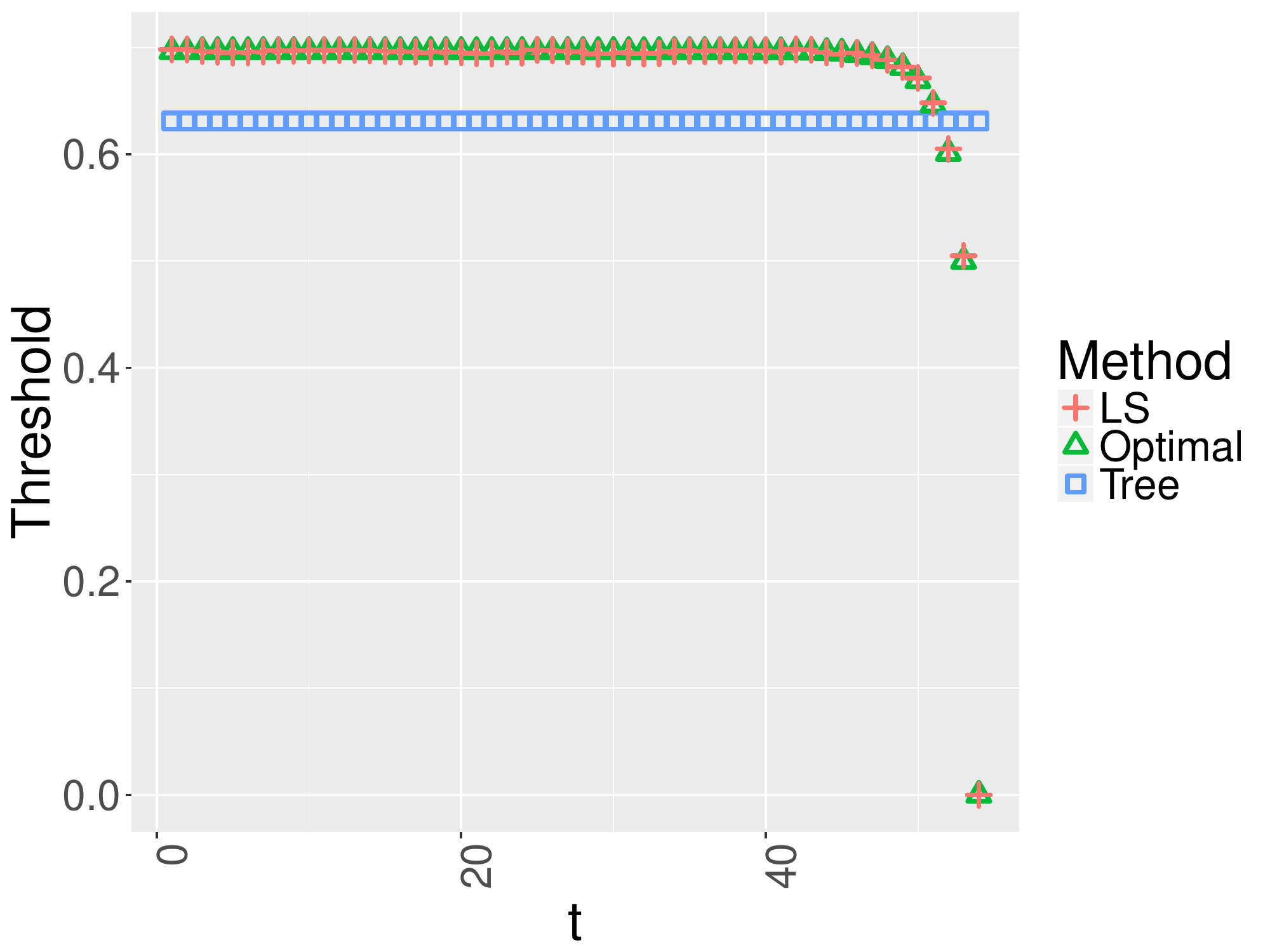}
                \caption{$\beta = 0.9$.}
        \end{subfigure}  \hfill
\begin{subfigure}[t]{0.5\textwidth}
                \centering
                \includegraphics[width = 0.85\textwidth]{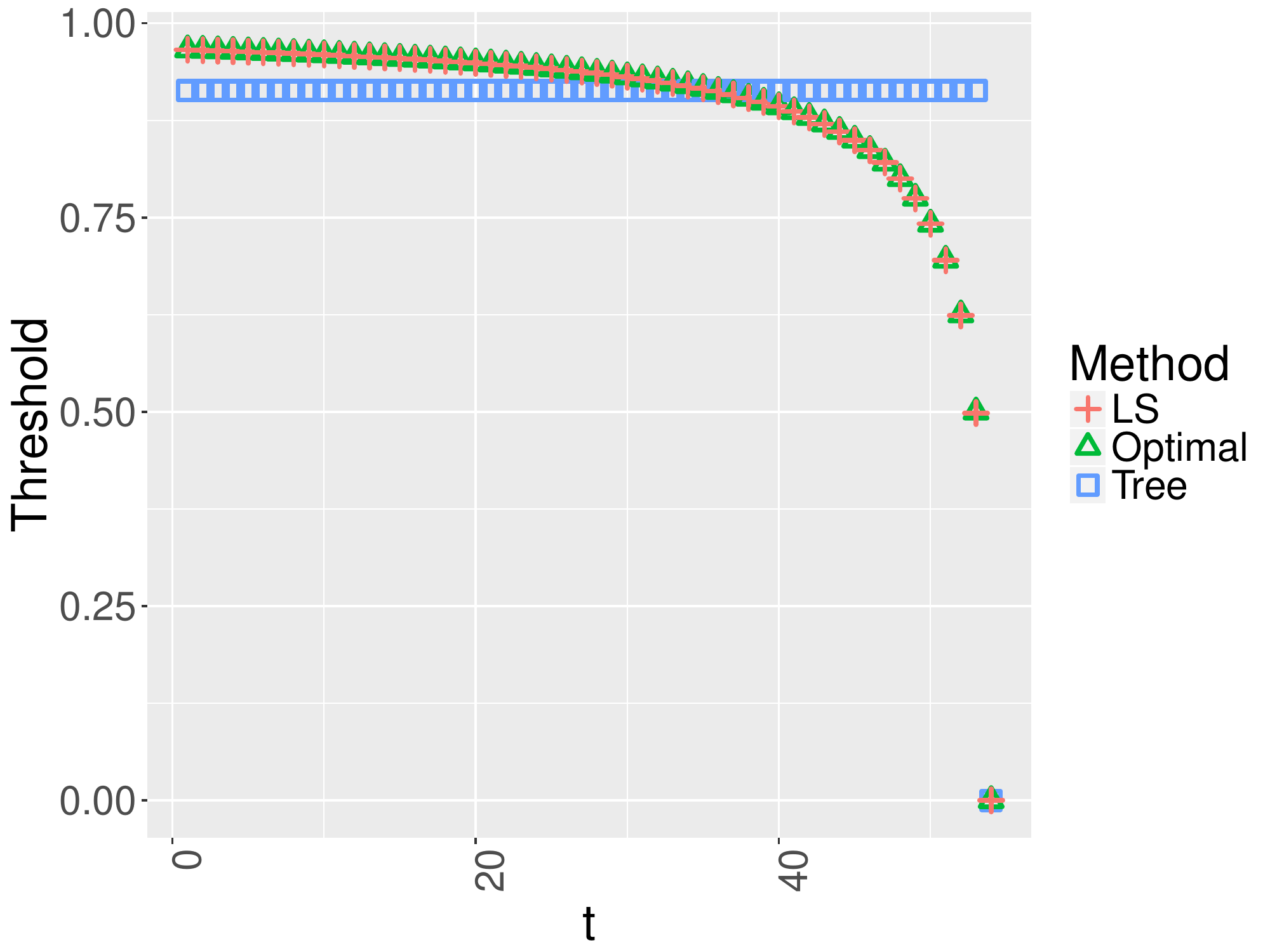}
                \caption{$\beta = 1.0$.}
        \end{subfigure} 
\caption{Comparison of thresholds used by the three different policies at $\beta = 0.9$ and $\beta = 1.0$. 
 (Note that in these experiments, the thresholds for ``LS'' and ``Optimal'' are almost identical, giving the appearance of a single curve.)
\label{figure:1D_uniform_thresholds}}
\end{figure}

Figure~\ref{figure:1D_uniform_example_tree} displays the tree policies for $\beta = 0.9, 0.95, 0.99, 1.0$. From this figure we can see that for $\beta =0.9$ and $\beta = 0.95$, the tree policy does not depend on time: at any $t$, we simply check whether the payoff $g$ is greater than or equal to some threshold value. For $\beta = 0.99$ and $\beta = 1.0$, we stop if either the payoff is greater than or equal to some threshold, or if we are in the last period. To compare these against the LS and optimal policies, we also plot in Figure~\ref{figure:1D_uniform_thresholds} the effective thresholds used by the three policies from $t = 1$ to $t = 54$.

From these two figures, we can see that the optimal thresholds are time-dependent, as we would expect, and the thresholds used by LS are essentially the same as those used by the optimal policy. In contrast, the thresholds used by the tree policies are significantly less time-dependent (for $\beta =0.9$, the threshold is constant, while for $\beta = 1.0$, the threshold only changes at $t = 54$, and is otherwise constant). 

Given the large differences in the policy structure, one would expect that the tree policies would be highly suboptimal. Surprisingly, however, this turns out not to be the case. To compare the policies, we show the out-of-sample performance of the tree, LS and optimal policies in Table~\ref{table:1D_uniform_example}. At low discount factors ($\beta \leq 0.99$), both the tree and LS policies are essentially optimal. For $\beta > 0.99$, LS displays a slight edge over the tree policies, and is essentially optimal; in the largest case, there is a difference of about 0.013 between the tree policies and the LS/optimal policies. 

\begin{table}[ht]
\centering
\small
\begin{tabular}{lcccc}
  \toprule $\beta$ & Tree & LS & Optimal (sim.) & Optimal (true) \\ 
\midrule 
0.9 & 0.6962 & 0.6961 & 0.6961 & 0.6964 \\ 
  0.95 & 0.7622 & 0.7622 & 0.7622 & 0.7620 \\ 
  0.97 & 0.8043 & 0.8043 & 0.8043 & 0.8044 \\ 
  0.98 & 0.8342 & 0.8342 & 0.8342 & 0.8340 \\ 
  0.99 & 0.8762 & 0.8763 & 0.8763 & 0.8763 \\ 
  0.995 & 0.9078 & 0.9086 & 0.9086 & 0.9087 \\ 
  0.999 & 0.9427 & 0.9507 & 0.9507 & 0.9507 \\ 
  0.9999 & 0.9528 & 0.9647 & 0.9647 & 0.9648 \\ 
  1 & 0.9532 & 0.9665 & 0.9665 & 0.9666 \\ 
   \bottomrule
\end{tabular}

\caption{Comparison of tree, LS and optimal policies for the 1D uniform problem. With the exception of ``Optimal (true)'' (the theoretical optimal reward for the problem), all values shown are out-of-sample rewards averaged over 5 independent replications. All standard errors are smaller than 0.0005. \label{table:1D_uniform_example}}
\end{table}

This experiment provides several insights. First, even for a simple problem such as this one, there may be near-optimal policies that are simpler than the optimal policy (for example, compare the constant threshold policy for $\beta = 0.9$ to the time-dependent threshold policy that is optimal) and that our construction algorithm can potentially discover such policies. Second, our approach is not a free lunch; indeed, for the higher discount rates, LS is able to recover the optimal policy, while our tree policies are suboptimal (albeit only by a small amount).

Lastly, we remark that in theory, the optimal policy for this problem could be represented as a tree policy, where one would have a sequence of splits on $t$, and each such split would be followed by a split on $g$. In this set of examples, our algorithm does not recover such a policy due to its greedy nature. One question in this direction is how to recognize splits that do not immediately yield an improvement in reward, but enable later splits that yield significant improvements. The answer to this question is not obvious, especially in light of the structure of the optimal stopping problem that is leveraged to efficiently optimize split points in our construction algorithm (Section~\ref{sec:construction}), and is an interesting direction for future research.

\section{Conclusion}
\label{sec:conclusion}

In this paper, we consider the problem of designing interpretable policies for optimal stopping, based on binary trees. We formulate the problem as an SAA problem, which we show to be theoretically intractable. Thus motivated, we develop a heuristic algorithm that greedily constructs a policy from the top down, by exploiting the stopping problem's structure in finding the optimal split point at each leaf. In numerical experiments on a standard option pricing problem, our algorithm attains better performance than state-of-the-art ADP methods for option pricing, while simultaneously producing policies that are significantly simpler and more transparent. In terms of future directions, besides those mentioned earlier, another valuable direction is to explore the application of this methodology in the healthcare domain, where optimal stopping appears in many problems \cite[for recent examples, see][]{iancu2018monitoring,cheng2019optimal}. %
We believe that this methodology represents an exciting starting point for future research at the intersection of stochastic control and interpretable machine learning.

\section*{Acknowledgments}
The authors sincerely thank the department editor Chung Piaw Teo, the associate editor and the three anonymous referees for their thoughtful comments that have helped to improve the paper. The authors also thank David Brown, Stephen Chick, Adam Elmachtoub, Paul Glasserman, Vishal Gupta, Francis Longstaff, Karthik Natarajan, Huseyin Topaloglu, Nikolaos Trichakis, Alexander Remorov and Spyros Zoumpoulis for helpful discussions, feedback and suggestions on earlier versions of this paper, and Nathan Kallus for a serendipitous suggestion that led to the development of the results in Section~\ref{subsec:epsilon_depth_bound_R2}.  %

\bibliographystyle{ormsv080}
\bibliography{IOS_literature}

\clearpage

\ECSwitch
\ECHead{Electronic companion for ``Interpretable Optimal Stopping''}

\section{Proofs and Additional Theoretical Results}
\label{appendix:proofs}

\subsection{Proofs and additional lemmas for Section \ref{subsec:problem_convergence}}
\label{sec:convergence_appendix}
In this section, we focus now on proving Theorem \ref{thm:pointwise_convergence} and Corollary \ref{cor:convergence}. Our proof relies on the fact that, via the strong law of large numbers, $J^{  \pi}(\bar \xb) \to \hat J^\pi(\bar \xb)$ almost surely for a \textit{fixed} policy $\pi$. On the other hand, Theorem \ref{thm:pointwise_convergence} requires almost sure convergence to hold simultaneously for \textit{all} policies in $\Pi_{\rm tree}(d)$, which is a possibly uncountable set. This implies that we cannot directly apply the previous fact which holds only for a fixed policy; instead, we construct a more intricate argument.

First, we deal with the issue of the uncountability of  $\Pi_{\rm tree}(d)$, where each policy $\pi$ in this space is induced by a tuple $(\mathcal{T}, \vb, \thetab, \ab)$ specifying a tree. Observe that since we have bounded the depth of the tree by $d$, the space of tree topologies, split variable indices and leaf actions $(\mathcal{T}, \vb, \ab)$ is finite, although the space of split points $\thetab$ remains uncountable. Thus, we need to prove that, fixing the tuple $(\mathcal{T}, \vb, \ab)$, we have almost sure convergence over all possible choice of split points $\thetab \in \mathcal{X}$. 

In order to make the parametrization of a tree policy $\pi$ in terms of the tree parameters explicit, we use the notation $\pi(\cdot; \mathcal{T}, \vb, \thetab, \ab)$. When unambiguous, for readability we use $\pi_{\thetab}$ to denote the tree policy $\pi(\cdot; \mathcal{T}, \vb, \thetab, \ab)$ and
$\tau_{\thetab}$ to denote $\tau_{\pi(\cdot; \mathcal{T}, \vb, \thetab, \ab)}$, the stopping time induced by this policy. 

\begin{lemma}
\label{lem:pointwise_convergence_splits}
Under Assumptions \ref{ass:tb}, \ref{ass:bounded_cost} and \ref{ass:state_smoothness}, for any fixed tree parameters  $(\mathcal{T}, \vb, \ab)$, starting state $\bar\xb$ and arbitrary $\epsilon > 0$, almost surely, there exists finite $\Omega_1 \in  \mathbb{N}^+$ such that for all $\Omega \geq \Omega_1$ and all split points $\thetab \in \mathcal{X}$,
\begin{equation*}
\left|
J^{  \pi(\cdot; \mathcal{T}, \vb, \thetab, \ab)}(\bar \xb) - \hat J^{\pi(\cdot; \mathcal{T}, \vb, \thetab, \ab)}(\bar \xb)
\right|
\leq \epsilon.
\end{equation*}
\end{lemma}

\proof{Proof.} We split the analysis into two cases.

Case 1: $\thetab$ is such that $\pi(\bar x; \mathcal{T}, \vb, \thetab, \ab) = \stop$. Then, clearly $\left|
J^{  \pi(\cdot; \mathcal{T}, \vb, \thetab, \ab)}(\bar \xb) - \hat J^{\pi(\cdot; \mathcal{T}, \vb, \thetab, \ab)}(\bar \xb)
\right| = 0$.

Case 2: $\thetab$ is such that $\pi(\bar x; \mathcal{T}, \vb, \thetab, \ab) = \go$. Now, let us consider the subset of $\Xcal$ defined as $\Xcal^\go \defeq \left\{\theta \in \Xcal \mid \pi(\bar x; \mathcal{T}, \vb, \thetab, \ab) = \go \right\}$. Since by Assumption \ref{ass:tb}, $\Xcal$ is totally bounded, it follows that $\Xcal^\go \subseteq \Xcal$ is also totally bounded. This implies that for some parameter $\delta > 0$ which we will set later, we can choose a set $Q_\delta = \left\{\thetab_1, \ldots, \thetab_K \right\}$ to be a $\delta$-cover of $\Xcal^\go$ in the $||\cdot||_\infty$ norm. That is, for any $\thetab \in \Xcal^\go$, there exists $\tilde\thetab \in Q_\delta$ such that $||\thetab - \tilde\thetab||_\infty \leq \delta$. 
For any $\thetab \in \Xcal^\go$ and some $\tilde\thetab \in Q_\delta$ with $||\thetab - \tilde\thetab||_\infty \leq \delta$, we then have
\begin{align*}
\left| J^{  \pi(\cdot; \mathcal{T}, \vb, \thetab, \ab)}(\bar \xb) - \hat J^{\pi(\cdot; \mathcal{T}, \vb, \thetab, \ab)}(\bar \xb) \right| 
&  = \left|J^{  \pi_{\thetab}}(\bar \xb) - \hat J^{\pi_{\thetab}}(\bar \xb)  \right|\\
&  \leq \underbrace{\left|J^{  \pi_{\thetab}}(\bar \xb) - J^{  \pi_{\tilde\thetab}}(\bar \xb) \right|}_{(a)} + \underbrace{\left|J^{  \pi_{\tilde\thetab}}(\bar \xb) - \hat J^{\pi_{\tilde\thetab}}(\bar \xb)  \right|}_{(b)} +
\underbrace{\left|\hat J^{\pi_{\tilde\thetab}}(\bar \xb) - \hat J^{\pi_{\thetab}}(\bar \xb)  \right|}_{(c)}.
\end{align*}
We now bound each of the three terms above. By Lemma \ref{lem:diff_true_j_theta} and using the fact that $|\theta_{s} - \tilde\theta_{s}| \leq \delta$ due to $Q_{\delta}$ being a $\delta$-cover, (a) is upper bounded by 
\[
GT \sum_{s\in \splits} f(|\theta_{s} - \tilde\theta_{s}|) \leq GT |\splits| f(\delta).
\]

For term (b), by the strong law of large numbers and using the fact that $Q_{\delta}$ is finite, there exists some $\Omega_b \in  \mathbb{N}^+$ such that for all $\Omega \geq \Omega_b$ and any $\tilde\thetab \in Q_\delta$,  
$\left|J^{  \pi_{\tilde\thetab}}(\bar \xb) - \hat J^{\pi_{\tilde\thetab}}(\bar \xb)  \right| \leq \epsilon/3$ almost surely.

Finally, for (c), by Lemma \ref{lem:diff_hat_j_theta},
\begin{align*}
\left|\hat J^{\pi_{\tilde\thetab}}(\bar \xb) - \hat J^{\pi_{\thetab}}(\bar \xb)  \right| 
& \leq
\frac{G}{\Omega} \sum_{\omega=1}^{\Omega} \sum_{t\in [T] \setminus \{1\}} \sum_{s\in \splits} \Ibb \left\{x_{v(s)}(\omega,t) \in [\min\left\{\theta_{s}, \tilde \theta_{s}\right\}, \max\left\{\theta_{s},  \tilde\theta_{s}\right\}]\right\}\\
& \leq
\frac{G}{\Omega} \sum_{\omega=1}^{\Omega} \sum_{t\in [T] \setminus \{1\}} \sum_{s\in \splits} \Ibb \left\{x_{v(s)}(\omega,t) \in [\tilde \theta_{s} - \delta, \tilde\theta_{s} + \delta]\right\}.
\end{align*}
Furthermore, since $Q_\delta$ is a finite set,  we can again invoke the strong law of large numbers to show that for some $\gamma>0$ to be chosen later, there exists some $\Omega_c$ such that for all $\Omega \geq \Omega_c$ and all $\tilde\thetab \in Q_\delta$,
\begin{align*}
\frac{G}{\Omega} \sum_{\omega=1}^{\Omega} \sum_{t\in [T]  \setminus \{1\}} \sum_{s\in \splits} \Ibb \left\{x_{v(s)}(\omega,t) \in [\tilde \theta_{s} - \delta, \tilde\theta_{s} + \delta]\right\} 
& \leq  
G \sum_{t\in [T] \setminus \{1\}} \sum_{s\in \splits} \Pr\left[x_{v(s)} \in [\tilde\theta_{s} - \delta,  \tilde\theta_{s} + \delta]\right] + \gamma. %
\end{align*}
We now define the set 
\begin{equation*}
A_s = \Xcal_1\times \ldots \times \Xcal_{v(s)-1} \times [\tilde\theta_{s} - \delta,  \tilde\theta_{s} + \delta] \times \Xcal_{v(s)+1} \times \ldots \times \Xcal_n.
\end{equation*} 
By Assumption \ref{ass:tb}, the Borel measure of $A_s$ is upper bounded by $|\tilde\theta_{s} + \delta -  (\tilde\theta_{s} - \delta)| \cdot ( \max_{\xb, \mathbf{y} \in \Xcal} ||\xb - \mathbf{y}||_\infty)^{n-1} \leq 2\delta$. Then, using Part 1 of Assumption \ref{ass:state_smoothness},
\begin{align*}
\Pr\left[x_{v(s)}(t) \in [\tilde\theta_{s} - \delta,  \tilde\theta_{s} + \delta]\right] = \Pr\left[\xb(t) \in A_s \right] \leq f(2\delta).
\end{align*}
Thus, for all $\tilde\thetab \in Q_\delta$,
\[
\left|\hat J^{\pi_{\tilde\thetab}}(\bar \xb) - \hat J^{\pi_{\thetab}}(\bar \xb)  \right|  \leq GT |\splits| f(2\delta) + \gamma.
\]
Setting $\gamma = \epsilon / 6$, and $\delta$ such that $GT |\splits| f(2\delta) \leq \epsilon / 6$, which we can do by Parts $2$ and $3$ of Assumption \ref{ass:state_smoothness}, we obtain that $GT |\splits| f(2\delta) + \gamma \leq \epsilon / 3$.

Putting everything together, we obtain that, almost surely, for any $\Omega \geq \max\left\{\Omega_b, \Omega_c \right\}$, 
\[
\left|J^{  \pi_{\thetab}}(\bar \xb) - \hat J^{\pi_{\thetab}}(\bar \xb) \right| \leq \epsilon, \ \ \textrm{for all $\thetab \in \mathcal{X}$}. \Halmos
\] 
\endproof

We now state and prove two lemmas which we have used in proving Lemma \ref{lem:pointwise_convergence_splits} above:
\begin{lemma}
\label{lem:diff_true_j_theta}
Under Assumptions \ref{ass:bounded_cost} and \ref{ass:state_smoothness}, for initial state $\bar\xb$ and any fixed tree parameters  $(\mathcal{T}, \vb, \ab)$  and $\thetab_1, \thetab_2 \in \mathcal{X}$ such that $ \pi(\bar x; \mathcal{T}, \vb, \thetab_1, \ab) =  \pi(\bar x; \mathcal{T}, \vb, \thetab_2, \ab) = \go$,
\begin{equation*}
\left|
J^{  \pi(\cdot; \mathcal{T}, \vb, \thetab_1, \ab)}(\bar \xb) - J^{  \pi(\cdot; \mathcal{T}, \vb, \thetab_2, \ab)}(\bar \xb)
\right|
\leq GT \sum_{s \in \splits} f(|\theta_{1,s} - \theta_{2,s}|).
\end{equation*}
\end{lemma}

\proof{Proof.}
We construct an upper bound on the difference in values of the two policies defined by $\thetab_1$ and $\thetab_2$. First, by Jensen's inequality, we have,
\begin{align*}
\left| J^{  \pi(\cdot; \mathcal{T}, \vb, \thetab_1, \ab)}(\bar \xb) - J^{  \pi(\cdot; \mathcal{T}, \vb, \thetab_2, \ab)}(\bar \xb)  \right| 
& = 
 \left|J^{  \pi_{\thetab_1}}(\bar \xb) - J^{  \pi_{\thetab_2}}(\bar \xb)  \right| \\
&  =
\left|\Exp\left[\beta^{\tau_{\thetab_1}-1} g(\tau_{\thetab_1}, \xb(\tau_{\thetab_1})) - \beta^{\tau_{\thetab_2}-1} g(\tau_{\thetab_2}, \xb(\tau_{\thetab_2})) \right]  \right| \\
&  \leq 
\Exp\left[\left|\beta^{\tau_{\thetab_1}-1} g(\tau_{\thetab_1}, \xb(\tau_{\thetab_1})) - \beta^{\tau_{\thetab_2}-1} g(\tau_{\thetab_2}, \xb(\tau_{\thetab_2})) \right| \right].
\end{align*}
We further bound the RHS of the last equation by conditioning on whether the stopping times of the two policies induced by 
$\thetab_1$ and $\thetab_2$ agree:
\begin{align*}
& \Exp\left[\left|\beta^{\tau_{\thetab_1}-1} g(\tau_{\thetab_1}, \xb(\tau_{\thetab_1})) - \beta^{\tau_{\thetab_2}-1} g(\tau_{\thetab_2}, \xb(\tau_{\thetab_2})) \right| \right] \\
& \quad =
\Exp\left[\left. \left|\beta^{\tau_{\thetab_1}-1} g(\tau_{\thetab_1}, \xb(\tau_{\thetab_1})) - \beta^{\tau_{\thetab_2}-1} g(\tau_{\thetab_2}, \xb(\tau_{\thetab_2})) \right|  \right| \tau_{\thetab_1} = \tau_{\thetab_2}\right]   \Pr\left[\tau_{\thetab_1} = \tau_{\thetab_2}\right]\\
& \quad
\quad + \Exp\left[\left. \left|\beta^{\tau_{\thetab_1}-1} g(\tau_{\thetab_1}, \xb(\tau_{\thetab_1})) - \beta^{\tau_{\thetab_2}-1} g(\tau_{\thetab_2}, \xb(\tau_{\thetab_2})) \right|  \right| \tau_{\thetab_1} \neq \tau_{\thetab_2}\right]   \Pr\left[\tau_{\thetab_1} \neq \tau_{\thetab_2}\right]\\
& \quad \leq
0 \cdot  \Pr\left[\tau_{\thetab_1} = \tau_{\thetab_2}\right] + G \cdot \Pr\left[\tau_{\thetab_1} \neq \tau_{\thetab_2}\right]\\
& \quad \leq 
G \cdot \Pr\left[ \exists t \in [T]\setminus \{1\} \ \st\  \pi_{\thetab_1}(t, \xb(t)) \neq \pi_{\thetab_2}(t, \xb(t))\right] \\
& \quad \leq
G \sum_{t\in [T] \setminus \{1\}} \Pr\left[\pi_{\thetab_1}(t, \xb(t)) \neq \pi_{\thetab_2}(t, \xb(t))\right] \\
& \quad \leq
G \sum_{t\in [T]\setminus \{1\}} \Pr\left[\exists s \in \splits\ \st\ x_{v(s)}(t) \in [\min\left\{\theta_{1s},  \theta_{2s}\right\}, \max\left\{\theta_{1s},  \theta_{2s}\right\}]\right] \\
& \quad \leq
G \sum_{t\in [T]\setminus \{1\}} \sum_{s\in \splits} \Pr\left[x_{v(s)}(t) \in [\min\left\{\theta_{1s},  \theta_{2s}\right\}, \max\left\{\theta_{1s},  \theta_{2s}\right\}]\right],%
\end{align*}
where in the first inequality, we have used Assumption \ref{ass:bounded_cost}, which guarantees that the difference in values of $\pi_{\theta_1}$ and $\pi_{\theta_2}$ cannot be larger than $G$; in the second, the fact that the stopping times being different implies that $\pi_{\thetab_1}$ and $\pi_{\thetab_2}$ differ in their action at some $t \in [T]$, which since $\pi(\bar x; \mathcal{T}, \vb, \thetab_1, \ab) =  \pi(\bar x; \mathcal{T}, \vb, \thetab_2, \ab) = \go$ cannot be $t=1$; in the third and fifth, a union bound over times and splits, respectively; in the fourth, the fact that $\pi_{\thetab_1}$ and $\pi_{\thetab_2}$ disagreeing at $t$ implies that there is some split at  which the policies disagree (i.e., the split variable is mapped to different child nodes in the two policies).%

Now to conclude, let us define the set 
\begin{equation*}
A_s = \Xcal_1\times \ldots \times \Xcal_{v(s)-1} \times  [\min\left\{\theta_{1s},  \theta_{2s}\right\}, \max\left\{\theta_{1s},  \theta_{2s}\right\}] \times \Xcal_{v(s)+1} \times \ldots \times \Xcal_n,
\end{equation*}
similarly to the proof of Lemma \ref{lem:pointwise_convergence_splits}. The Borel measure of $A_s$ is upper bounded by $|\theta_{1s}- \theta_{2s}| (\max_{\xb, \mathbf{y} \in \Xcal} ||\xb - \mathbf{y}||_\infty)^{n-1} \leq |\theta_{1s}- \theta_{2s}|$, where the last inequality follows from Assumption \ref{ass:tb}. %
We thus have, using Part 1 of Assumption \ref{ass:state_smoothness},
\begin{align*}
\Pr\left[x_{v(s)}(t) \in [\min\left\{\theta_{1s},  \theta_{2s}\right\}, \max\left\{\theta_{1s},  \theta_{2s}\right\}]\right] 
& = \Pr\left[\xb(t) \in A_s\right] \leq f(|\theta_{1s}- \theta_{2s}|)
\end{align*}
which completes the proof. \Halmos
\endproof

\begin{lemma}
\label{lem:diff_hat_j_theta}
Under Assumption \ref{ass:bounded_cost}, for starting state $\bar\xb$ and any fixed tree parameters  $(\mathcal{T}, \vb, \ab)$  and $\thetab_1, \thetab_2 \in \mathcal{X}$ such that $\pi(\bar x; \mathcal{T}, \vb, \thetab_1, \ab) =  \pi(\bar x; \mathcal{T}, \vb, \thetab_2, \ab) = \go$,
\begin{align*}
&\left|\hat J^{\pi(\cdot; \mathcal{T}, \vb, \thetab_1, \ab)}(\bar \xb) - \hat J^{\pi(\cdot; \mathcal{T}, \vb, \thetab_2, \ab)}(\bar \xb) \right| \\
&\quad \quad \quad \leq 
\frac{G}{\Omega} \sum_{\omega=1}^{\Omega} \sum_{t \in [T] \setminus \{1\}} \sum_{s \in \splits} 
\left[ \Ibb \left\{ x_{v(s)} \in [\min\left\{\theta_{1s},  \theta_{2s}\right\}, \max\left\{\theta_{1s},  \theta_{2s}\right\}] \right\} \right].
\end{align*}
\end{lemma}

\proof{Proof.}
Similarly to the argument for Lemma \ref{lem:diff_true_j_theta},
\begin{align*}
& \left|\hat J^{\pi(\cdot; \mathcal{T}_0, \vb_0, \thetab_1, \ab_0)}(\bar \xb) - \hat J^{\pi(\cdot; \mathcal{T}_0, \vb_0, \thetab_2, \ab_0)}(\bar \xb)  \right| \\
& \quad \quad =
\left|\frac{1}{\Omega} \sum_{\omega=1}^{\Omega}  \left[\beta^{\tau_{\thetab_1,\omega}-1} g(\tau_{\thetab_1,\omega}, \xb(\tau_{\thetab_1,\omega})) - \frac{1}{\Omega} \sum_{\omega=1}^{\Omega}  \beta^{\tau_{\thetab_2,\omega}-1} g(\tau_{\thetab_2,\omega}, \xb(\tau_{\thetab_2,\omega})) \right]  \right| \\
& \quad \quad \leq
\frac{1}{\Omega} \sum_{\omega=1}^{\Omega} \left|\beta^{\tau_{\thetab_1,\omega}-1} g(\tau_{\thetab_1,\omega}, \xb(\tau_{\thetab_1,\omega})) -  \beta^{\tau_{\thetab_2,\omega}-1} g(\tau_{\thetab_2,\omega}, \xb(\tau_{\thetab_2,\omega})) \right|\\
& \quad \quad \leq
\frac{1}{\Omega} \sum_{\omega=1}^{\Omega} G \Ibb \left\{\tau_{\thetab_1,\omega} \neq \tau_{\thetab_2,\omega} \right\}\\
& \quad \quad \leq
\frac{1}{\Omega} \sum_{\omega=1}^{\Omega} G \sum_{t\in [T]  \setminus \{1\}} \Ibb \left\{\pi_{\thetab_1}(t, \xb(\omega,t)) \neq \pi_{\thetab_2}(t, \xb(\omega, t))\right\}\\
& \quad \quad \leq
\frac{G}{\Omega} \sum_{\omega=1}^{\Omega} \sum_{t\in [T]  \setminus \{1\}} \sum_{s\in \splits} \Ibb \left\{x_{v(s)}(\omega,t) \in [\min\left\{\theta_{1s},  \theta_{2s}\right\}, \max\left\{\theta_{1s},  \theta_{2s}\right\}]\right\},
\end{align*}
where in the second inequality we used Assumption \ref{ass:bounded_cost}, and the third the assumption that $\pi(\bar x; \mathcal{T}, \vb, \thetab_1, \ab) =  \pi(\bar x; \mathcal{T}, \vb, \thetab_2, \ab) = \go$ to exclude $t=1$.\Halmos \\
\endproof

We end this section with the proofs of Theorem \ref{thm:pointwise_convergence} and Corollary \ref{cor:convergence} from Section \ref{subsec:problem_convergence}.\\

\noindent \proof{Proof of Theorem \ref{thm:pointwise_convergence}.}
Since we are restricting our analysis to policies in $\Pi_{\rm tree}(d)$, it follows that there is a finite number of tuples $(\mathcal{T}, \vb, \ab)$. For each choice of such tuple, Lemma~\ref{lem:pointwise_convergence_splits} guarantees the existence of a $\Omega(\Tcal, \vb, \ab)$ such that almost surely,
$\left| J^{\pi(\cdot; \mathcal{T}, \vb, \thetab, \ab)}(\bar \xb) - \hat J^{\pi(\cdot; \mathcal{T}, \vb, \thetab, \ab)}(\bar \xb) \right| \leq \epsilon$ for all $\thetab \in \Xcal$ and $\Omega \geq \Omega(\mathcal{T}, \vb, \ab)$. The result follows by taking $\Omega_0 = \max_{(\mathcal{T}, \vb, \ab)}  \{ \Omega(\mathcal{T}, \vb, \ab) \}$. \Halmos \\
\endproof 

 \proof{Proof of Corollary \ref{cor:convergence}.}
Consider the case that $\sup_{\pi \in \Pi_{\rm tree}(d)} J^{  \pi}(\bar \xb) \geq \sup_{\pi \in \Pi_{\rm tree}(d)} \hat J^\pi(\bar \xb)$. This is without loss of generality since the other case is symmetric. We know that there must exist some policy $\underline\pi \in \Pi_{\rm tree}(d)$ such that
\begin{equation}
\label{eq:j_underbar_pi}
J^{  \underline\pi}(\bar \xb) \geq \sup_{\pi \in \Pi_{\rm tree}(d)} J^{  \pi}(\bar \xb)  - \frac{\epsilon}{2}.
\end{equation} 
Additionally, by Theorem \ref{thm:pointwise_convergence}, almost surely there exists $\Omega_0$ such that for all $\Omega \geq \Omega_0$,  
$\left| J^{  \pi}(\bar \xb) - \hat J^\pi(\bar \xb) \right| \leq \epsilon/2$ for all policies $\pi \in \Pi_{\rm tree}(d)$ and thus, 
\begin{equation}
\label{eq:underbar_pi_convergence}
\hat J^{\underline\pi}(\bar \xb) \geq J^{  \underline\pi}(\bar \xb) - \frac{\epsilon}{2} \ \as
\end{equation}

Together, equations \eqref{eq:j_underbar_pi} and \eqref{eq:underbar_pi_convergence} imply that
\begin{align*}
\left|
\sup_{\pi \in \Pi_{\rm tree}(d)} J^{  \pi}(\bar \xb) - \sup_{\pi \in \Pi_{\rm tree}(d)} \hat J^\pi(\bar \xb)
\right|
& = 
\sup_{\pi \in \Pi_{\rm tree}(d)} J^\pi(\bar \xb) - \sup_{\pi \in \Pi_{\rm tree}(d)} \hat J^{  \pi}(\bar \xb)  \\
& \leq
\sup_{\pi \in \Pi_{\rm tree}(d)} J^\pi(\bar \xb) - \hat J^{ \underline  \pi}(\bar \xb)  \\
& \leq 
J^{  \underline\pi}(\bar \xb) - \hat J^{ \underline  \pi}(\bar \xb)  + \frac{\epsilon}{2} \\
& \leq 
\epsilon \ \as \Halmos
\end{align*}
\endproof

\subsection{An example where SAA convergence fails with unbounded depth}
\label{sec:convergence_can_fail}

In this section, we show that the SAA convergence result does not hold if the depth restriction on the trees is lifted. In particular, we will construct an example where the optimal sample-based reward converges to a quantity that in the limit is strictly higher than the reward obtained by the optimal policy. \\

\noindent \emph{Analogy to binary classification}. To intuitively understand why the convergence result does not hold when the depth of the tree can be arbitrarily large, it is helpful to think of the analogous result in the binary classification setting. Consider a binary classification problem where $\xb \in \mathbb{R}^n$ is the vector of features and $y \in \{-1, +1\}$ is the class label. Assume that the joint distribution of $(\xb, y)$ is such that $\xb$ is uniformly distributed in a box $[L,U]^n \subseteq \mathbb{R}^n$, and there exists an $\epsilon > 0$ such that for every $\bar{\xb} \in [L,U]^n$, we have that $\epsilon < P(y = +1 \mid \xb = \bar{\xb}) < 1-\epsilon$. In other words, if we fix a point $\xb \in [L,U]^n$, the label can be either +1 or -1 with positive probability (i.e., it is not deterministically either +1 or -1). In this setup, it is known that the optimal classifier $h$, which minimizes the expected 0-1 error with respect to the true probability distribution of $(\xb,y)$, is
\begin{equation}
h(\bar{\xb}) =\left\{ \begin{array}{ll} +1 & \text{if}\ P(y = +1 \mid \xb = \bar{\xb} ) > 1/2, \\ -1 & \text{otherwise}, \end{array}\right. 
\end{equation}
i.e., at any point $\bar{\xb}$, it selects the class label $\bar{y}$ with the highest conditional probability. Moreover, the expected 0-1 error of this classifier is known to be
\begin{equation}
\Exp_{(\xb,y)}[ \Ibb\{ h(\xb) \neq y\} ] = \Exp_\xb [  \min\{ P(y = +1 \mid \xb = \bar{\xb}),  P(y = -1 \mid \xb = \bar{\xb}) \} ]
\end{equation}
which by our assumption on $P(y = +1 \mid \xb = \bar{\xb})$ is bounded below by $\epsilon > 0$. Therefore, any tree of arbitrary depth, will obtain a 0-1 error of at least $\epsilon$ -- with respect to the true joint distribution of $\xb$ and $y$. 

Now let us consider a sample of points $(\xb^i, y^i)$, where $\xb^i \in [L,U]^n$ and $y^i \in \{-1,+1\}$, for $i = 1,\dots, \Omega$. We wish to find a classifier $h$ that minimizes the sample-based 0-1 error:
\begin{equation}
\frac{1}{\Omega} \sum_{i = 1}^{\Omega} \Ibb\{ h(\xb^i) \neq y^i \}. \label{eq:sample_01_error}
\end{equation}
Observe that by our assumption on the distribution of $\xb$ we will have that, almost surely, $\xb^i \neq \xb^{i'}$ for any $i \neq i'$. Thus, without any restriction on the depth of the tree, one can find a tree that places each point $(\xb^i, y^i)$ in its own leaf, and classifies that point according to its exact label (this is $y^i$); in this way, the tree is such that the sample-based error~\eqref{eq:sample_01_error} is zero. Therefore, letting $H_{trees}$ be the set of all classifiers represented as a binary tree of some arbitrary depth, it follows that 
\begin{equation*}
\min_{h \in H_{trees}} \  \frac{1}{\Omega} \sum_{i = 1}^{\Omega} \Ibb\{ h(\xb^i) \neq y^i \}  \to 0, \text{a.s.},
\end{equation*}
as $\Omega \to \infty$, whereas $\min_{h \in H_{trees}} \Exp_{(\xb,y)}[ \Ibb\{ h(\xb) \neq y\} ] \geq \epsilon$. \\

\noindent \emph{Example of divergence in optimal stopping.} We now describe our example in the optimal stopping setting. Let $T = 3$, and consider a simple stochastic system with only one state variable $x(t)$, where each $x(t)$ is independently drawn from a $\text{Uniform}(0,1)$ distribution.\footnote{For simplicity, our counterexample does not start from a deterministic state. This can be easily rectified by adding an additional $0$-th period where $x(0) = 0 \ \as$ to the original system.} We set the reward as $g(t, x(t)) = x(t)$ and set the discount factor $\beta = 1$. (This is actually the same setup that we consider in Section~6, with the exception of fixing $\beta$ to 1, which we do here to simplify some calculations.) For $t = 1,2,3$, we compute the optimal value function $J_t(x)$, where $x$ is the state at time $t$ (i.e., $x(t)$):
\begin{align*}
J_3(x) & = \max\{ x, 0\}  \\
& = x, \\
J_2(x) & = \max\{ x,  \beta \cdot \Exp_{\tilde{x}} [ J_3(\tilde{x}) \mid x] \} \\
& = \max\{ x, \Exp_{\tilde{x}} [ \tilde{x} ] \} \\
& = \max\{ x, \frac{1}{2} \} \\
J_1(x) & = \max\{ x, \beta \cdot \Exp_{\tilde{x}} [ J_2(\tilde{x}) \mid x] \} \\
& = \max\{ x, \Exp_{\tilde{x}}[ \max\{\tilde{x}, 1/2\} ] \} \\
& = \max\{ x,  \ \Exp_{\tilde{x}}[ \max\{\tilde{x}, 1/2\} \mid \tilde{x} \leq 1/2 ] \cdot P(\tilde{x} \leq 1/2) \\
& \quad \qquad \qquad +  \Exp_{\tilde{x}}[ \max\{\tilde{x}, 1/2\} \mid \tilde{x} > 1/2 ] \cdot P(\tilde{x} > 1/2) \} \\
& = \max\{ x,  \ (1/2)  P(\tilde{x} \leq 1/2) + \Exp_{\tilde{x}}[ \tilde{x}  \mid \tilde{x} > 1/2 ]  P(\tilde{x} > 1/2) \} \\
& = \max\{ x, (1/2) (1/2) + (3/4) (1/2) \} \\
& = \max\{ x, 5/8\} 
\end{align*}
Since the initial state $x(1)$ is random, we compute the expected value of $J_1(x)$ over this initial state:
\begin{align*}
\Exp_x [J_1(x)] & = \Exp_x [ \max\{ x, 5/8\} ] \\
& = \Exp_x [ \max\{ x, 5/8\} \mid x \geq 5/8] \cdot P(x \geq 5/8) + \Exp_x [ \max\{ x, 5/8\} \mid x < 5/8] \cdot P(x < 5/8) \\
& = \Exp_x [ x \mid x \geq 5/8] \cdot P(x \geq 5/8) + (5/8) \cdot P(x < 5/8) \\
& = (13/16) \cdot (3/8) + (5/8) \cdot (5/8) \\
& = 89 / 128 \\
& = 0.6953125
\end{align*}
In other words, the maximum reward of any policy is approximately 0.6953. Therefore, we know that 
\begin{equation*}
\max_{\pi \in \Pi_{tree}} J^{\pi} \leq 0.6953125,
\end{equation*}
where $\Pi_{tree}$ is the set of all policies that are allowed to split on the state $x$ and on the period $t$. 

Now, we turn our attention to the objective value of the SAA problem, $\max_{\pi \in \Pi_{tree}} \hat{J}^{\pi}$. Suppose that we have $\Omega$ trajectories of this stochastic system. Thus, for each $\omega \in [\Omega]$ and $t \in \{1,2,3\}$, $x(\omega, t)$ is drawn independently from the $\text{Uniform}(0,1)$ distribution. Given some $\Omega \in \mathbb{N}$, we have that almost surely, the state for each trajectory-period pair will be different -- that is, $x(\omega, t) \neq x(\omega', t')$, for any $(\omega, t) \neq (\omega',t')$. Because of this, SAA can construct a tree policy $\pi$ by splitting on both $x$ and $t$ such that each trajectory stops at the period in which $x(t)$ is highest, i.e., the stopping time  $\tau_{\pi,\omega}$ satisfies $\tau_{\pi,\omega} = \arg \max_{1 \leq t' \leq 3} x(\omega, t')$. Intuitively, this allows the SAA optimization to create an (anticipatory) policy that perfectly matches what a clairvoyant will do. It therefore follows that the sample-based reward of this policy is 
\begin{equation*}
\frac{1}{\Omega} \sum_{\omega = 1}^{\Omega} \max\{ x(\omega,1), x(\omega,2), x(\omega,3) \}.
\end{equation*}
Since each random variable $Z_\omega = \max\{ x(\omega,1), x(\omega,2), x(\omega,3) \}$ is independent and identically distributed, we can invoke the strong law of large numbers to assert that, almost surely,  
\begin{equation*}
\frac{1}{\Omega} \sum_{\omega = 1}^{\Omega} Z_{\omega} \to \Exp[Z]
\end{equation*}
as $\Omega \to \infty$, where $Z = \max\{ X_1, X_2, X_3\}$ and $X_1, X_2, X_3$ are independent $\text{Uniform}(0,1)$ random variables. By definition, $Z$ follows a $\text{Beta}(3,1)$ distribution. Therefore, $E[Z] = 3/(3 + 1) = 3/4 = 0.75$. We emphasize that this quantity corresponds to the sample-based optimization problem $\max_{\pi \in \Pi_{tree}} \hat{J}^{\pi}$; thus, while $\Pi_{tree}$ contains non-anticipatory policies, the sample-based nature of the optimization problem and the lack of depth restriction on the tree policy allows SAA to do as well as clairvoyant policies.

We observe that the value of $0.75$ is higher than the bound of $89/128  = 0.6953125$ that we obtained before; thus, in this example, the SAA convergence result does not hold if the depth restriction is lifted.

\subsection{Proofs and additional results for Section \ref{subsec:problem_complexity}}

\subsubsection{Proof of Proposition~\ref{proposition:leaf_decision_NPHard}}

\label{proof:leaf_decision_NPHard}

We will prove this by showing that the minimum vertex cover problem reduces to the leaf action SAA problem~\eqref{prob:optimal_stopping_SAA_leafdecision}. The minimum vertex cover problem is stated as follows:
\begin{quote}
\textsc{Minimum Vertex Cover}: Given a graph $(V, E)$, where $V$ is the set of nodes and $E$ is the set of edges, find the smallest subset of nodes $S \subseteq V$ such that every edge is covered, i.e., for all $e \in E$, there exists a node $i$ in $S$ such that $e$ is incident to $i$. 
\end{quote}
Given a graph $(V,E)$, we will now define an instance of the leaf action SAA problem. To do this, we need to specify a sample of trajectories, a tree topology, split variable indices and split points. We start by defining the trajectories. We assume that the nodes in $V$ are indexed from 1 to $|V|$, and the edges in $E$ are indexed from 1 to $|E|$. \\

\noindent \emph{Trajectories}. We will assume that the ambient space of our stochastic process is $\Xcal = \mathbb{R}^{|V|+1}$. We assume that we have $\Omega = |V| + |E|$ trajectories, each with $T = 3$ periods. We will consider two types of trajectories:
\begin{enumerate}
\item \emph{Edge trajectories}: These trajectories are indexed by $\omega = 1,\dots, |E|$, and each such trajectory corresponds to an edge in the graph. For each edge $e$, let $e_1$ and $e_2$ be the two vertices to which edge $e$ is incident. For each trajectory $\omega$ that corresponds to an edge $e \in E$, we assume that it takes the following values:
\begin{align*}
x_i(\omega,1) & = 0, \quad \forall\ i \in \{1,\dots, |V|+1\}, \\
x_i(\omega,2) & = \Ibb\{ i = e_1 \}, \quad \forall\ i \in \{1,\dots, |V|+1\},  \\
x_i(\omega,3) & = \Ibb\{ i = e_2 \}, \quad \forall\ i \in \{1,\dots, |V|+1\}. 
\end{align*}
In words, a trajectory for an edge $e$ is zero at all coordinates, except at the vertices to which $e$ is incident, for which the corresponding coordinates tick up to 1 at periods $t = 2$ and $t = 3$. 
\item \emph{Vertex trajectories}: These trajectories are indexed by $\omega = |E|+1,\dots, |E| + |V|$, and each such trajectory corresponds to a vertex in the graph. For each such trajectory $\omega$ that corresponds to a vertex $v \in V$, we assume that it takes the following values:
\begin{align*}
x_i(\omega,1) & = \Ibb\{ i = v\}, \quad \forall\ i \in \{1,\dots, |V|+1\}, \\
x_i(\omega,2) & = 0, \quad \forall\ i \in \{1,\dots, |V|+1 \}, \\
x_i(\omega,3) & = \Ibb\{ i = |V|+1 \}, \quad \forall \ i \in \{1,\dots, |V|+1 \}.
\end{align*}
In words, this trajectory is zero in all coordinates, except for coordinates $v$ and $|V|+1$: for coordinate $v$, it starts at 1 at $t = 1$, and comes down to 0 at $t = 2$, and for coordinate $|V|+1$, it starts at 0 at $t = 1$ and ticks up to 1 at $t = 3$. 
\end{enumerate}
We define the payoff function $g$ as follows:
\begin{equation*}
g(t,\xb) = \left\{  \begin{array}{ll} 1 & \text{if}\ t \in \{2,3\} \ \text{and} \  x_i = 1\ \text{for any}\ i \in \{1,\dots,|V|\}, \\
1/(|V| + 1) & \text{if}\ t = 3 \ \text{and} \ x_{|V|+1} = 1, \\
0 & \text{otherwise.} 
\end{array}  \right.
\end{equation*}
To understand the payoff function, we receive a reward of 1 if we stop at $t = 2$ or $t = 3$ in a state $\xb$ where any coordinate $i \in \{1,\dots, |V|\}$ is one; we receive a reward of $1/(|V| + 1)$ if we stop at $t = 3$ in a state $\xb$ where coordinate $|V|+1$ is one; and for any other state-time pair, we receive zero. For each edge trajectory $\omega = 1,\dots, |E|$, this means that if we stop at $t = 2$ or $t = 3$ of any edge trajectory $\omega = 1,\dots, |E|$, we will get a reward of 1; otherwise, if we stop at $t = 1$, we will get a reward of 0. For each vertex trajectory $\omega = |E|+1, \dots, |E| + |V|$, stopping at $t = 1$ or $t = 2$ gives us a reward of zero, and stopping at $t = 3$ gives us a reward of $1 / (|V| + 1)$. \\

\begin{figure}
\caption{Visualization of topology of tree policy for proof of Proposition~\ref{proposition:leaf_decision_NPHard}.}
\begin{subfigure}[t]{0.4\textwidth}
                \centering
                \includegraphics[width = \textwidth]{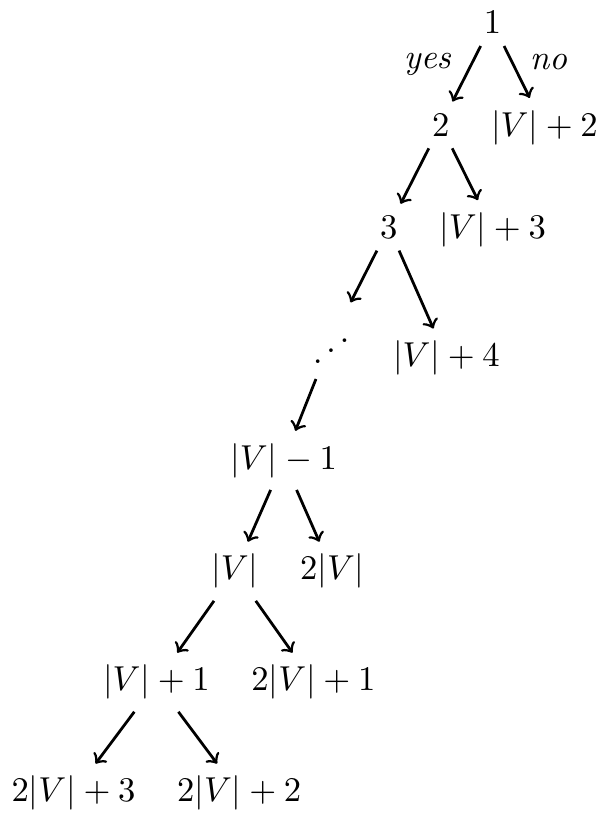}
                \caption{Topology with split node and leaf node indices for Proposition~\ref{proposition:leaf_decision_NPHard} proof.}
                \label{figure:nphardness_topology_only}
        \end{subfigure}  \hfill
\begin{subfigure}[t]{0.45\textwidth}
                \centering
                \includegraphics[width = \textwidth]{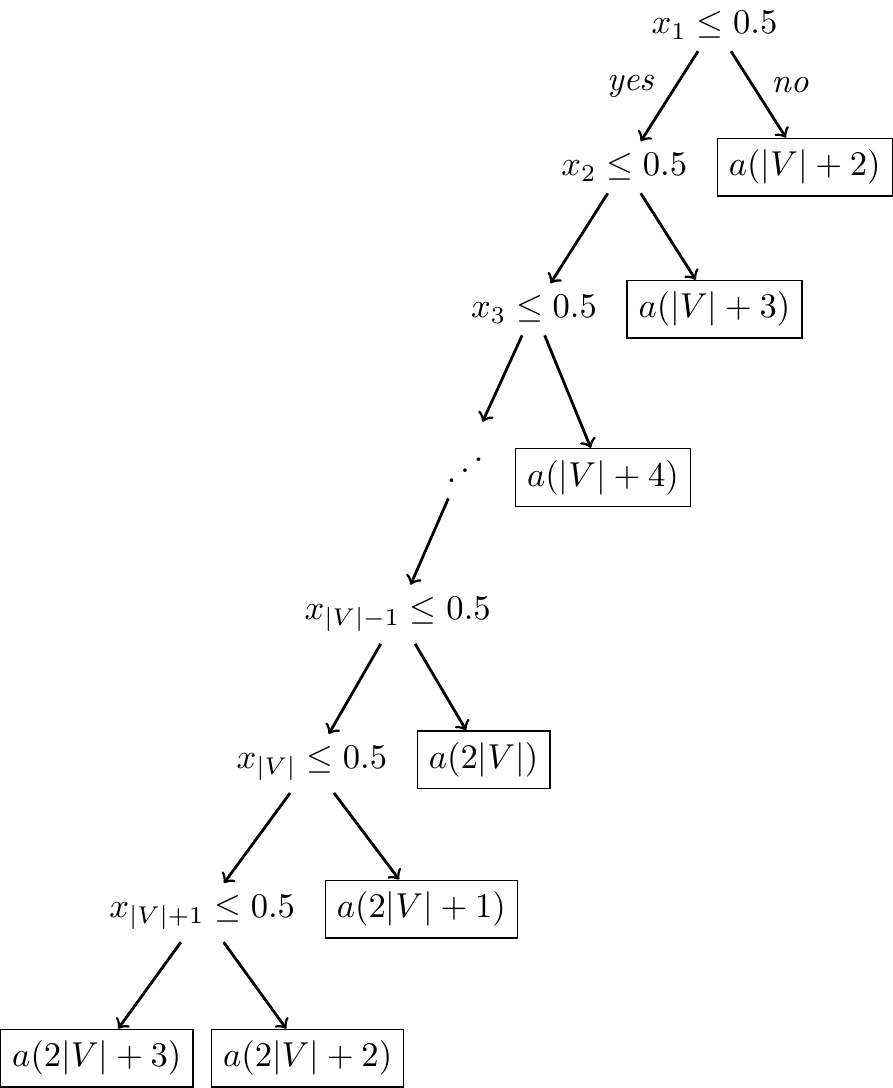} 
                \caption{Topology with split variables and split points for Proposition~\ref{proposition:leaf_decision_NPHard} proof.}
                \label{figure:nphardness_splitVarsPoints}
        \end{subfigure}

\end{figure}

\noindent \emph{Tree topology}. Our tree consists of $|V|+1$ splits and $|V|+2$ leaves. We number the split nodes from 1 to $|V|+1$, and we number the leaf nodes from $|V|+2$ to $2|V| + 3$. We define the left and right children of the split nodes as follows:
\begin{align*}
& \leftchild(s) = s+1, \quad s = 1,\dots, |V|, \\
& \rightchild(s) = |V|+1+s, \quad s = 1,\dots, |V|, \\
& \leftchild(|V|+1) = 2|V| + 3, \\
& \rightchild(|V|+1) = 2|V| + 2.
\end{align*}
Figure~\ref{figure:nphardness_topology_only} visualizes this topology. \\

\noindent \emph{Split variable indices and split points}. For each split $s = 1,\dots, |V|+1$, we define the variable $v(s) = s$, and we define the split point $\theta(s) = 0.5$. Figure~\ref{figure:nphardness_splitVarsPoints} shows the tree topology together with the split variable indices and split points. \\

\noindent Note that with this tree topology, split points and split variable indices, the policy will behave as follows:
\begin{enumerate}
\item For an edge trajectory corresponding to an edge $e$:
\begin{itemize}
\item At $t = 1$, the policy will take action $a(2|V|+3)$;
\item At $t = 2$, the policy will take action $a(|V| + e_1)$; and 
\item At $t = 3$, the policy will take action $a(|V| + e_2)$.
\end{itemize}
\item For a vertex trajectory corresponding to a vertex $v$:
\begin{itemize}
\item At $t = 1$, the policy will take action $a(|V| + 1 + v)$;
\item At $t = 2$, the policy will take action $a(2|V| + 3)$; and
\item At $t = 3$, the policy will take action $a(2|V| + 2)$. 
\end{itemize}
\end{enumerate}

\noindent \emph{Reduction of vertex cover to leaf action problem.} Consider now the problem of deciding the leaf actions. We will show that the optimal solution of this problem provides an optimal solution of the vertex cover problem. Without loss of generality, and for ease of exposition, we omit the $(1/\Omega)$ factor in the objective function of the leaf action SAA problem~\eqref{prob:optimal_stopping_SAA_leafdecision}.

Let $\ab^*$ be an optimal collection of leaf actions (i.e., an optimal solution of the leaf action SAA problem). Observe that the leaf actions for $\ell = |V|+2, \dots, 2|V|+1$ can be interpreted as a collection of vertices: if the leaf action is $\stop$, this indicates that the vertex is in the collection, whereas if the action is $\go$, this indicates that the vertex is not in the collection. We therefore specify the set of vertices as $S = \{ v \in V \ | \ a^*(|V|+1+v) = \stop \}$. 

We consider two cases:\\

\noindent \textbf{Case 1}: $E = \emptyset$, i.e., there are no edges. In this degenerate case, an upper bound on the reward we can achieve is given by $|V|/(|V| + 1)$. (This is because the highest possible reward for each of the vertex trajectories is $1/(|V|+1)$, and there are $|V|$ such trajectories; recall that there are no edge trajectories.) This upper bound can only be attained by setting $a^*(|V|+1+v) = \go$ for $v = 1,\dots, |V|$, setting $a^*(2|V|+2) = \stop$ and setting $a^*(2|V|+3) = \go$. Note that in this case, the corresponding set of vertices $S$ is exactly the empty set, which is exactly the optimal solution of the vertex cover problem in this case. \\

\noindent \textbf{Case 2}: $E \neq \emptyset$, i.e., there is at least one edge.

We first argue that the set $S$ must be a bona fide cover. To see this, we argue by contradiction. Suppose that $S$ is not a cover of $E$; then there exists an edge $e$ such that the corresponding vertices $i$ and $j$ are not in the cover. This would mean that both $a^*(|V|+1+i)$ and $a^*(|V|+1+j)$ are set to $\go$. This, in turn means that the reward of 1 of the edge trajectory $e$ cannot be earned by the tree policy. (To obtain the reward, we must stop at $t = 2$ or $t = 3$ in this trajectory; given the trajectory and the topology, split variable indices and split points, this can only happen if at least one of the two actions is $\stop$.) This means that the reward that the policy can obtain is upper bounded by $|E|-1 + |V|/(|V|+1)$. However, if we set $\ab$ so that the corresponding vertex set $S$ of $\ab$ covers $E$, and we additionally set $a(2|V|+2) = \stop$ and $a(2|V|+3) = \go$, our reward will be at least $|E|$. Therefore, it must be that $S$ covers all edges.

We now argue that $S$ is a minimal cover. To see this, we again argue by contradiction. Suppose that there is a cover $S'$ such that $|S'| < |S|$. Consider the set of leaf actions $\ab'$, which is defined as follows:
\begin{align*}
a'(|V|+1+v) & = \left\{ \begin{array}{ll} \stop & \text{if}\ v \in S', \\ \go & \text{if}\ v \notin S', \end{array} \right.  \\
a'(2|V|+2) & = \stop, \\
a'(2|V|+3) & = \go.
\end{align*}
The corresponding tree policy with $\ab'$ has the following behavior. For all edge trajectories, it stops at $t = 2$ or $t = 3$, because $S'$ is a cover; as a result, it accrues a reward of 1 from all edge trajectories. For those vertex trajectories corresponding to vertices in $S'$, it stops at $t = 1$, accruing a reward of zero. For those vertex trajectories corresponding to vertices \emph{not} in $S'$, it stops at $t = 3$, accruing a reward of $1/(|V|+1)$. Therefore, the total reward is $|E| + (|V| - |S'|)/(|V|+1)$. 

Now, we show that for $\ab^*$, the reward is upper bounded by $|E| + (|V| - |S|)/(|V|+1)$. Observe that if this is the case, we immediately have our contradiction because the reward $|E| + (|V| - |S|)/(|V|+1)$ will be smaller than $|E| + (|V| - |S'|)/(|V|+1)$ (this follows from $|S| > |S'|$), which will contradict the fact that $\ab^*$ is an optimal solution of the leaf action SAA problem. To see this, observe that for each edge trajectory, since $S$ is a cover, the policy can garner at most a reward of 1, resulting in a total reward of $|E|$ from the edge trajectories. For the vertex trajectories, observe that for all vertices $v \in S$, the reward for the corresponding vertex trajectory must be zero, because the tree policy must stop at $t = 1$ (the reward for any vertex trajectory at $t = 1$ is zero). For each vertex trajectory corresponding to a vertex $v \notin S$, the most reward that can be garnered is $1/(|V|+1)$ (obtained by stopping at $t = 3$). Therefore, an upper bound for the reward is 
\begin{align*}
& |E| + |V \setminus S| \,  / \, (|V| +1)  \\
= & |E| + (|V| - |S|)/(|V|+1).
\end{align*}
This gives us our contradiction, and establishes that $S$ is a minimal cover.

We have thus shown that any instance $(V,E)$ of the minimum vertex cover problem can be transformed to an instance of the leaf action problem~\eqref{prob:optimal_stopping_SAA_leafdecision}; moreover, the instance is polynomially-sized in terms of $V$ and $E$. Since the minimum vertex cover problem is known to be NP-Complete \citep{garey2002computers}, it follows that the leaf action problem is NP-Hard. \Halmos

\subsubsection{Solving the leaf action SAA problem for fixed $\Omega$}
\label{appendix:leaf_action_problem_fixed_Omega}

For each trajectory, let $\ell_{\omega,1}, \dots, \ell_{\omega, m_{\omega}}$ be the sequence of leaves that the trajectory reaches, in the order that they are reached. Let $t_{\omega,1}, \dots, t_{\omega,m_{\omega}}$ be the times at which those leaves are reached. We therefore have that 
\begin{equation}
t_{\omega,1} < t_{\omega,2} < \dots < t_{\omega, m_{\omega}} \leq T.
\end{equation}

Observe that no matter how we specify the leaf action vector $\ab$, in each trajectory the policy will stop at some time $t_{\omega,i}$, or it will never stop. We can thus encode the time at which the policy stops in a trajectory by the integer $i \in [ m_{\omega} + 1] $, where values between 1 and $m_{\omega}$ inclusive correspond to an actual time in $[T]$, while the integer $m_{\omega}+1$ corresponds to never stopping. \\

Before we can formulate the leaf action problem, we require some additional definitions. First, for each trajectory $\omega$ and each $i \in [m_{\omega}+1]$, let us define $\Lcal_{\stop}(\omega,i)$ as 
\begin{equation*}
\Lcal_{\stop}(\omega,i) = \left\{ \begin{array}{ll} \{ \ell_{\omega,i} \}, & \text{if}\ i \in [m_{\omega}], \\
\emptyset, & \text{if}\ i = m_{\omega}+1. \end{array} \right.
\end{equation*}
In words, this is the set of leaves $\ell$ for which $a(\ell)$ must be set to $\stop$ in order for the trajectory to stop at $t_{\omega,i}$. Similarly, we define the set $\Lcal_{\go}(\omega,i)$ as 
\begin{equation*}
\Lcal_{\go}(\omega,i) = \{ \ell_{\omega,1}, \dots, \ell_{\omega,i-1} \}.
\end{equation*}
In words, $\Lcal_{\go}(\omega,i)$ is the set of leaves that must be set to $\go$ in order for the trajectory to stop at $t_{\omega,i}$; these leaves are precisely those leaves that are reached before $\ell_{\omega,i}$. Lastly, we define $J(\ib)$ as 
\begin{equation}
J(\ib) = \frac{1}{\Omega} \sum_{\omega = 1}^{\Omega} \sum_{i=1}^{m_{\omega}} \Ibb\{i_{\omega} = i\} \cdot \beta^{t_{\omega,i_{\omega}}-1} \cdot g(t_{\omega,i_{\omega}}, \xb(\omega, t_{\omega,i_{\omega}})),
\end{equation}
i.e., it is the reward from stopping at the times indexed by $\ib$. 
With these definitions, we now claim that one can solve the leaf action SAA problem by solving the following problem:
\begin{subequations}
\begin{alignat}{2}
& \underset{\ib \in \prod_{\omega = 1}^{\Omega} [m_{\omega}+1 ] }{ \text{maximize} } & \quad & J(\ib) \\
& \text{subject to} & &  \left(  \bigcup_{\omega = 1}^{\Omega} \Lcal_{\stop}(\omega, i_{\omega})  \right) \cap \left(  \bigcup_{\omega = 1}^{\Omega} \Lcal_{\go}(\omega, i_{\omega}) \right) = \emptyset. \label{prob:index_optimization_achievable}
\end{alignat}
\label{prob:index_optimization}%
\end{subequations}
The optimization problem involves searching over all vectors $\ib$ of indices that indicate which leaf each trajectory will stop at. The constraint~\eqref{prob:index_optimization_achievable} requires that the vector $\ib$ of indices is actually achievable. Specifically, the set $\bigcup_{\omega = 1}^{\Omega} \Lcal_{\stop}(\omega, i_{\omega})$ is the set of leaves for which $a(\ell)$ must be set to $\stop$, while the set $\bigcup_{\omega = 1}^{\Omega} \Lcal_{\go}(\omega, i_{\omega})$ is the set of leaves for which $a(\ell)$ must be set to $\go$. If the two sets intersect, then there is a leaf $\ell$ for which $a(\ell)$ must be set to both $\stop$ and $\go$, which is impossible. After solving the optimization problem, we can obtain an optimal $\ab$ as follows:
\begin{enumerate}
\item For all $\ell \in \bigcup_{\omega = 1}^{\Omega} \Lcal_{\stop}(\omega, i_{\omega})$, set $a(\ell) = \stop$;
\item For all $\ell \in \bigcup_{\omega = 1}^{\Omega} \Lcal_{\go}(\omega, i_{\omega})$, set $a(\ell) = \go$; and
\item For all other leaves $\ell$, set $a(\ell)$ to either $\stop$ or $\go$. 
\end{enumerate}

Now, with regard to the complexity of solving problem~\eqref{prob:index_optimization}, we make the following observations:
\begin{itemize}
\item The set of index vectors $\prod_{\omega = 1}^{\Omega} [m_{\omega}+1]$ is of size $O( T^{\Omega} )$;
\item For each $\omega$, one can compute the sets $\Lcal_{\stop}(\omega,i)$ and $\Lcal_{\go}(\omega,i)$ with $O(T)$ computations; thus, constructing these sets for all $\Omega$ trajectories amounts to $O(\Omega T)$ computations; and
\item For each index vector $\ib \in \prod_{\omega = 1}^{\Omega} [m_{\omega}+1]$, checking whether the intersection $\left(  \bigcup_{\omega = 1}^{\Omega} \Lcal_{\stop}(\omega, i_{\omega})  \right) \cap \left(  \bigcup_{\omega = 1}^{\Omega} \Lcal_{\go}(\omega, i_{\omega}) \right)$ is empty requires $O(\Omega |\leaves|)$ computations. (Each set $\Lcal_{\go}(\cdot,\cdot)$ and $\Lcal_{\stop}(\cdot, \cdot)$ can be viewed as a binary vector indicating whether leaf $\ell$ is in the set or not.) 
\end{itemize}
Thus, the overall time complexity is $O(\Omega T + T^{\Omega} \cdot \Omega |\leaves|) = O(T^{\Omega} \cdot \Omega |\leaves|)$. If $\Omega$ is fixed, then this quantity scales polynomially in the problem parameters, and one can solve the leaf action SAA problem in polynomial time. 

\subsubsection{NP-Hardness of leaf action SAA problem when $n = 1$}
\label{appendix:leaf_action_problem_NPHard_neq1}

Consider the leaf action problem~\eqref{prob:optimal_stopping_SAA_leafdecision} when $n$ is restricted to 1. 

\begin{proposition}
The leaf action SAA problem~\eqref{prob:optimal_stopping_SAA_leafdecision} with $n = 1$ is NP-Hard. 
\end{proposition}

\proof{Proof.}
We again prove this result by reducing the minimum vertex cover problem to the leaf action SAA problem with $n = 1$. Consider a vertex cover instance with a set of vertices $V$, indexed from 1 to $|V|$, and a set of edges $E$, indexed from 1 to $|E|$. We set the time horizon to $T = 3$. For convenience, since there is only one state variable, we will drop the subscript $v$ and refer to this state variable simply by $x$, i.e., $x(\omega,t)$ is the value of the state variable in trajectory $\omega$ and at time $t$. We also set the discount factor $\beta = 1$.  

Consider a set of $\Omega = |V| + |E|$ trajectories, defined as follows:
\begin{itemize}
\item \emph{Edge trajectories}: for each edge $e = 1,\dots, |E|$, which is a pair of nodes $e = (e_1, e_2) \in V \times V$, we define the trajectory $\omega = e$ as
\begin{align}
x(e,1) & = 0, \\
x(e,2) & = e_1 + 1/2, \\
x(e,3) & = e_2 + 1/2.
\end{align}
\item \emph{Vertex trajectories}: for each vertex $v = 1,\dots, |V|$, we define the trajectory $\omega = |E| + v$ as 
\begin{align}
x(|E| + v, 1) & = v + 1/2, \\
x(|E| + v, 2) & = 0, \\
x(|E| + v, 3) & = 1/2. 
\end{align}
\end{itemize}
We define the reward function as:
\begin{equation*}
g(t, x) = \left\{ \begin{array}{ll} 1 & \text{if}\ x > 1 \ \text{and}\ t \in \{2,3\}, \\
1/(|V| + 1) & \text{if} \ x \in (0,1) \ \text{and}\ t = 3, \\
0 & \text{otherwise}. 
\end{array} \right.
\end{equation*}

\begin{figure}
\centering 
\begin{tabular}{ccc}
$\omega = e$ for $e \in E$: & &$\omega = |E| + v$ for $v \in V$: \\[0.5em]
\includegraphics[width=0.3\textwidth]{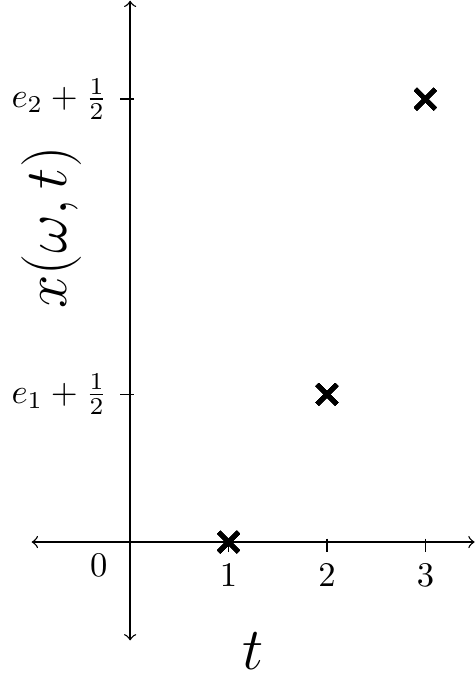} & \quad & \hspace{2em} \includegraphics[width=0.3\textwidth]{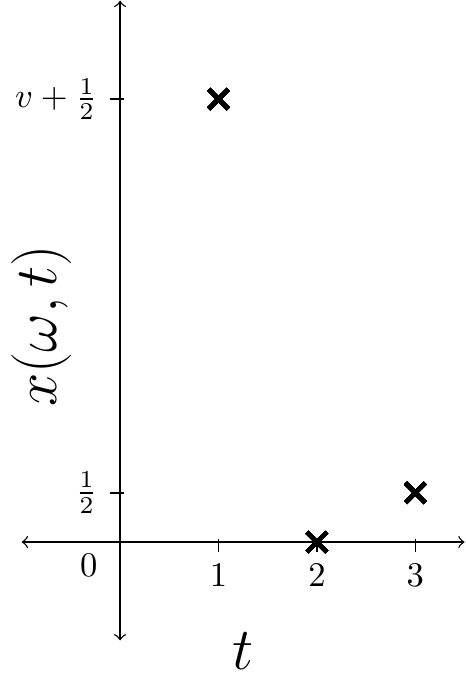} \\[1em]
\includegraphics[width=0.3\textwidth]{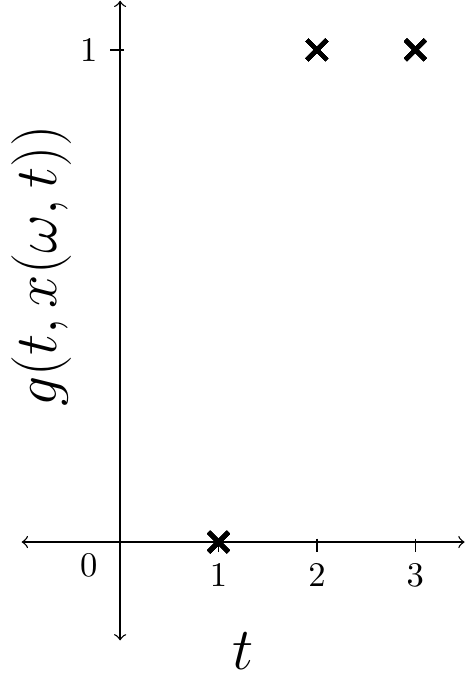} & \quad &\includegraphics[width=0.35\textwidth]{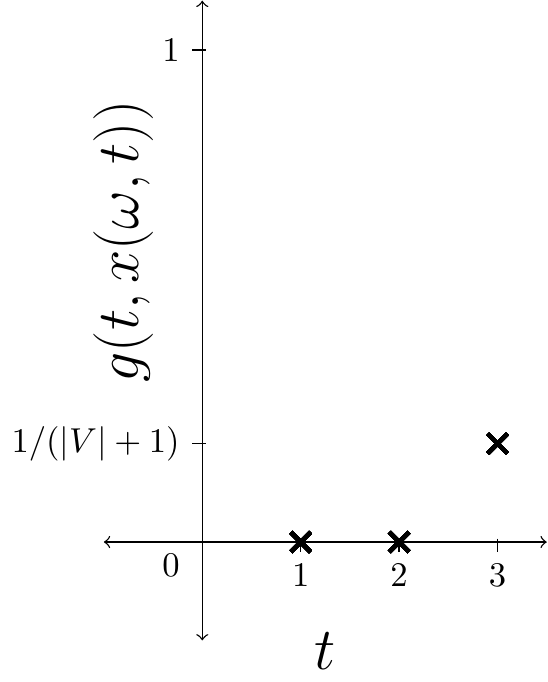} \\
\end{tabular}

\caption{Visualization of state variable $x(\omega, t)$ (top subfigures) and reward $g(t, x(\omega,t))$ (bottom subfigures) for edge trajectories ($\omega = e$ for $e \in E$; left subfigures) and vertex trajectories ($\omega = |E| + v$ for $v \in V$; right subfigures) for $n = 1$ leaf action SAA reduction. \label{figure:leaf_action_SAA_neq1_trajectories}}
\end{figure}

The two types of trajectories and their rewards are plotted in Figure~\ref{figure:leaf_action_SAA_neq1_trajectories}.

\begin{figure}
\centering
\includegraphics[width=0.4\textwidth]{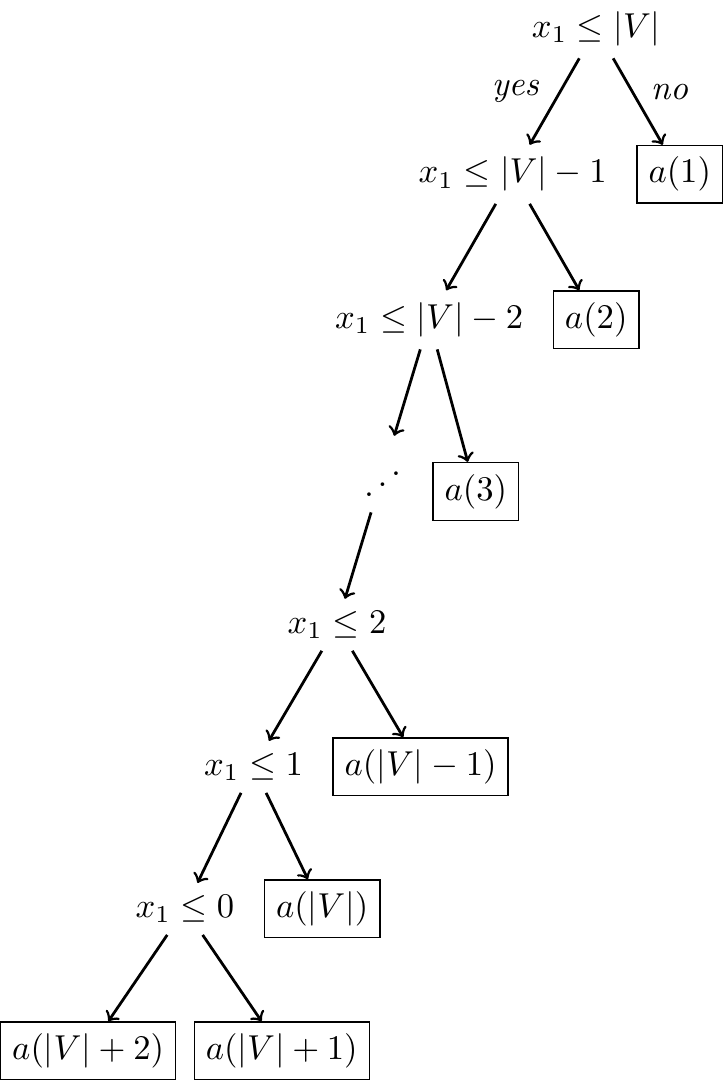}

\caption{Tree topology and split points for $n = 1$ leaf action SAA problem reduction. \label{figure:leaf_action_SAA_neq1_tree}}
\end{figure}

The tree policy we will consider is defined by the topology $\Tcal$ and split points $\thetab$ as shown in Figure~\ref{figure:leaf_action_SAA_neq1_tree}; the leaves are indexed from 1 to $|V|+2$, in the manner shown in the figure.

The leaf action SAA problem~\eqref{prob:optimal_stopping_SAA_leafdecision} amounts to deciding the actions $a(1), a(2), \dots, a(|V| +2)$. For ease of exposition, in the reduction that follows, we drop the $1/\Omega$ factor in the objective function of the leaf SAA problem.

Let $\ab$ be the optimal solution of the problem. At optimality, we note that we will have $a(|V| + 2) = \go$, because any solution in which $a(|V| + 2) = \stop$ will stop at $t = 1$ in all edge trajectories and garner a reward of at most $|V|/(|V|+1)$ , whereas a policy that sets $a(v) = \stop$ for all $v = 1,\dots, |V|+1$ and sets $a(|V| + 2) = \go$ will garner a reward of at least $|E|$ which is strictly greater than $|V|/(|V|+1)$, implying that any optimal policy must therefore achieve a reward that is strictly greater than $|V|/(|V|+1)$. 

The rest of the proof follows along the same lines as our original proof of the NP-Hardness of the leaf action SAA problem. For completeness, we provide the details here. Given the optimal solution $\ab$, consider the set of vertices $S$ defined as
\begin{equation*}
S = \left\{ v \in V \mid a(v) = \stop \right\}.
\end{equation*}
We now argue that $S$ is both a feasible cover, and is minimal.

\emph{Feasible cover}. To see why $S$ is feasible, suppose that $S$ is not feasible. If this is the case, there must exist an edge $e = (e_1, e_2)$, with $e_1, e_2 \in V$, such that $e_1 \notin S$ and $e_2 \notin S$, implying that $a(e_1) = a(e_2) = \go$. In this case, it must be that the policy defined by $\ab$ does not stop at period 2 or period 3 of trajectory $e$, which means that the reward of $\ab$ is at most 
\begin{equation}
(|E| - 1) + |V| / (|V| + 1), \label{eq:suboptimal_reward}
\end{equation}
where the first term represents the most reward that can be garnered from the edge trajectories, while the second term represents the most reward that can be garnered from the vertex trajectories. However, observe that just by setting $a(v) = \stop$ for $v = 1,\dots, |V| + 1$ and $a(|V| +2) = \go$, the policy will stop at period 2 or 3 in all edge trajectories, and garner a reward of at least $|E|$, which is strictly greater than \eqref{eq:suboptimal_reward}, contradicting the fact that $\ab$ is optimal for the leaf action SAA problem. Therefore, $S$ must be a feasible cover. 

\emph{Minimal cover}. To see why $S$ must be minimal, suppose that $S$ is not minimal, in which case there exists a $S'$ that is a feasible cover of $V$ with $|S'| < |S|$. In this case, define a solution $\ab'$  to the leaf action SAA problem such that 
\begin{equation}
a'(v) = \left\{ \begin{array}{ll} \stop & \text{if}\ v \in S', \\
\go & \text{if} \ v \in V \setminus S', \\
\stop & \text{if}\ v = |V| + 1,\\
\go & \text{if}\ v = |V| + 2.
\end{array} \right.
\end{equation}
By the construction of this new leaf action vector $\ab'$, it is clear that the corresponding tree policy will stop at period 2 or 3 for each edge trajectory $\omega \in \{1,\dots, |E|\}$, at period 1 for every trajectory $\omega = |E| + v$ for each $v \in S'$, and at period 3 for every trajectory $\omega = |E| + v$ for each $v \in V \setminus S'$. As a result, the reward of this policy is $|E| + (|V| - |S'|)/(|V| + 1)$. However, the reward of $\ab$ is at most $|E| + (|V| - |S|)/(|V| + 1)$. Since $|S'| < |S|$, it follows that 
\begin{equation}
|E| + (|V| - |S'|)/(|V| + 1) < |E| + (|V| - |S|)/(|V| + 1),
\end{equation}
which implies that $\ab$ is not optimal, and thus gives rise to a contradiction. It therefore must be that $S$ is a minimal cover.

Thus, we have shown that if one can solve the $n = 1$ leaf action SAA problem to optimality, then one can solve the minimum vertex cover problem. Since the instance of the leaf action SAA problem that we have constructed is polynomially sized in $|V|$ and $|E|$, and since the minimum vertex cover problem is NP-Hard, it follows that the leaf action SAA problem with $n = 1$ is also NP-Hard. $\square$
\endproof

\subsubsection{Proof of Proposition~\ref{proposition:split_variable_NPHard}}
\label{proof:split_variable_NPHard}

To prove that the split variable index SAA problem is NP-Hard, we will again show that the minimum vertex cover problem reduces to it. We will prove the result by considering an instance $(V, E)$ of the minimum vertex cover problem, and considering the same set of trajectories, payoff function and tree topology from the proof of Proposition~\ref{proposition:leaf_decision_NPHard}. 

\begin{figure}
  \centering
                \includegraphics[width = 0.4\textwidth]{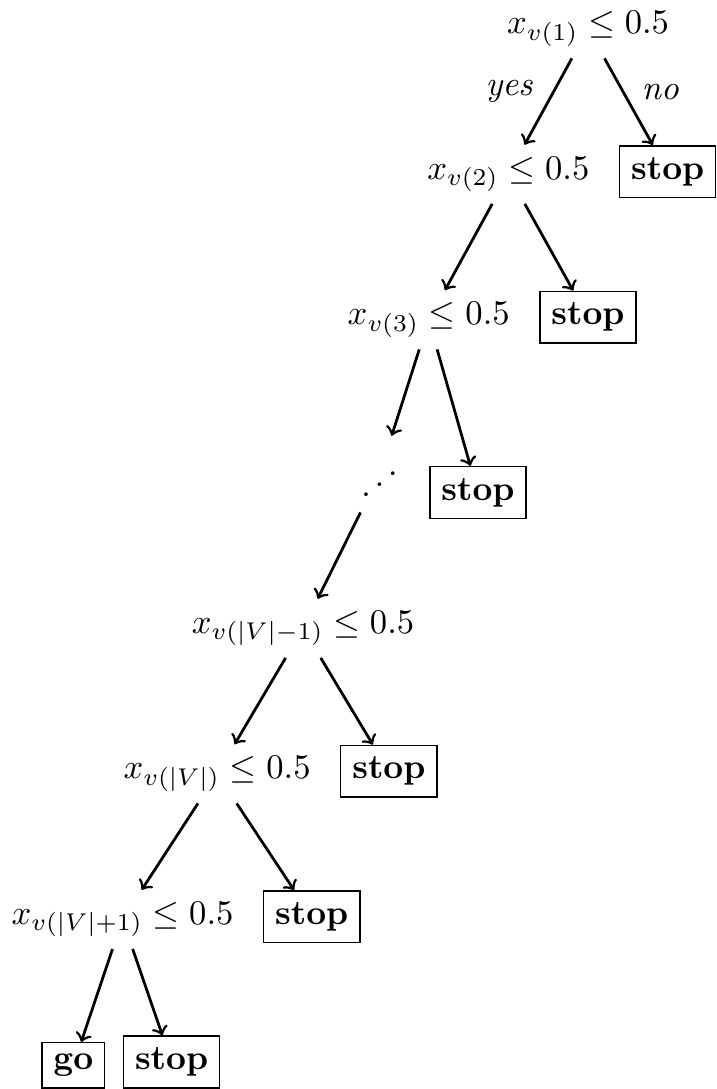} 
                \caption{Tree policy for Proposition~\ref{proposition:split_variable_NPHard} proof.}
                \label{figure:nphardness_SplitVarOpt_policy}
\end{figure}

In this version of the overall tree optimization problem, the leaf actions $\ab$ are now fixed, and the split variable indices must be chosen. For the leaf actions, we set them as follows:
\begin{align*}
& a(|V|+1+i)  = \stop, \quad \forall \ i \in V, \\
& a(2|V|+2)  = \stop, \\
& a(2|V|+3)  = \go.
\end{align*}
We also use $v(1), \dots, v(|V|+1)$ to denote the split variable indices of splits $s = 1,\dots, |V|+1$. Figure~\ref{figure:nphardness_SplitVarOpt_policy} visualizes the tree policy structure.

Note that with this tree topology, split points and leaf actions, the policy will behave as follows:
\begin{enumerate}
\item For an edge trajectory corresponding to an edge $e$:
\begin{itemize}
\item At $t = 1$, the policy will take the action $\go$; 
\item At $t = 2$, the policy will take the action $\stop$ if and only if $v(s) = e_1$ for some $s \in \{1, \dots, |V|+1\}$; and
\item At $t = 3$, the policy will take the action $\stop$ if and only if $v(s) = e_2$ for some $s \in \{1, \dots, |V|+1\}$. 
\end{itemize}
\item For a vertex trajectory corresponding to a vertex $v$:
\begin{itemize}
\item At $t = 1$, the policy will take the action $\stop$ if and only if $v(s) = v$ for some $s \in \{1,\dots, |V| + 1\}$; 
\item At $t = 2$, the policy will take the action $\go$; and 
\item At $t = 3$, the policy will take the action $\stop$ if and only if $v(s) = |V| + 1$, for some $s \in \{1,\dots, |V| + 1\}$. 
\end{itemize}
\end{enumerate}
As with the leaf action SAA problem, we again drop the $(1/\Omega)$ factor in the objective function of problem~\eqref{prob:optimal_stopping_SAA_splitvariable} to ease the exposition.

Let $\vb^*$ be an optimal solution of the corresponding split variable index SAA problem. Let $S = \{ i \in V \mid v(s) = i \ \text{for some} \ s \in \splits\}$ be the set of split variable indices. Effectively, $\vb^*$ encodes a subset of the vertices of $V$, and we denote this subset by $S$. We show that $S$ is an optimal solution of the vertex cover problem. We split our analysis in two cases: \\

\noindent  \textbf{Case 1}: $E = \emptyset$. In the case that there are no edges, an upper bound on the objective value of the split variable index SAA problem is $|V| / (|V| + 1)$, which follows because the most reward that can be obtained from any of the vertex trajectories is $1 / (|V| +1)$ and there are $|V|$ such trajectories. This upper bound can only be attained by setting $v(1) = \dots = v(|V| + 1) = |V| +1$, which corresponds to a set of vertices $S = \emptyset$, which is precisely the optimal solution in this case. \\

\noindent \textbf{Case 2}: $E \neq \emptyset$. In this case, we show that $S$ is both a feasible cover and a minimal cover.

To see that $S$ is feasible, suppose that $S$ were not a feasible cover. This would mean that there exists an edge $e$ such that the two vertices that are incident to this edge, $e_1$ and $e_2$, are not chosen as split variable indices -- in other words, for all $s \in \splits$, $v(s) \notin \{e_1, e_2\}$. This means that the policy cannot garner the reward of 1 from the edge trajectory corresponding to $e$. As a result, the reward of the tree policy that corresponds to $\vb^*$ is upper bounded by $|E| - 1 + |V| / (|V| + 1)$. However, note that one can attain a reward of at least $|E|$, which is higher, by simply setting $v(1) = 1, v(2) = 2, \dots, v(|V|) = |V|$ in the tree policy (this policy is guaranteed to obtain a reward of 1 from every edge trajectory). This leads to a contradiction, because $\vb^*$ is assumed to be optimal. Therefore, it must be the case that $S$ is a feasible cover.

To see that $S$ is a minimal cover, suppose that we could find $S' \subseteq V$ that covers $E$ such that $|S'| < |S|$. Enumerate the nodes in $S'$ as $S' = \{i_1, \dots, i_{|S'|}\}$, and consider the solution $\vb'$ to the split variable index SAA problem obtained by setting
\begin{align*}
& v'(1)  = i_1, \\
& v'(2)  = i_2, \\
& \vdots  \\
& v'(|S'|)  = i_{|S'|},\\
& v'(|S'|+1)  = |V| + 1, \\
& \vdots  \\
& v'(|V|)  = |V| + 1, \\
& v'(|V| + 1) = |V| + 1.
\end{align*}
The resulting tree policy will stop at either $t = 2$ or $t = 3$ in every edge trajectory (because $S'$ is assumed to be a cover), and accrue a reward of 1 from each such trajectory. For each vertex $i \in S'$, it will stop at $t = 1$ in the corresponding vertex trajectory, and obtain a reward of 0. For each vertex $i \notin S'$, it will stop at $t = 3$ in the $i$th vertex trajectory and obtain a reward of $1 / (|V| +1)$. The total reward is therefore $|E| + (|V| - |S'|) / (|V| + 1)$.

Observe that for $\vb^*$, the total reward is at most $|E| + (|V| - |S| ) / (|V| + 1)$. Since $|S'| < |S|$, it follows that
\begin{equation*}
|E| + (|V| - |S|) / (|V| + 1) < |E| + (|V| - |S'|) / (|V| + 1),
\end{equation*}
which implies that $\vb'$ attains a higher objective than $\vb^*$ in the split variable index SAA problem. This immediately results in a contradiction, because $\vb^*$ is assumed to be an optimal solution of that problem. It therefore follows that the cover $S$ we defined above in terms of $\vb^*$ is a minimal cover. 

Thus, by solving the split variable index problem, we are able to solve the minimum vertex cover problem. Since the minimum vertex cover problem is NP-Complete and the instance we have described is polynomially sized in terms of $E$ and $V$, it follows that the split variable index problem must be NP-Hard. \Halmos

\subsubsection{Proof of Proposition~\ref{proposition:split_point_NPHard}}
\label{proof:split_point_NPHard}

The proof that the split point SAA problem is NP-Hard follows along similar lines as the proofs of Propositions~\ref{proposition:leaf_decision_NPHard} and \ref{proposition:split_variable_NPHard} -- in particular, by fixing a specific tree topology, split variable indices and leaf actions, and a specific set of trajectories, one can show that the split point SAA problem can be used to solve the minimum vertex cover problem.

\begin{figure}
 \centering
                \includegraphics[width = 0.4\textwidth]{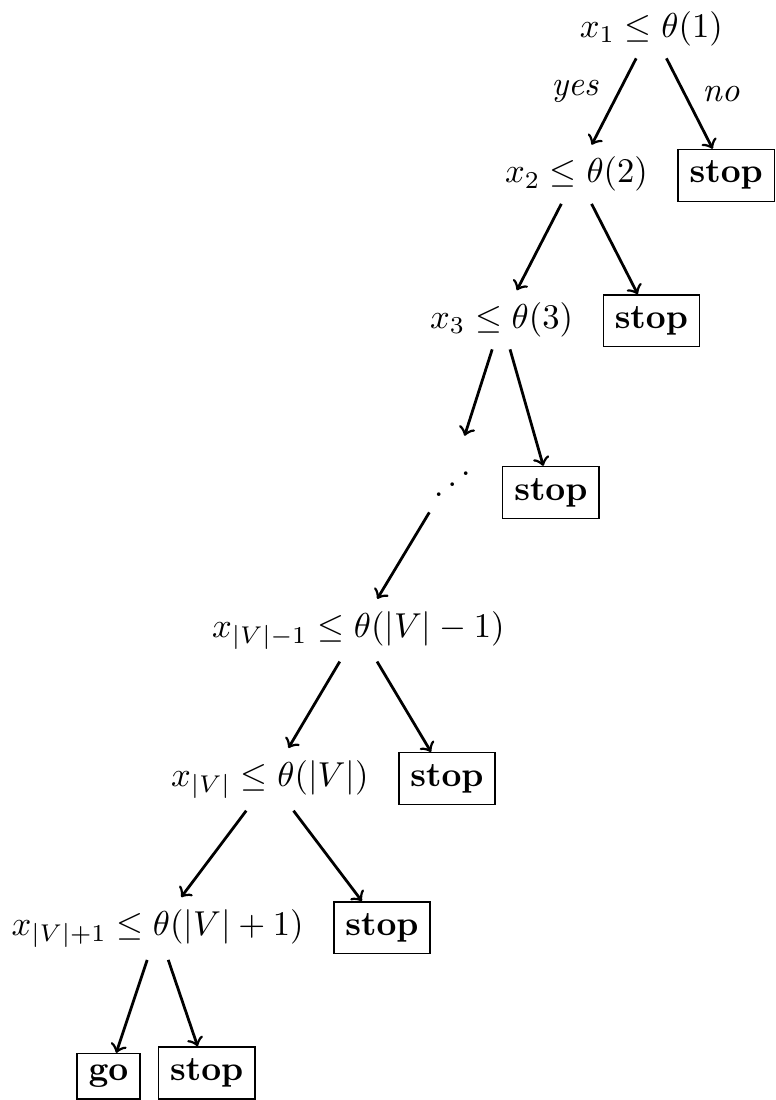} 
                \caption{Tree policy for Proposition~\ref{proposition:split_point_NPHard} proof.}
                \label{figure:nphardness_SplitPointOpt_policy}
\end{figure}

Given an instance $(V,E)$ of the minimum vertex cover problem, consider the same set of trajectories and payoffs as in the proofs of Propositions~\ref{proposition:split_variable_NPHard}, and the same tree topology and split variable indices. Define the leaf actions as 
\begin{align*}
& a(|V|+2) = \stop, \\
& a(|V| + 3) = \stop, \\
& \vdots \\
& a(2|V| + 1) = \stop, \\
& a(2|V| + 2) = \stop, \\
& a(2|V| + 3) = \go.
\end{align*}
Figure~\ref{figure:nphardness_SplitPointOpt_policy} visualizes the tree policy structure.

Note that with this tree topology, split variable indices and leaf actions, the policy will behave as follows:
\begin{enumerate}
\item For an edge trajectory corresponding to an edge $e$:
\begin{itemize}
\item At $t = 1$, the policy will take the action $\go$ if and only if $\theta(s) \geq 0$ for all $s$; 
\item At $t = 2$, the policy will take the action $\go$ if and only if $\theta(s) \geq 0 $ for all $s \neq e_1$ and $\theta(e_1) \geq 1$; and 
\item At $t = 3$, the policy will take the action $\go$ if and only if $\theta(s) \geq 0 $ for all $s \neq e_2$ and $\theta(e_2) \geq 1$.
\end{itemize}
\item For a vertex trajectory corresponding to a vertex $v$:
\begin{itemize}
\item At $t = 1$, the policy will take the action $\go$ if and only if $\theta(s) \geq 0$ for all $s \neq v$ and $\theta(v) \geq 1$; 
\item At $t = 2$, the policy will take the action $\go$ if and only if $\theta(s) \geq 0$ for all $s$; and 
\item At $t = 3$, the policy will take the action $\go$ if and only if $\theta(s) \geq 0$ for all $s \neq |V| + 1$ and $\theta(|V| + 1) \geq 1$.
\end{itemize}
\end{enumerate}
As with the leaf action and split variable index SAA problems, we again drop the $(1/\Omega)$ factor in the objective function of problem~\eqref{prob:optimal_stopping_SAA_splitpoint} to ease the exposition. 

Let $\thetab^* = (\theta^*(1), \theta^*(2), \dots, \theta^*(|V|+1))$ be an optimal solution to the split point SAA problem. Define the set of vertices $S$ as 
\begin{equation*}
S = \{ i \in V \, | \, \theta(i) \in [0, 1) \}.
\end{equation*}
We will now show that $S$ is an optimal solution of the vertex cover problem. 

Before we show this, we first argue that in any optimal solution $\thetab^*$ of the split point SAA problem, it must be that $\theta^*(i) \geq 0$ for all $i = 1,\dots, |V|+1$. To see this, observe that if there exists an $i$ such that $\theta^*(i) < 0$, then at $t = 1$, we will stop in every trajectory, because every coordinate of every trajectory is greater than or equal to 0 at $t = 1$. As a result, we will accrue a total reward of zero from such a policy. However, observe that if we set every $\theta(i)$ to 0.5, then we are guaranteed to stop at $t = 2$ or $t = 3$ in every edge trajectory and accrue a reward of 1 from each such trajectory; as a result, the total reward will be at least $|E|$, which would contradict the optimality of $\thetab^*$. 

We now prove that $S$ solves the minimum vertex cover problem. We consider two cases:\\

\noindent \textbf{Case 1}: $E = \emptyset$. In this case, observe that an upper bound on the reward that can be obtained is $|V| / (|V| + 1)$. This reward can only be attained by setting $\theta(v) \geq 1$ for every $v = 1,\dots, |V|$, and $\theta(|V|+1) \in [0, 1)$. Note that the set of vertices $S$ that corresponds to $\thetab$ in this case is the empty set, which is exactly the optimal solution of the minimum vertex cover problem in this case. \\

\noindent \textbf{Case 2}: $E \neq \emptyset$. In this case, we show that $S$ is both a feasible cover and a minimal cover.

To see that $S$ is a feasible cover, suppose that $S$ were not a feasible cover. This would mean that there is an edge $e$ that is not covered, i.e., that $e_1, e_2 \notin S$. By the definition of $S$, this would mean that $\theta(e_1) \notin [0, 1)$ and $\theta(e_2) \notin [0,1)$. By our earlier observation that $\theta^*(i) \geq 0$ for any optimal $\theta^*$, this means that $\theta(e_1) \geq 1$ and $\theta(e_2) \geq 1$. Note that by the definition of the tree policy and the trajectories, this would mean that the policy defined by $\thetab^*$ either stops at $t = 1$ for the edge trajectory corresponding to $e$ (resulting in a reward of zero from that trajectory), or it does not stop at any $t$ (again, resulting in a reward of zero from the trajectory). As a result, the reward obtained by $\thetab^*$ would be upper bounded by $|E| - 1 + |V| / (|V| + 1)$. However, observe that just by setting all $\thetab(s) = 0.5$ for $s = 1,\dots, |V|+1$, we would obtain a policy that stops each edge trajectory at either $t = 2$ or $t= 3$, which would guarantee a reward of at least $|E|$, which is greater. Since this would contradict the optimality of $\thetab^*$, it follows that $S$ must be a feasible cover.

To see that $S$ is a minimal cover, suppose this were not the case. This would mean that there exists a set $S'$ that covers $E$ and is smaller, i..e, $|S'| < |S|$. For this cover $S'$, define a new split point vector $\thetab'$ as
\begin{align*}
& \theta'( s ) = 0.5, \quad \forall s \in S', \\
& \theta'( s ) = 1, \quad \forall s \in V \setminus S', \\
& \theta'( |V| + 1 ) = 0.5. 
\end{align*}
Since $S'$ covers $E$, this policy is guaranteed to stop at $t = 2$ or $t = 3$ for every edge trajectory, and thus will garner a reward of 1 from each such trajectory. For each vertex $i \in S$, it will stop at $t = 1$ in the corresponding vertex trajectory, and garner a reward of zero. For each vertex $i \notin S$, it will stop at $t = 3$ in the corresponding vertex trajectory, and earn a reward of $1 / (|V| + 1)$. Thus, the reward of the policy defined by $\thetab'$ will be $|E| + (|V| - |S'|) / (|V| + 1)$.

Now, observe that for the optimal policy $\thetab^*$, since $S$ corresponds to a feasible cover, the total reward from all edge trajectories for $\thetab^*$ is upper bounded by $|E|$. For $s \in S$, by definition of $S$, the reward from the vertex trajectory corresponding to $s$ must be zero. For $s \in V \setminus S$, the reward from the vertex trajectory corresponding to $s$ is upper bounded by $1 / (|V|+1)$. Therefore, the reward of $\thetab^*$ is upper bounded by $|E| + (|V| - |S|)/(|V|+1)$. Since $|S'| < |S|$, this would imply that the reward of $\thetab'$ is greater than the reward of $\thetab^*$; however, this would contradict the fact that $\thetab^*$ is an optimal solution to the split point SAA problem. Therefore, it must be that $S$ as defined above is a minimal cover. 

This establishes that the split point SAA problem is NP-Hard, as required. \Halmos

\subsubsection{NP-Hardness of optimization over thresholds}
\label{appendix:threshold_opt_nphard}

\begin{figure}[!ht]
\centering
\includegraphics[width=0.7\textwidth]{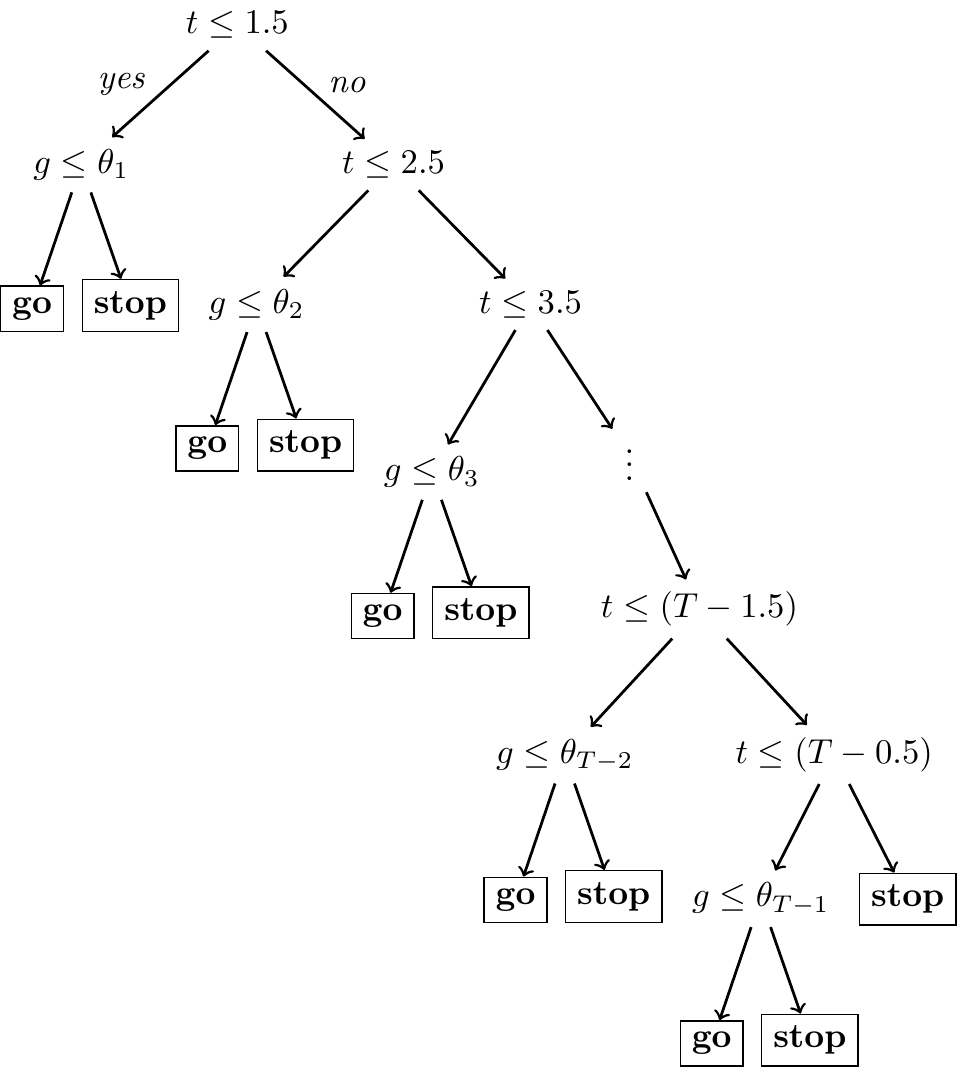}
\caption{Threshold tree policy structure. \label{figure:threshold_policy_tree} }
\end{figure}

In this section, we consider a specific type of tree policy optimization problem, where one optimizes a set of time-dependent thresholds. Suppose that the time horizon is $T$, and that the specification of the state variable $\xb(t)$ includes the time $t$ and the reward $g(t, \xb(t))$. We assume that $g(t, \xb(t)) \geq 0$ for all $t$. Suppose that the tree policy follows the form given in Figure~\ref{figure:threshold_policy_tree}. %

Let $\Pi_{\text{Threshold}}$ be the set of policies obtained by varying the thresholds $\theta_1,\dots, \theta_{T-1}$. Let $\Omega$ be the number of trajectories in the training set, and let $\xb(\omega,t)$ be the state variable at time $t$ in trajectory $\omega$. The time-dependent threshold SAA problem is defined as
\begin{equation}
\underset{\pi \in \Pi_{\text{Threshold}}}{\text{maximize}} \ \frac{1}{\Omega} \sum_{\omega =1}^{\Omega}  \beta^{\tau_{\pi,\omega}-1} \cdot g(\tau_{\pi, \omega}, \xb(\omega, \tau_{\pi,\omega})). \label{prob:timedependentthresholdSAA}
\end{equation}

We then have the following result.
\begin{proposition}
The time-dependent threshold SAA problem~\eqref{prob:timedependentthresholdSAA} is NP-Hard.
\end{proposition}

\proof{Proof.} We will prove this result by reducing the minimum vertex cover problem to problem~\eqref{prob:timedependentthresholdSAA}. Given a vertex cover instance defined by a set of vertices $V$ and a set of edges $E$, we will construct a corresponding instance for problem~\eqref{prob:timedependentthresholdSAA} as follows. We set the time horizon $T = |V| + 1$. For convenience, we assume that the vertices in $V$ are indexed from 1 to $|V|$. 

We create two types of trajectories:
\begin{enumerate}
\item \emph{Edge trajectories}: For each edge $e \in E$, we create $2|V| + 1$ copies of the same trajectory, which has the following reward structure:
\begin{align*}
g(t, \xb(\omega,t)) & = 1, \quad \text{if}\ t \in \{e_1, e_2\}, \\
g(t, \xb(\omega,t)) & = 0, \quad \text{if}\ t \notin \{e_1,e_2\},
\end{align*}
where $e_1$ and $e_2$ are the first and second vertex, respectively, to which edge $e$ is incident. 
\item \emph{Vertex trajectories}: For each vertex $v \in V$, we create a single trajectory, which has the following reward structure:
\begin{align*}
g(t, \xb(\omega,t)) & = 1, \quad \text{if}\ t = v, \\
g(t, \xb(\omega,t)) & = 2, \quad \text{if}\ t = |V|+1, \\
g(t, \xb(\omega,t)) & = 0, \quad \text{if}\ t \notin \{v, |V|+1\}.
\end{align*}
\end{enumerate}

We now proceed with the reduction. Note that unlike the previous reductions, we no longer suppress the $(1/\Omega)$ factor in the objective of the SAA problem.

Let $(\theta_1,\dots, \theta_{T-1})$ be the thresholds of the optimal solution to problem~\eqref{prob:timedependentthresholdSAA}. Since $T - 1 = |V|$, we have exactly $|V|$ thresholds. Given this solution, let us construct a candidate solution to the vertex cover problem by specifying the set of vertices $S$ as
\begin{equation}
S = \{ i \in V \ \vline \ \theta_i < 1 \}. \label{eq:cover_defn}
\end{equation}
We now argue that $S$ is a cover and that $S$ is the minimal cover.

\noindent \emph{Feasibility}. To see why $S$ is a cover, let us proceed by contradiction and suppose that it is not. If it is not a cover, then there exists an edge $e \in E$ such that $e_1$ and $e_2$ are not contained in $S$. If $e_1$ and $e_2$ are not contained in $S$, then this means that $\theta_{e_1} \geq 1$ and $\theta_{e_2} \geq 1$. If this is the case, then observe that for any of the $2|V|+1$ copies of the edge $e$ trajectory, the policy will not stop at $t = e_1$ or at $t = e_2$ for any of those trajectories. Thus, for any of those trajectories, the reward will be zero. We therefore have that the reward of the policy satisfies the following bound
\begin{equation}
(1/\Omega) \cdot (2|V|+1) \cdot (|E| - 1) + (1/\Omega) \cdot 2|V|  \label{eq:obj_bound}
\end{equation}
where the first term corresponds to the highest possible reward from the edge trajectories (i.e., garnering a reward of 1 in every edge trajectory except those corresponding to edge $e$; there are $|E|-1$ other edges beside $e$, and $2|V|+1$ copies of each edge's trajectory), and the second term corresponds to the highest possible reward from the vertex trajectories (i.e., garnering a reward of 2 in every vertex trajectory by stopping at $t = T = |V|+1$). 

We now argue that the policy defined by $(\theta_1,\dots, \theta_{|V|})$ is not optimal, thereby leading to a contradiction. Consider a different policy where we set $\theta_i = 0.5$ for each $i = 1,\dots, |V|$. The corresponding policy will stop at either $t = e'_1$ or $e'_2$ for every trajectory corresponding to each edge $e'$. Thus, the reward of this policy will be at least $(1/\Omega) \cdot (2|V|+1) \cdot |E|$, which is strictly higher than the bound~\eqref{eq:obj_bound}. This leads to an immediate contradiction, because the policy defined by $(\theta_1,\dots, \theta_{|V|})$ was assumed to be optimal. This establishes that $S$, as defined in \eqref{eq:cover_defn}, is a bona fide cover.

\noindent \emph{Optimality}. To see why $S$ is a minimal cover, let us suppose there exists a cover $S'$ such that $|S'| < |S|$. Using $S'$, let us construct a new policy for problem~\eqref{prob:timedependentthresholdSAA} that is defined by the following thresholds:
\begin{equation}
\theta'_i = \left\{ \begin{array}{ll} 0.5 &\text{if}\ i \in S, \\
1.5 & \text{if}\ i \notin S. \end{array} \right. 
\end{equation}
For this policy, it can be verified that the total reward is 
\begin{align}
& \frac{1}{\Omega} \cdot (2 |V|+1) \cdot |E| + \frac{1}{\Omega} \cdot ( |S'| + 2(|V| - |S'|) ) \\
& = \frac{1}{\Omega} \cdot (2 |V|+1) \cdot |E| + \frac{1}{\Omega} \cdot ( 2|V| - |S'| )
\end{align}
Similarly, it can be verified that for the original policy defined by $(\theta_1, \dots, \theta_{|V|})$, the total reward is at most
\begin{align}
\frac{1}{\Omega}  \cdot (2 |V|+1) \cdot |E| + \frac{1}{\Omega}  \cdot ( 2|V| - |S| ).
\end{align}
Since $|S'| < |S|$, it follows that $(2|V| - |S'|) > (2|V| - |S|)$, which implies that the policy defined by $(\theta'_1,\dots, \theta'_{|V|})$ achieves a higher reward than the policy defined by $(\theta_1,\dots, \theta_{|V|})$. This contradicts the fact that $(\theta_1,\dots, \theta_{|V|})$ corresponds to an optimal policy for problem~\eqref{prob:timedependentthresholdSAA}. Therefore, $S$ must be a minimal cover.

Since we have proven that the solution $S$ derived from a solution of problem~\eqref{prob:timedependentthresholdSAA} is a solution to the minimum vertex cover problem, and since the optimal stopping instance that we derive is polynomially sized in $(V, E)$, it follows that problem~\eqref{prob:timedependentthresholdSAA} is NP-Hard. $\square$

\endproof

\subsection{Proofs and additional lemmas for Section \ref{subsec:tree_optimality}}
\label{sec:tree_optimality_appendix}

\begin{figure}
\begin{center}
\includegraphics[width=0.5\textwidth]{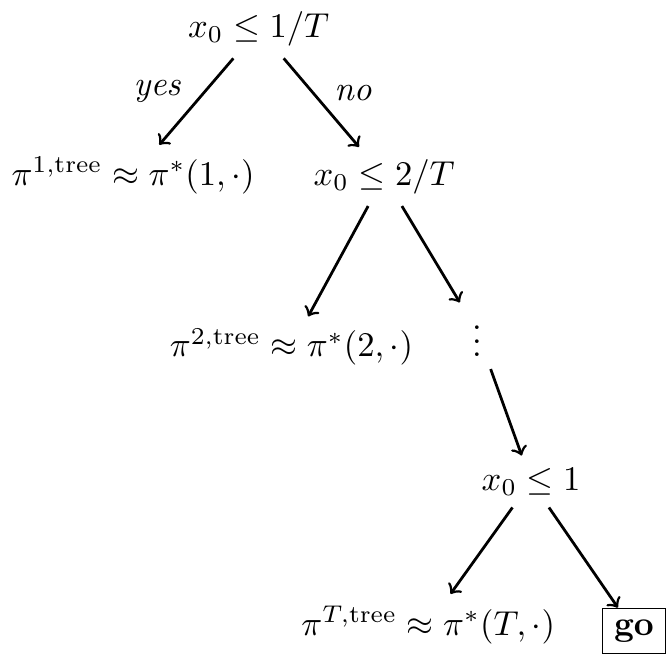}
\end{center}
\caption{The approximate tree. \label{fig:approximate_tree_time_splits}}
\end{figure}

\textbf{Proof of Theorem~\ref{thm:tree_policy_appr_opt}}.
Our proof proceeds in the following way. We note that, as per Assumption \ref{ass:compactness} Part 1, we have augmented the state with a variable $x_0$ which is set to $t/T$ if and only if the time period is $t$. We build a tree which first checks whether $x_0$ corresponds to period $t$ as in Figure \ref{fig:approximate_tree_time_splits}; the left subtree of the split node $\{x_0 \leq t/T\}$ will approximate $\pi^*(t, \cdot)$.

It thus remains to show that for some time period $t$, there exists a subtree which can approximate $\pi^*(t, \cdot)$. This is trivial for $t=1$ (the left subtree of the split $x_0 \leq 1/T$ is set to a leaf prescribing the action that $\pi^*(1, \bar\xb)$ takes) and thus we restrict our attention to stopping problems whose optimal policies do not stop at the starting state, i.e. policies for which $\pi^*(1,\bar \xb) = \go$.

Proceeding to this case, we show that $\pi^*(t,\cdot)$ for $t > 1$ can be approximated by a simpler policy %
which takes the action $\stop$ if and only if the current state belongs to some finite union of half-open boxes, i.e. $\xb(t) \in \bigcup_{j \in K(t)} (a^{j,t}_1, b^{j,t}_1] \times \ldots \times (a^{j,t}_n, b^{j,t}_n]$. We then argue (via Lemma \ref{le:from_boxes_to_trees}) that such a box policy can be alternatively represented as a tree policy %
for all $t \in [T]$. Finally, we show that this policy disagrees with $\pi^*(t, \cdot)$ on $\xb(t)$ with small probability, and thus achieves roughly the same value as $\pi^*$.

We proceed with the first step. Using Assumption~\ref{ass:state_smoothness}, let $\delta > 0$ be such that 
\begin{equation*}
f(\delta) \leq \frac{\epsilon}{GT}
\end{equation*}
Now, consider any $t \in [T] \setminus \{1\}$. Since $\Xcal^{\stop}_t$ is compact, it is also closed and therefore a Borel set. By the definition of Borel measure, there exists a countable collection of sets $\{ \Bcal_{j,t} \}_{j \in J(t)}$, where each set $\Bcal_{j,t}$ is an open box of the form $(a^{j,t}_1, b^{j,t}_1) \times \dots \times (a^{j,t}_n, b^{j,t}_n)$, such that:
\begin{equation}
\mu( \Xcal^{\stop}_t ) \leq \mu \left( \bigcup_{j \in J(t)} \Bcal_{j,t} \right) \leq \mu(\Xcal^{\stop}_t) + \delta. 
\end{equation}
Since $\Xcal^{\stop}_t$ is compact, it follows that there exists a finite subcollection of boxes $\{ \Bcal_{j,t} \}_{j \in K(t)}$ from the collection $\{ \Bcal_{j,t} \}_{j \in J(t)}$ that cover $\Xcal^{\stop}_t$. Since $\Xcal^{\stop}_t \subseteq \bigcup_{j \in K(t)} \Bcal_{j, t} \subseteq \bigcup_{j \in J(t)} \Bcal_{j,t}$, we have that 
\begin{equation}
\mu(\Xcal^{\stop}_\go) \leq \mu\left( \bigcup_{j \in K(t)} \Bcal_{j, t} \right) \leq \mu( \Xcal^{\stop}_\go) + \delta.
\end{equation}
Lastly, if each box $\Bcal_{j,t}$ is replaced by the half-closed boxed $\bar \Bcal_{j,t} = (a^{j,t}_1, b^{j,t}_1 ] \times \dots \times (a^{j,t}_n, b^{j,t}_n]$, we still have that 
\begin{equation}
\mu\left( \bigcup_{j \in K(t)} \bar\Bcal_{j, t} \right) \leq \mu( \Xcal^{\stop}_\go) + \delta.
\end{equation}
Rearranging the above inequality, we therefore obtain that
\begin{equation}
\mu \left( \bigcup_{j \in K(t)}  \bar\Bcal_{j, t} \setminus \Xcal^{\stop}_{\go} \right) \leq \delta. \label{eq:disagreement_volume}
\end{equation}

Using the finite collection of boxes $\{\bar \Bcal_{j,t}\}_{j \in K(t)}$, we can define a policy $\pi^{t, \textrm{boxes}}$ as follows: 
\begin{equation}
\pi^{t, \textrm{boxes}}(\xb)  = \left\{ \begin{array}{ll} \stop & \text{if} \ \xb \in \bigcup_{j \in K(t)} \bar\Bcal_{j, t}, \\
\go & \text{otherwise},
\end{array} \right.
\end{equation} 
i.e., a policy that stops at state $\xb$ if and only if $\xb \in \bigcup_{j \in K(t)} \bar \Bcal_{j,t}$. By Lemma \ref{le:from_boxes_to_trees}, we can represent $\pi^{t, \textrm{boxes}}$  as a tree policy, whose corresponding tree we append as the left sub-tree of the $x_0(t) \leq t/T$ split to form $\pi^{\textrm{tree}}$. Then,
\begin{align*}
& J^*(\bar\xb) - J^{\pi^{\textrm{tree}}}(\bar\xb) \\
& \quad =
\Exp\left[ \left. \beta^{\tau_{\pi^*} - 1} g\left(\tau_{\pi^*}, \xb\left(\tau_{\pi^*} \right) \right)-  \beta^{\tau_{\pi^{\textrm{tree}}} - 1} g\left(\tau_{\pi^{\textrm{tree}}}, \xb\left(\tau_{\pi^{\textrm{tree}}} \right) \right) \right| \xb(1) = \bar\xb \right] \\
& \quad = \Exp\left[ \left. \beta^{\tau_{\pi^*} - 1} g\left(\tau_{\pi^*}, \xb\left(\tau_{\pi^*} \right) \right)-  \beta^{\tau_{\pi^{\textrm{tree}}} - 1} g\left(\tau_{\pi^{\textrm{tree}}}, \xb\left(\tau_{\pi^{\textrm{tree}}} \right) \right) \right| \tau_{\pi^*} \neq \tau_{\pi^{\textrm{tree}}},  \xb(1) = \bar\xb \right] \cdot \Pr\left[ \left.  \tau_{\pi^*} \neq \tau_{\pi^{\textrm{tree}}} \right|\xb(1) = \bar\xb  \right] \\
& \quad\quad\quad + 0 \cdot  \Pr\left[ \left.  \tau_{\pi^*} = \tau_{\pi^{\textrm{tree}}} \right|\xb(1)= \bar\xb  \right] \\
& \quad \leq G \cdot \Pr\left[ \left.  \tau_{\pi^*} \neq \tau_{\pi^{\textrm{tree}}} \right|\xb(1) = \bar\xb  \right] \\
& \quad \leq G \cdot \Pr\left[ \left.  \exists t \in [T] \setminus \{1\} \st \pi^*(t, \xb(t)) \neq \pi^{\textrm{tree}}(t, \xb(t)) \right|\xb(1) = \bar\xb  \right] \\
& \quad \leq G \cdot \sum_{t \in [T] \setminus \{1\}} \Pr\left[ \left. \pi^*(t, \xb(t)) \neq \pi^{t,\textrm{boxes}}(t, \xb(t)) \right|\xb(1) = \bar\xb  \right] \\
& \quad = G \cdot \sum_{t \in [T] \setminus \{1\}} \Pr\left[ \xb(t) \in \bigcup_{j \in K(t)} \bar \Bcal_{j,t} \setminus \Xcal_t^{\stop} \ \vline \ \xb(1) = \bar\xb  \right] \\
& \quad \leq G \cdot T \cdot f(\delta) \\
& \quad = \epsilon.
\end{align*}
In the above, the first inequality follows from Assumption \ref{ass:bounded_cost} and the second inequality follows by construction, since the two policies can only disagree after $t>1$. The penultimate equality follows because the states $\xb$ where $\pi^{\text{boxes}}$ and $\pi^*$ disagree are exactly those where $\xb$ is inside the finite subcover $\{ \bar \Bcal_{j, t} \}_{j \in K(t)}$, but not inside the optimal stopping region $\Xcal^{\stop}_t$. The last inequality then follows from Assumption \ref{ass:state_smoothness} and equation \eqref{eq:disagreement_volume}, ignoring the augmented state variable $x_0$. \Halmos

\begin{lemma}
\label{le:from_boxes_to_trees}
Consider a collection of $K$ boxes $\{ \bar \Bcal_{k} \}_{k=1}^K$, where $\bar \Bcal_k = (a^k_1, b^k_1] \times \dots (a^k_n, b^k_n]$, and a policy $\pi^{\text{boxes}}$ of the form
\begin{equation}
\pi^{\text{boxes}}(\xb) = \left\{ \begin{array}{ll} \stop & \text{if} \ \xb \in \bigcup_{k=1}^K \bar \Bcal_k, \\
\go & \text{otherwise}. \end{array} \right.
\end{equation}
Then there exists a tree policy $\pi^{\text{tree}}$ such that $\pi^{\text{tree}}(\xb) = \pi^{\text{boxes}}(\xb)$ for all $\xb \in \Xcal$. 
\end{lemma}

\proof{Proof.}

\begin{figure}
\centering
\includegraphics[width=0.2\textwidth]{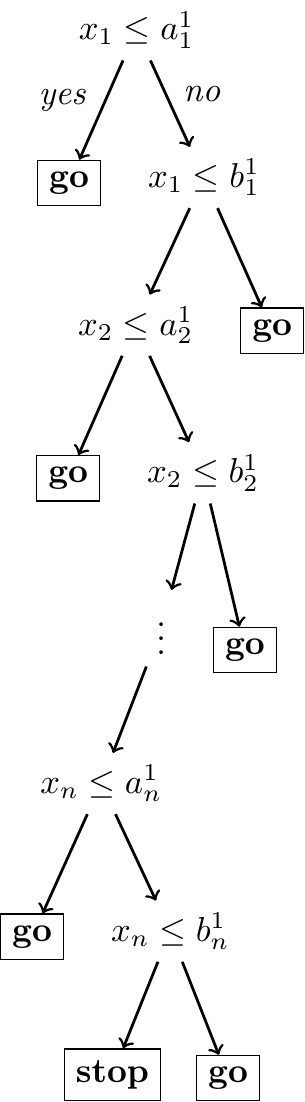}
\caption{Tree representation of $\pi^{\text{boxes}}$ when the number of boxes $K = 1$ and $\bar \Bcal_1 = (a^1_1, b^1_1] \times \dots (a^1_n, b^1_n]$. \label{figure:policy_boxes_Keq1}}
\end{figure}

We will prove this statement by induction over $K$, the number of boxes. First, observe that when $K = 1$, there is only one box $\bar \Bcal_1 = (a^1_1, b^1_1] \times \dots \times (a^1_n, b^1_n]$. The corresponding policy can be exactly represented by a tree, as shown in Figure~\ref{figure:policy_boxes_Keq1}. 

\begin{figure}
\begin{subfigure}[t]{0.25\textwidth}
                \centering
                \includegraphics[width = 0.6\textwidth]{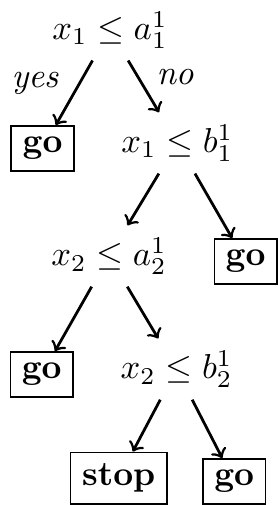}
                \caption{Tree policy corresponding to $k = 1$ box $\bar \Bcal_1 = (a^1_1, b^1_1] \times (a^1_2, b^1_2]$.}
                \label{figure:box_tree_induct1}
        \end{subfigure}  \hfill
\begin{subfigure}[t]{0.25\textwidth}
                \centering
                \includegraphics[width = 0.6\textwidth]{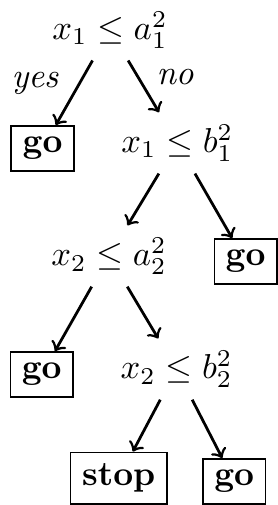} 
                \caption{Tree policy corresponding to $k = 2$ box $\bar \Bcal_2 = (a^2_1, b^2_1] \times (a^2_2, b^2_2]$.}
                \label{figure:box_tree_induct2}
        \end{subfigure} 
        \hfill
\begin{subfigure}[t]{0.4\textwidth}
                \centering
                \includegraphics[width = 1.1\textwidth]{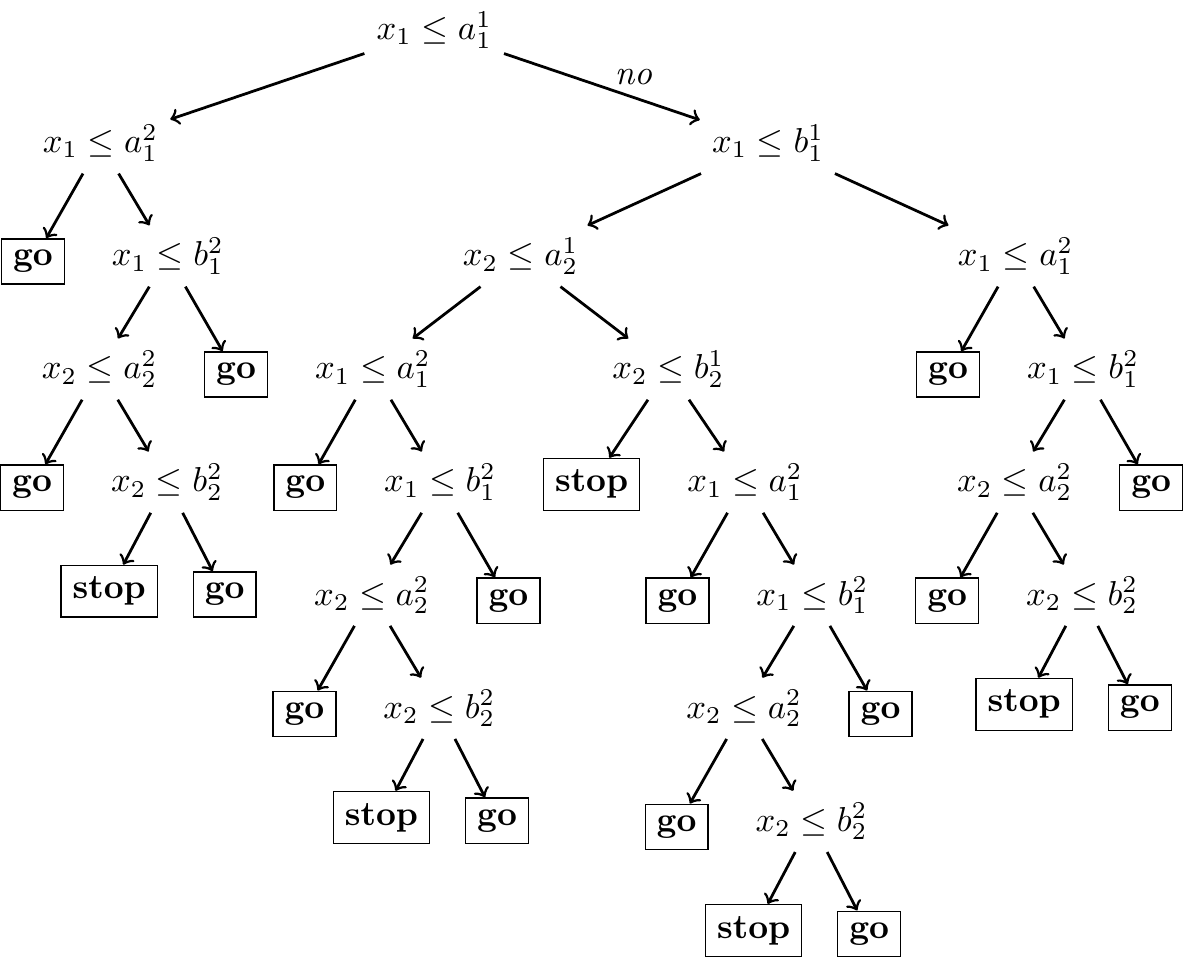} 
                \caption{Tree policy corresponding to $\pi^{\text{boxes}}$, obtained by replacing each $\go$ leaf in (a) by the tree in (b).}
                \label{figure:box_tree_induct1and2}
        \end{subfigure}

\caption{Visualization of induction step for Lemma~\ref{le:from_boxes_to_trees} with $K = 2$ boxes, $n = 2$ state variables. \label{figure:box_induction_example}}

\end{figure}

Now, suppose that the statement holds for $K - 1$ boxes, and we need to prove it for $K$ boxes. For the box $\bar \Bcal_{K} = (a^K_1, b^K_1] \times \dots \times (a^K_n, b^K_n]$, one can construct a tree in the same way as for the base case, as shown in Figure~\ref{figure:policy_boxes_Keq1}. Then, for the tree corresponding to $K - 1$ boxes, replace each leaf node with the $\go$ action by the tree corresponding to the box $\bar \Bcal_K$. Figure~\ref{figure:box_induction_example} illustrates the procedure for $K = 2$ boxes, $n = 2$ state variables. 

The new policy constructed in this way will have the following behavior: it will stop if $\xb \in \bigcup_{k=1}^{K-1} \bar \Bcal_k$; otherwise, if $\xb \notin \bigcup_{k=1}^{K-1} \bar \Bcal_k$, it will check if $\xb \in \bar \Bcal_K$; if it is, it will stop; otherwise, it will continue. It is straightforward to see that this policy behaves in exactly the same way as $\pi^{\text{boxes}}$, which establishes the lemma. \Halmos

\endproof

\clearpage
\pagebreak

\subsection{A constructive procedure for characterizing the tree depth necessary for a certain level of error}

\label{subsec:epsilon_depth_bound_R2}

In this section, we present a theoretical construction of a class of tree policies for which (1) the tree policy achieve an arbitrary small gap with respect to the optimal policy and (2) the depth can be explicitly characterized, and scales gracefully with the optimality gap.

The theoretical argument that we present here is different from the one used in Section~\ref{subsec:tree_optimality}. To establish Theorem~\ref{thm:tree_policy_appr_opt}, we made relatively weak assumptions about the compactness of the stopping region and the smoothness of the marginal probability distributions of the stochastic process $\{\xb(t)\}_{t=1}^T$; we then used the measurability and the compactness of the stopping region to generate a collection of boxes, which could be easily used to form a tree policy that is within $\epsilon$ of the optimal policy. In that argument, the compactness of the stopping region is used in a ``black-box'' fashion, and does not come with any guarantee on the size of the subcover; since the depth of the tree directly scales with the size of the subcover, providing an accompanying guarantee on the depth of the tree policy is challenging.

Thus, in this section, we take a different approach, relying on a different set of assumptions; we outline these assumptions now. Our first assumption is on the state space $\Xcal$. 
\begin{assumption}
$\Xcal = [0,1]^n$. 
\end{assumption}

Our second assumption is on the stochastic process $\{ \xb(t) \}_{t=1}^T$.
\begin{assumption} \label{ass:Markovian}
The stochastic process $\{ \xb(t) \}_{t=1}^T$ is Markovian, i.e., for any $t \in \{1,\dots, T-1\}$ and any Borel set $A \subseteq \Xcal$,
\begin{equation*}
\Pr( \xb(t+1) \in A \mid \Fcal_t ) = \Pr( \xb(t+1) \in A \mid \xb(t)),
\end{equation*} 
where $\Fcal \defeq \{\Fcal_t\}_{t=1}^T$ is the natural filtration of the stochastic process $\{\xb(t)\}_{t=1}^T$. %
\end{assumption}
Under the assumption that the stochastic process is Markovian, we can define the (optimal) value function and continuation value function using the Bellman recursion. In particular, we define the value function $J^*_t: \Xcal \to \Rbb$ and the continuation value function $C^*_t: \Xcal \to \Rbb$ as 
\begin{align*}
C^*_T(\xb) & = 0, \quad \forall\ \xb \in \Xcal, \\
C^*_t(\xb) & = \Exp[ J^*_{t+1}( \xb(t+1) ) \mid \xb(t) = \xb], \quad \forall \ t \in \{1, \dots, T-1\}, \ \xb \in \Xcal, \\
J^*_t(\xb) & = \max\{ g(t, \xb), \beta \cdot C^*_t(\xb) \}, \quad \forall t \in \{1,\dots, T\}, \ \xb \in \Xcal.
\end{align*}
Our next assumption concerns the smoothness of $g$ and $C^*_t$. 
\begin{assumption}
For each $t$, the functions $g(t,\cdot): \Xcal \to \Rbb$ and $C^*_t(\cdot): \Xcal \to \Rbb$ are Lipschitz continuous with common Lipschitz constant $L$. 
\label{ass:g_C_Lipschitz}
\end{assumption}

Our last assumption is that the transition kernel of $\xb(t)$ is bounded.
\begin{assumption} \label{ass:bounded_density}
For any Borel set $A$ and any $t \in \{1,\dots,T-1\}$, we have, almost surely,
\begin{equation}
\Pr( \xb(t+1) \in A \mid \xb(t)) \leq \lambda \cdot \mu(A).
\end{equation}
\end{assumption}

For a fixed integer $m$, we let $H_m$ denote the set of all $m^n$ hypercubes with side length $1 / m$. In addition, let $h(\xb)$ be the hypercube in $H_m$ containing the state $\xb \in \Xcal$.

We now define the approximating policy $\hat{\pi}$ with respect to $H_m$. The policy $\hat{\pi}$ is defined inductively as follows:
\begin{enumerate}
\item For $t = 1$, we define $\hat{\pi}(1,\cdot)$ as 
\begin{equation}
\hat{\pi}(1, \xb) = \left\{ \begin{array}{ll} \stop & \text{if} \quad \Exp[ g(1, \xb(1) ) \mid \xb(1) \in h(\xb) ] > \beta \Exp[  C^*_1(\xb(1)) \mid \xb(1) \in h(\xb) ], \\
\go & \text{if} \quad \Exp[ g(1, \xb(1) ) \mid \xb(1) \in h(\xb) ] \leq \beta \Exp[  C^*_1(\xb(1)) \mid \xb(1) \in h(\xb) ]. \end{array} \right. \label{eq:hatpi_teq1_defn}
\end{equation}
\item For $t = 2, \dots, T$, first define the event $A_t$ as 
\begin{equation}
A_t = \{ \hat{\pi}(t', \xb(t')) = \go \ \text{for all} \ t' < t \}.
\end{equation}
We then define $\hat{\pi}(t, \cdot)$ as 
\begin{equation}
\hat{\pi}(t, \xb) = \left\{ \begin{array}{ll} \stop & \text{if} \quad \Exp[ g(t, \xb(t) ) \mid \xb(t) \in h(\xb), A_t ] > \beta \Exp[  C^*_t(\xb(t)) \mid \xb(t) \in h(\xb), A_t], \\
\go & \text{if} \quad \Exp[ g(t, \xb(t) ) \mid \xb(t) \in h(\xb), A_t ] \leq \beta \Exp[  C^*_t (\xb(t)) \mid \xb(t) \in h(\xb), A_t ]. \end{array} \right. \label{eq:hatpi_tgeq2_defn}
\end{equation}
\end{enumerate}

(Note that in \eqref{eq:hatpi_teq1_defn} and \eqref{eq:hatpi_tgeq2_defn}, the conditioning event in the conditional expectations may have probability zero, in which case the conditional expectation may not be well defined. In the definition of $\hat{\pi}$ and in the analysis that follows, we by default set such conditional expectations to be zero when they occur.)

Two comments are in order regarding $\hat{\pi}$. First, note that the policy $\hat{\pi}$ is piecewise constant on each hypercube. Since each hypercube can be represented with a suitably large number of axis-aligned splits, it should be clear that one can construct a tree that represents $\hat{\pi}$ by splitting on both time $t$ and the state variables $x_1(t), \dots, x_n(t)$. 

Second, it is helpful to understand the rationale behind the policy. For each hypercube, $\hat{\pi}$ must select the same action for all states within that hypercube. Given the constraint that the policy must be constant, an intuitive choice for the action on each hypercube is the constant action that is greedy with respect to the optimal policy. 

As an example, if we consider $t = 1$ and a given state $\xb$, the quantity $\Exp[ \beta C^*_1(\xb(1)) \mid \xb(1) \in h(\xb) ]$ is the expected reward we would get for choosing $\go$ for every state in $h(\xb)$, conditional on $\xb(1)$ being in $h(\xb)$, and then following the optimal policy from that point on. The quantity $\Exp[ g(t, \xb(1)) \mid \xb(1) \in h(\xb) ]$ is the expected reward we get for choosing $\stop$ for every state in $h(\xb)$, conditional on $\xb(1)$ being in $h(\xb)$. From the definition of $\hat{\pi}$, the policy $\hat{\pi}$ chooses the action which gives the best of the two rewards. Conditional on $\xb(1) \in h(\xb)$, if we were to follow $\hat{\pi}$ in period $t = 1$ and then the optimal policy from $t = 2$ on, the reward we would get would be exactly
\begin{equation}
\max\{  \Exp[ g(1, \xb(1) ) \mid \xb(1) \in h(\xb) ],  \Exp[ \beta C^*_1(\xb(1)) \mid \xb(1) \in h(\xb) ] \}.
\end{equation}
The definition of $\hat{\pi}$ for $t > 1$ exhibits similar behavior. In particular, conditional on $\hat{\pi}$ not having stopped before period $t$ and $\xb(t)$ being in $h(\xb)$, the reward we get from period $t$ onwards by taking $\hat{\pi}$ at $t$ and then the optimal policy for $t' > t$ is 
\begin{equation}
\max\{ \Exp[ g(t, \xb(t) ) \mid \xb(t) \in h(\xb), A_t ],  \Exp[ \beta C^*_1(\xb(t)) \mid \xb(t) \in h(\xb), A_t ] \}.
\end{equation}
(Note that the above reward is discounted to period $t$.)

In addition to $\hat{\pi}$, we will also need to define a collection of \emph{intermediate policies}, denoted by $\tipi^0, \tipi^1,\dots, \tipi^T$. The intermediate policy $\tipi^i$ is the policy that follows $\hat{\pi}$ from period 1 to $i$, and then follows an optimal policy $\pi^*$ from $i+1$ onwards. %

The main result that we will prove is the following theorem.
\begin{theorem} \label{theorem:bound_Jstar_Jhatpi}
Let $m$ be any integer and let $\bar{\xb} \in \Xcal$ be the initial state of the system at $t = 1$. For $\hat{\pi}$ as defined above with respect to $H_m$, we have 
\begin{equation}
J^*(\bar{\xb}) - J^{\hat{\pi}}(\bar{\xb})  \leq \frac{2 T \lambda L \sqrt{n}}{m}.
\end{equation}

\end{theorem}
From this theorem, it will be possible to prove Theorem~\ref{theorem:ultimate_depth_bound}, which asserts the existence of a tree with an explicitly characterized depth that approximates the optimal policy to a desired accuracy $\epsilon$. We now provide a sketch of the strategy we will use to prove Theorem~\ref{theorem:bound_Jstar_Jhatpi}. 

\proof{Proof Sketch.}
First, we can write the difference $J^*(\bar{\xb}) - J^{\hat{\pi}}(\bar{\xb})$ using the intermediate policies $\tipi^0, \dots, \tipi^T$: 
\begin{equation*}
J^*(\bar{\xb}) - J^{\hat{\pi}}(\bar{\xb}) = (J^{\tipi_0}(\bar{\xb}) - J^{\tipi_1}(\bar{\xb})) + (J^{\tipi_1}(\bar{\xb}) -  J^{\tipi_2}(\bar{\xb})) + (J^{\tipi_2}(\bar{\xb}) - J^{\tipi_3}(\bar{\xb})) + \dots + (J^{\tipi_{T-1}}(\bar{\xb}) - J^{\tipi_T}(\bar{\xb})).
\end{equation*}
To bound $J^*(\bar{\xb}) - J^{\hat{\pi}}(\bar{\xb})$, we will bound each difference of successive intermediate policies, i.e., each difference $J^{\tipi_i}(\bar{\xb}) - J^{\tipi_{i+1}}(\bar{\xb})$. 

Each difference $J^{\tipi_i}(\bar{\xb}) - J^{\tipi_{i+1}}(\bar{\xb})$ involves two policies: $\tipi_i$ and $\tipi_{i+1}$. These two policies are similar: they follow $\hat{\pi}$ at first, and then follow an optimal policy. These policies differ by only one time period: the policy $\tipi_i$ follows $\hat{\pi}$ for periods 1 to $i$, and then an optimal policy  from period $i+1$ onward, whereas $\tipi_{i+1}$ follows $\hat{\pi}$ for one more period ($i+1$) before switching to an optimal policy. Along any sample path of $\{\xb(t)\}_{t=1}^T$, these two policies take the same actions for periods 1 to $i$, so we only need to understand what happens at period $i+1$, where $\tipi_i$ will take an optimal action, while $\tipi_{i+1}$ will take the same action as $\hat{\pi}$. 

Conditional on not stopping before $i+1$ and the state $\xb(i+1)$ being in a particular hypercube $h \in H_m$, the difference between the expected rewards garnered by $\tipi_{i+1}$ and $\tipi_i$ can be bounded by leveraging the Lipschitz smoothness of $g$ and $C^*_t$ (Lemma~\ref{lemma:Lipschitz_bound_for_single_box}). This bound on the reward difference, together with a simple bound on the conditional distribution of $\xb(i+1)$ (Lemma~\ref{lemma:density_bound_for_pihat}), can be aggregated across all $m^n$ hypercubes to obtain a bound on $J^{\tipi_i}(\bar{\xb}) - J^{\tipi_{i+1}}(\bar{\xb})$ (Lemma~\ref{lemma:bound_of_consecutive_tipi}). Adding up all of these bounds completes the proof. \Halmos
\endproof

The first lemma that we will prove ensures that the distribution of $\xb(t)$, the state at time $t$, remains bounded conditional on the policy $\hat{\pi}$ not having stopped by time $t$. The lemma follows in a straightforward manner from Assumptions~\ref{ass:Markovian} and \ref{ass:bounded_density}. 
\begin{lemma} \label{lemma:density_bound_for_pihat}
Suppose that $A_t = \{ \hat{\pi}(t', \xb(t')) = \go \ \text{for all}\ t' < t\}$. Then for any Borel set $S$, we have
\begin{equation}
\Pr( \xb(t) \in S \mid A_t) \leq \lambda \mu(S).
\end{equation}
\end{lemma}
\proof{Proof.}
We have:
\begin{align}
\Pr( \xb(t) \in S \mid A_t) & = \Exp_{ \xb(t-1) \mid A_t} \left[ \Pr( \xb(t) \in S \mid \xb(t-1), A_t) \right] \\
& = \Exp_{ \xb(t-1) \mid A_t} \left[ \Pr( \xb(t) \in S \mid \xb(t-1) ) \right] \\
& \leq \Exp_{ \xb(t-1) \mid A_t} \left[  \lambda \mu(S) \right] \\
&=  \lambda \mu(S),
\end{align}
where the first step follows by the tower property of conditional expectation; the second follows by Assumption~\ref{ass:Markovian}; and the inequality follows by Assumption~\ref{ass:bounded_density}. \Halmos
\endproof

Our second lemma provides a bound between the optimal policy and a constant policy via the Lipschitz smoothness of the reward function $g(t,\cdot)$ and the continuation value function $C^*_t(\cdot)$. 
\begin{lemma} \label{lemma:Lipschitz_bound_for_single_box}
Let $t \in [T]$ and let $\zb, \zb'$ be random variables in $h \in H_m$. Then
\begin{equation}
\Exp_\zb \left[ \max\{ g(t, \zb), \beta C^*_t(\zb)\} - \max\{ \Exp_{\zb'} [g(t,\zb')], \Exp_{\zb'}[ \beta C^*_t(\zb')] \} \right] \leq \frac{2L \sqrt{n}}{m}. 
\end{equation}
\end{lemma}

\proof{Proof.}
For notational convenience, let $V = \max\{ \Exp_{\zb'} [g(t,\zb')], \Exp_{\zb'} [ \beta C^*_t(\zb')] \}$, and let $B$ denote the event that stopping is optimal at the random state $\zb$:
\begin{equation}
B = \{ g(t, \zb) \geq \beta C^*_t(\zb) \} .
\end{equation}
Note that $B^C$ denotes the event that stopping is not optimal (i.e., continuing is optimal). We then have:
\begin{align*}
& \Exp_\zb[ \max\{ g(t, \zb), \beta C^*_t(\zb) \} - V ] \\[0.5em]
& = \Exp_\zb[ \max\{ g(t, \zb), \beta C^*_t(\zb) \} - V \mid B ] \cdot \Pr( B ) + \Exp_\zb[ \max\{ g(t, \zb), \beta C^*_t(\zb) \} - V \mid B^C ] \cdot \Pr( B^C ) \\
& = \Exp_\zb[ g(t, \zb) - V \mid B ] \cdot \Pr( B) + \Exp_\zb[ \beta C^*_t(\zb) - V \mid B^C ] \cdot \Pr( B^C ) \\[0.25em]
& \leq \Exp_\zb[ g(t, \zb) - \Exp_{\zb'}[g(t,\zb')] \mid B ]  \cdot \Pr( B)  \ + \ \Exp_{\zb}[ \beta C^*_t(\zb) - \Exp_{\zb'}[ \beta C^*_t(\zb') ] \mid B^C  ] \cdot \Pr( B^C)  \\
& \leq \Exp_\zb[ | g(t, \zb) - \Exp_{\zb'}[g(t,\zb')]|  \mid B ] \cdot \Pr( B)  \ + \ \Exp_{\zb}[  |\beta C^*_t(\zb) - \Exp_{\zb'}[ \beta C^*_t(\zb') ]|  \mid B^C ] \cdot \Pr( B^C )  \\
& \leq \Exp_{\zb, \zb'}[ | g(t, \zb) - g(t,\zb')| \mid B]  \ + \ \Exp_{\zb, \zb'}[  |\beta C^*_t(\zb) - \beta C^*_t(\zb') | \mid  B^C ] \\
& \leq \Exp_{\zb, \zb'}[ L \| \zb - \zb' \|  \mid B] +  \Exp_{\zb, \zb'}[ \beta L \| \zb - \zb' \| \mid B^C] \\ 
& \leq L \cdot \frac{\sqrt{n}}{m} \ + \ \beta L \cdot \frac{\sqrt{n}}{m}  \\
& \leq \frac{2L \sqrt{n}}{m},
\end{align*}
where the first equality follows by the definition of conditional expectation; the second equality follows by the definition of the events $B$ and $B^C$; the first inequality follows by the definition of $V$ as a maximum of $\Exp_{\zb'}[g(t,\zb')]$ and $\Exp_{\zb'}[ \beta C^*_t(\zb')]$; the second inequality follows by the fact that $a \leq |a|$ for any real $a$; the third inequality by Jensen's inequality; the fourth inequality by Assumption~\ref{ass:g_C_Lipschitz}; the fifth inequality by the fact that the maximum Euclidean distance of any two points in the hypercube $h$ is bounded by $\sqrt{n}/m$; and the final inequality by the fact that $\beta \leq 1$.  \Halmos %
\endproof

Our next result uses Lemmas~\ref{lemma:density_bound_for_pihat} and \ref{lemma:Lipschitz_bound_for_single_box} to obtain a bound on the difference in reward between two successive intermediate policies.
\begin{lemma}
Fix $\bar{\xb} \in \Xcal$. For every $j \in \{0,1, \dots, T -1\}$, let $\tipi^j$ be the policy that follows $\hat{\pi}$ from period 1 to period $j$, and then follows $\pi^*$ from $j+1$ to $T$. For any $i \in \{0,1,\dots, T-1\}$, we have
\begin{equation}
J^{\tipi_i}(\bar{\xb}) - J^{\tipi_{i+1}}(\bar{\xb}) \leq \frac{2 L \lambda \sqrt{n}}{m}.
\end{equation}

\label{lemma:bound_of_consecutive_tipi}
\end{lemma}

\proof{Proof.}
For notational convenience, let $\Delta$ be the random variable of the difference in reward between the two policies:
\begin{equation}
\Delta = \beta^{\tau_{\tipi_i} - 1} \cdot g( \tau_{\tilde{\pi}_i}, \xb( \tau_{\tipi_i}) ) - \beta^{\tau_{\tipi_{i+1} } - 1 } \cdot g( \tau_{\tipi_{i+1}}, \xb( \tau_{\tipi_{i+1} }) ).
\end{equation}

Let us also define the event $A_{i+1}$ as 
\begin{equation}
A_{i+1} = \{ \hat{\pi}(t', \xb(t')) = \go \ \text{for all}\ t' < i+1\}.
\end{equation}
We then have
\begin{align*}
J^{\tipi_i}(\bar{\xb}) - J^{\tipi_{i+1}}(\bar{\xb}) & = \Exp[ \Delta ] \\
& = \Exp[ \Delta \mid A_{i+1} ] \cdot \Pr( A_{i+1}) + \Exp[ \Delta \mid A^C_{i+1} ] \cdot \Pr( A^C_{i+1}) \\
& = \Exp[ \Delta \mid A_{i+1} ] \cdot \Pr( A_{i+1}) \\
& \leq \Exp[ \Delta \mid A_{i+1} ]
\end{align*}
In the above sequence of steps, the last equality follows by the fact that when the event $A_{i+1}$ does not occur, then the policy $\hat{\pi}$ stopped in a period between 1 and $i$; since both $\tipi_i$ and $\tipi_{i+1}$ follow policy $\hat{\pi}$ from for periods up to and including $i$, the difference in rewards $\Delta$ must be zero in this case. The inequality follows because $J^{\tipi_i}(\bar{\xb}) - J^{\tipi_{i+1}}(\bar{\xb})$, by the definition of $\tipi_i$ and $\tipi_{i+1}$, is nonnegative ($\tipi_i$ follows the optimal policy from $i+1$ on, whereas $\tipi_{i+1}$ only follows it from $i+2$ on; $\tipi_i$ must attain a higher reward than $\tipi_{i+1}$). Since $J^{\tipi_i}(\bar{\xb}) - J^{\tipi_{i+1}}(\bar{\xb})$ is nonnegative, $\Exp[ \Delta \mid A_{i+1} ]$ must be nonnegative. 

We shall now bound $\Exp[ \Delta \mid A_{i+1} ]$. Before doing so, let us define for notational convenience the random variable $W$ as 
\begin{align}
W & = \Ibb\{ \hat{\pi}(i+1, \xb(i+1)) = \stop\} \cdot g(i+1, \xb(i+1))  \nonumber \\
& \phantom{=}  +  \Ibb\{ \hat{\pi}(i+1, \xb(i+1)) = \go\} \cdot \beta C^*_{i+1}( \xb(i+1))
\end{align}
The random variable $W$ represents the expected reward of $\tipi_{i+1}$, conditional on the state of the system at time $i+1$ (the random variable $\xb(i+1)$). Note that if the policy $\hat{\pi}$ selects $\stop$, then the expected reward garnered is just $g(i+1, \xb(i + 1))$; if the policy selects $\go$, then the expected reward is $\beta C^*_{i+1}( \xb(i+1))$, which is the expected reward from following the optimal policy from period $i+2$ on. Note that $W$ is the expected reward discounted to period $i+1$; the quantity $\beta^{i} W$ represents this expected reward discounted to the first period $t = 1$. 

With $W$ defined, we can write $\Exp[ \Delta \mid A_{i+1} ]$ in terms of the hypercube set $H_m$: 
\begin{align}
& \phantom{=} \Exp[ \Delta \mid A_{i+1} ]  \nonumber \\
& = \Exp[ \beta^{i} \max\{ g(i+1, \xb(i+1)), \beta C^*_{i+1}(\xb(i+1))\}  - \beta^{i}  W \mid A_{i+1}] \nonumber \\
& =  \beta^{i} \Exp[  \max\{ g(i+1, \xb(i+1)), \beta C^*_{i+1}(\xb(i+1))\}  - W \mid A_{i+1}] \nonumber \\
& =  \beta^{i}  \sum_{h \in H_m} \Exp[  \max\{ g(i+1, \xb(i+1) ), \beta C^*_{i+1}(\xb(i+1))\}   - W \mid A_{i+1}, \xb(i+1) \in h]  \cdot \Pr(\xb(i+1) \in h \mid A_{i+1}) \label{eq:sum_over_boxes} %
\end{align}

We now bound the conditional expectations in the sum. Note that, by the definition of $\hat{\pi}(i+1, \xb)$ as picking the best action in each hypercube that maximizes the expected reward conditional on reaching that hypercube, we will have that 
\begin{equation}
\Exp[ W \mid A_{i+1}, \xb(i+1) \in h] = \max\{ \Exp[g( i+1, \yb )], \beta \Exp[ C^*_{i+1}(\yb) ] \},
\end{equation}
where $\yb$ is a random variable that follows the conditional distribution of $\xb(i+1)$ given the events $A_{i+1}$ and $\xb(i+1) \in h$. With this definition, we have:
\begin{align*}
& \Exp[  \max\{ g(i+1, \xb(i+1)\}, \beta C^*_{i+1}(\xb(i+1))\}   - W \mid A_{i+1}, \xb(i+1) \in h]  \\
& = \Exp[  \max\{ g(i+1, \xb(i+1)\}, \beta C^*_{i+1}(\xb(i+1))\}   - \max\{ \Exp[g( i+1, \yb )], \beta \Exp[ C^*_{i+1}(\yb) ] \} \mid A_{i+1}, \xb(i+1) \in h] \\
& \leq \frac{2L \sqrt{n}}{m},
\end{align*}
where the inequality follows by Lemma~\ref{lemma:Lipschitz_bound_for_single_box}. Returning to \eqref{eq:sum_over_boxes}, we thus have 
\begin{align*}
\Exp[ \Delta \mid A_{i+1} ] & \leq \beta^{i}  \sum_{h \in H_m} \frac{2L \sqrt{n}}{m} \cdot \Pr(\xb(i+1) \in h \mid A_{i+1}) \\
& \leq \beta^{i}  \sum_{h \in H_m} \frac{2L \sqrt{n}}{m} \cdot \frac{\lambda}{ m^n }  \\
& = \beta^{i} m^n \cdot \frac{2L \sqrt{n}}{m} \cdot \frac{\lambda}{ m^n }  \\
& \leq \frac{2 \lambda L \sqrt{n}}{m}
\end{align*}
where the second inequality follows by Lemma~\ref{lemma:density_bound_for_pihat}; the second equality follows by the fact that $H_m$ contains $m^n$ boxes; and the last equality by the fact that $\beta^i < 1$. Thus, we have established that $J^{\tipi_i}(\bar{\xb}) - J^{\tipi_{i+1}}(\bar{\xb}) \leq 2 \lambda L \sqrt{n} / m$, which is the required result. \Halmos
\endproof

With Lemma~\ref{lemma:bound_of_consecutive_tipi} in hand, we are ready to prove Theorem~\ref{theorem:bound_Jstar_Jhatpi}. 
\proof{Proof of Theorem~\ref{theorem:bound_Jstar_Jhatpi}.}
Write the difference in the policy performances $J^*(\bar{\xb}) - J^{\hat{\pi}}(\bar{\xb})$ as the following telescoping sum:
\begin{align*}
J^*(\bar{\xb}) - J^{\hat{\pi}}(\bar{\xb}) & = J^{\tipi_0}(\bar{\xb}) - J^{\tipi_T}(\bar{\xb}) \\
& = J^{\tipi_0}(\bar{\xb}) - J^{\tipi_1}(\bar{\xb}) +  J^{\tipi_1}(\bar{\xb}) -  J^{\tipi_2}(\bar{\xb}) + J^{\tipi_2}(\bar{\xb}) - J^{\tipi_3}(\bar{\xb}) + \dots + J^{\tipi_{T-1}}(\bar{\xb}) - J^{\tipi_T}(\bar{\xb}) \\
& \leq \frac{2 \lambda L \sqrt{n}}{m} + \frac{2 \lambda L \sqrt{n}}{m} + \frac{2 \lambda L \sqrt{n}}{m} + \dots + \frac{2 \lambda L \sqrt{n}}{m} \\
& = \frac{2 T \lambda L \sqrt{n}}{m},
\end{align*}
where the inequality follows by Lemma~\ref{lemma:bound_of_consecutive_tipi}.  \Halmos
\endproof

Our final step in the analysis is to establish Theorem~\ref{theorem:ultimate_depth_bound}.

\begin{theorem}
For any $\epsilon > 0$ and state $\bar{\xb} \in \Xcal$, there exists a tree policy $\hat{\pi}$ such that $J^*(\bar{\xb}) - J^{\pi}(\bar{\xb}) \leq \epsilon$ and the depth of the corresponding tree is at most $T + n \lceil \log_2( \frac{2TL \lambda \sqrt{n}}{\epsilon} )\rceil$. 
\label{theorem:ultimate_depth_bound}
\end{theorem}

\proof{Proof.}
We will show that the policy $\hat{\pi}$ exactly satisfies the requirement of the theorem. 
For a fixed integer $m$, the policy $\hat{\pi}$ as described above requires a tree of depth $T + n \lceil \log_2(m)\rceil$. (This follows because one requires $T$ splits on $t$ to ensure that the policy correctly discerns the period $t$; then, for the set of hypercubes $H_m$, one needs $n \lceil \log_2(m) \rceil$ levels of splits. The deepest leaf will be at a depth of at most $T + n \log_2(m)$.) Substituting $m = 2TL \lambda \sqrt{n} / \epsilon$ into $T + n \lceil \log_2(m)\rceil$ and the bound of Theorem~\ref{theorem:bound_Jstar_Jhatpi}, the result follows.  \Halmos
\endproof

\clearpage

\pagebreak

\section{Numeric examples for construction algorithm}

\label{sec:examples}

\subsection{Example of complete construction algorithm}

\label{sec:toy_example}

In this section, we provide a small demonstration of our construction algorithm (Algorithm~\ref{algorithm:construction}) in Section~\ref{sec:construction}. 

\begin{figure}
\centering
\includegraphics[width=0.9\textwidth]{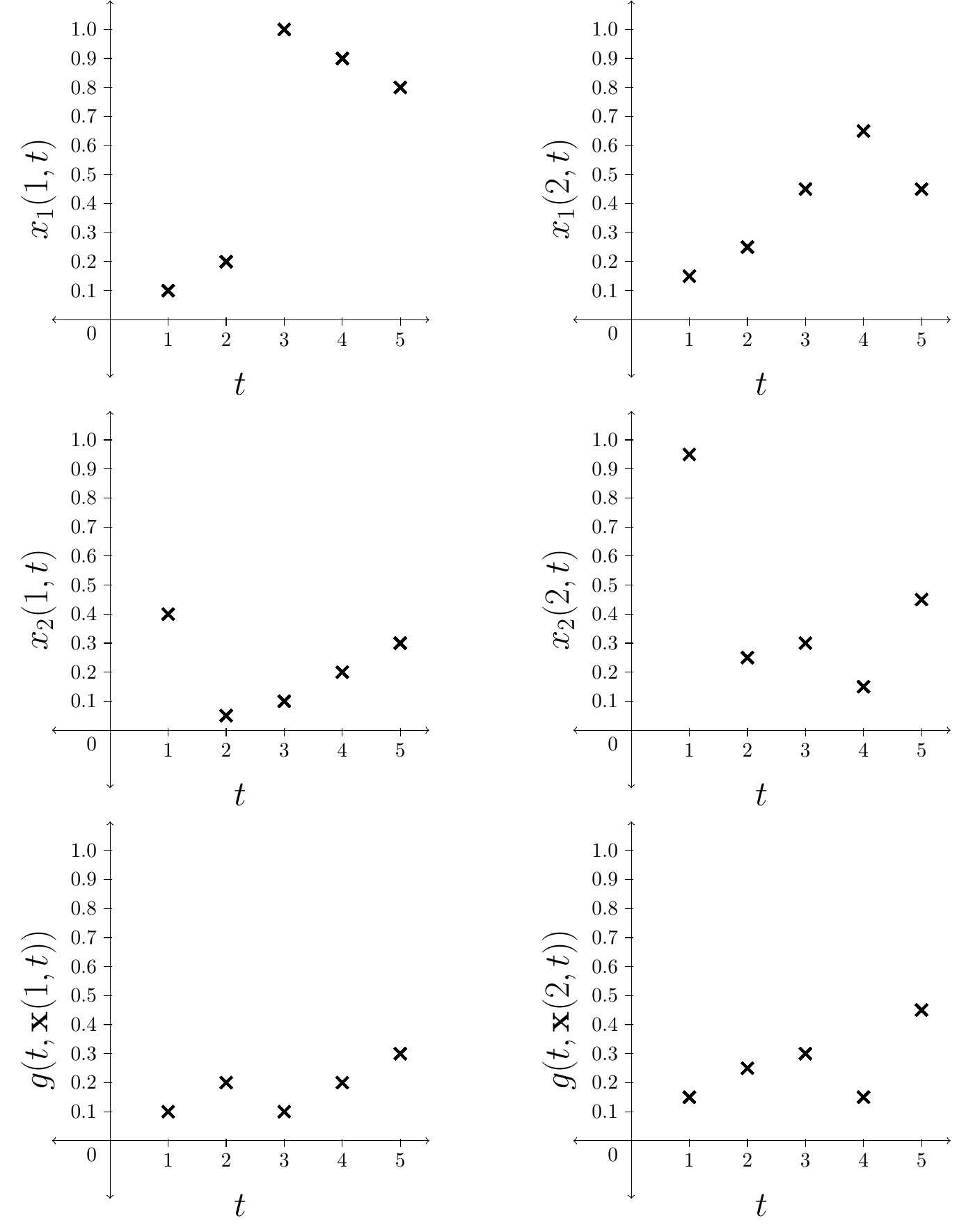}
\caption{Trajectories for Example~\ref{example:ToyExample} of overall construction algorithm. Each column corresponds to a different trajectory ($\omega = 1, 2$); the top two plots show the two state variables, while the bottom plot shows the reward in each trajectory. \label{figure:toy_example_trajectories} }
\end{figure}

\begin{figure}
\centering
\includegraphics[width=0.8\textwidth]{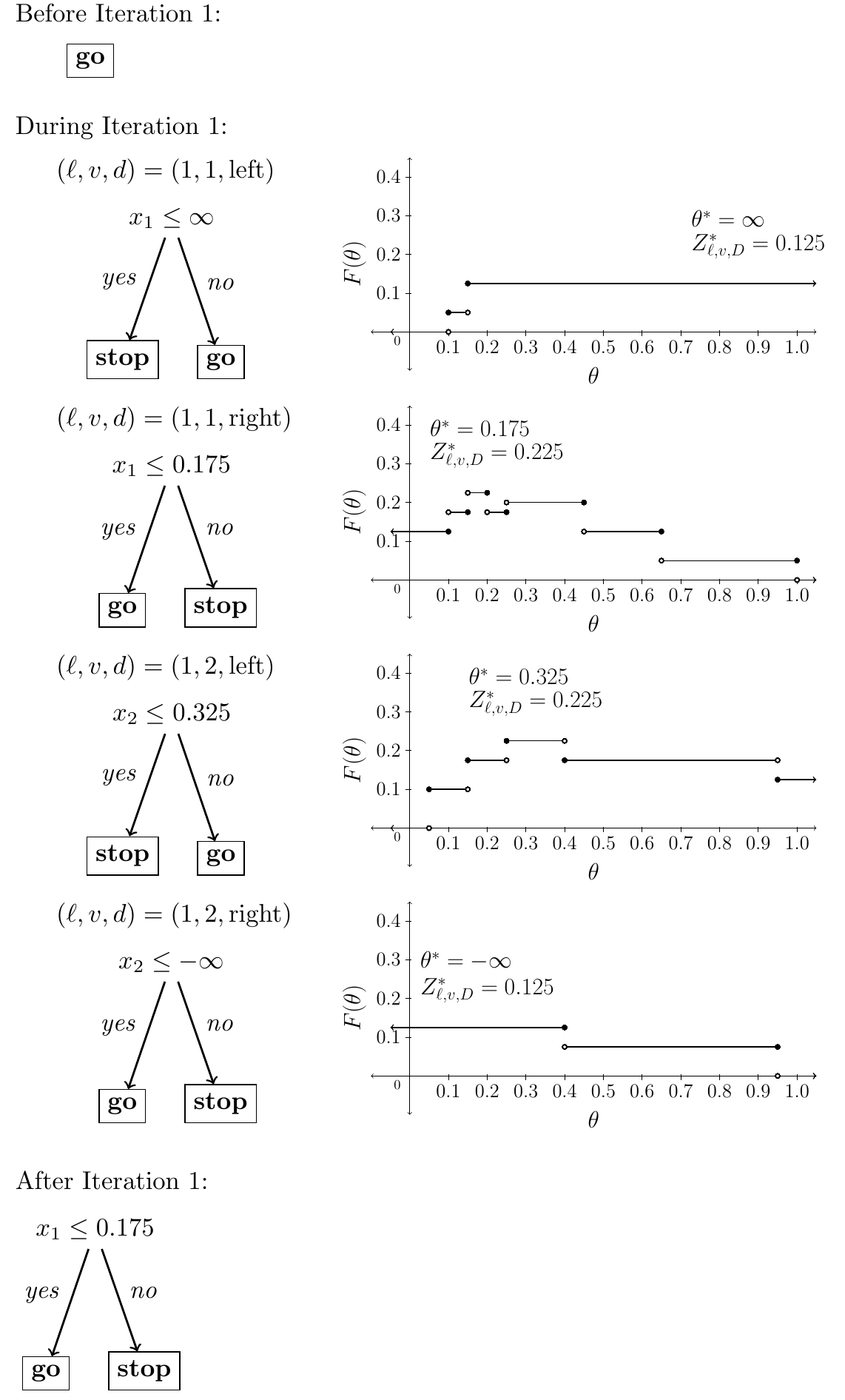}
\caption{Iteration 1 in Example~\ref{example:ToyExample} of overall construction algorithm. \label{figure:toy_example_iteration1} }
\end{figure}

\begin{sidewaysfigure}
\centering
\vspace{45em}
\includegraphics[width = 0.9\textwidth]{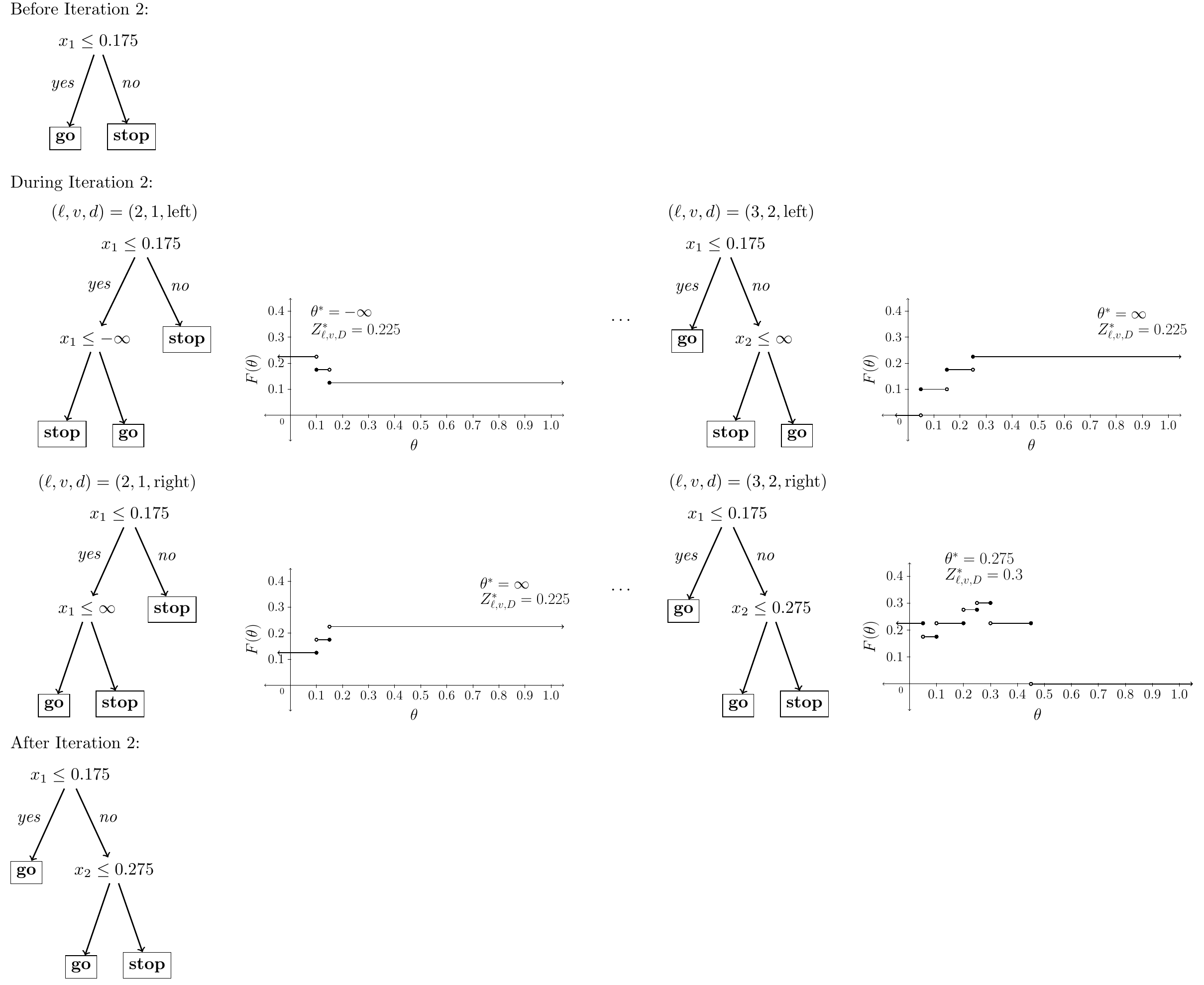}
\caption{Iteration 2 in Example~\ref{example:ToyExample} of overall construction algorithm. Due to space limitations, trees for five out of the total eight $(\ell,v,D)$ combinations are omitted. Note: leaf nodes 2 and 3 are the left child and right child, respectively, of the tree under ``Before Iteration 2''. \label{figure:toy_example_iteration2}}
\end{sidewaysfigure}

\begin{example}
\label{example:ToyExample}
Consider a small example with $n = 2$ state variables, $T = 5$ time periods, and $\Omega = 2$ trajectories. The two trajectories, together with their rewards, are displayed in Figure~\ref{figure:toy_example_trajectories}. We assume the discount factor $\beta$ is set to 1 and that $\gamma$ (the relative improvement parameter for the construction algorithm) is set to zero. 

Figure~\ref{figure:toy_example_iteration1} shows the first iteration of the construction algorithm. At the top, the initial tree is shown, which is simply a degenerate tree that always chooses $\go$. As discussed in Section~\ref{subsec:construction_description}, the construction algorithm considers replacing every leaf in the tree with a split. This results in a collection of candidate trees. Each such candidate tree corresponds to a tuple $(\ell, v, D)$, where $\ell$ is a leaf, $v$ is a split variable index and $D$ is a subtree orientation (left or right, to indicate which child leaf corresponds to the $\stop$ action). Here, because there is only one leaf, two split variables and two subtree orientations, there are exactly four candidate trees, which are shown in the middle of the figure. For each tree shown, the split point chosen is the optimal split point, using the procedure in Section~\ref{subsec:construction_OptimizeSplitPoint} (Algorithm~\ref{algorithm:OptimizeSplitPoint}). To the right of each candidate tree, a plot of $F(\theta)$ is provided, which completely describes the objective value of the tree as a function of the split point $\theta$ chosen for the split; the optimal $\theta^*$ and the objective value $Z^*_{\ell, v,D}$ for this choice of $\theta$ are also shown in the plot. After all of these candidate trees are evaluated, we replace the initial tree (which was only the leaf with the $\go$ action), with the subtree corresponding to $(\ell, v, D) = (1,1,\text{right})$. Note that in this particular example, the candidate tree for $(1,2,\text{left})$ achieves the same objective value of 0.225 and is also optimal, and thus we could have chosen this subtree as well. In such cases, we break ties arbitrarily. At the conclusion of the iteration, the new tree policy is one where the split variable is $x_1$, the split point is 0.175 and the right child is the $\stop$ action (shown under ``After Iteration 1''). 

In the next iteration, starting from this tree, we repeat the same procedure, and we again consider replacing every leaf with a split. Figure~\ref{figure:toy_example_iteration2} shows the (partial) steps involved in this iteration. In this iteration, because we started from a tree with two leaves (shown under ``Before Iteration 2''), there are eight possible candidate trees; due to space limitations, we show only four candidate trees, for the tuples $(2,1,\text{left})$, $(2,1,\text{right})$, $(3,2,\text{left})$ and $(3,2,\text{right})$. The optimal tree turns out to be this last tree, and becomes our tree for the third iteration. 

In the third iteration (not shown), all candidate trees achieve the same objective value as the tree following iteration 2. After this iteration, the algorithm terminates with the tree at the bottom of Figure~\ref{figure:toy_example_iteration2}. \Halmos
\end{example}
\clearpage

\pagebreak
\clearpage

\pagebreak
\clearpage

\clearpage
\pagebreak

\subsection{Example of \textsc{OptimizeSplitPoint} function}

\label{sec:example_optimizesplitpoint}

In this section, we provide a demonstration of the calculations involved in the \textsc{OptimizeSplitPoint} function on a simple example.

\begin{example}
\label{example:2}
In this example, let us suppose that we are at the first iteration of the construction algorithm, where the tree policy is simply a single leaf with the action $\go$. We want to find the optimal split point with respect to variable 1 for a split on this leaf, assuming that the subtree we place will be a right-stop subtree. 

To illustrate, we fix a single trajectory $\omega$. Figure~\ref{figure:CEX1_data} shows the relevant data for this trajectory, which are the values $x_1(\omega,t)$ and the rewards $g(t, \xb(\omega,t) )$. We assume that $T = 18$ for all trajectories.

\begin{figure}
\centering
\begin{tabular}{c}
\includegraphics[width=0.6\textwidth]{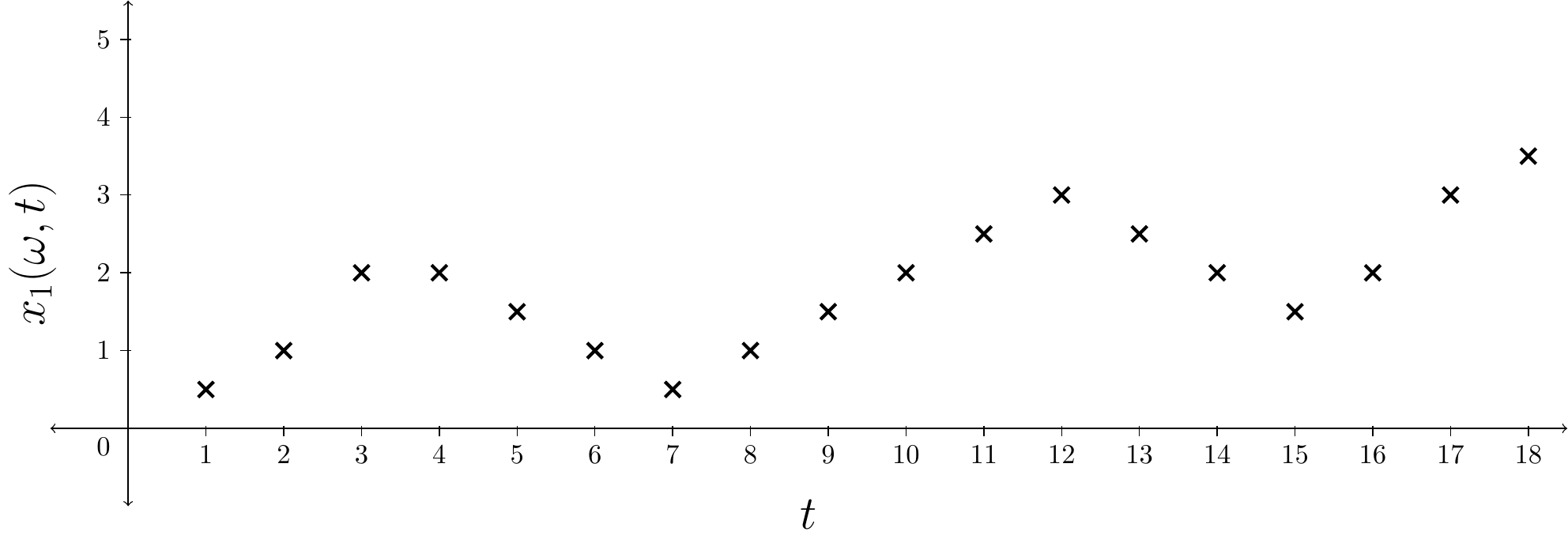}\\
\includegraphics[width=0.6\textwidth]{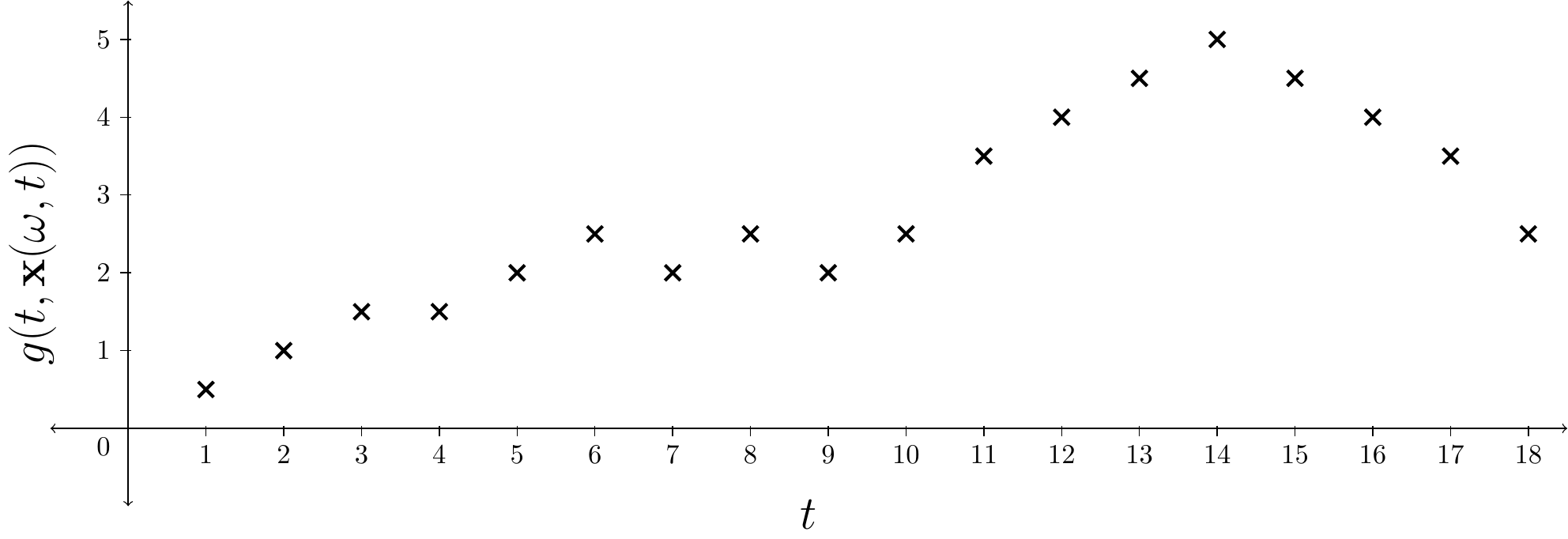}
\end{tabular}
\caption{Data for Example \ref{example:2}. The top plot shows coordinate 1 ($x_1$) of a single trajectory $\omega$, and the bottom plot shows the corresponding reward ($g$). \label{figure:CEX1_data} }
\end{figure}

Since we are at the first iteration, we obtain that for the root node $\ell$, the no-stop time is $\tau_{-\ell,\omega} = +\infty$, because there are no other leaves in the tree, and thus we have that $S_{\omega} = \{1,\dots, 18\}$, i.e., every period is a valid in-leaf period for the root node. Observe also that the no-stop value $f_{\omega, \text{ns}} = 0$, because there is no other leaf in which the trajectory is stopped.

Having determined the in-leaf periods, we determine the set of permissible stop periods $P_{\omega}$. To do this, since we are placing a right-stop subtree, we follow the computation in equation~\eqref{eq:permissible_stop_rightstop_subtree}. The left-hand side of Figure~\ref{figure:CEX1_PS_fomega} shows the same data as in Figure~\ref{figure:CEX1_data}, but with the permissible stop periods indicated with red squares. 

To intuitively understand the computation in equation~\eqref{eq:permissible_stop_rightstop_subtree}, one can imagine a horizontal line, superimposed on the top plot of Figure~\ref{figure:CEX1_data}, that starts at $-\infty$ and slowly moves up towards $+\infty$. The vertical height of this horizontal line is the split point $\theta$. As we move this line, we track the first period in time at which the trajectory exceeds this horizontal line: this is where the trajectory would stop, if we fixed the split point to the height of that line. Notice that when we start, the first period is $t = 1$. As soon as we exceed $x_1(\omega,1) = 0.5$, we will no longer stop at $t = 1$, but at $t = 2$. As soon as our line goes above $x_1(\omega,2)$, the new period at which we stop will be $t = 3$. Once our line goes above $x_1(\omega,3)$, we will stop at $t = 11$, as this is the earliest that the trajectory passes above our line. We keep going in this way, until the line exceeds the value $\max_{t} x_1(\omega,t)$ and we have determined all of the permissible stop periods. %

To further elaborate on this process, notice that $t = 1, 2, 3$ are permissible stop periods, but period $t = 5$ is not. In order to stop at $t = 5$, we would have to set the threshold $\theta$ to a value lower than $x_1(\omega,5) = 1.5$, and $t = 5$ would have to be the first period at which the trajectory exceeds $\theta$, i.e., $x_1(\omega,5) > \theta$. This is impossible with the trajectory of Figure~\ref{figure:CEX1_data}. For this reason, the values of $x_1$ and $g$ at $t = 5$ are not relevant in determining the reward garnered from the trajectory as $\theta$ is varied.

\begin{figure}
\centering
\includegraphics[width=0.8\textwidth]{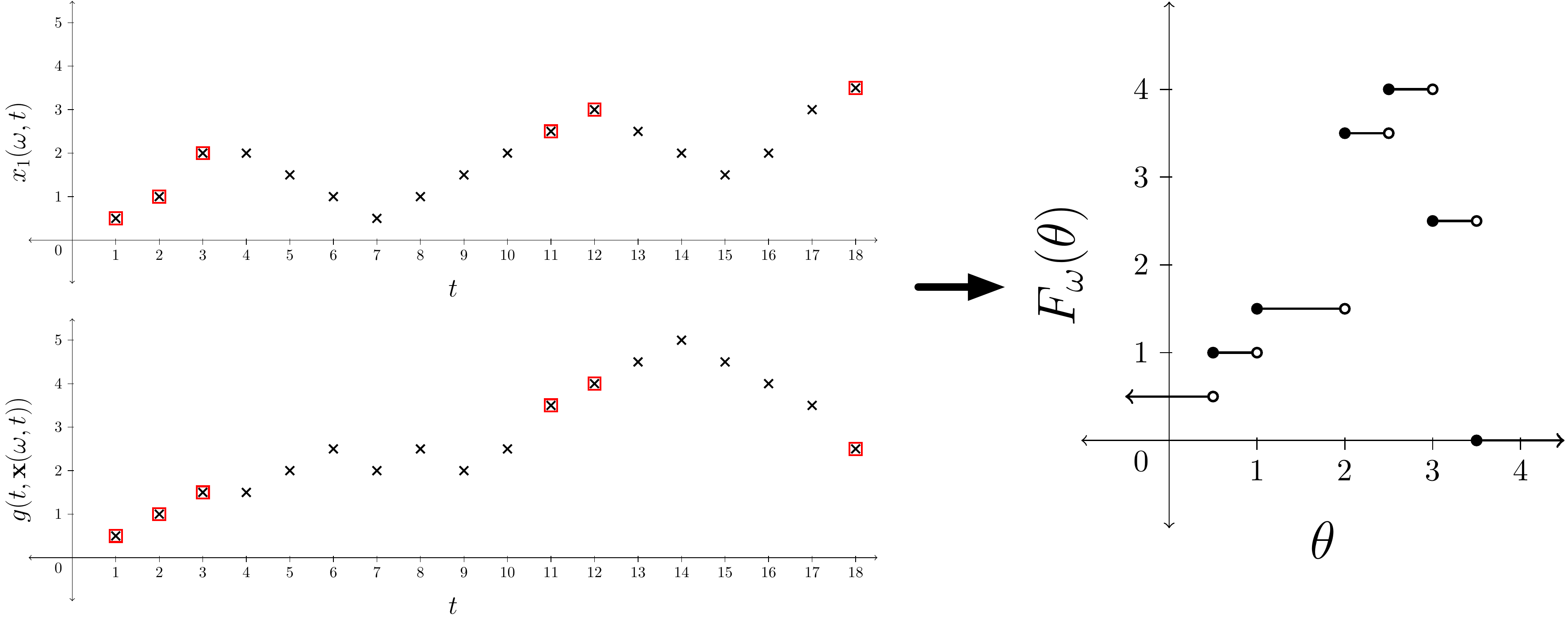}
\caption{Process for creating the function $F_{\omega}$. The left-hand side shows $x_1(\omega,t)$ and $g(t, \xb(\omega,t))$ with the permissible stop periods indicated by red squares, while the right-hand side shows the function $F_{\omega}(\theta)$ that corresponds to this trajectory. \label{figure:CEX1_PS_fomega} }
\end{figure}

The permissible stop periods $P_{\omega} = \{t_1, t_2, \dots, t_{|P_{\omega}|} \}$ allow us to define the breakpoints $b_{\omega,1},\dots, b_{\omega, |P_{\omega}|}$ and the values $f_{\omega,1}, \dots, f_{\omega, |P_{\omega}|}$ and $f_{\omega,\text{ns}}$ of our piecewise constant function. The corresponding piecewise constant function is shown on the right-hand side of Figure~\ref{figure:CEX1_PS_fomega}. We then repeat this construction for all of the trajectories in our training set, and average them to obtain the overall function $F(\cdot)$. Figure~\ref{figure:CEX1_F_schematic} visualizes this averaging process, assuming that there are three trajectories in total (i.e., $\Omega = 3$), including the one trajectory displayed in Figures~\ref{figure:CEX1_data} and \ref{figure:CEX1_PS_fomega} (which we denote by $\omega = 1$).

\begin{figure}
\centering
\begin{tabular}{cccc}
\includegraphics[width=0.24\textwidth]{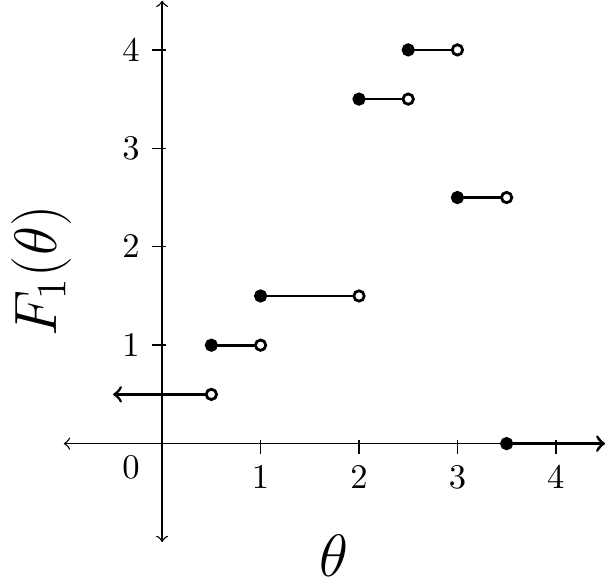} & \includegraphics[width=0.24\textwidth]{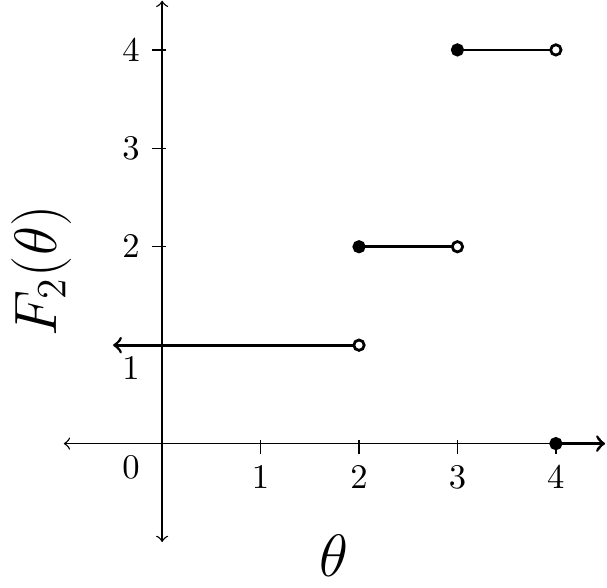} & \includegraphics[width=0.24\textwidth]{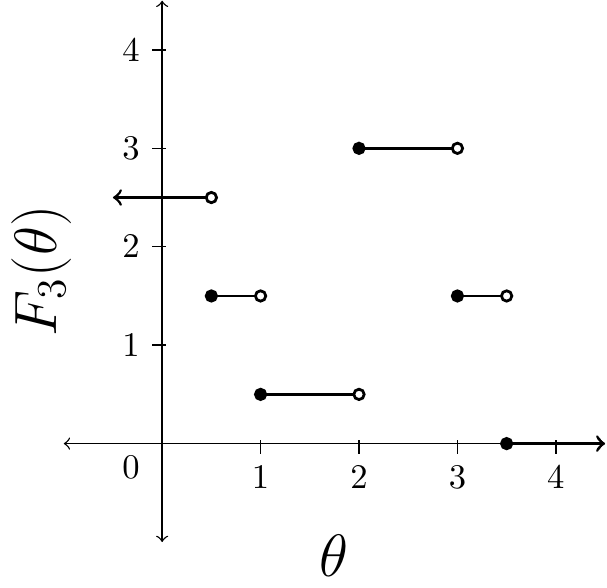} & \includegraphics[width=0.24\textwidth]{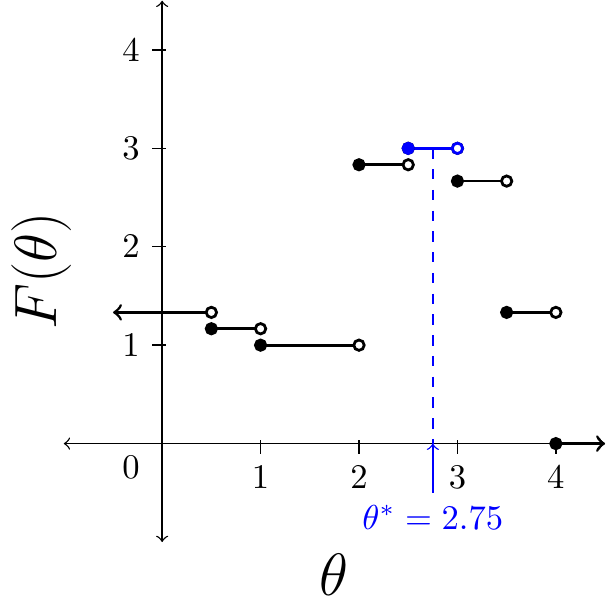} \\
\end{tabular}

\caption{Process for creating the function $F$ from $F_{1}, F_2, F_3$. The left three graphs show the trajectory-specific functions $F_1, F_2, F_3$ while the right-most graph shows the overall sample average function $F(\theta)$, obtained by averaging $F_1,F_2, F_3$. On the plot of $F$, the blue piece of the function is the one on which $F$ is maximized; the midpoint of that interval is $\theta^* = 2.75$, which is used as the final split point. \label{figure:CEX1_F_schematic} }
\end{figure}

After $F_1, F_2, F_3$ are averaged, we determine $\arg \max_{\theta \in \mathbb{R}} F(\theta)$ to be $[2.5, 3)$. Since this is a bounded interval, we take the midpoint of this interval, which is 2.75, to be our final split point $\theta^*$.  \Halmos

\end{example}

\pagebreak
\clearpage

\section{Procedure for selecting $\gamma$ via $k$-fold cross-validation}
\label{sec:kfold_cross_validation}

In this section, we present a procedure for selecting $\gamma$ using $k$-fold cross-validation, where $k$ is a positive integer. For simplicity, we assume that the number of trajectories $\Omega$ is such that $\Omega / k = \Omega'$ is a positive integer, and that the folds comprise disjoint subsets of $[\Omega]$.

A simple way that one might perform cross-validation is to pick a grid of $\gamma$ values, estimate a tree policy using the $k - 1$ training folds with each $\gamma$ value, evaluate the tree policy's objective on the hold-out fold, and repeat for each of the $k$ folds. One would then average the hold-out fold objective over all choices of hold-out fold, and then pick the grid value of $\gamma$ that gives the best average performance on the hold-out fold. While one can certainly take this approach, it may be unattractive for a couple of reasons. First, it requires one to pre-select the grid of $\gamma$ values; it is not clear how one should do this so as to not omit values of $\gamma$ that may lead to good cross-validated performance. Second, it is computationally wasteful. Suppose that we estimate a tree from a data set with $\gamma = 0.1$, and we then estimate a tree from the same data set with $\gamma = 0.05$; in this case, the second run of the construction algorithm with  $\gamma = 0.05$ will actually repeat the same iterations as the first run with $\gamma = 0.1$ until the relative improvement goes below $\gamma = 0.1$. 

These two observations suggest a different way that we might run the construction algorithm. We present the pseudocode of this procedure below as Algorithm~\ref{algorithm:CalculateKFoldCVBreakpoints}. We remark that in the definition of this algorithm, for a given fold $i$ out of $k$, the function $Z_i(\Tcal, \vb, \thetab, \ab)$ represents the sample average objective of the tree policy defined by $(\Tcal, \vb, \thetab, \ab)$ evaluated on the hold-out fold $i$. We also note that the function $\textsc{OptimizeSplitPoint}$ has an additional input parameter $i$, indicating that it is executed using the trajectories not in the $i$th fold. 

\begin{algorithm}
{\SingleSpacedXI
\begin{algorithmic}
\REQUIRE $\gamma_{\min} > 0$. 
\STATE Define $S_i = \emptyset$ for $i \in [k]$.
\FOR{$i \in [k]$ }
    \STATE $\bar{\gamma} = \infty$. 
    \STATE  $\Tcal \leftarrow ( \{1\}, \{1\}, \emptyset, \leftchild, \rightchild )$, $\vb \leftarrow \emptyset$, $\thetab \leftarrow \emptyset$, $a(1) \leftarrow \go$
    \STATE  $Z \leftarrow 0$
    \STATE  $\textbf{existsImprovement} \leftarrow \textbf{true}$ \\[1em]
    \WHILE{ \textbf{existsImprovement} }
    	\FOR{$\ell \in \leaves$, $v \in [n]$, $\dir \in \{ \text{left}, \text{right} \}$ }
    				\STATE $Z^*_{\ell,v, \dir }, \theta^*_{\ell,v, \dir} \leftarrow \textsc{OptimizeSplitPoint}(i, \ell, v, \dir; \Tcal, \vb, \thetab, \ab)$
    	\ENDFOR\\[1em]
    	$\textbf{existsImprovement} \leftarrow \Ibb \left[\max_{\ell, v, \dir} Z^*_{\ell,v, \dir} \geq (1 + \gamma_{\min})  Z \right]$\\[1em]

    	\IF{ $\max_{\ell, v, \dir} Z^*_{\ell,v, \dir} >  Z$}
    		\STATE ($\ell^*, v^*, \dir^*) \leftarrow \arg \max_{\ell, v, D} Z^*_{\ell,v, \dir }$
    		\STATE $\textsc{GrowTree}(\Tcal, \ell^*)$
    		\STATE $v(\ell^*) \leftarrow v^*$
    		\STATE $\theta(\ell^*) \leftarrow \theta^*_{\ell^*, v^*, \dir^*}$ \\[0.5em]
    		\IF{ $\dir^* = \text{left}$ }
    		 	\STATE $a( \leftchild(\ell^*) ) \leftarrow \stop$,  $a( \rightchild(\ell^*) ) \leftarrow \go$
    		\ELSE
    		 	\STATE $a( \leftchild(\ell^*) ) \leftarrow \go$, $a( \rightchild(\ell^*) ) \leftarrow \stop$
    		\ENDIF\\[1em]
		
		\IF{ $(1 + \bar{\gamma}) > (\max_{\ell, v, \dir} Z^*_{\ell,v, \dir} ) / Z$}
			\STATE $\bar{\gamma} \leftarrow (\max_{\ell, v, \dir} Z^*_{\ell,v, \dir} ) / Z - 1$. 
			\STATE $Z_h \leftarrow Z_{i}(\Tcal, \vb, \thetab, \ab)$ 
			\STATE $S_i \leftarrow S_i \cup \{ (\bar{\gamma}, Z_h) \}$
		\ENDIF
    		\STATE $Z \leftarrow Z^*_{\ell^*, v^*, \dir^*}$ 
    	\ENDIF
    \ENDWHILE
\ENDFOR
\RETURN Breakpoint sets $S_1, \dots, S_k$.  
\end{algorithmic}
}
\caption{ \textsc{CalculateKFoldCVBreakpoints} function. \label{algorithm:CalculateKFoldCVBreakpoints}}
\end{algorithm}

In this procedure we specify a minimum $\gamma$ denoted by $\gamma_{\min}$, which is the lowest value of $\gamma$ that we will consider. For each fold $i$ of the $k$ folds, we then execute the construction algorithm until reaching the relative improvement tolerance of $\gamma_{\min}$. Through each run of the construction algorithm, we maintain a variable $\bar{\gamma}$ which tracks the lowest relative improvement seen thus far. With each iteration of the construction algorithm, we check whether the relative improvement was below $\bar{\gamma}$; if so, we update $\bar{\gamma}$, we run the current tree on the hold-out fold and calculate the hold-out fold objective $Z_{h}$, and we record the $(\bar{\gamma}, Z_{h})$ pair in a set $S_i$ corresponding to the fold $i$. Upon termination, Algorithm~\ref{algorithm:CalculateKFoldCVBreakpoints} returns the sets $S_1,\dots, S_k$. 

The idea of the set $S_i$ is that each tuple in $S_i$ corresponds to a relative improvement tolerance $\gamma$ at which the construction algorithm may terminate for the current fold. Each point in $S_i$ therefore corresponds to a point on a piecewise constant function that represents the hold-out objective for fold $i$ as a function of $\gamma$. For each $i$, let us write the points in $S_i$ as 
\begin{equation}
S_i = \{ (\bar{\gamma}_1, Z_{h,1}), \dots, (\bar{\gamma}_{M_i}, Z_{h,M_i}) \}
\end{equation}
such that $\bar{\gamma}_1 > \bar{\gamma}_2 > \dots > \bar{\gamma}_{M_i}$. We then define the function $\nu_i$ as 
\begin{equation}
\nu_i(\gamma) = \left\{ \begin{array}{ll} 
Z_{h,1} & \text{if} \ \gamma > \bar{\gamma}_1, \\
Z_{h,2} & \text{if} \ \bar{\gamma}_1 \geq \gamma > \bar{\gamma}_2, \\
Z_{h,3} & \text{if} \ \bar{\gamma}_2 \geq \gamma > \bar{\gamma}_3, \\
\vdots & \vdots  \\
Z_{h,M_i-1} & \text{if} \ \bar{\gamma}_{ {M_i}-2} \geq \gamma > \bar{\gamma}_{ {M_i}-1}, \\
Z_{h,M_i} & \text{if} \ \bar{\gamma}_{M_i-1} \geq \gamma > \bar{\gamma}_{M_i}.
\end{array} \right. 
\end{equation}
For each of the $k$ folds, we will have one such function $\nu_i$ which indicates the objective on the hold-out fold $i$ as a function of $\gamma$. We then average these functions to obtain the cross-validated estimate of the out-of-sample objective:
\begin{equation}
\nu(\gamma) = \frac{1}{k} \sum_{i=1}^k \nu_i(\gamma).
\end{equation}
To set $\gamma$, we now find the optimum of this function for $\gamma \geq \gamma_{\min}$: 
\begin{equation}
\gamma^* = \arg \max_{\gamma \in [\gamma_{\min}, +\infty)} \nu(\gamma). \label{eq:kfoldcv_argmax}
\end{equation}

\clearpage

\section{Additional numerical results}
\label{appendix:results}

\subsection{Additional results on sensitivity to $\gamma$} 
\label{appendix:results_sensitivity}

In our previous experiments, we estimated the tree policies with a fixed $\gamma$ value of 0.005, which corresponds to a requirement of a 0.5\% relative improvement in the objective with each iteration of the construction algorithm. In this section, we consider how the tree policies and their performance changes as $\gamma$ varies. 

To understand how the performance changes as a function of $\gamma$, we follow a similar setup as in Section~\ref{subsec:results_problemdefinition} and randomly generate training and testing sets consisting of 2000 and 100,000 trajectories, respectively, for values of $n \in \{4,8,16\}$ and $\bar{p} \in \{90, 100, 110\}$, for 10 replications. We use a common correlation of $\bar{\rho} = 0$. For each replication, we run our construction algorithm on the training set and at each iteration, we compute the performance of the current tree on the training set and the testing set; we stop the construction algorithm at a $\gamma$ value of 0.0001. (We do not test smaller values of $\gamma$, due to the prohibitively large number of iterations that such smaller tolerances result in.) We test three different sets of basis functions: \textsc{prices}; \textsc{prices}, \textsc{time}; and \textsc{prices}, \textsc{time}, \textsc{payoff}, \textsc{KOind}. 

Table~\ref{table:sensitivity_gamma_neq8_p0eq90} shows the in-sample and out-of-sample performance of the tree policies at different values of $\gamma$ for the three sets of basis functions. The values are averaged over the ten replications. To simplify the presentation, we focus on $n = 8$ and $\bar{p} = 90$, as the results for other values of $n$ and $\bar{p}$ are qualitatively similar. While the in-sample performance improves as $\gamma$ decreases, the out-of-sample performance improves up to a point, after which it begins to decrease. For example, with \textsc{prices} only, the out-of-sample performance begins to decrease for $\gamma$ lower than 0.001. However, in all cases the deterioration appears to be mild (relative to the best average out-of-sample objective, the lowest value of $\gamma$ is only a few percent lower). This suggests that some care needs to be taken in selecting the right value of $\gamma$, as the tree policy may overfit the available training data. One possibility, as alluded to earlier, is to use part of the training data for building the tree policy and to use the remainder as a validation set, and to select the value of $\gamma$ that gives the highest objective on the validation set; then, the whole training set would be used to re-estimate the tree policy.

\begin{table}[ht]
\caption{In-sample and out-of-sample rewards for tree policies as a function of $\gamma$ for $n = 8$ assets and $\bar{p} = 90$ and common correlation $\bar{\rho} = 0$. \label{table:sensitivity_gamma_neq8_p0eq90} }
\centering
\begin{tabular}{lcccc}
\toprule
State variables & $\gamma$ & Training Obj. & Test. Obj. & Num. Iter. \\  \midrule
\textsc{prices} & 0.1000 & 20.38 & 19.52 & 3.2 \\ 
   & 0.0500 & 33.37 & 32.66 & 10.8 \\ 
   & 0.0100 & 35.48 & 34.62 & 14.2 \\ 
   & 0.0050 & 36.74 & 35.60 & 19.7 \\ 
   & 0.0010 & 39.93 & 37.56 & 53.0 \\ 
   & 0.0005 & 40.89 & 37.48 & 83.8 \\ 
   & 0.0001 & 42.60 & 36.60 & 239.0 \\[0.5em]
  \textsc{prices}, \textsc{time} & 0.1000 & 27.68 & 27.10 & 2.0 \\ 
 & 0.0500 & 28.17 & 27.52 & 2.4 \\ 
  & 0.0100 & 36.10 & 35.47 & 12.5 \\ 
   & 0.0050 & 39.76 & 38.89 & 23.8 \\ 
   & 0.0010 & 44.81 & 43.22 & 59.1 \\ 
   & 0.0005 & 45.65 & 43.53 & 82.5 \\ 
   & 0.0001 & 46.96 & 42.36 & 209.4 \\[0.5em]
  \textsc{prices}, \textsc{time},  & 0.1000 & 44.65 & 44.68 & 2.9 \\ 
 \textsc{payoff}, \textsc{KOind}  & 0.0500 & 44.81 & 44.84 & 3.0 \\ 
   & 0.0100 & 45.40 & 45.35 & 4.6 \\ 
   & 0.0050 & 45.47 & 45.38 & 5.0 \\ 
   & 0.0010 & 45.64 & 45.41 & 7.5 \\ 
   & 0.0005 & 45.78 & 45.39 & 11.6 \\ 
   & 0.0001 & 46.05 & 45.31 & 36.2 \\ 
\bottomrule
\end{tabular}

\end{table}

\clearpage

\subsection{Additional results for $k$-fold cross-validation procedure}
\label{appendix:results_kfoldcv}

\begin{figure}
\centering
\includegraphics[width = 0.7\textwidth]{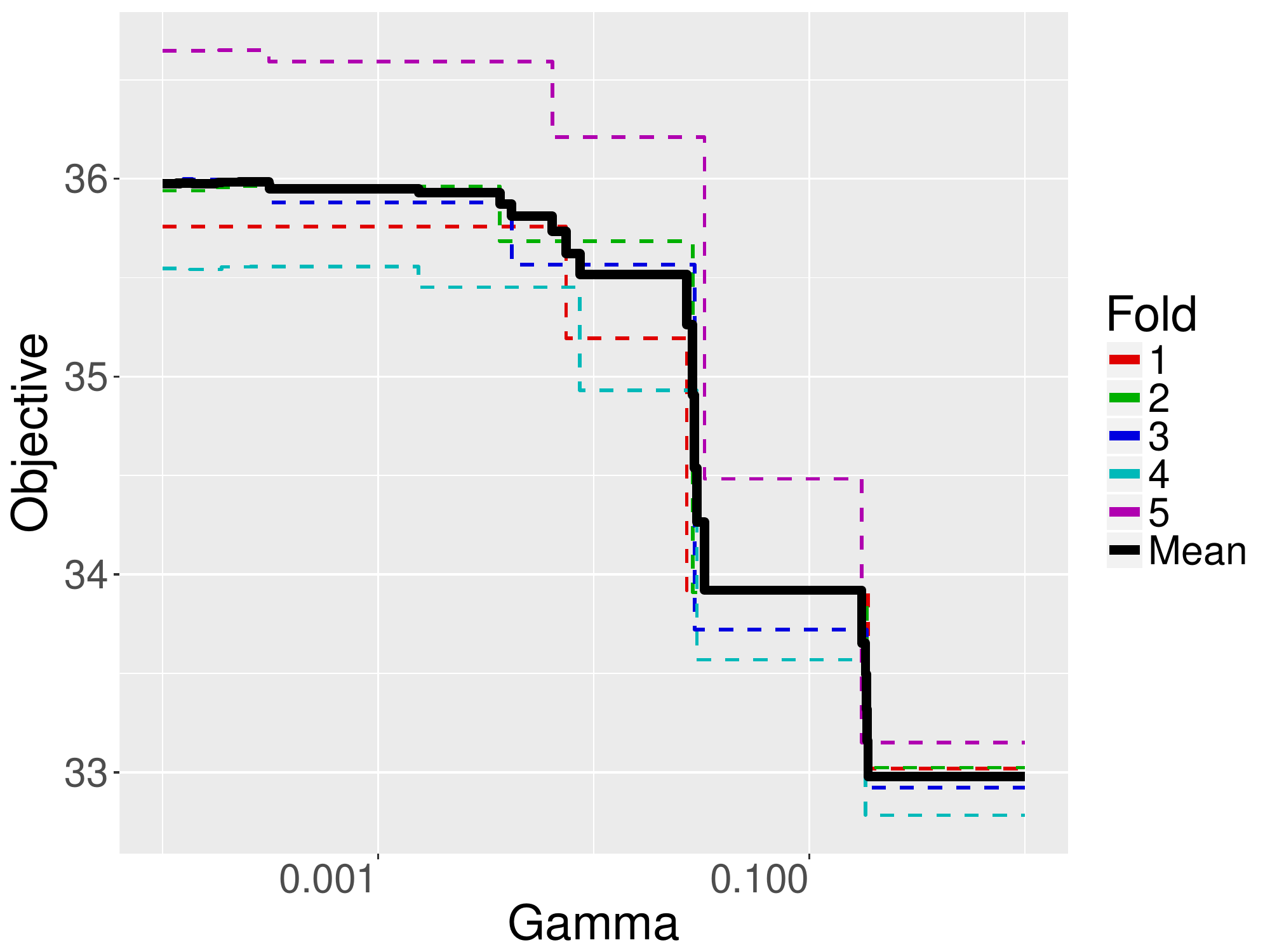}
\caption{$k$-fold cross-validation curves for all folds and the mean of all folds, for an $n = 8$, $\Omega = 20,000$ instance with \textsc{Prices}, \textsc{Payoff}, \textsc{KOind} and \textsc{Time} as state variables, and $\gamma_{\min} = 10^{-4}$. For clarity, the $x$ axis is plotted on a log scale. Note also that the right-most segment of each curve extends to $\gamma = +\infty$; for simplicity, we end each curve at $\gamma = 1$.  \label{figure:kfoldcv_allcurves} }

\end{figure}

We now show how the $k$-fold cross-validation procedure of Section~\ref{sec:kfold_cross_validation} works with a numerical example. In this example, we follow the same option pricing example as in Section~\ref{subsec:results_problemdefinition}. We generate 20,000 trajectories of $n = 8$ stocks with common correlation $\bar{\rho} = 0$. We run our $k$-fold cross-validation procedure (Algorithm~\ref{algorithm:CalculateKFoldCVBreakpoints}) with $k = 5$. Figure~\ref{figure:kfoldcv_allcurves} shows the functions $\nu_1, \dots, \nu_5$, as well as the mean function $\nu$. From Figure~\ref{figure:kfoldcv_allcurves}, the set of maximizers of $\nu(\cdot)$ is 
$(2.55 \times 10^{-4}, 3.12 \times 10^{-4}]$. 
However, it appears that a large range of $\gamma$ values will give near-optimal cross-validated performance. Thus, one might consider selecting $\gamma$ by different means -- for example, selecting the largest $\gamma$ that is within some percentage of optimal. (This is akin to the ``1 standard error'' rule used for cross-validating other types of machine learning models, such as LASSO regression models; in LASSO, for example, there is often a range of values of the regularization parameter $\lambda$ that are near optimal, and rather than picking the exact $\arg \max$ of the cross-validated error, one picks the largest $\lambda$ that is within one standard error of optimal, as this leads to the smallest number of non-zero coefficients in the model.)

As a further comparison, we re-run our comparison of our construction algorithm and LS with $\Omega = 20,000$ training trajectories, 100,000 test trajectories, $n = 8$ stocks and common correlation $\bar{\rho} = 0$. We test our construction algorithm with 5-fold cross-validation and without cross-validation (using the same value of $\gamma = 0.005$ as in previous experiments). For the cross-validated policy, we use $\gamma_{\min} = 10^{-4}$, and since the $\arg \max$ in equation~\eqref{eq:kfoldcv_argmax} will always be an interval, we pick the midpoint as our optimal value of $\gamma$. Table~\ref{table:kfoldcv_results} reports the performance of the three methods (LS, and the construction algorithm with and without cross-validation). From this table, we can see that 5-fold cross-validation leads to a slight improvement in performance over the choice of $\gamma = 0.005$. %
The main takeaway of this section is that the main input parameter of our construction algorithm -- $\gamma$ -- can be specified in a fully data-driven way, and the specification of this parameter is not an obstacle to the deployment of our algorithm in practice. %

\begin{table}

\begin{tabular}{lllccc} \toprule
$\Omega$ & Method & State variables / Basis Functions & $\bar{p} = 90$ & $\bar{p} = 100$ & $\bar{p} = 110$ \\ \midrule
20,000 & LS & \textsc{pricesKO}, \textsc{KOind}, \textsc{payoff} & 43.79\enskip (0.016) & 49.84\enskip (0.017) & 53.06\enskip (0.018) \\[0.25em] 
   & Tree & \textsc{prices}, \textsc{time}, \textsc{payoff}, \textsc{KOind} & 45.36\enskip (0.016) & 51.27\enskip (0.019) & 54.52\enskip (0.013) \\[0.25em] 
   & Tree-CV & \textsc{prices}, \textsc{time}, \textsc{payoff}, \textsc{KOind} & 45.43\enskip (0.018) & 51.34\enskip (0.013) & 54.52\enskip (0.013) \\ \bottomrule
 \end{tabular}
 
 \caption{Comparison of LS, tree policies with $\gamma = 0.005$ and tree policies with $\gamma$ chosen by 5-fold cross-validation ($\gamma_{\min} = 10^{-4}$) \label{table:kfoldcv_results}}
 \end{table}

\subsection{Performance results for $\bar{\rho} = 0$}
\label{appendix:results_rho0_performance}

Table~\ref{table:OOS_LS_vs_tree_neq4_16} provides additional performance results in the $\bar{\rho} = 0$ case, analogous to those in Table~\ref{table:OOS_LS_vs_tree_neq8}, for $n = 4$ and $n = 16$.

\begin{table}[!ht]
\caption{Comparison of out-of-sample performance between Longstaff-Schwartz and tree policies for $n = 4$ and $n = 16$ assets, for different initial prices $\bar{p}$ and common correlation $\bar{\rho} = 0$. In each column, the best performance is indicated in bold.}
 \label{table:OOS_LS_vs_tree_neq4_16}
 \scriptsize
  \centering
\begin{tabular}{lllccc}
  \toprule
$n$  & Method & State variables / Basis functions &  \multicolumn{3}{c}{Initial Price} \\
 &  & & $\bar{p} = 90$ & $\bar{p} = 100$ & $\bar{p} = 110$ \\ 
  \midrule
  4 & LS & one & 24.68\enskip (0.021) & 31.77\enskip (0.020) & 37.47\enskip (0.017) \\ 
    4 & LS & prices & 25.76\enskip (0.021) & 32.06\enskip (0.024) & 37.40\enskip (0.023) \\ 
    4 & LS & pricesKO & 28.52\enskip (0.018) & 38.31\enskip (0.027) & 46.58\enskip (0.027) \\ 
    4 & LS & pricesKO KOind & 30.23\enskip (0.016) & 39.05\enskip (0.021) & 46.58\enskip (0.023) \\ 
    4 & LS & pricesKO KOind payoff & 32.73\enskip (0.029) & 41.22\enskip (0.022) & 47.75\enskip (0.015) \\ 
    4 & LS & pricesKO KOind payoff  & 33.04\enskip (0.023) & 41.31\enskip (0.020) & 47.76\enskip (0.014) \\ 
    & & maxpriceKO \\
    4 & LS & pricesKO KOind payoff & 32.98\enskip (0.022) & 41.35\enskip (0.022) & 47.81\enskip (0.012) \\ 
    & & maxpriceKO max2priceKO  & \\
    4 & LS & pricesKO payoff & 33.47\enskip (0.021) & 41.70\enskip (0.015) & 47.72\enskip (0.007) \\ 
    4 & LS & pricesKO prices2KO KOind payoff & 33.43\enskip (0.022) & 41.81\enskip (0.022) & 48.03\enskip (0.013) \\[0.5em]
  4 & PO & prices & 31.43\enskip (0.017) & 38.91\enskip (0.018) & 43.40\enskip (0.016) \\ 
    4 & PO & pricesKO KOind payoff & 31.39\enskip (0.037) & 40.61\enskip (0.044) & 48.46\enskip (0.015) \\ 
    4 & PO & pricesKO KOind payoff & 32.20\enskip (0.032) & 41.11\enskip (0.037) & 48.49\enskip (0.015) \\ 
    & & maxpriceKO max2priceKO \\
    4 & PO & pricesKO prices2KO KOind payoff & 33.66\enskip (0.026) & 42.49\enskip (0.013) & 48.77\enskip (0.013) \\[0.5em]
      4 & Tree & payoff time & \bfseries 34.30\enskip (0.028) & \bfseries 43.08\enskip (0.022) & \bfseries 49.38\enskip (0.019) \\ 
    4 & Tree & prices & 27.14\enskip (0.048) & 36.91\enskip (0.021) & 45.15\enskip (0.043) \\ 
    4 & Tree & prices payoff & 27.34\enskip (0.022) & 37.13\enskip (0.015) & 45.77\enskip (0.014) \\ 
    4 & Tree & prices time & 33.76\enskip (0.126) & 39.51\enskip (0.470) & 40.94\enskip (0.370) \\ 
    4 & Tree & prices time payoff & \bfseries 34.30\enskip (0.028) & \bfseries 43.08\enskip (0.022) & \bfseries 49.38\enskip (0.019) \\ 
    4 & Tree & prices time payoff KOind & \bfseries 34.30\enskip (0.028) & \bfseries 43.08\enskip (0.022) & \bfseries 49.38\enskip (0.019) \\ \midrule
 16 & LS & one & 39.09\enskip (0.021) & 43.19\enskip (0.016) & 47.12\enskip (0.025) \\ 
   16 & LS & prices & 38.99\enskip (0.022) & 43.11\enskip (0.015) & 47.04\enskip (0.027) \\ 
   16 & LS & pricesKO & 50.34\enskip (0.024) & 53.40\enskip (0.009) & 54.68\enskip (0.007) \\ 
   16 & LS & pricesKO KOind & 50.41\enskip (0.026) & 53.72\enskip (0.008) & 54.97\enskip (0.009) \\ 
   16 & LS & pricesKO KOind payoff & 50.51\enskip (0.020) & 53.30\enskip (0.009) & 54.78\enskip (0.012) \\ 
   16 & LS & pricesKO KOind payoff  & 50.51\enskip (0.020) & 53.30\enskip (0.009) & 54.78\enskip (0.012) \\ 
   & & maxpriceKO \\
   16 & LS & pricesKO KOind payoff  & 50.50\enskip (0.020) & 53.29\enskip (0.009) & 54.78\enskip (0.013) \\ 
   & & maxpriceKO max2priceKO \\
   16 & LS & pricesKO payoff & 50.30\enskip (0.018) & 52.93\enskip (0.012) & 54.45\enskip (0.012) \\ 
   16 & LS & pricesKO prices2KO KOind payoff & 50.27\enskip (0.016) & 53.06\enskip (0.011) & 54.59\enskip (0.012) \\[0.5em]
 16 & PO & prices & 45.58\enskip (0.025) & 48.05\enskip (0.019) & 50.34\enskip (0.017) \\ 
   16 & PO & pricesKO KOind payoff & 51.26\enskip (0.012) & 53.92\enskip (0.007) & 55.28\enskip (0.009) \\ 
   16 & PO & pricesKO KOind payoff  & 51.22\enskip (0.015) & 53.92\enskip (0.006) & 55.27\enskip (0.008) \\ 
   & & maxpriceKO max2priceKO \\[0.5em]
    16 & Tree & payoff time & \bfseries 51.85\enskip (0.015) &\bfseries  54.62\enskip (0.008) & \bfseries 56.00\enskip (0.010) \\ 
   16 & Tree & prices & 39.70\enskip (0.155) & 42.60\enskip (0.125) & 43.91\enskip (0.144) \\ 
   16 & Tree & prices payoff & 49.35\enskip (0.012) & 54.17\enskip (0.008) & 55.96\enskip (0.008) \\ 
   16 & Tree & prices time & 39.51\enskip (0.089) & 43.10\enskip (0.020) & 46.17\enskip (0.502) \\ 
   16 & Tree & prices time payoff & \bfseries 51.85\enskip (0.015) &\bfseries  54.62\enskip (0.007) & \bfseries 56.00\enskip (0.010) \\ 
   16 & Tree & prices time payoff KOind & \bfseries 51.85\enskip (0.015) & \bfseries 54.62\enskip (0.007) & \bfseries 56.00\enskip (0.010) \\ \bottomrule
\end{tabular}
\end{table}

\clearpage 

\subsection{Performance results for $\bar{\rho} \in \{-0.05, +0.05, +0.10, +0.20\}$}
\label{appendix:results_rhoneq0_performance}

Tables~\ref{table:OOS_LS_vs_tree_neq4_8_16_rho-0.05}, \ref{table:OOS_LS_vs_tree_neq4_8_16_rho+0.05}, \ref{table:OOS_LS_vs_tree_neq4_8_16_rho+0.10} and \ref{table:OOS_LS_vs_tree_neq4_8_16_rho+0.20} report the performance of the tree and Longstaff-Schwartz policies for $\bar{\rho}$ values of $-0.05$, $+0.05$, $+0.10$ and $+0.20$, respectively. The instantaneous correlation matrix of the asset price process is set so that $\rho_{ii} = 1$ and $\rho_{ij} = \bar{\rho}$ for $i \neq j$. The same experimental set-up as for $\bar{\rho} = 0$ is followed (the number of assets $n$ varies in $\{4,8,16\}$, the initial price $\bar{p}$ varies in $\{90, 100, 110\}$, and each value reported for each method is averaged over ten replications).

\begin{table}[!ht]
\caption{Comparison of out-of-sample performance between Longstaff-Schwartz and tree policies for $n \in \{4, 8, 16\}$ assets, for different initial prices $\bar{p}$ and common correlation $\bar{\rho} = -0.05$. In each column, the best performance is indicated in bold.}
 \label{table:OOS_LS_vs_tree_neq4_8_16_rho-0.05}
 \scriptsize
  \centering
\begin{tabular}{lllccc}
  \toprule
$n$  & Method & State variables / Basis functions &  \multicolumn{3}{c}{Initial Price} \\
 &  & & $\bar{p} = 90$ & $\bar{p} = 100$ & $\bar{p} = 110$ \\ 
  \midrule
  4 & LS & one & 25.61\enskip (0.021) & 32.57\enskip (0.010) & 38.11\enskip (0.014) \\ 
    4 & LS & prices & 26.54\enskip (0.019) & 32.79\enskip (0.016) & 38.00\enskip (0.016) \\ 
    4 & LS & pricesKO & 29.48\enskip (0.027) & 39.39\enskip (0.023) & 47.48\enskip (0.023) \\ 
    4 & LS & pricesKO KOind & 31.07\enskip (0.012) & 40.04\enskip (0.022) & 47.48\enskip (0.017) \\ 
    4 & LS & pricesKO KOind payoff & 33.64\enskip (0.022) & 42.12\enskip (0.026) & 48.47\enskip (0.014) \\ 
    4 & LS & pricesKO KOind payoff & 33.88\enskip (0.024) & 42.20\enskip (0.026) & 48.47\enskip (0.016) \\ 
     &  & maxpriceKO & & & \\
    4 & LS & pricesKO KOind payoff & 33.85\enskip (0.024) & 42.23\enskip (0.021) & 48.51\enskip (0.013) \\ 
     &  & maxpriceKO max2priceKO  & & & \\
    4 & LS & pricesKO payoff & 34.36\enskip (0.026) & 42.56\enskip (0.020) & 48.41\enskip (0.013) \\ 
    4 & LS & pricesKO prices2KO KOind payoff & 34.26\enskip (0.022) & 42.65\enskip (0.018) & 48.70\enskip (0.019) \\ [0.5em]
   4 & Tree & payoff time & \bfseries 35.20\enskip (0.064) & \bfseries 43.86\enskip (0.022) & \bfseries 50.02\enskip (0.020) \\ 
    4 & Tree & prices & 27.83\enskip (0.024) & 37.76\enskip (0.026) & 46.08\enskip (0.037) \\ 
    4 & Tree & prices payoff & 28.06\enskip (0.028) & 37.97\enskip (0.026) & 46.61\enskip (0.019) \\ 
    4 & Tree & prices time & 34.74\enskip (0.062) & 40.33\enskip (0.341) & 41.61\enskip (0.347) \\ 
    4 & Tree & prices time payoff & \bfseries 35.20\enskip (0.064) & \bfseries 43.86\enskip (0.021) & \bfseries 50.03\enskip (0.020) \\ 
    4 & Tree & prices time payoff KOind & \bfseries 35.20\enskip (0.064) & \bfseries 43.86\enskip (0.021) & \bfseries 50.03\enskip (0.020) \\ \midrule 
  8 & LS & one & 34.99\enskip (0.018) & 39.57\enskip (0.014) & 44.02\enskip (0.017) \\ 
    8 & LS & prices & 34.95\enskip (0.021) & 39.44\enskip (0.016) & 43.92\enskip (0.021) \\ 
    8 & LS & pricesKO & 43.00\enskip (0.021) & 50.51\enskip (0.017) & 53.76\enskip (0.012) \\ 
    8 & LS & pricesKO KOind & 43.40\enskip (0.013) & 50.55\enskip (0.016) & 54.03\enskip (0.009) \\ 
    8 & LS & pricesKO KOind payoff & 45.15\enskip (0.017) & 50.76\enskip (0.018) & 53.58\enskip (0.006) \\ 
    8 & LS & pricesKO KOind payoff & 45.15\enskip (0.017) & 50.75\enskip (0.018) & 53.58\enskip (0.006) \\
    &  & maxpriceKO & & & \\ 
    8 & LS & pricesKO KOind payoff & 45.17\enskip (0.016) & 50.74\enskip (0.019) & 53.57\enskip (0.005) \\ 
    &  & maxpriceKO max2priceKO & & & \\
    8 & LS & pricesKO payoff & 45.35\enskip (0.012) & 50.56\enskip (0.017) & 53.25\enskip (0.008) \\ 
    8 & LS & pricesKO prices2KO KOind payoff & 45.27\enskip (0.014) & 50.77\enskip (0.018) & 53.58\enskip (0.009) \\ [0.5em]
 8 & Tree & payoff time & \bfseries 46.58\enskip (0.020) & \bfseries 52.11\enskip (0.024) & \bfseries 54.92\enskip (0.011) \\
    8 & Tree & prices & 36.73\enskip (0.101) & 44.54\enskip (0.091) & 47.83\enskip (0.089) \\ 
    8 & Tree & prices payoff & 40.56\enskip (0.018) & 49.66\enskip (0.020) & 54.35\enskip (0.010) \\ 
    8 & Tree & prices time & 39.13\enskip (0.196) & 40.81\enskip (0.503) & 43.41\enskip (0.123) \\ 
    8 & Tree & prices time payoff & \bfseries 46.58\enskip (0.020) & \bfseries 52.11\enskip (0.024) & \bfseries 54.92\enskip (0.010) \\ 
    8 & Tree & prices time payoff KOind & \bfseries 46.58\enskip (0.020) &\bfseries  52.11\enskip (0.024) & \bfseries 54.92\enskip (0.010) \\ \midrule
16 & LS & one & 40.34\enskip (0.021) & 44.27\enskip (0.035) & 48.10\enskip (0.022) \\ 
   16 & LS & payoff time & 43.78\enskip (0.035) & 45.83\enskip (0.045) & 48.19\enskip (0.018) \\ 
   16 & LS & prices & 40.25\enskip (0.022) & 44.21\enskip (0.036) & 48.05\enskip (0.020) \\ 
   16 & LS & pricesKO & 51.78\enskip (0.015) & 54.13\enskip (0.007) & 55.06\enskip (0.008) \\ 
   16 & LS & pricesKO KOind & 51.76\enskip (0.012) & 54.16\enskip (0.009) & 55.16\enskip (0.011) \\ 
   16 & LS & pricesKO KOind payoff & 51.56\enskip (0.013) & 53.78\enskip (0.010) & 55.12\enskip (0.008) \\ 
   16 & LS & pricesKO KOind payoff & 51.56\enskip (0.013) & 53.78\enskip (0.010) & 55.12\enskip (0.008) \\
   &  & maxpriceKO & & & \\  
   16 & LS & pricesKO KOind payoff & 51.54\enskip (0.013) & 53.76\enskip (0.009) & 55.11\enskip (0.010) \\ 
   &  & maxpriceKO max2priceKO & & & \\ 
   16 & LS & pricesKO payoff & 51.51\enskip (0.014) & 53.65\enskip (0.010) & 55.01\enskip (0.010) \\ 
   16 & LS & pricesKO prices2KO KOind payoff & 51.30\enskip (0.015) & 53.53\enskip (0.015) & 54.93\enskip (0.008) \\
[0.5em]
    16 & Tree & payoff time & \bfseries 52.75\enskip (0.011) & \bfseries 54.94\enskip (0.006) & \bfseries 56.14\enskip (0.009) \\ 
   16 & Tree & prices & 40.81\enskip (0.089) & 43.11\enskip (0.167) & 44.06\enskip (0.118) \\ 
   16 & Tree & prices payoff & 51.09\enskip (0.011) & 54.81\enskip (0.007) & 56.13\enskip (0.009) \\ 
   16 & Tree & prices time & 41.16\enskip (0.076) & 44.23\enskip (0.095) & 47.28\enskip (0.509) \\ 
   16 & Tree & prices time payoff & 52.74\enskip (0.011) & \bfseries 54.94\enskip (0.006) &  56.13\enskip (0.009) \\ 
   16 & Tree & prices time payoff KOind & 52.74\enskip (0.011) & \bfseries 54.94\enskip (0.006) &  56.13\enskip (0.009) \\ \bottomrule
\end{tabular}
\end{table}

\begin{table}
\caption{Comparison of out-of-sample performance between Longstaff-Schwartz and tree policies for $n \in \{4, 8, 16\}$ assets, for different initial prices $\bar{p}$ and common correlation $\bar{\rho} = +0.05$. In each column, the best performance is indicated in bold.}
 \label{table:OOS_LS_vs_tree_neq4_8_16_rho+0.05}
 \scriptsize
  \centering
\begin{tabular}{lllccc}
  \toprule
$n$  & Method & State variables / Basis functions &  \multicolumn{3}{c}{Initial Price} \\
 &  & & $\bar{p} = 90$ & $\bar{p} = 100$ & $\bar{p} = 110$ \\ 
  \midrule 
  4 & LS & one & 23.84\enskip (0.022) & 31.00\enskip (0.020) & 36.77\enskip (0.013) \\ 
    4 & LS & payoff time & 32.07\enskip (0.031) & 39.32\enskip (0.015) & 42.96\enskip (0.073) \\ 
    4 & LS & prices & 25.00\enskip (0.023) & 31.37\enskip (0.019) & 36.72\enskip (0.012) \\ 
    4 & LS & pricesKO & 27.68\enskip (0.029) & 37.42\enskip (0.033) & 45.68\enskip (0.029) \\ 
    4 & LS & pricesKO KOind & 29.50\enskip (0.023) & 38.21\enskip (0.021) & 45.73\enskip (0.026) \\ 
    4 & LS & pricesKO KOind payoff & 31.93\enskip (0.028) & 40.41\enskip (0.028) & 47.03\enskip (0.019) \\ 
    4 & LS & pricesKO KOind payoff  & 32.30\enskip (0.027) & 40.55\enskip (0.029) & 47.04\enskip (0.018) \\
       &  & maxpriceKO  & & & \\ 
    4 & LS & pricesKO KOind payoff  & 32.28\enskip (0.027) & 40.59\enskip (0.029) & 47.08\enskip (0.017) \\ 
        &  & maxpriceKO max2priceKO & & & \\
    4 & LS & pricesKO payoff & 32.69\enskip (0.024) & 40.88\enskip (0.020) & 47.03\enskip (0.014) \\ 
    4 & LS & pricesKO prices2KO KOind payoff & 32.71\enskip (0.024) & 41.01\enskip (0.023) & 47.34\enskip (0.008) \\ 
 [0.5em]
4 & Tree & payoff time & \bfseries 33.68\enskip (0.048) & \bfseries 42.26\enskip (0.033) & \bfseries 48.70\enskip (0.024) \\ 
    4 & Tree & prices & 26.47\enskip (0.021) & 36.00\enskip (0.044) & 44.29\enskip (0.054) \\ 
    4 & Tree & prices payoff & 26.74\enskip (0.023) & 36.29\enskip (0.027) & 44.95\enskip (0.019) \\ 
    4 & Tree & prices time & 32.97\enskip (0.088) & 38.90\enskip (0.423) & 40.48\enskip (0.287) \\ 
    4 & Tree & prices time payoff & \bfseries 33.68\enskip (0.048) & \bfseries 42.26\enskip (0.033) & \bfseries 48.70\enskip (0.024) \\ 
    4 & Tree & prices time payoff KOind & \bfseries 33.68\enskip (0.048) & \bfseries 42.26\enskip (0.033) & \bfseries 48.70\enskip (0.024) \\  \midrule
   8 & LS & one & 32.65\enskip (0.022) & 37.73\enskip (0.020) & 42.31\enskip (0.018) \\ 
    8 & LS & payoff time & 40.40\enskip (0.044) & 42.95\enskip (0.030) & 44.46\enskip (0.025) \\ 
    8 & LS & prices & 32.81\enskip (0.023) & 37.61\enskip (0.019) & 42.19\enskip (0.019) \\ 
    8 & LS & pricesKO & 39.95\enskip (0.034) & 48.18\enskip (0.017) & 52.40\enskip (0.009) \\ 
    8 & LS & pricesKO KOind & 40.40\enskip (0.027) & 48.15\enskip (0.022) & 52.74\enskip (0.017) \\ 
    8 & LS & pricesKO KOind payoff & 42.48\enskip (0.033) & 48.93\enskip (0.015) & 52.53\enskip (0.014) \\ 
    8 & LS & pricesKO KOind payoff  & 42.57\enskip (0.032) & 48.93\enskip (0.016) & 52.52\enskip (0.014) \\ 
    &  & maxpriceKO  & & & \\ 
    8 & LS & pricesKO KOind payoff  & 42.60\enskip (0.033) & 48.96\enskip (0.013) & 52.52\enskip (0.015) \\ 
    &  & maxpriceKO  max2priceKO & & & \\ 
    8 & LS & pricesKO payoff & 42.79\enskip (0.027) & 48.69\enskip (0.013) & 52.06\enskip (0.012) \\ 
    8 & LS & pricesKO prices2KO KOind payoff & 42.88\enskip (0.030) & 49.07\enskip (0.011) & 52.61\enskip (0.014) \\
 [0.5em]
8 & Tree & payoff time & \bfseries44.24\enskip (0.021) & \bfseries 50.43\enskip (0.008) & \bfseries 54.08\enskip (0.009) \\ 
    8 & Tree & prices & 34.17\enskip (0.217) & 42.65\enskip (0.105) & 46.41\enskip (0.104) \\ 
    8 & Tree & prices payoff & 37.75\enskip (0.025) & 47.08\enskip (0.023) & 52.84\enskip (0.012) \\ 
    8 & Tree & prices time & 37.74\enskip (0.259) & 40.02\enskip (0.337) & 41.81\enskip (0.109) \\ 
    8 & Tree & prices time payoff & \bfseries 44.24\enskip (0.022) & \bfseries 50.43\enskip (0.008) & \bfseries 54.08\enskip (0.010) \\ 
    8 & Tree & prices time payoff KOind & \bfseries 44.24\enskip (0.022) & \bfseries 50.43\enskip (0.008) & \bfseries 54.08\enskip (0.010) \\ \midrule
   16 & LS & one & 37.96\enskip (0.023) & 42.12\enskip (0.023) & 46.19\enskip (0.019) \\ 
   16 & LS & payoff time & 42.70\enskip (0.041) & 44.48\enskip (0.020) & 46.73\enskip (0.031) \\ 
   16 & LS & prices & 37.84\enskip (0.026) & 42.01\enskip (0.029) & 46.09\enskip (0.019) \\ 
   16 & LS & pricesKO & 48.95\enskip (0.011) & 52.60\enskip (0.010) & 54.33\enskip (0.011) \\ 
   16 & LS & pricesKO KOind & 49.00\enskip (0.015) & 53.08\enskip (0.011) & 54.72\enskip (0.010) \\ 
   16 & LS & pricesKO KOind payoff & 49.44\enskip (0.014) & 52.73\enskip (0.009) & 54.45\enskip (0.012) \\ 
   16 & LS & pricesKO KOind payoff & 49.44\enskip (0.014) & 52.73\enskip (0.009) & 54.45\enskip (0.012) \\
    &  & maxpriceKO  & & & \\  
   16 & LS & pricesKO KOind payoff  & 49.44\enskip (0.015) & 52.72\enskip (0.008) & 54.44\enskip (0.011) \\ 
    &  & maxpriceKO  max2priceKO & & & \\ 
   16 & LS & pricesKO payoff & 49.14\enskip (0.010) & 52.21\enskip (0.013) & 53.96\enskip (0.013) \\ 
   16 & LS & pricesKO prices2KO KOind payoff & 49.26\enskip (0.012) & 52.51\enskip (0.011) & 54.29\enskip (0.012) \\ 
 [0.5em]
   16 & Tree & payoff time & \bfseries 50.86\enskip (0.019) & \bfseries 54.17\enskip (0.007) & \bfseries 55.86\enskip (0.010) \\ 
   16 & Tree & prices & 38.76\enskip (0.147) & 41.83\enskip (0.114) & 43.49\enskip (0.140) \\ 
   16 & Tree & prices payoff & 47.62\enskip (0.021) & 53.29\enskip (0.012) & 55.73\enskip (0.012) \\ 
   16 & Tree & prices time & 38.74\enskip (0.151) & 42.07\enskip (0.028) & 44.82\enskip (0.642) \\ 
   16 & Tree & prices time payoff & \bfseries 50.86\enskip (0.019) & \bfseries 54.17\enskip (0.007) & \bfseries 55.86\enskip (0.010) \\ 
   16 & Tree & prices time payoff KOind & \bfseries 50.86\enskip (0.019) & \bfseries 54.17\enskip (0.007) & \bfseries 55.86\enskip (0.010) \\ \bottomrule
  \end{tabular}
  \end{table}

\begin{table}
\caption{Comparison of out-of-sample performance between Longstaff-Schwartz and tree policies for $n \in \{4, 8, 16\}$ assets, for different initial prices $\bar{p}$ and common correlation $\bar{\rho} = +0.10$. In each column, the best performance is indicated in bold.}
 \label{table:OOS_LS_vs_tree_neq4_8_16_rho+0.10}
 \scriptsize
  \centering
\begin{tabular}{lllccc}
  \toprule
$n$  & Method & State variables / Basis functions &  \multicolumn{3}{c}{Initial Price} \\
 &  & & $\bar{p} = 90$ & $\bar{p} = 100$ & $\bar{p} = 110$ \\ 
  \midrule 
4 & LS & one & 23.03\enskip (0.023) & 30.22\enskip (0.018) & 36.13\enskip (0.024) \\ 
    4 & LS & payoff time & 31.24\enskip (0.019) & 38.71\enskip (0.028) & 42.70\enskip (0.058) \\ 
    4 & LS & prices & 24.30\enskip (0.022) & 30.70\enskip (0.018) & 36.14\enskip (0.019) \\ 
    4 & LS & pricesKO & 26.88\enskip (0.026) & 36.43\enskip (0.030) & 44.74\enskip (0.016) \\ 
    4 & LS & pricesKO KOind & 28.77\enskip (0.026) & 37.36\enskip (0.019) & 44.85\enskip (0.019) \\ 
    4 & LS & pricesKO KOind payoff & 31.07\enskip (0.020) & 39.58\enskip (0.019) & 46.29\enskip (0.015) \\ 
    4 & LS & pricesKO KOind payoff & 31.51\enskip (0.019) & 39.77\enskip (0.019) & 46.32\enskip (0.015) \\
       &  & maxpriceKO  & & & \\   
    4 & LS & pricesKO KOind payoff & 31.49\enskip (0.026) & 39.79\enskip (0.019) & 46.36\enskip (0.021) \\ 
        &  & maxpriceKO  max2priceKO & & & \\  
    4 & LS & pricesKO payoff & 31.89\enskip (0.014) & 40.10\enskip (0.020) & 46.36\enskip (0.013) \\ 
    4 & LS & pricesKO prices2KO KOind payoff  & 31.92\enskip (0.018) & 40.26\enskip (0.017) & 46.64\enskip (0.016) \\ 
     [0.5em]
 4 & Tree & payoff time & \bfseries 32.81\enskip (0.030) & \bfseries 41.53\enskip (0.023) & \bfseries 47.99\enskip (0.020) \\ 
    4 & Tree & prices & 25.83\enskip (0.021) & 35.22\enskip (0.029) & 43.39\enskip (0.042) \\ 
    4 & Tree & prices payoff & 26.07\enskip (0.014) & 35.53\enskip (0.022) & 44.11\enskip (0.029) \\ 
    4 & Tree & prices time & 32.37\enskip (0.112) & 38.31\enskip (0.322) & 39.90\enskip (0.253) \\ 
    4 & Tree & prices time payoff & \bfseries 32.81\enskip (0.030) & 41.52\enskip (0.021) & \bfseries 47.99\enskip (0.020) \\ 
    4 & Tree & prices time payoff KOind & \bfseries 32.81\enskip (0.030) &  41.52\enskip (0.021) & \bfseries 47.99\enskip (0.020) \\  \midrule 
8 & LS & one & 31.46\enskip (0.020) & 36.76\enskip (0.022) & 41.48\enskip (0.018) \\ 
    8 & LS & payoff time & 39.52\enskip (0.042) & 42.58\enskip (0.035) & 44.00\enskip (0.034) \\ 
    8 & LS & prices & 31.75\enskip (0.015) & 36.68\enskip (0.021) & 41.36\enskip (0.023) \\ 
    8 & LS & pricesKO & 38.45\enskip (0.035) & 46.95\enskip (0.017) & 51.68\enskip (0.015) \\ 
    8 & LS & pricesKO KOind & 38.97\enskip (0.026) & 46.85\enskip (0.018) & 52.00\enskip (0.022) \\ 
    8 & LS & pricesKO KOind payoff & 41.16\enskip (0.027) & 47.93\enskip (0.020) & 51.95\enskip (0.017) \\ 
    8 & LS & pricesKO KOind payoff & 41.29\enskip (0.025) & 47.94\enskip (0.021) & 51.94\enskip (0.016) \\
       &  & maxpriceKO  & & & \\   
    8 & LS & pricesKO KOind payoff & 41.33\enskip (0.024) & 47.98\enskip (0.019) & 51.95\enskip (0.014) \\ 
       &  & maxpriceKO max2priceKO & & & \\  
    8 & LS & pricesKO payoff & 41.52\enskip (0.017) & 47.74\enskip (0.017) & 51.46\enskip (0.013) \\ 
    8 & LS & pricesKO prices2KO KOind payoff & 41.67\enskip (0.021) & 48.15\enskip (0.010) & 52.06\enskip (0.016) \\ 
   [0.5em]
8 & Tree & payoff time & \bfseries 43.09\enskip (0.022) & \bfseries 49.58\enskip (0.024) & \bfseries 53.53\enskip (0.017) \\ 
    8 & Tree & prices & 33.80\enskip (0.196) & 41.40\enskip (0.094) & 45.77\enskip (0.117) \\ 
    8 & Tree & prices payoff & 36.46\enskip (0.015) & 45.85\enskip (0.022) & 51.98\enskip (0.011) \\ 
    8 & Tree & prices time & 36.55\enskip (0.218) & 39.51\enskip (0.238) & 41.23\enskip (0.139) \\ 
    8 & Tree & prices time payoff & \bfseries 43.09\enskip (0.022) & \bfseries 49.58\enskip (0.024) & \bfseries 53.53\enskip (0.017) \\ 
    8 & Tree & prices time payoff KOind & \bfseries 43.09\enskip (0.022) & \bfseries 49.58\enskip (0.024) & \bfseries 53.53\enskip (0.017) \\  \midrule
16 & LS & one & 36.74\enskip (0.017) & 41.06\enskip (0.016) & 45.25\enskip (0.020) \\ 
   16 & LS & payoff time & 42.25\enskip (0.032) & 43.85\enskip (0.032) & 46.07\enskip (0.024) \\ 
   16 & LS & prices & 36.62\enskip (0.019) & 40.93\enskip (0.019) & 45.12\enskip (0.024) \\ 
   16 & LS & pricesKO & 47.47\enskip (0.015) & 51.77\enskip (0.013) & 53.94\enskip (0.015) \\ 
   16 & LS & pricesKO KOind & 47.42\enskip (0.017) & 52.28\enskip (0.017) & 54.40\enskip (0.016) \\ 
   16 & LS & pricesKO KOind payoff & 48.22\enskip (0.016) & 52.06\enskip (0.016) & 54.09\enskip (0.011) \\ 
   16 & LS & pricesKO KOind payoff & 48.23\enskip (0.016) & 52.06\enskip (0.016) & 54.09\enskip (0.012) \\ 
    &  & maxpriceKO & & & \\  
   16 & LS & pricesKO KOind payoff & 48.24\enskip (0.015) & 52.05\enskip (0.016) & 54.09\enskip (0.012) \\ 
    &  & maxpriceKO max2priceKO & & & \\  
   16 & LS & pricesKO payoff & 47.93\enskip (0.013) & 51.47\enskip (0.011) & 53.53\enskip (0.013) \\ 
   16 & LS & pricesKO prices2KO KOind payoff & 48.11\enskip (0.012) & 51.89\enskip (0.017) & 53.96\enskip (0.013) \\
  [0.5em] 
16 & Tree & payoff time & \bfseries 49.78\enskip (0.021) & \bfseries 53.61\enskip (0.009) & \bfseries 55.63\enskip (0.013) \\ 
   16 & Tree & prices & 37.11\enskip (0.103) & 41.12\enskip (0.083) & 43.10\enskip (0.092) \\ 
   16 & Tree & prices payoff & 45.92\enskip (0.023) & 52.31\enskip (0.015) & 55.35\enskip (0.012) \\ 
   16 & Tree & prices time & 37.54\enskip (0.111) & 41.14\enskip (0.076) & 44.45\enskip (0.270) \\ 
   16 & Tree & prices time payoff & \bfseries 49.78\enskip (0.021) & \bfseries 53.61\enskip (0.009) & \bfseries 55.63\enskip (0.013) \\ 
   16 & Tree & prices time payoff KOind & \bfseries 49.78\enskip (0.021) & \bfseries 53.61\enskip (0.009) & \bfseries 55.63\enskip (0.013) \\ \bottomrule
  \end{tabular}
  \end{table}

  \begin{table}
\caption{Comparison of out-of-sample performance between Longstaff-Schwartz and tree policies for $n \in \{4, 8, 16\}$ assets, for different initial prices $\bar{p}$ and common correlation $\bar{\rho} = +0.20$. In each column, the best performance is indicated in bold.}
 \label{table:OOS_LS_vs_tree_neq4_8_16_rho+0.20}
 \scriptsize
  \centering
\begin{tabular}{lllccc}
  \toprule
$n$  & Method & State variables / Basis functions &  \multicolumn{3}{c}{Initial Price} \\
 &  & & $\bar{p} = 90$ & $\bar{p} = 100$ & $\bar{p} = 110$ \\ 
  \midrule 
 4 & LS & one & 21.41\enskip (0.025) & 28.60\enskip (0.023) & 34.77\enskip (0.018) \\ 
    4 & LS & payoff time & 29.70\enskip (0.026) & 37.27\enskip (0.034) & 42.07\enskip (0.038) \\ 
    4 & LS & prices & 22.85\enskip (0.024) & 29.28\enskip (0.027) & 34.84\enskip (0.017) \\ 
    4 & LS & pricesKO & 25.22\enskip (0.031) & 34.46\enskip (0.026) & 42.83\enskip (0.019) \\ 
    4 & LS & pricesKO KOind & 27.36\enskip (0.018) & 35.66\enskip (0.016) & 43.17\enskip (0.023) \\ 
    4 & LS & pricesKO KOind payoff & 29.48\enskip (0.027) & 37.80\enskip (0.017) & 44.79\enskip (0.023) \\ 
    4 & LS & pricesKO KOind payoff & 30.04\enskip (0.023) & 38.09\enskip (0.016) & 44.85\enskip (0.022) \\
        &  & maxpriceKO & & & \\   
    4 & LS & pricesKO KOind payoff & 30.02\enskip (0.026) & 38.12\enskip (0.018) & 44.90\enskip (0.020) \\
        &  & maxpriceKO max2priceKO & & & \\   
    4 & LS & pricesKO payoff & 30.34\enskip (0.017) & 38.42\enskip (0.021) & 44.95\enskip (0.014) \\ 
    4 & LS & pricesKO prices2KO KOind payoff & 30.42\enskip (0.024) & 38.59\enskip (0.021) & 45.16\enskip (0.025) \\[0.5em]
4 & Tree & payoff time & \bfseries 31.40\enskip (0.047) &\bfseries 39.90\enskip (0.025) & \bfseries 46.58\enskip (0.014) \\ 
    4 & Tree & prices & 24.59\enskip (0.039) & 33.47\enskip (0.030) & 41.66\enskip (0.039) \\ \bfseries
    4 & Tree & prices payoff & 24.79\enskip (0.014) & 33.86\enskip (0.024) & 42.44\enskip (0.014) \\ 
    4 & Tree & prices time & 30.82\enskip (0.124) & 37.06\enskip (0.227) & 38.58\enskip (0.188) \\ 
    4 & Tree & prices time payoff & \bfseries 31.40\enskip (0.047) & \bfseries 39.90\enskip (0.025) & \bfseries 46.58\enskip (0.014) \\ 
    4 & Tree & prices time payoff KOind & \bfseries 31.40\enskip (0.047) & \bfseries 39.90\enskip (0.025) &\bfseries  46.58\enskip (0.014) \\ \midrule
 8 & LS & one & 29.13\enskip (0.021) & 34.94\enskip (0.028) & 39.94\enskip (0.013) \\ 
    8 & LS & payoff time & 37.69\enskip (0.024) & 41.99\enskip (0.039) & 43.32\enskip (0.031) \\ 
    8 & LS & prices & 29.73\enskip (0.026) & 34.94\enskip (0.033) & 39.81\enskip (0.016) \\ 
    8 & LS & pricesKO & 35.56\enskip (0.025) & 44.48\enskip (0.028) & 50.16\enskip (0.014) \\ 
    8 & LS & pricesKO KOind & 36.38\enskip (0.027) & 44.42\enskip (0.032) & 50.30\enskip (0.020) \\ 
    8 & LS & pricesKO KOind payoff & 38.65\enskip (0.030) & 45.92\enskip (0.017) & 50.64\enskip (0.019) \\ 
    8 & LS & pricesKO KOind payoff & 38.93\enskip (0.031) & 45.96\enskip (0.014) & 50.63\enskip (0.018) \\
        &  & maxpriceKO & & & \\   
    8 & LS & pricesKO KOind payoff & 38.96\enskip (0.026) & 46.01\enskip (0.015) & 50.65\enskip (0.017) \\ 
        &  & maxpriceKO max2priceKO & & & \\  
    8 & LS & pricesKO payoff & 39.13\enskip (0.026) & 45.83\enskip (0.014) & 50.19\enskip (0.013) \\ 
    8 & LS & pricesKO prices2KO KOind payoff & 39.37\enskip (0.025) & 46.27\enskip (0.011) & 50.78\enskip (0.016) \\[0.5em]
8 & Tree & payoff time & \bfseries 40.73\enskip (0.014) &\bfseries 47.74\enskip (0.025) &\bfseries 52.30\enskip (0.013) \\ 
    8 & Tree & prices & 32.11\enskip (0.170) & 39.28\enskip (0.318) & 44.34\enskip (0.177) \\ 
    8 & Tree & prices payoff & 33.98\enskip (0.033) & 43.36\enskip (0.021) & 50.16\enskip (0.015) \\ 
    8 & Tree & prices time & 36.14\enskip (0.172) & 37.78\enskip (0.334) & 40.17\enskip (0.231) \\ 
    8 & Tree & prices time payoff & \bfseries 40.73\enskip (0.014) & \bfseries 47.74\enskip (0.025) & \bfseries 52.30\enskip (0.014) \\ 
    8 & Tree & prices time payoff KOind & \bfseries 40.73\enskip (0.014) & \bfseries 47.74\enskip (0.025) & \bfseries 52.30\enskip (0.014) \\ \midrule
16 & LS & one & 34.45\enskip (0.020) & 39.00\enskip (0.024) & 43.47\enskip (0.017) \\ 
   16 & LS & payoff time & 41.39\enskip (0.042) & 42.94\enskip (0.037) & 44.84\enskip (0.025) \\ 
   16 & LS & prices & 34.41\enskip (0.017) & 38.82\enskip (0.025) & 43.30\enskip (0.019) \\ 
   16 & LS & pricesKO & 44.41\enskip (0.027) & 50.02\enskip (0.010) & 52.94\enskip (0.016) \\ 
   16 & LS & pricesKO KOind & 44.24\enskip (0.024) & 50.35\enskip (0.019) & 53.47\enskip (0.015) \\ 
   16 & LS & pricesKO KOind payoff & 45.64\enskip (0.021) & 50.54\enskip (0.012) & 53.21\enskip (0.013) \\ 
   16 & LS & pricesKO KOind payoff & 45.72\enskip (0.024) & 50.54\enskip (0.012) & 53.20\enskip (0.014) \\
       &  & maxpriceKO & & & \\   
   16 & LS & pricesKO KOind payoff & 45.76\enskip (0.020) & 50.55\enskip (0.012) & 53.20\enskip (0.015) \\ 
       &  & maxpriceKO max2priceKO & & & \\  
   16 & LS & pricesKO payoff & 45.48\enskip (0.017) & 49.95\enskip (0.011) & 52.59\enskip (0.011) \\ 
   16 & LS & pricesKO prices2KO KOind payoff & 45.73\enskip (0.019) & 50.44\enskip (0.011) & 53.17\enskip (0.013) \\[0.5em]
16 & Tree & payoff time & \bfseries 47.50\enskip (0.026) & \bfseries 52.25\enskip (0.012) & \bfseries 54.91\enskip (0.010) \\ 
   16 & Tree & prices & 34.95\enskip (0.128) & 39.89\enskip (0.127) & 42.47\enskip (0.070) \\ 
   16 & Tree & prices payoff & 42.62\enskip (0.018) & 50.05\enskip (0.022) & 54.19\enskip (0.014) \\ 
   16 & Tree & prices time & 35.21\enskip (0.122) & 38.94\enskip (0.120) & 41.96\enskip (0.199) \\ 
   16 & Tree & prices time payoff & \bfseries 47.50\enskip (0.026) & \bfseries 52.25\enskip (0.012) &  54.90\enskip (0.010) \\ 
   16 & Tree & prices time payoff KOind & \bfseries 47.50\enskip (0.026) & \bfseries 52.25\enskip (0.012) &  54.90\enskip (0.010) \\ \bottomrule
  \end{tabular}
  \end{table}
  
\clearpage 
\subsection{Timing results for $n \in \{4,16\}$ and $\bar{\rho} = 0$}
\label{appendix:results_rho0_timing}

Table~\ref{table:OOS_LS_vs_tree_neq4_16_rho0} reports the computation time for the tree and Longstaff-Schwartz policies for $n \in \{4,16\}$ and $\bar{p} \in \{90, 100, 110\}$, for the uncorrelated ($\bar{\rho} = 0$) case. The computation times are averaged over the ten replications for each combination of $n$ and $\bar{p}$.

\begin{table}[!ht]
\caption{Comparison of estimation time between Longstaff-Schwartz, pathwise optimization and tree policies for $n \in \{4, 16\}$ assets, for different initial prices $\bar{p}$ and common correlation $\bar{\rho} = 0$.}
 \label{table:OOS_LS_vs_tree_neq4_16_rho0}
 \scriptsize
 \centering
\begin{tabular}{lllrrr}
  \toprule
$n$  & Method & State variables / Basis functions &  \multicolumn{3}{c}{Initial Price} \\
 &  & & $\bar{p} = 90$ & $\bar{p} = 100$ & $\bar{p} = 110$ \\ 
  \midrule
  4 & LS & one & 1.3\enskip (0.1) & 1.3\enskip (0.1) & 1.3\enskip (0.0) \\ 
    4 & LS & prices & 1.3\enskip (0.1) & 1.2\enskip (0.1) & 1.2\enskip (0.0) \\ 
    4 & LS & pricesKO & 1.6\enskip (0.1) & 1.6\enskip (0.1) & 1.4\enskip (0.1) \\ 
    4 & LS & pricesKO KOind & 1.7\enskip (0.2) & 1.9\enskip (0.2) & 1.5\enskip (0.2) \\ 
    4 & LS & pricesKO KOind payoff & 1.7\enskip (0.1) & 1.7\enskip (0.1) & 1.5\enskip (0.2) \\ 
    4 & LS & pricesKO KOind payoff & 2.0\enskip (0.1) & 2.2\enskip (0.2) & 2.2\enskip (0.2) \\ 
    & & maxpriceKO \\
    4 & LS & pricesKO KOind payoff & 2.0\enskip (0.2) & 1.9\enskip (0.1) & 1.9\enskip (0.2) \\ 
    & & maxpriceKO max2priceKO \\
    4 & LS & pricesKO payoff & 2.0\enskip (0.1) & 1.9\enskip (0.1) & 1.5\enskip (0.1) \\ 
    4 & LS & pricesKO prices2KO KOind payoff & 3.2\enskip (0.3) & 4.0\enskip (0.3) & 3.2\enskip (0.2) \\[0.5em]
  4 & PO & prices & 25.9\enskip (0.8) & 27.6\enskip (0.9) & 27.4\enskip (1.1) \\ 
    4 & PO & pricesKO KOind payoff & 65.0\enskip (5.0) & 79.7\enskip (9.0) & 60.9\enskip (3.7) \\ 
    4 & PO & pricesKO KOind payoff & 86.0\enskip (4.3) & 89.9\enskip (4.9) & 82.3\enskip (4.2) \\ 
    & & maxpriceKO max2priceKO \\
    4 & PO & pricesKO prices2KO KOind payoff & 165.9\enskip (10.1) & 157.5\enskip (7.2) & 125.8\enskip (4.4) \\[0.5em]
     4 & Tree & payoff time & 10.8\enskip (1.5) & 6.4\enskip (0.6) & 4.9\enskip (0.4) \\ 
    4 & Tree & prices & 22.3\enskip (2.1) & 39.2\enskip (4.6) & 30.8\enskip (2.4) \\ 
    4 & Tree & prices payoff & 4.3\enskip (0.5) & 3.1\enskip (0.3) & 4.0\enskip (0.4) \\ 
    4 & Tree & prices time & 96.1\enskip (7.3) & 141.8\enskip (21.9) & 75.9\enskip (10.6) \\ 
    4 & Tree & prices time payoff & 18.6\enskip (2.8) & 11.9\enskip (1.5) & 8.4\enskip (0.6) \\ 
    4 & Tree & prices time payoff KOind & 19.6\enskip (2.0) & 11.5\enskip (0.8) & 9.4\enskip (1.1) \\ \midrule
 16 & LS & one & 1.1\enskip (0.0) & 1.1\enskip (0.0) & 1.1\enskip (0.0) \\ 
   16 & LS & prices & 2.1\enskip (0.0) & 2.1\enskip (0.0) & 2.0\enskip (0.0) \\ 
   16 & LS & pricesKO & 2.4\enskip (0.2) & 2.1\enskip (0.2) & 2.3\enskip (0.2) \\ 
   16 & LS & pricesKO KOind & 2.2\enskip (0.1) & 2.4\enskip (0.2) & 2.1\enskip (0.1) \\ 
   16 & LS & pricesKO KOind payoff & 2.2\enskip (0.1) & 1.9\enskip (0.1) & 2.3\enskip (0.2) \\ 
   16 & LS & pricesKO KOind payoff  & 2.3\enskip (0.2) & 2.2\enskip (0.1) & 1.9\enskip (0.1) \\ 
   & & maxpriceKO \\
   16 & LS & pricesKO KOind payoff  & 2.2\enskip (0.2) & 2.1\enskip (0.3) & 2.1\enskip (0.1) \\ 
   & & maxpriceKO max2priceKO \\
   16 & LS & pricesKO payoff & 1.9\enskip (0.2) & 2.0\enskip (0.2) & 2.1\enskip (0.2) \\ 
   16 & LS & pricesKO prices2KO KOind payoff & 36.0\enskip (1.8) & 35.7\enskip (1.8) & 33.5\enskip (1.0) \\[0.5em]
 16 & PO & prices & 59.4\enskip (0.9) & 59.4\enskip (1.7) & 53.5\enskip (1.3) \\ 
   16 & PO & pricesKO KOind payoff & 114.7\enskip (6.0) & 85.6\enskip (4.3) & 60.6\enskip (2.0) \\ 
   16 & PO & pricesKO KOind payoff  & 96.6\enskip (2.9) & 85.2\enskip (3.2) & 67.3\enskip (2.9) \\
& &    maxpriceKO max2priceKO \\ [0.5em]
 16 & Tree & payoff time & 4.0\enskip (0.1) & 3.5\enskip (0.1) & 2.2\enskip (0.1) \\ 
   16 & Tree & prices & 281.1\enskip (11.3) & 260.8\enskip (8.4) & 202.9\enskip (4.1) \\ 
   16 & Tree & prices payoff & 10.0\enskip (0.1) & 9.4\enskip (0.2) & 9.1\enskip (0.2) \\ 
   16 & Tree & prices time & 195.0\enskip (2.2) & 175.3\enskip (3.2) & 150.8\enskip (13.6) \\ 
   16 & Tree & prices time payoff & 17.5\enskip (0.2) & 16.0\enskip (0.3) & 10.0\enskip (0.2) \\ 
   16 & Tree & prices time payoff KOind & 17.8\enskip (0.2) & 16.0\enskip (0.2) & 10.3\enskip (0.1) \\ \bottomrule
\end{tabular}
\end{table}

\clearpage

\subsection{Example of LS policy}
\label{appendix:results_LS_policy_example}

Table~\ref{table:LS_policy_example} below provides an example of a policy produced by LS (namely, the regression coefficients for predicting the continuation value at each $t$) for the basis function architecture consisting of \textsc{pricesKO}, \textsc{KOind} and \textsc{payoff}.

\begin{table}[ht]
\scriptsize
\centering
\begin{tabular}{rrrrrrrrrrr} \toprule
$t$ & $r_{p_1(t) y(t)}$ & $r_{p_2(t) y(t)}$ & $r_{p_3(t) y(t)}$ & $r_{p_4(t) y(t)}$ & $r_{p_5(t) y(t)}$ & $r_{p_6(t) y(t)}$ & $r_{p_7(t) y(t)}$ & $r_{p_8(t) y(t)}$ & $r_{y(t)}$ & $r_{g(t)}$ \\  \midrule
1 & 0.187 & 0.157 & 0.093 & 0.141 & 0.153 & 0.096 & 0.153 & 0.116 & -53.234 & 1.244 \\ 
  2 & 0.160 & 0.137 & 0.117 & 0.105 & 0.137 & 0.103 & 0.138 & 0.108 & -46.019 & 0.938 \\ 
  3 & 0.166 & 0.136 & 0.122 & 0.117 & 0.139 & 0.091 & 0.151 & 0.122 & -49.177 & 0.253 \\ 
  4 & 0.163 & 0.126 & 0.130 & 0.113 & 0.137 & 0.091 & 0.151 & 0.132 & -49.359 & 0.136 \\ 
  5 & 0.148 & 0.127 & 0.122 & 0.114 & 0.121 & 0.092 & 0.140 & 0.134 & -45.788 & 0.223 \\ 
  6 & 0.145 & 0.122 & 0.129 & 0.124 & 0.118 & 0.088 & 0.146 & 0.132 & -46.506 & 0.183 \\ 
  7 & 0.137 & 0.118 & 0.126 & 0.121 & 0.119 & 0.095 & 0.139 & 0.130 & -45.206 & 0.188 \\ 
  8 & 0.131 & 0.114 & 0.123 & 0.117 & 0.115 & 0.096 & 0.132 & 0.130 & -43.315 & 0.189 \\ 
  9 & 0.128 & 0.108 & 0.125 & 0.114 & 0.111 & 0.100 & 0.128 & 0.125 & -42.122 & 0.178 \\ 
  10 & 0.130 & 0.106 & 0.125 & 0.110 & 0.111 & 0.098 & 0.127 & 0.124 & -41.591 & 0.169 \\ 
  11 & 0.126 & 0.101 & 0.121 & 0.114 & 0.109 & 0.094 & 0.126 & 0.123 & -40.556 & 0.178 \\ 
  12 & 0.118 & 0.096 & 0.121 & 0.108 & 0.108 & 0.095 & 0.121 & 0.114 & -38.273 & 0.201 \\ 
  13 & 0.120 & 0.093 & 0.119 & 0.107 & 0.107 & 0.094 & 0.120 & 0.112 & -37.695 & 0.208 \\ 
  14 & 0.111 & 0.094 & 0.117 & 0.103 & 0.102 & 0.095 & 0.114 & 0.112 & -36.255 & 0.222 \\ 
  15 & 0.110 & 0.095 & 0.109 & 0.098 & 0.098 & 0.090 & 0.108 & 0.109 & -34.205 & 0.241 \\ 
  16 & 0.107 & 0.092 & 0.106 & 0.097 & 0.096 & 0.089 & 0.107 & 0.107 & -33.227 & 0.248 \\ 
  17 & 0.109 & 0.086 & 0.106 & 0.095 & 0.089 & 0.087 & 0.097 & 0.103 & -31.397 & 0.267 \\ 
  18 & 0.109 & 0.086 & 0.107 & 0.094 & 0.084 & 0.085 & 0.098 & 0.100 & -30.839 & 0.264 \\ 
  19 & 0.108 & 0.082 & 0.105 & 0.095 & 0.083 & 0.085 & 0.094 & 0.097 & -29.996 & 0.263 \\ 
  20 & 0.102 & 0.075 & 0.101 & 0.098 & 0.082 & 0.084 & 0.092 & 0.096 & -28.817 & 0.272 \\ 
  21 & 0.101 & 0.074 & 0.092 & 0.093 & 0.077 & 0.077 & 0.091 & 0.095 & -27.092 & 0.291 \\ 
  22 & 0.099 & 0.068 & 0.084 & 0.098 & 0.077 & 0.077 & 0.087 & 0.093 & -26.011 & 0.296 \\ 
  23 & 0.091 & 0.069 & 0.077 & 0.094 & 0.071 & 0.076 & 0.083 & 0.089 & -23.579 & 0.308 \\ 
  24 & 0.088 & 0.062 & 0.083 & 0.090 & 0.066 & 0.073 & 0.072 & 0.087 & -21.895 & 0.333 \\ 
  25 & 0.086 & 0.067 & 0.078 & 0.089 & 0.070 & 0.071 & 0.075 & 0.083 & -22.130 & 0.332 \\ 
  26 & 0.079 & 0.068 & 0.073 & 0.091 & 0.060 & 0.072 & 0.076 & 0.087 & -21.551 & 0.339 \\ 
  27 & 0.082 & 0.065 & 0.071 & 0.086 & 0.060 & 0.061 & 0.076 & 0.081 & -19.951 & 0.340 \\ 
  28 & 0.079 & 0.064 & 0.069 & 0.081 & 0.056 & 0.060 & 0.076 & 0.076 & -19.046 & 0.373 \\ 
  29 & 0.078 & 0.062 & 0.074 & 0.084 & 0.056 & 0.064 & 0.073 & 0.078 & -19.757 & 0.354 \\ 
  30 & 0.074 & 0.056 & 0.069 & 0.079 & 0.060 & 0.060 & 0.067 & 0.074 & -18.156 & 0.376 \\ 
  31 & 0.074 & 0.056 & 0.066 & 0.081 & 0.055 & 0.058 & 0.069 & 0.072 & -18.033 & 0.384 \\ 
  32 & 0.068 & 0.051 & 0.059 & 0.075 & 0.056 & 0.059 & 0.069 & 0.065 & -16.544 & 0.414 \\ 
  33 & 0.070 & 0.051 & 0.067 & 0.072 & 0.055 & 0.059 & 0.070 & 0.067 & -17.626 & 0.402 \\ 
  34 & 0.066 & 0.052 & 0.064 & 0.075 & 0.047 & 0.059 & 0.064 & 0.064 & -16.432 & 0.408 \\ 
  35 & 0.068 & 0.045 & 0.059 & 0.071 & 0.049 & 0.051 & 0.059 & 0.056 & -14.880 & 0.447 \\ 
  36 & 0.062 & 0.047 & 0.059 & 0.066 & 0.058 & 0.057 & 0.056 & 0.058 & -15.804 & 0.451 \\ 
  37 & 0.066 & 0.047 & 0.057 & 0.063 & 0.054 & 0.056 & 0.055 & 0.055 & -15.697 & 0.457 \\ 
  38 & 0.062 & 0.044 & 0.059 & 0.064 & 0.057 & 0.052 & 0.056 & 0.060 & -16.549 & 0.467 \\ 
  39 & 0.061 & 0.047 & 0.056 & 0.059 & 0.052 & 0.045 & 0.056 & 0.059 & -16.029 & 0.499 \\ 
  40 & 0.063 & 0.042 & 0.056 & 0.063 & 0.047 & 0.044 & 0.052 & 0.062 & -16.101 & 0.502 \\ 
  41 & 0.059 & 0.045 & 0.048 & 0.058 & 0.049 & 0.048 & 0.053 & 0.054 & -15.502 & 0.511 \\ 
  42 & 0.056 & 0.041 & 0.048 & 0.051 & 0.050 & 0.053 & 0.046 & 0.048 & -15.099 & 0.545 \\ 
  43 & 0.052 & 0.042 & 0.045 & 0.058 & 0.050 & 0.052 & 0.045 & 0.050 & -15.850 & 0.550 \\ 
  44 & 0.053 & 0.039 & 0.046 & 0.050 & 0.049 & 0.041 & 0.047 & 0.037 & -14.233 & 0.579 \\ 
  45 & 0.046 & 0.038 & 0.039 & 0.042 & 0.041 & 0.039 & 0.033 & 0.040 & -11.716 & 0.603 \\ 
  46 & 0.050 & 0.032 & 0.043 & 0.042 & 0.045 & 0.037 & 0.039 & 0.041 & -13.478 & 0.613 \\ 
  47 & 0.039 & 0.034 & 0.039 & 0.039 & 0.041 & 0.038 & 0.036 & 0.036 & -12.441 & 0.639 \\ 
  48 & 0.040 & 0.034 & 0.036 & 0.036 & 0.027 & 0.034 & 0.026 & 0.023 & -9.532 & 0.652 \\ 
  49 & 0.031 & 0.025 & 0.029 & 0.028 & 0.033 & 0.040 & 0.033 & 0.029 & -9.949 & 0.672 \\ 
  50 & 0.023 & 0.018 & 0.027 & 0.026 & 0.027 & 0.027 & 0.024 & 0.031 & -7.687 & 0.716 \\ 
  51 & 0.027 & 0.019 & 0.019 & 0.027 & 0.032 & 0.029 & 0.023 & 0.025 & -8.271 & 0.712 \\ 
  52 & 0.008 & 0.017 & 0.015 & 0.016 & 0.023 & 0.019 & 0.021 & 0.016 & -4.681 & 0.774 \\ 
  53 & 0.011 & 0.014 & 0.019 & 0.011 & 0.008 & 0.014 & 0.009 & 0.016 & -3.458 & 0.808 \\ 
  54 & 0.000 & 0.000 & 0.000 & 0.000 & 0.000 & 0.000 & 0.000 & 0.000 & 0.000 & 0.000 \\  \bottomrule
\end{tabular}
\caption{Example of LS regression coefficients for $n = 8$, with \textsc{pricesKO}, \textsc{KOind} and \textsc{payoff} basis functions. \label{table:LS_policy_example}}
\end{table}

\subsection{Performance in the high-dimensional regime}
\label{appendix:results_highdim}
Thus far in our option pricing experiments, we have studied instances of the problem where the number of stocks, which is the main driver of the size of the state space, is small. In addition, our prior experiments studied a setting where the number of trajectories is large relative to the number of stocks. 

 In this section, we seek to understand how our tree policies as well as Longstaff-Schwartz perform in the high-dimensional regime, where (1) the number of stocks $n$ is large and (2) the number of trajectories $\Omega$ may be quite small relative to the number of stocks $n$. 
 
Towards this goal, we consider an experiment using the same option pricing problem as in Section~\ref{subsec:results_problemdefinition}. We fix the initial price $\bar{p} = 100$, common correlation $\bar{\rho} = 0$ and the number of stocks $n = 1000$. We generate 100,000 trajectories for the test set. We vary the number of trajectories generated for the training set $\Omega$ in the set $\{1000, 1500, 2000, 5000, 10000, 20000\}$. For each value of $\Omega$, we perform 10 replications. In each replication, we run LS with \textsc{pricesKO}, \textsc{payoff} and \textsc{KOind} as the basis functions, and we run our construction algorithm with the complete state variable specification (\textsc{time}, \textsc{prices}, \textsc{payoff}, \textsc{KOind}), with a relative improvement tolerance $\gamma = 0.005$. 

\begin{figure}
\centering
\includegraphics[width=0.7\textwidth]{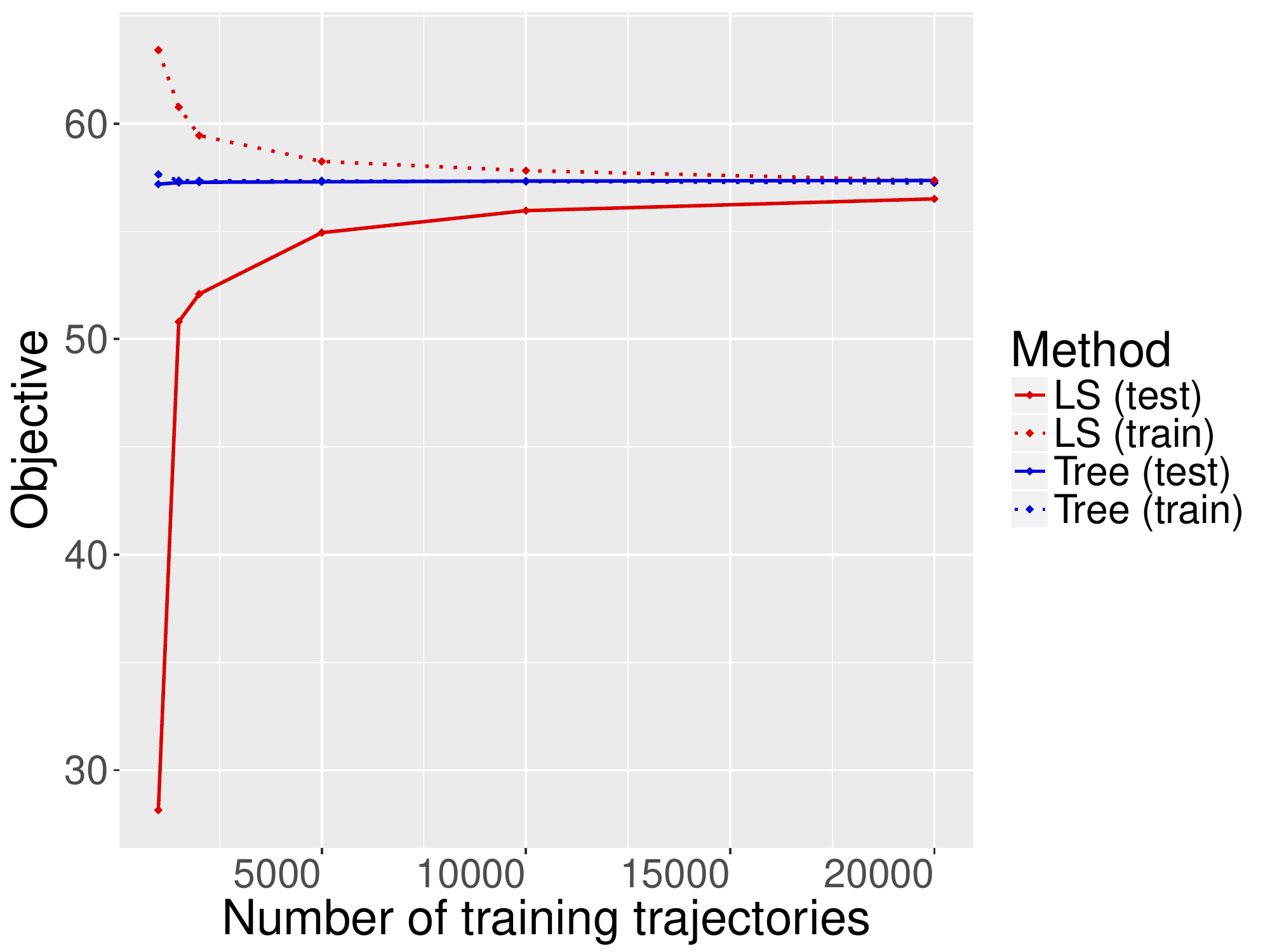}
\caption{Training set and test set performance of LS and tree policies in high-dimensional example as a function of number of training trajectories $\Omega$. \label{figure:highdim_plot} }
\end{figure}

Figure~\ref{figure:highdim_plot} compares the performance of the two methods as a function of $\Omega$. For each method, the in-sample (training set) reward is shown by the dashed curve, with the out-of-sample/test set reward shown by the solid curve. 

From this plot, we obtain several insights. First, even when the number of training trajectories $\Omega = 20,000$, the tree policies exhibit an edge over the LS policies (mean reward of 57.36 compared to 56.51 for LS). Second, as the number of training trajectories decreases, the out-of-sample performance of the tree policy remains very stable. In contrast, the out-of-sample performance of LS deteriorates, dropping to about 50\% of its performance with $\Omega = 20,000$. 

Interestingly, the training set performance of the two methods also behaves differently. Whereas the training set performance of the tree policies is close to the test set performance and mirrors its stability, the training set performance of LS increases as the number of trajectories decreases. The reason for this is due to the use of least squares in the LS algorithm. The LS algorithm works by running a least squares regression at each period starting from $t = T-1$ to $t = 1$, where the ``observations'' in this regression correspond to trajectories, and the dependent variable is the continuation value (specifically, the reward from following the LS policy in future periods). As the number of trajectories decreases, the regression model estimated in LS at each period will be more and more overfitted to the continuation values. When the number of trajectories is the same as or lower than the number of basis functions, the regression model at each period will be able to perfectly predict the continuation value; it is not difficult to see that such a policy will essentially behave in the same way as the (unattainable) perfect foresight policy on the training set. In contrast, the tree policies appear to be resistant to this type of overfitting.

\subsection{Experiment with high effective dimension}
\label{appendix:results_sparseKO}

In our previous experiments, we considered a max-call option with a knock-out barrier. This setup leads to a high-dimensional option pricing problem, where the number of state variables scales with the number of underlying assets. However, the option pricing instances were structured in a way that both the payoff of this option and the knock-out barrier are essentially dependent on the maximum price. This is perhaps one reason why a tree policy that judiciously splits on the time and payoff can perform effectively. One could thus argue that, while the nominal dimension of the optimal stopping problem may be large because of the large number of assets, the ``effective'' dimension is small, and the problem is thus an ``easy'' stopping problem. 

In this section, we consider another family of option pricing instances where the effective dimension of the problem is large. We define these instances as follows. Given $n$ assets, the first $m$ assets drive the knock-out behavior. In particular, the stock is knocked out at time $t$ and the payoff becomes zero for all $t' \geq t$ if the price $p_i$ of any asset $i \in \{1,\dots,m\}$ exceeds the barrier price $B$. The corresponding state variable $y(t)$, which indicates whether the option has not been knocked out by time $t$, is defined as
\begin{equation*}
y(t) = \Ibb \left\{ \max_{1 \leq j \leq m, 1 \leq t' \leq t} p_j(t) < B \right\}
\end{equation*}
The payoff of the option is defined using the maximum price of the remaining $n - m$ assets, i.e.,
\begin{equation*}
g(t) = \max \left\{ 0, \max_{m+1 \leq j \leq n} p_j(t) - K \right\} \cdot y(t).
\end{equation*}
The prices of the assets, $p_1(t), \dots, p_n(t)$, follow the same geometric Brownian motion dynamics as in Section~\ref{subsec:results_problemdefinition}. We test $n = 16$ and $m \in \{4, 8\}$, and initial prices $\bar{p} \in \{90, 100, 110\}$. 

To understand the rationale behind these instances, suppose that we fix a particular $t$ and a particular value of $g(t)$, and consider the following two situations. In the first situation, suppose that the asset prices $p_1(t), \dots, p_m(t)$ are far from the knock-out barrier $B$. For this situation, if there is enough time remaining in the horizon and assuming a positive drift in the asset prices, we would expect the optimal action to be to continue, as the asset prices $p_{m+1}(t), \dots, p_{n}(t)$ will continue to grow, and we are not in danger of $p_1(t), \dots, p_m(t)$ exceeding $B$ and eliminating the payoff completely. In the second situation, suppose one of the first $m$ assets may have a price that is close to the barrier $B$. For this situation, it is reasonable to expect that the optimal action will be to stop, because there is a risk that one of the $m$ assets will cross the barrier and force the payoff to zero. 

Thus, in these instances, we expect that a tree policy that splits solely on the time and the payoff should not perform well, precisely because it cannot recognize those states where there is a risk of a knock-out occurring. In contrast, a tree policy that can additionally split on the asset prices should be able to recognize these risky states, and achieve better performance.

In Table~\ref{table:sparseKO_performance}, we display the performance of the tree algorithm with two different choices of state variables: \textsc{payoff} and \textsc{time}; and \textsc{prices}, \textsc{time}, \textsc{payoff} and \textsc{KOind}. We additionally compare the performance of these two policies against LS with \textsc{pricesKO}, \textsc{KOind} and \textsc{payoff} as the basis functions. From this table, we can see that in general, the tree policy that uses only \textsc{payoff} and \textsc{time} achieves appreciably lower rewards than the policy that also uses \textsc{prices}, as we would expect. 

\begin{table}
\small
\centering
\begin{tabular}{lllllll} \toprule 
$n$ & $m$ & Method & State variables / basis functions & $\bar{p} = 90$ & $\bar{p} = 100$ & $\bar{p} = 110$ \\ \midrule
   16 &   4 & LS & \textsc{pricesKO}, \textsc{KOind}, \textsc{payoff} & 51.62\enskip (0.034) & 56.65\enskip (0.042) & 58.82\enskip (0.027) \\ 
   16 &   4 & Tree & \textsc{payoff}, \textsc{time}  & 48.21\enskip (0.021) & 50.64\enskip (0.045) & 52.10\enskip (0.053) \\ 
   16 &   4 & Tree & \textsc{prices}, \textsc{time}, \textsc{payoff}, \textsc{KOind}  & \bfseries 56.74\enskip (0.072) & \bfseries 59.46\enskip (0.218) & \bfseries 64.19\enskip (0.124) \\[0.5em]
   16 &   8 & LS & \textsc{pricesKO}, \textsc{KOind}, \textsc{payoff} & 37.21\enskip (0.022) & 40.70\enskip (0.020) & 43.79\enskip (0.034) \\ 
   16 &   8 & Tree & \textsc{payoff}, \textsc{time}  & 34.73\enskip (0.026) & 37.73\enskip (0.019) & 41.03\enskip (0.043) \\ 
   16 &   8 & Tree & \textsc{prices}, \textsc{time}, \textsc{payoff}, \textsc{KOind} & \bfseries 39.06\enskip (0.080) & \bfseries 43.52\enskip (0.060) & \bfseries 46.30\enskip (0.069) \\[0.5em]
   16 &  12 & LS & \textsc{pricesKO}, \textsc{KOind}, \textsc{payoff} & 23.65\enskip (0.029) & 27.73\enskip (0.014) & \bfseries 32.47\enskip (0.019) \\ 
   16 &  12 & Tree & \textsc{payoff}, \textsc{time} & 22.57\enskip (0.023) & 26.25\enskip (0.033) & 31.05\enskip (0.022) \\ 
   16 &  12 & Tree & \textsc{prices}, \textsc{time}, \textsc{payoff}, \textsc{KOind} & \bfseries 24.83\enskip (0.049) & \bfseries 29.38\enskip (0.043) & 31.72\enskip (0.276) \\ 
   \bottomrule
\end{tabular}
\caption{Performance for effective dimension option pricing instances. \label{table:sparseKO_performance}}
\end{table}

\begin{figure}
\begin{center}
\includegraphics[width=0.7\textwidth]{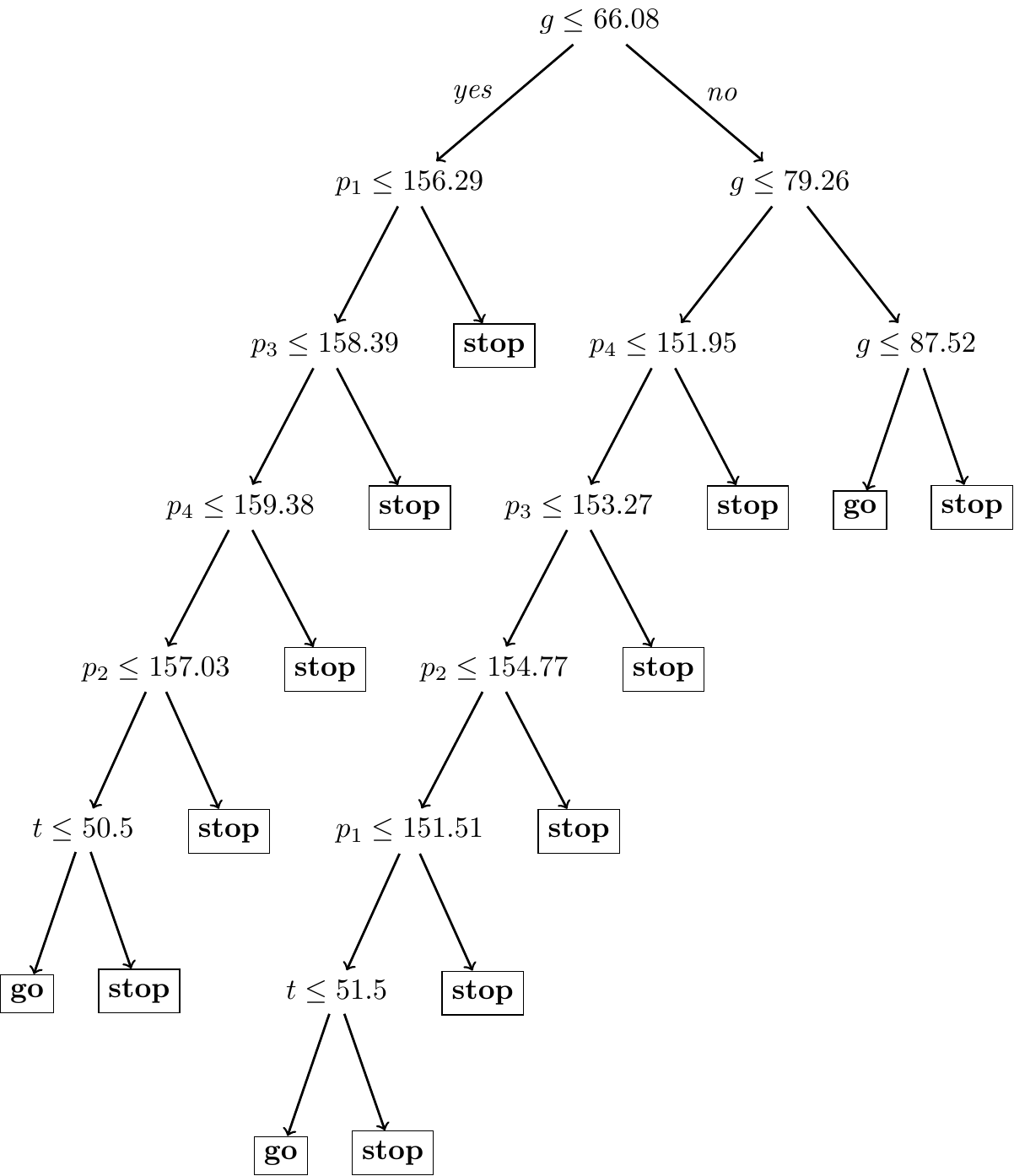}
\end{center}
\caption{Example of tree policy for an instance with $n = 16$, $m = 4$, $\bar{p} = 110$ and with \textsc{prices}, \textsc{time}, \textsc{payoff} and \textsc{KOind} as state variables. \label{figure:sparseKO_example_tree}}
\end{figure}

In addition, it is also helpful to examine the tree policies that emerge for these instances. Figure~\ref{figure:sparseKO_example_tree} displays an example of a tree policy for one instance with $n = 16$ and $m = 4$. Observe that while the tree retains some of the salient features of the tree policies shown earlier for the simpler option pricing instances (namely splitting on $g$ and $t$; see Figure~\ref{figure:example_tree_policies_peq90}), the tree also includes a chain of splits on the first 4 asset prices, as one would expect for these instances. Thus, the tree algorithm was able to recognize that certain values of these state variables are predictive of future payoffs. We believe that this is indicative of the potential of our greedy construction procedure to scale to instances where optimal and near optimal policies may exhibit a more complicated structure than those studied in Section~\ref{subsec:results_policy_structure}.

\end{document}